\documentclass{article}

\usepackage[utf8]{inputenc}
\usepackage[T1]{fontenc}
\usepackage[british]{babel}
\usepackage[style=british]{csquotes}
\usepackage[en-GB]{datetime2}
\DTMlangsetup[en-GB]{ord=omit}
\usepackage{xcolor}

\usepackage{amssymb}
\usepackage{amsmath}
\usepackage[amsmath,hyperref,thmmarks]{ntheorem}
\usepackage{mathtools}
\usepackage[bookmarksnumbered]{hyperref}

\usepackage{tikz}
\usepackage{tikz-cd}

\usepackage[lining,scaled=.8]{FiraMono}
\usepackage[tt=false,semibold]{libertine}
\usepackage[libertine,smallerops,vvarbb]{newtxmath}
\usepackage[cal=euler,scr=boondoxo,frak=euler]{mathalpha}
\tikzcdset{
    arrow style=tikz,
    arrows={/tikz/line width=.5pt},
    diagrams={>={Straight Barb[scale=0.8]}}
}

\usepackage{bm}
\usepackage{cases}
\usepackage{accents}

\usepackage{setspace}
\usepackage{subdepth} 

\newcommand{\initlengths}{%
    \setlength{\abovedisplayshortskip}{3pt plus 9pt minus 3pt}%
    \setlength{\belowdisplayshortskip}{9pt plus 9pt minus 9pt}%
    \setlength{\abovedisplayskip}{9pt plus 9pt minus 9pt}%
    \setlength{\belowdisplayskip}{9pt plus 9pt minus 9pt}%
}

\usepackage{titling}

\usepackage[immediate]{silence}
\usepackage{sectsty}
\DeactivateWarningFilters[temp]
\allsectionsfont{\bfseries\boldmath}

\makeatletter
\def\theorem@optheaderfont{\normalfont}
\renewtheoremstyle{plain}%
    {\item[\theorem@headerfont\hskip\labelsep ##1\ ##2\theorem@separator]}%
    {\item[\theorem@headerfont\hskip\labelsep ##1\ ##2\ {\theorem@optheaderfont(##3)}\theorem@separator]}
\newtheoremstyle{prooflike}%
    {\item[\theorem@headerfont\hskip\labelsep ##1\theorem@separator]}%
    {\item[\theorem@headerfont\hskip\labelsep ##3\theorem@separator]}
\makeatother

\theoremseparator{.}
\newtheorem{theorem}{Theorem}[section]
\newtheorem{lemma}[theorem]{Lemma}

\newtheorem{proposition}[theorem]{Proposition}
\newtheorem{conjecture}[theorem]{Conjecture}
\theoremsymbol{\ensuremath{\Box}}
\newtheorem{theoremq}[theorem]{Theorem}

\theorembodyfont{\normalfont}
\theoremsymbol{\ensuremath{\lhd}}
\newtheorem{definition}[theorem]{Definition}
\newtheorem{remark}[theorem]{Remark}
\newtheorem{example}[theorem]{Example}
\newtheorem{construction}[theorem]{Construction}
\newtheorem{situation}[theorem]{Situation}
\newtheorem{notation}[theorem]{Notation}
\newtheorem{assumption}[theorem]{Assumption}
\theoremstyle{prooflike}
\theoremheaderfont{\normalfont\itshape}
\theoremsymbol{\ensuremath{\Box}}
\newtheorem{proof}{Proof}
\theoremheaderfont{\bfseries}
\newtheorem{bproof}{Proof}

\numberwithin{equation}{section}
\newcounter{subequation}
\numberwithin{subequation}{equation}

\usepackage{enumitem}
\setlist[enumerate]{label={\upshape(\roman*)}}

\usepackage[labelsep=period, labelfont={bf}]{caption}

\usepackage[
    backend=biber,
    style=oxnum,
    giveninits,
    scnames,
    maxbibnames=5,
    sorting=nyt
]{biblatex}

\DeclareFontFamily{OMX}{MnSymbolE}{}
\DeclareSymbolFont{MnLargeSymbols}{OMX}{MnSymbolE}{m}{n}
\SetSymbolFont{MnLargeSymbols}{bold}{OMX}{MnSymbolE}{b}{n}
\DeclareFontShape{OMX}{MnSymbolE}{m}{n}{
    <-6>  MnSymbolE5
   <6-7>  MnSymbolE6
   <7-8>  MnSymbolE7
   <8-9>  MnSymbolE8
   <9-10> MnSymbolE9
  <10-12> MnSymbolE10
  <12->   MnSymbolE12
}{}
\DeclareFontShape{OMX}{MnSymbolE}{b}{n}{
    <-6>  MnSymbolE-Bold5
   <6-7>  MnSymbolE-Bold6
   <7-8>  MnSymbolE-Bold7
   <8-9>  MnSymbolE-Bold8
   <9-10> MnSymbolE-Bold9
  <10-12> MnSymbolE-Bold10
  <12->   MnSymbolE-Bold12
}{}
\DeclareMathDelimiter{[}    {\mathopen} {MnLargeSymbols}{'000}{MnLargeSymbols}{'000}
\DeclareMathDelimiter{]}    {\mathclose}{MnLargeSymbols}{'005}{MnLargeSymbols}{'005}
\DeclareMathDelimiter{\llbr}{\mathopen} {MnLargeSymbols}{'102}{MnLargeSymbols}{'102}
\DeclareMathDelimiter{\rrbr}{\mathclose}{MnLargeSymbols}{'107}{MnLargeSymbols}{'107}
                     
\makeatletter
\def\big#1{{\hbox{$\left#1\vbox to8.5\p@{}\right.\n@space$}}}
\makeatother

\newcommand{\Aone}{$\mathbb{A}^{\<1}$}
\newcommand{\Aut}{\opname{Aut}}
\newcommand{\bbA}{\mathbb{A}}
\newcommand{\bbC}{\mathbb{C}}

\newcommand{\bbK}{\mathbb{K}}
\newcommand{\bbL}{\mathbb{L}}
\newcommand{\bbN}{\mathbb{N}}
\newcommand{\bbQ}{\mathbb{Q}}
\newcommand{\bbR}{\mathbb{R}}
\newcommand{\bbT}{\mathbb{T}}
\newcommand{\bbZ}{\mathbb{Z}}
\newcommand{\BO}{\mathrm{BO}}
\newcommand{\BSp}{\mathrm{BSp}}
\newcommand{\BU}{\mathrm{BU}}
\newcommand{\calA}{\mathcal{A}}

\newcommand{\calC}{\mathcal{C}}

\newcommand{\calExt}[1][]{{\smash{#1{\mathcal{E}}}\mspace{-2mu}\mathit{xt}}}
\newcommand{\calHom}[1][]{{\smash{#1{\mathcal{H}}}\mspace{-2mu}\mathit{om}}}

\newcommand{\calM}{\mathcal{M}}
\newcommand{\calMul}{\mathord{\mspace{1.6mu}\underaccent{\bar}{\mspace{-1.6mu}{\calM}}}}
\newcommand{\calO}{\mathcal{O}}
\newcommand{\calP}{\mathcal{P}}
\newcommand{\calU}{\mathcal{U}}
\newcommand{\calV}{\mathcal{V}}
\newcommand{\calX}{\mathcal{X}}

\newcommand{\cat}[1]{\mathord{\mathsf{#1}}}
\newcommand{\CH}{\mathrm{CH}}
\newcommand{\ch}{\opname{ch}}
\newcommand{\citestacks}[1]{\cite[Tag~\href{https://stacks.math.columbia.edu/tag/#1}{\texttt{#1}}]{Stacks}}
\newcommand{\cl}{{\smash{\mathrm{cl}}}}

\newcommand{\Db}{\cat{D}^{\mathrm{b}}}
\newcommand{\dprime}{\mathord{\mathchoice{\prime \< \prime}{\prime \< \prime}{\prime \< \prime}{\> \prime \prime}}}
\newcommand{\dslash}{\mathord{/\mspace{-5mu}/}}
\newcommand{\dslashsub}[1]{\mathord{/\mspace{-5mu}/\mspace{-4mu}_{#1}\mspace{3mu}}}

\newcommand{\ev}{\opname{ev}}
\newcommand{\Ext}{\opname{Ext}}
\newcommand{\frg}{\mathfrak{g}}
\newcommand{\frM}{\mathfrak{M}}
\newcommand{\frm}{\mathfrak{m}}

\newcommand{\frS}{\mathfrak{S}}
\newcommand{\frt}{\mathfrak{t}}
\newcommand{\frX}{\mathfrak{X}}

\newcommand{\fund}[1]{[#1]_{\mathrm{fund}}}
\newcommand{\GL}{\mathrm{GL}}
\newcommand{\Gm}{\mathbb{G}_{\mathrm{m}}}
\newcommand{\heart}{\mathbin{\heartsuit}}

\newcommand{\Hom}{\opname{Hom}}
\newcommand{\id}{\mathrm{id}}
\newcommand{\inftyCAT}{\infty\mathhyphen\textsf{\textup{\textsc{Cat}}}}
\newcommand{\inftyCat}{\infty\mathhyphen\cat{Cat}}
\newcommand{\inftyGpd}{\infty\mathhyphen\cat{Gpd}}

\newcommand{\inv}[1]{[#1]_{\mathrm{inv}}}
\newcommand{\itExt}[1][]{{\smash{#1{E}}\mspace{-1.5mu}\mathit{xt}}}
\newcommand{\llparen}{\mathopen{(\mspace{-3.2mu}(}}
\newcommand{\longhookrightarrow}{\lhook\joinrel\longrightarrow}
\newcommand{\longsimto}{\mathrel{\overset{\smash{\raisebox{-.8ex}{$\sim$}}\mspace{3mu}}{\longrightarrow}}}
\newcommand{\mathhyphen}{\textnormal{-}}
\newcommand{\mathscalebox}[2]{\scalebox{#1}{$#2$}}
\newcommand{\numberthis}{\addtocounter{equation}{1}\tag{\theequation}}
\newcommand{\Omegaul}{\mathord{\mspace{1.6mu}\underaccent{\bar}{\mspace{-1.6mu}{\Omega}}}}
\newcommand{\otimeshat}{\mathbin{\hat{\otimes}}}
\newcommand{\op}{{\smash{\mathrm{op}}}}
\newcommand{\OPerf}{\cat{OPerf}}
\newcommand{\opname}[1]{\operatorname{\smash[b]{\mathrm{#1}}}}
\newcommand{\Perf}{\cat{Perf}}
\newcommand{\pl}{{\smash{\mathrm{pl}}}}
\newcommand{\plss}{{\smash{\mathrm{pl},\mspace{1mu}\mathrm{ss}}}}
\newcommand{\pr}{\opname{pr}}
\renewcommand{\qed}{\hfill$\Box$}
\newcommand{\rank}{\opname{rank}}
\newcommand{\res}{\opname{res}}

\newcommand{\RcalHom}{\operatorname{\mathbb{R}\calHom}}
\renewcommand{\ring}[1]{\accentset{\smash{\raisebox{-0.12ex}{$\scriptscriptstyle\circ$}}}{#1}}
\newcommand{\rrparen}{\mathclose{)\mspace{-2.5mu})}}
\newcommand{\sd}{{\smash{\mathrm{sd}}}}
\newcommand{\sdss}{{\smash{\mathrm{sd},\mspace{1mu}\mathrm{ss}}}}
\newcommand{\sdst}{{\smash{\mathrm{sd},\mspace{1mu}\mathrm{st}}}}
\newcommand{\second}{^{\dprime}}
\newcommand{\sharps}{{\smash{\sharp}}}
\newcommand{\sharpsd}{{\smash{\sharp, \> \sd}}}
\newcommand{\simto}{\mathrel{\overset{\smash{\raisebox{-.8ex}{$\sim$}}\mspace{3mu}}{\to}}}
\newcommand{\simTo}{\mathrel{\overset{\smash{\raisebox{-.7ex}{$\sim$}}\mspace{3mu}}{\Rightarrow}}}
\newcommand{\SO}{\mathrm{SO}}
\newcommand{\Sp}{\mathrm{Sp}}
\newcommand{\Spec}{\opname{Spec}}
\newcommand{\SpPerf}{\cat{SpPerf}}
\renewcommand{\ss}{{\smash{\mathrm{ss}}}}

\newcommand{\sss}[2]{_{\smash{#1}}^{\smash{#2}}}
\newcommand{\supp}{\opname{supp}}
\newcommand{\Sym}{\opname{Sym}}

\newcommand{\td}{\opname{td}}
\newcommand{\Thetaul}[1][]{\mathord{\mspace{1.6mu}\underaccent{\bar}{\mspace{-1.6mu}#1{\Theta}}}}
\newcommand{\tria}{{\mathbin{\vartriangle}}}
\newcommand{\tw}{{\smash{\mathrm{tw}}}}
\newcommand{\upB}{\mathrm{B}}

\newcommand{\upe}{\mathrm{e}}
\newcommand{\upi}{\mathrm{i}}
\newcommand{\upJ}{\mathrm{J}}
\newcommand{\upO}{\mathrm{O}}
\newcommand{\upU}{\mathrm{U}}
\newcommand{\vdim}{\opname{vdim}}
\newcommand{\virt}[1]{[#1]_{\mathrm{virt}}}
\newcommand*{\<}{\mspace{-2mu}}
\renewcommand*{\>}{\mspace{2mu}}

\newcommand{\nicesubstack}[1]{\substack{\makebox{\scriptsize$\begin{gathered}#1\end{gathered}$}}}
\newcommand{\leftsubstack}[2][6em]{\substack{\makebox[#1][l]{\scriptsize$\begin{aligned}#2\end{aligned}$}}}

\newcommand{\LSs}{1.5ex}
\newcommand{\LS}{2.5ex}
\newcommand{\LSl}{3.5ex}
\newcommand{\LAP}[2]{\mathord{\hspace{#1}\mathclap{#2}\hspace{#1}}}

\title{Enumerative invariants in self-dual categories.\\II. Homological invariants}
\author{Chenjing Bu}
\date{}

\begin{document}

\initlengths

\maketitle

\begin{abstract}
    This is the second paper in a series
on enumerative invariants counting
self-dual objects in self-dual categories.
Ordinary enumerative invariants in abelian categories
can be seen as invariants for the structure group $\GL (n)$,
and our theory is an extension of this
to structure groups $\upO (n)$ and $\Sp (2n)$.
Examples of our theory include
counting principal orthogonal or symplectic bundles,
and counting self-dual quiver representations.

In the present paper, we propose a conjectural picture
on homological enumerative invariants counting
self-dual objects in self-dual categories,
which are homology classes lying in the ordinary homology of moduli stacks.
This is a self-dual analogue of the conjectures of
Gross--Joyce--Tanaka~\cite{GrossJoyceTanaka}.
We study algebraic structures arising from the homology of these moduli stacks,
including vertex algebras and twisted modules,
and we formulate wall-crossing formulae for the invariants
using these algebraic structures.
We also provide a partial proof of our conjecture
in the case of self-dual quiver representations.

\end{abstract}

\tableofcontents

\clearpage
\section{Introduction}

This is the second paper in a series
studying enumerative invariants
counting \emph{self-dual objects in self-dual categories},
and is a sequel to the author's work \cite{Bu2023}.
\nopagebreak
\vskip\topsep

One major branch of enumerative geometry
is the study of moduli stacks of objects in an \emph{abelian category},
such as the category of coherent sheaves
on a smooth projective variety, or the category of
representations of a quiver.
Numerous types of enumerative invariants of this kind
have been studied in the literature,
including \emph{Donaldson--Thomas invariants}
studied by Thomas~\cite{Thomas2000}
and generalized by Joyce--Song~\cite{JoyceSong2012}
and Kontsevich--Soibelman~\cite{KontsevichSoibelman2008},
\emph{Donaldson invariants} studied by Mochizuki~\cite{Mochizuki2009},
Joyce's recent work~\cite{Joyce2021} on homological enumerative invariants,
and so on.

On the other hand, these are all part of a more general problem,
which is that of counting \emph{principal $G$-bundles} on a variety,
for an algebraic group $G$,
or some coherent sheaf version of them.
One could also consider quiver analogues of this problem,
using a suitable notion of \emph{$G$-quivers},
such as the notion defined in Derksen--Weyman~\cite{DerksenWeyman2002},
or a modification of it.
We shall refer to such invariants as \emph{$G$-invariants}.

Theories of invariants in abelian categories, as discussed above,
address the cases when $G = \mathrm{GL} (n)$, or sometimes~$\mathrm{SL} (n)$,
which belong to the A family of algebraic groups.
One advantage in this case is that the objects in question
often form an \emph{abelian category},
which provides rich structures one could work with.

The present series of papers address the cases when
$G = \upO (n)$~or $\Sp (2n)$, which cover the B, C, and D families of algebraic groups.
In this case, the objects in question are often
\emph{self-dual objects in a self-dual category}.
\nopagebreak
\vskip\topsep

The basic setting for studying these Type~B/C/D enumerative invariants
involves the following data:

\begin{itemize}
    \item 
        A \emph{self-dual category} $\calA$,
        with moduli stack $\calM$.
        Here, \emph{self-dual} means that there is an equivalence
        $(-)^\vee \colon \calA \simto \calA^\op$,
        which squares to the identity of~$\calA$.
    \item 
        A \emph{category of self-dual objects} $\calA^\sd$,
        with moduli stack $\calM^\sd$.
        Here, a \emph{self-dual object} is an object $E \in \calA$
        with an isomorphism $E \simto E^\vee$
        that squares to the identity of~$E$.
\end{itemize}
Our two main examples are as follows:

\begin{itemize}
    \item 
        $\calA$ is the category of vector bundles on a variety $X$,
        and $\calA^\sd$ is the category of orthogonal or symplectic
        principal bundles on $X$.
    \item
        $\calA$ is the category of representations of a \emph{self-dual quiver} $Q$,
        and $\calA^\sd$ is the category of \emph{self-dual representations}.
\end{itemize}
Here, \emph{self-dual quivers} are a natural quiver analogue
of Type B/C/D principal bundles.
They were defined by Derksen--Weyman~\cite{DerksenWeyman2002}
under the name of \emph{symmetric quivers},
and studied in the context of enumerative geometry
by Young~\cite{Young2015,Young2020}.
\nopagebreak
\vskip\topsep

The first paper in the series~\cite{Bu2023}
established a theory of \emph{motivic invariants}
in self-dual categories,
which is a Type B/C/D analogue of
Joyce's series of work~\cite{Joyce2006I,Joyce2007II,Joyce2007III,Joyce2008IV}
on motivic invariants in abelian categories.
We defined enumerative invariants as elements in a ring of motives,
and we established \emph{wall-crossing formulae} for these invariants
when changing the stability condition.
\nopagebreak
\vskip\topsep

In the present paper, we consider a different flavour of enumerative invariants,
which we call \emph{homological invariants},
defined as elements in the ordinary homology of moduli stacks.

For Type~A invariants, this approach was taken by
Gross--Joyce--Tanaka~\cite{GrossJoyceTanaka},
who formulated conjectures on universal structures of these invariants.
These conjectures were subsequently proved by Joyce~\cite{Joyce2021}
in several cases, including counting sheaves on curves or surfaces, 
counting quiver representations, and several others.

Joyce~\cite{Joyce2021,JoyceHall} showed that for a moduli stack $\calM$
of objects in an abelian category, under certain assumptions,
the homology $H_* (\calM; \bbQ)$ carries the structure
of a \emph{vertex algebra}.
See, for example, Frenkel--Ben-Zvi~\cite{FrenkelBenZvi2004}
for an account of vertex algebras.
These are complicated algebraic structures
used to describe local structures of a
\emph{conformal field theory} in physics,
and it is a surprising phenomenon that
such a structure arises in the context of moduli theory.

Each vertex algebra $V$ can be associated with a Lie algebra $V / D (V)$,
where $D$ is the translation operator.
In our context, this Lie algebra coincides with the homology
of the \emph{$\Gm$-rigidification} $\calM^\pl = \calM / [*/\Gm]$,
obtained from the moduli stack $\calM$ by removing scalar automorphisms.

Joyce's enumerative invariants are elements
\begin{equation}
    \label{eq-intro-inv-mplss}
    \inv{\calM^\plss_\alpha (\tau)} \in H_* (\calM^\pl; \bbQ)
\end{equation}
in the Lie algebra $H_* (\calM^\pl; \bbQ)$,
where $\alpha$ is a $K$-theory class,
$\tau$ is the stability condition being used,
and $\calM^\plss_\alpha (\tau) \subset \calM^\pl$ is the open substack
of $\tau$-semistable objects of class~$\alpha$.
The invariants coincide with fundamental classes
$\fund{\calM^\plss_\alpha (\tau)}$
or \emph{virtual fundamental classes}
$\virt{\calM^\plss_\alpha (\tau)}$ 
in the sense of Behrend--Fantechi~\cite{BehrendFantechi1997},
whenever the fundamental classes or the virtual fundamental classes are defined.

Joyce~\cite[\S5]{Joyce2021} also proved \emph{wall-crossing formulae}
for these invariants when varying the stability condition.
These formulae are of the form
\begin{multline}
    \label{eq-intro-wcf-gl}
    \inv{ \calM^\plss_\alpha (\tilde{\tau}) } =
    \sum_{ \leftsubstack[8em]{
        \\[-3ex]
        & n > 0;
        \alpha_1, \dotsc, \alpha_n \in C (\calA) \colon \\[-.5ex]
        & \alpha = \alpha_1 + \cdots + \alpha_n
    } } {} 
    \tilde{U} (\alpha_1, \dotsc, \alpha_n; \tau, \tilde{\tau}) \cdot {} 
    \\[.5ex] 
    \bigl[ \bigl[ \bigl[ 
    \inv{ \calM^\plss_{\alpha_1} (\tau) } \, , 
    \inv{ \calM^\plss_{\alpha_2} (\tau) } \bigr] , \dotsc \bigr] ,
    \inv{ \calM^\plss_{\alpha_n} (\tau) } \bigr] ,
\end{multline}
where $\tau, \tilde{\tau}$ are stability conditions,
$C (\calA)$ is the set of classes $\alpha$
realized by non-zero objects of $\calA$,
$\tilde{U} ({\cdots})$ are certain combinatorial coefficients,
and the Lie brackets come from the Lie algebra structure described above.
The coefficients $\tilde{U} ({\cdots})$
coincide with the ones used in the wall-crossing formulae for 
motivic invariants discussed earlier.
\nopagebreak
\vskip\topsep

The goal of the present paper is to propose a theory
of homological enumerative invariants in the Type B/C/D case,
parallel to the theory described above,
and formulate conjectures on universal structures of these invariants.
We also give a partial proof for our conjectures
in the case of quiver representations.
This will be an analogue of the work of Gross--Joyce--Tanaka~\cite{GrossJoyceTanaka}.

We show that the homology $H_* (\calM^\sd; \bbQ)$
of the moduli stack $\calM^\sd$ of self-dual objects
carries the structure of a \emph{twisted module}
over the vertex algebra $H_* (\calM; \bbQ)$,
in the sense of~\S\ref{sect-va-mod} below.
This is similar to a vertex algebra module,
as in~\cite[\S5.1]{FrenkelBenZvi2004},
but allows more singularities in the vertex operator products.

This also equips $H_* (\calM^\sd; \bbQ)$
with the structure of a \emph{twisted module}
over the Lie algebra $H_* (\calM^\pl; \bbQ)$.
Such an algebraic structure already appeared
in the context of motivic invariants, 
in the author's work \cite[\S7.3]{Bu2023}.

Our main conjecture, Conjecture~\ref{conj-main}, states that,
under certain assumptions, one can construct
\emph{enumerative invariants}, which are homology classes
\begin{equation}
    \label{eq-intro-inv-msdss}
    \inv{\calM^\sdss_\theta (\tau)} \in H_* (\calM^\sd; \bbQ),
\end{equation}
where $\theta$ is a $K$-theory class,
$\tau$ is the stability condition being used,
and $\calM^\sdss_\theta (\tau) \subset \calM^\sd$ is the open substack
of $\tau$-semistable self-dual objects of class~$\theta$.
The invariants should coincide with the fundamental classes
$\fund{\calM^\sdss_\theta (\tau)}$
or the Behrend--Fantechi virtual fundamental classes
$\virt{\calM^\sdss_\theta (\tau)}$ 
whenever the latter are defined.

In addition, when one varies the stability condition,
we conjecture that these invariants should satisfy
the \emph{wall-crossing formula}
\begin{align*}
    \numberthis
    \label{eq-intro-wcf-sd}
    \hspace{2em} & \hspace{-2em}
    \inv{ \calM^\sdss_\theta (\tilde{\tau}) } = {}
    \\[.5ex] &
    \sum_{ \leftsubstack[8em]{
        \\[-3ex]
        & n \geq 0; \, m_1, \dotsc, m_n > 0; \\[-.5ex]
        & \alpha_{1,1}, \dotsc, \alpha_{1,m_1}; \dotsc;
        \alpha_{n,1}, \dotsc, \alpha_{n,m_n} \in C (\calA); \,
        \rho \in \smash{C^\sd (\calA)} \colon \\[-.5ex]
        & \theta = (\bar{\alpha}_{1,1} + \cdots + \bar{\alpha}_{1,m_1})
        + \cdots + (\bar{\alpha}_{n,1} + \cdots + \bar{\alpha}_{n,m_n}) + \rho
    } } {} 
    \tilde{U}^\sd (\alpha_{1,1}, \dotsc, \alpha_{1,m_1}; \dotsc;
    \alpha_{n,1}, \dotsc, \alpha_{n,m_n}; \tau, \tilde{\tau}) \cdot {} 
    \\[.5ex] &
    \bigl[ \bigl[ \bigl[ 
    \inv{ \calM^\plss_{\alpha_{1,1}} (\tau) } \, , 
    \inv{ \calM^\plss_{\alpha_{1,2}} (\tau) } \bigr] , \dotsc \bigr] ,
    \inv{ \calM^\plss_{\alpha_{1,m_1}} (\tau) } \bigr] \heart \cdots \heart {}
    \\ &
    \bigl[ \bigl[ \bigl[ 
    \inv{ \calM^\plss_{\alpha_{n,1}} (\tau) } \, , 
    \inv{ \calM^\plss_{\alpha_{n,2}} (\tau) } \bigr] , \dotsc \bigr] ,
    \inv{ \calM^\plss_{\alpha_{n,m_n}} (\tau) } \bigr] \heart 
    \inv{ \calM^\sdss_{\rho} (\tau) } \ ,
\end{align*}
where $\tau, \tilde{\tau}$ are stability conditions;
$C (\calA)$ is as in~\eqref{eq-intro-wcf-gl};
$C^\sd (\calA)$ is the set of self-dual $K$-theory classes
realized by self-dual objects, not necessarily non-zero;
for $\alpha \in C (\calA)$,
$\bar{\alpha}$ roughly means $\alpha + \alpha^\vee$;
$\tilde{U}^\sd ({\cdots})$ are certain combinatorial coefficients;
the Lie brackets come from the Lie algebra structure on
$H_* (\calM^\pl; \bbQ)$ discussed above;
and $\heart$ is the action of the Lie algebra $H_* (\calM^\pl; \bbQ)$
on the twisted module $H_* (\calM^\sd; \bbQ)$.
Note that the coefficients $\tilde{U}^\sd ({\cdots})$
coincide with the ones used in the wall-crossing formulae for 
motivic invariants, in \cite[Theorem~8.16]{Bu2023}.

We also provide a partial proof for our main conjecture
in the special case of counting self-dual representations
of self-dual quivers, as Theorem~\ref{thm-quiver-main} below.
Following an approach parallel to that of
Gross--Joyce--Tanaka~\cite{GrossJoyceTanaka},
we construct the invariants~\eqref{eq-intro-inv-msdss}
for self-dual quivers, and show that they satisfy the
wall-crossing formula~\eqref{eq-intro-wcf-sd}.
We also show that these invariants agree with fundamental classes
$\fund{\calM^\sdss_\theta (\tau)}$ when $\theta$ is \emph{tame},
in the sense of \S\ref{sect-quiv-inv} below.
We hope to generalize the last statement to all classes~$\theta$
such that the fundamental class is defined.
\nopagebreak
\vskip\topsep

This paper is organized as follows.
In~\S\ref{sect-background}, we provide background material
on formal power series, algebraic topology, and algebraic geometry.
In~\S\ref{sect-cat}, we discuss \emph{self-dual categories}
and their moduli stacks,
largely following ideas from the previous paper~\cite{Bu2023}.
We also discuss stability conditions on self-dual categories,
and we discuss a few general properties of stable and semistable self-dual objects.
In~\S\ref{sect-vertex},
we provide background on \emph{vertex algebras},
and introduce a notion of \emph{twisted modules} for vertex algebras.
In~\S\ref{sect-vertex-moduli},
we study algebraic structures
that arise from moduli stacks of objects in 
$\bbC$-linear categories,
and moduli of self-dual objects in self-dual $\bbC$-linear categories.
This includes Joyce's vertex algebra construction,
which we will discuss in \S\ref{sect-cons-va},
and our construction of twisted modules in \S\ref{sect-cons-tw-mod}.

In~\S\ref{sect-inv},
we formulate our main conjecture
on the universal structure of homological enumerative invariants
counting self-dual objects in self-dual $\bbC$-linear categories.
In~\S\ref{sect-quiver},
we study one of our main examples,
the moduli of self-dual representations of self-dual quivers.
We state Theorem~\ref{thm-quiver-main},
which partially verifies our main conjecture.
We also obtain explicit formulae for the vertex algebra
and twisted module structures obtained in this case.
In~\S\ref{sect-principal},
we study another main example,
that is the moduli of principal orthogonal or symplectic bundles
on a smooth projective variety.
We have not yet been able to obtain any results on enumerative invariants in this case,
but the constructions of algebraic structures, including vertex algebras and twisted modules,
do apply in this case, and we spell out these constructions in detail.
Finally, in \S\ref{sect-proofs-va} and \S\ref{sect-proof-quiver},
we prove the main results in \S\ref{sect-vertex-moduli} and \S\ref{sect-quiver}.

\subsection*{Acknowledgements}

The author is grateful to his supervisor,
Dominic Joyce, for his support and guidance throughout this project.
The author would also like to thank
Arkadij Bojko,
Alexei Latyntsev,
and Si Li for helpful discussions.

This work was done during the author's PhD programme supported by the Engineering
and Physical Sciences Research Council.

\subsection*{Conventions}

\begin{itemize}
    \item 
        A \emph{topological space} is always assumed to be paracompact and Hausdorff,
        and have the homotopy type of a CW complex.
        We consider topological spaces up to homotopy equivalence, unless otherwise stated.

    \item
        An \emph{$\infty$-category}, or an \emph{$(\infty, 1)$-category},
        means a quasi-category in the sense of Lurie~\cite{LurieHTT},
        considered up to categorical equivalence.
        To circumvent size issues, we assume the existence of
        strongly inaccessible cardinals $\kappa < \kappa'$,
        and all our $\infty$-categories are $\kappa$-small, unless otherwise stated.
\end{itemize}

\clearpage
\section{Background material}

\label{sect-background}

\subsection{Formal power series}
\label{sect-pow-ser}

\begin{definition}
    Let $R$ be an integral domain, and let $n \geq 0$ be an integer. Let
    \[
        R \llbr z_1, \dotsc, z_n \rrbr
    \]
    be the ring of formal power series in $n$ variables, with coefficients in $R$.
    Its elements are of the form
    \begin{equation}
        a (z_1, \dotsc, z_n) = \sum_{k_1, \dotsc, k_n \geq 0}
        a_{k_1, \dotsc, k_n} z\sss{1}{k_1} \cdots z\sss{n}{k_n},
    \end{equation}
    with $a_{k_1, \dotsc, k_n} \in R$ for all $k_1, \dotsc, k_n \geq 0$,
    and there can be infinitely many non-zero terms.
    One can check, by induction on $n$, that $R \llbr z_1, \dotsc, z_n \rrbr$
    is again an integral domain,
    and that an element $a (z_1, \dotsc, z_n)$ is invertible
    precisely when $a_{0, \dotsc, 0} \in R$ is invertible.
    
    Now, let $\bbK$ be a field. Define
    \[
        \bbK \llparen z_1, \dotsc, z_n \rrparen
    \]
    to be the fraction field of $\bbK \llbr z_1, \dotsc, z_n \rrbr$.

    For a $\bbK$-vector space $V$, define a $\bbK$-vector space
    \begin{equation}
        V \llbr z_1, \dotsc, z_n \rrbr
        = V \underset{\bbK}{\otimeshat} \bbK \llbr z_1, \dotsc, z_n \rrbr,
    \end{equation}
    whose elements are of the form
    \begin{equation}
        v (z_1, \dotsc, z_n) = \sum_{k_1, \dotsc, k_n \geq 0}
        v_{k_1, \dotsc, k_n} z\sss{1}{k_1} \cdots z\sss{n}{k_n},
    \end{equation}
    with $v_{k_1, \dotsc, k_n} \in V$ for all $k_1, \dotsc, k_n \geq 0$,
    and there can be infinitely many non-zero terms.
    One can naturally regard $V \llbr z_1, \dotsc, z_n \rrbr$
    as a $\bbK \llbr z_1, \dotsc, z_n \rrbr$-module.
    
    Define a $\bbK$-vector space
    \begin{equation}
        V \llparen z_1, \dotsc, z_n \rrparen
        = V \llbr z_1, \dotsc, z_n \rrbr
        \underset{\bbK \llbr z_1, \dotsc, z_n \rrbr}{\otimes}
        \bbK \llparen z_1, \dotsc, z_n \rrparen.
    \end{equation}
    Its elements are of the form
    \[
        \frac{v (z_1, \dotsc, z_n)}{a (z_1, \dotsc, z_n)} \, ,
    \]
    where $v (z_1, \dotsc, z_n) \in V \llbr z_1, \dotsc, z_n \rrbr$
    and $a (z_1, \dotsc, z_n) \in \bbK \llbr z_1, \dotsc, z_n \rrbr \setminus \{ 0 \}$.
\end{definition}

\begin{definition}
    \label{def-iota}
    Let $\bbK$ be a field, and let $n \geq 0$ be an integer.
    We have an inclusion map
    \begin{equation}
        \iota_{z_1, \> \dotsc \> , \> z_n} \colon
        \bbK \llparen z_1, \dotsc, z_n \rrparen \longhookrightarrow
        \bbK \llparen z_1 \rrparen \cdots \llparen z_n \rrparen,
    \end{equation}
    induced by the inclusion
    \[
        \bbK \llbr z_1, \dotsc, z_n \rrbr \simeq
        \bbK \llbr z_1 \rrbr \cdots \llbr z_n \rrbr \longhookrightarrow
        \bbK \llparen z_1 \rrparen \cdots \llparen z_n \rrparen,
    \]
    as one can verify that non-zero elements on the left-hand side
    are sent to invertible elements on the right-hand side.
    
    More generally, we have inclusion maps of the form
    \begin{multline}
        \iota_{\{ z_{1,1}, \> \dotsc \> , \> z_{1, m_1} \}, \> \dotsc \> ,
            \{ z_{n,1}, \> \dotsc \> , \> z_{n, m_n} \}} \colon
        \bbK \llparen z_{1,1}, \dotsc, z_{1, m_1}, \dotsc,
            z_{n,1}, \dotsc, z_{n, m_n} \rrparen 
        \\ \longhookrightarrow
        \bbK \llparen z_{1,1}, \dotsc, z_{1, m_1} \rrparen \cdots
        \llparen z_{n,1}, \dotsc, z_{n, m_n} \rrparen,
    \end{multline}
    where $m_1, \dotsc, m_n \geq 0$ are integers. 
    This is induced by the inclusion
    \[
        \bbK \llbr z_{1,1}, \dotsc, z_{1, m_1}, \dotsc,
            z_{n,1}, \dotsc, z_{n, m_n} \rrbr \longhookrightarrow
        \bbK \llparen z_{1,1}, \dotsc, z_{1, m_1} \rrparen \cdots
            \llparen z_{n,1}, \dotsc, z_{n, m_n} \rrparen.
    \]

    Furthermore, for a $\bbK$-vector space $V$, we have
    inclusion maps
    \begin{multline}
        \iota_{\{ z_{1,1}, \> \dotsc \> , \> z_{1, m_1} \}, \> \dotsc \> ,
            \{ z_{n,1}, \> \dotsc \> , \> z_{n, m_n} \}} \colon
        V \llparen z_{1,1}, \dotsc, z_{1, m_1}, \dotsc,
            z_{n,1}, \dotsc, z_{n, m_n} \rrparen 
        \\ \longhookrightarrow
        V \llparen z_{1,1}, \dotsc, z_{1, m_1} \rrparen \cdots
        \llparen z_{n,1}, \dotsc, z_{n, m_n} \rrparen.
    \end{multline}
    This is induced by the inclusions
    \begin{align*}
        V \llbr z_{1,1}, \dotsc, z_{1, m_1}, \dotsc,
            z_{n,1}, \dotsc, z_{n, m_n} \rrbr 
        & \longhookrightarrow
        V \llparen z_{1,1}, \dotsc, z_{1, m_1} \rrparen \cdots
            \llparen z_{n,1}, \dotsc, z_{n, m_n} \rrparen, \\
        \bbK \llparen z_{1,1}, \dotsc, z_{1, m_1}, \dotsc,
            z_{n,1}, \dotsc, z_{n, m_n} \rrparen 
        & \longhookrightarrow
        V \llparen z_{1,1}, \dotsc, z_{1, m_1} \rrparen \cdots
            \llparen z_{n,1}, \dotsc, z_{n, m_n} \rrparen,
    \end{align*}
    and the universal property of the tensor product.
    
    In applications,
    we sometimes restrict these inclusion maps to a smaller subspace.
    For example, in the context of vertex algebras,
    one often encounters the map
    \begin{equation*}
        \iota_{w, \> z} \colon
        \bbK \llbr z, w \rrbr [z^{-1}, w^{-1}, (z-w)^{-1}]
        \longrightarrow
        \bbK \llparen w \rrparen \llparen z \rrparen.
    \end{equation*}
    Explicitly, the effect of this map is that it expands $(z - w)^n$ as
    \begin{equation}
        \iota_{w, \> z} \bigl( (z - w)^n \bigr) = \begin{cases}
            (z - w)^n, & n \geq 0, \\
            \displaystyle
            (-1)^n \cdot \sum_{i \geq 0} {}
            \binom{i-n-1}{i} \,
            \frac{z^i}{w^{i-n}} \, , & n < 0.
        \end{cases}
    \end{equation}
\end{definition}

\subsection{Background on algebraic topology}

We review some basic concepts in algebraic topology,
fixing some terminology that will be used later.

In the following, a \emph{topological space}
is always assumed to be paracompact and Hausdorff,
and have the homotopy type of a CW complex.

For a topological space $X$ and a commutative ring $R$, write
\[
    H_* (X; R) = \bigoplus_{n \geq 0} H_n (X; R), \qquad 
    H^* (X; R) = \bigoplus_{n \geq 0} H^n (X; R)
\]
for the ordinary homology and cohomology of $X$,
as $\bbZ$-graded $R$-modules.

\begin{definition}
    Let $X$ be a topological space.
    The \emph{topological $K$-theory} of $X$ is the commutative ring
    \begin{equation}
        \label{eq-def-ktop}
        K (X) = [X, \BU \times \bbZ],
    \end{equation}
    where the right-hand side denotes
    the set of homotopy classes of continuous maps,
    as in May~\cite[\S24.1]{May1999}.
    The ring structure is inherited from a ring structure on $\BU \times \bbZ$.
    When $X$ is compact, $K (X)$ can be identified with the abelian group
    generated by elements $[E]$ for complex vector bundles $E$ on $X$,
    with relations $[E \oplus F] = [E] + [F]$
    for complex vector bundles $E, F$ on $X$.
    The multiplication in the ring $K (X)$ is given by
    taking the tensor product of vector bundles.

    For classes $E \in K (X)$, one can define
    \emph{Chern classes} $c_i (E) \in H^{2i} (X; \bbZ)$,
    as well as the \emph{Chern character} $\ch_i (E) \in H^{2i} (X; \bbQ)$.
    See, for example, \cite[Chapters~23--24]{May1999}.

    Some operations of vector bundles induce operations on topological $K$-theory.
    For example, we have maps
    \[
        (-)^\vee, \ \Sym^n, \ {\wedge}^n \colon K (X)
        \longrightarrow K (X),
    \]
    where $n \in \bbZ_{\geq 0}$\,.
    These maps send a vector bundle to its dual vector bundle,
    its $n$-th symmetric product,
    and its $n$-th exterior product, respectively.
    For a class $[E] - [F] \in K (X)$,
    where $E, F$ are vector bundles over $X$, these operations are given by
    \begin{align}
        ([E] - [F])^\vee & =
        [E^\vee] - [F^\vee], \\
        \Sym^n ([E] - [F]) & =
        \sum_{k=0}^{n} {} (-1)^{k} \cdot
        [\Sym^{n-k} (E) \otimes {\wedge}^{k} (F)], \\
        {\wedge}^n ([E] - [F]) & =
        \sum_{k=0}^{n} {} (-1)^{k} \cdot
        [{\wedge}^{n-k} (E) \otimes \Sym^{k} (F)].
    \end{align}
    The map $(-)^\vee$ is a ring isomorphism,
    but $\Sym^n$ and ${\wedge}^n$
    do not respect addition or multiplication in general.
    See also Lemma~\ref{lem-ch-sym} below.

    There are also the \emph{Adams operations}
    \[
        \psi^n \colon K (X) \longrightarrow K (X),
    \]
    as in~\cite[\S5]{Adams1962},
    where $n \in \bbZ$, which are ring homomorphisms
    uniquely determined by the property that
    $\psi^n ([L]) = [L^{\otimes n}]$
    for any line bundle $L$ on $X$.
\end{definition}

\begin{definition}
    \label{def-bu1-wt}
    Let $X$ be a topological space, acted on by the group $\BU (1)$.
    We say that a class $E \in K (X)$
    has \emph{$\upU (1)$-weight $k$}, where $k \in \bbZ$,
    if we have
    \begin{equation}
        a^* (E) = L^{\otimes k} \boxtimes E
    \end{equation}
    in $K (\BU (1) \times X)$,
    where $a \colon \BU (1) \times X \to X$ is the action map,
    and $L \in K (\BU (1))$ is the class of the universal complex line bundle.

    Similarly, if $X$ is acted on by $\BU (1)^n$,
    then one can speak of classes in $K (X)$
    of $\upU (1)$-weight $(k_1, \dotsc, k_n)$,
    where $k_1, \dotsc, k_n \in \bbZ$.
    These are classes $E \in K (X)$ with
    \begin{equation}
        a^* (E) = 
        L^{\otimes k_1} \boxtimes \cdots \boxtimes L^{\otimes k_n}
        \boxtimes E,
    \end{equation}
    where $a \colon \BU (1)^n \times X \to X$ is the action map.
\end{definition}

Next, we define the Chern series of a $K$-theory class,
and we state some of its properties that will become useful later.

\begin{definition}
    Let $X$ be a topological space, and let $E \in K (X)$. Write
    \begin{equation}
        c_z (E) = \sum_{i \geq 0} z^i c_i (E) 
        \quad \in H^* (X; \bbZ) \llbr z \rrbr .
    \end{equation}
    We have the relations
    \begin{align}
        c_z (E + F) & = c_z (E) \, c_z (F), \\
        c_z (E^\vee) & = c_{-z} (E)
    \end{align}
    for $E, F \in K (X)$.
\end{definition}

\begin{lemma}
    \label{lem-ch-to-c}
    Let $X$ be a topological space, and let $E \in K (X)$. Then
    \begin{equation}
        \label{eq-ch-to-c}
        c_z (E) = \exp \biggl(
            \sum_{i > 0} {}
            (-1)^{i-1} (i-1)! \, z^i \ch_i (E)
        \biggr)
    \end{equation}
    as formal power series in $z$.
\end{lemma}

\begin{proof}
    Note that~\eqref{eq-ch-to-c} holds when $E = 0$,
    and both sides of~\eqref{eq-ch-to-c} are multiplicative in~$E$.
    Therefore, if \eqref{eq-ch-to-c} holds for~$E$,
    then it holds for~$-E$.
    This means that we only need to prove~\eqref{eq-ch-to-c}
    in the case when $E$ is the class of a vector bundle,
    as every element in~$K (X)$ is the difference of two vector bundles.

    Using the splitting principle, we may further assume that
    $E = [L]$ is the class of a line bundle $L \to X$.
    Setting $c = c_1 (L)$, \eqref{eq-ch-to-c} becomes
    \begin{equation}
        1 + z c =
        \exp \circ \log (1 + z c),
    \end{equation}
    which is true.
\end{proof}

\begin{lemma}
    \label{lem-ch-sym}
    Let $X$ be a topological space, and let $E \in K (X)$. Then 
    \begin{align}
        \label{eq-sym-to-adams}
        \Sym^2 (E) & =
        \frac{1}{2} (E^2 + \psi^2 (E)), \\
        \label{eq-wedge-to-adams}
        {\wedge}^2 (E) & =
        \frac{1}{2} (E^2 - \psi^2 (E))
    \end{align}
    in $K (X),$
    where $\psi^n$ is the $n$-th Adams operation.
    In particular,
    \begin{align}
        \label{eq-ch-sym-2}
        \ch (\Sym^2 (E)) & =
        \frac{1}{2} \bigl( \ch (E)^2 + \ch_{(2)} (E) \bigr), \\
        \label{eq-ch-wedge-2}
        \ch ({\wedge}^2 (E)) & =
        \frac{1}{2} \bigl( \ch (E)^2 - \ch_{(2)} (E) \bigr),
    \end{align}
    where
    \[
        \ch_{(2)} (E) = \sum_{k = 0}^\infty 2^k \ch_k (E).
    \]
\end{lemma}

\begin{proof}
    Using the splitting principle,
    we see that \eqref{eq-sym-to-adams}--\eqref{eq-wedge-to-adams}
    hold when $E$ is the class of a vector bundle,
    and the general case follows from an elementary computation,
    writing an arbitrary element of $K (X)$
    as the difference of two vector bundles.
    The splitting principle gives
    \begin{equation}
        \ch (\psi^2 (E)) = \ch_{(2)} (E),
    \end{equation}
    so that \eqref{eq-ch-sym-2}--\eqref{eq-ch-wedge-2} follow
    from \eqref{eq-sym-to-adams}--\eqref{eq-wedge-to-adams}.
\end{proof}

\subsection{Derived algebraic geometry}

We briefly introduce basic notions in derived algebraic geometry,
following Toën--Vezzosi~\cite{ToenVezzosi2008}
and Lurie~\cite{Lurie2004}.
This is mainly to clarify definitions while fixing some terminology,
and we do not attempt to give a self-contained account of the subject.

\begin{definition}
    Let $\cat{Aff}_{\bbC}$ be the category of affine $\bbC$-schemes,
    and let $\inftyGpd$ be the $\infty$-category of small $\infty$-groupoids.

    An \emph{$\infty$-stack} over $\bbC$ is an $\infty$-functor
    \[
        \frX \colon \cat{Aff}_{\bbC}^\op
        \longrightarrow \inftyGpd,
    \]
    satisfying descent for the étale topology on $\cat{Aff}_{\bbC}$\,.

    Define the $\infty$-category $\cat{dAff}_{\bbC}$
    of \emph{derived affine $\bbC$-schemes}
    to be the opposite category of
    the $\infty$-category of simplicial $\bbC$-algebras.
    As in~\cite[Definition~2.2.2.12]{ToenVezzosi2008},
    one can define the \emph{étale topology} on $\cat{dAff}_{\bbC}$\,.

    A \emph{derived stack} over $\bbC$ is an $\infty$-functor
    \[
        \frX \colon \cat{dAff}_{\bbC}^\op
        \longrightarrow \inftyGpd,
    \]
    satisfying descent for the étale topology on $\cat{dAff}_{\bbC}$\,.

    The \emph{classical truncation} of a derived stack $\frX$,
    denoted by $\frX_\cl$\,,
    is the $\infty$-stack obtained as
    the restriction of $\frX$ to the subcategory $\cat{Aff}_{\bbC}^\op$.
    In this case, $\frX$ is also called
    a \emph{derived enhancement} of $\frX_\cl$\,.
\end{definition}

An important piece of information of a derived stack
is encoded in its \emph{cotangent complex}.

\begin{definition}
    \label{def-cot-cx}
    Let $\frX$ be a derived stack over $\bbC$.
    As in Toën--Vezzosi~\cite[\S1.4]{ToenVezzosi2008}
    or Lurie~\cite[\S3.2]{Lurie2004},
    when $\frX$ is \emph{geometric}, it admits a \emph{cotangent complex}
    \[
        \bbL_{\frX} \in \cat{QCoh} (\frX),
    \]
    where $\cat{QCoh} (\frX)$ denotes the stable $\infty$-category of
    quasi-coherent complexes on $\frX$.

    If $\frX$ is locally finitely presented,
    then the cotangent complex $\bbL_{\frX}$ is a perfect complex~%
    \cite[Proposition~3.2.14]{Lurie2004}.
    The \emph{tangent complex} of $\frX$ is defined
    as the dual of the cotangent complex,
    $\bbT_{\frX} = \bbL \sss{\frX}{\vee}$.

    We say that $\frX$ is a \emph{smooth}
    (resp.~\emph{quasi-smooth}) \emph{derived Artin stack},
    if the classical truncation $\frX_\cl$ of $\frX$
    is an Artin $1$-stack, i.e.~an algebraic stack,
    and $\bbL_{\frX}$ is a perfect complex with
    amplitude~$\geq 0$ (resp.~$\geq -1$).
    In fact, $\bbL_{\frX}$ will have amplitude in
    $[0, 1]$, resp.~$[-1, 1]$.

    Similarly, we say that $\frX$ is a \emph{smooth}
    (resp.~\emph{quasi-smooth}) \emph{derived Deligne--Mumford stack},
    if its classical truncation is a Deligne--Mumford $1$-stack,
    and its cotangent complex is perfect with
    amplitude~$\geq 0$ (resp.~$\geq -1$).
    In fact, $\bbL_{\frX}$ will have amplitude on
    $\{ 0 \}$, resp.~in $[-1, 0]$.
\end{definition}

We also discuss \emph{virtual fundamental classes},
introduced by Behrend--Fantechi~\cite[\S5]{BehrendFantechi1997},
which can reinterpreted using derived algebraic geometry,
as is often done in the literature. See, for example,
\cite[Remark~2.14]{Joyce2021}.

\begin{definition}
    \label{def-virt}
    Let $\frX$ be a quasi-smooth derived Deligne--Mumford stack over $\bbC$,
    and let $\calX = \frX_\cl$ be its classical truncation,
    which is a Deligne--Mumford $1$-stack.
    Denote by $i_\cl \colon \calX \to \frX$ the natural inclusion.

    The map $i_\cl$ induces a morphism $i_\cl^* (\bbL_{\frX}) \to \bbL_{\calX}$
    of complexes on $\calX$, 
    and this is a \emph{perfect obstruction theory}
    in the sense of Behrend--Fantechi~\cite[\S5]{BehrendFantechi1997}.
    Then \cite[\S5]{BehrendFantechi1997} produces
    a \emph{virtual fundamental class}
    \[
        \virt{\calX} \in \CH_{\vdim \calX} (\calX; \bbQ),
    \]
    where $\vdim \calX = \rank \bbL_{\frX}$
    is the \emph{virtual dimension} of $\calX$,
    and is assumed to be constant on $\calX$.
    This virtual class depends on the choice of the derived enhancement $\frX$.
    Note that Behrend--Fantechi assumed that the perfect obstruction theory
    has a global resolution,
    but it was later shown by Kresch~\cite[\S5.2]{Kresch1999}
    that this was unnecessary.

    Using the \emph{cycle map}
    $\CH_* (\calX; \bbQ) \to H \sss{2*}{\mathrm{BM}} (\calX; \bbQ)$
    into the \emph{Borel--Moore homology} of~$\calX$,
    as in Khan~\cite[\S1.2]{Khan2022},
    one obtains a virtual fundamental class
    \[
        \virt{\calX} \in H \sss{2 \vdim \calX}{\mathrm{BM}} (\calX; \bbQ).
    \]
    When $\calX$ is proper, one can identify
    the Borel--Moore homology $H \sss{*}{\mathrm{BM}} (\calX; \bbQ)$
    with the usual homology $H_* (\calX; \bbQ)$
    defined in \S\ref{sect-top-real} below.
\end{definition}

\subsection{Topological realization of stacks}
\label{sect-top-real}

For a stack or an $\infty$-stack $\frX$ over $\bbC$,
one can define a topological space $|\frX|$,
called the \emph{topological realization} of $\frX$,
defined by Simpson~\cite{Simpson1996},
also discussed in Morel--Voevodsky~\cite[\S3]{MorelVoevodsky1999}
and Blanc~\cite{Blanc2016}.

Using the realization functor, one can define the \emph{homology} and \emph{cohomology}
\[
    H_* (\frX; R) = H_* (|\frX|; R), \qquad
    H^* (\frX; R) = H^* (|\frX|; R)
\]
of a stack $\frX$,
where $R$ is a commutative ring.
In the context of enumerative geometry,
one considers the cases when $\frX$ is the moduli stack of objects
in an abelian category or a dg-category,
as discussed by Joyce~\cite{JoyceHall} and Gross~\cite{Gross2019}.

For our purposes, in addition to the above,
we will also be interested in the case when $\frX$
is the moduli stack of self-dual objects in a self-dual category
(see \S\ref{sect-cat} below).

Following Blanc~\cite[\S3]{Blanc2016},
we define the topological realization functor as follows.

\begin{definition}
    Let $|{-}| \colon \cat{Aff}_{\bbC} \to \cat{S}$
    denote the topological realization of affine $\bbC$-schemes,
    given by the underlying topological space of the analytification.
    Here, $\cat{S}$ is the $\infty$-category of spaces,
    i.e.~$\infty$-groupoids.

    Let $\cat{hSt}_{\bbC}$ be the \Aone-homotopy category
    of $\infty$-stacks over~$\bbC$,
    which is an $\infty$-category obtained from
    the $\infty$-category of $\infty$-stacks over~$\bbC$
    by localizing along the projections $X \times \bbA^1 \to X$ for all $X$.
    For a precise formulation of this, see, for example,
    Blanc~\cite[\S3.2]{Blanc2016}.
    There is an inclusion functor
    $i \colon \cat{Aff}_{\bbC} \to \cat{hSt}_{\bbC}$\,.

    Define the \emph{topological realization} functor
    \[
        |{-}| \colon \cat{hSt}_{\bbC} \longrightarrow \cat{S}
    \]
    to be the left Kan extension of the previously defined
    functor $|{-}|$ for affine schemes along $i$, as in the diagram
    \[ \begin{tikzcd}
        \cat{Aff}_{\bbC} \ar[d, "i"'] \ar[dr, "|{-}|"]
        \\ \cat{hSt}_{\bbC} \ar[r, "|{-}|"'] 
        & \cat{S} \rlap{ .}
    \end{tikzcd} \]
    The diagram commutes up to a natural isomorphism.

    In particular, by this definition,
    the topological realization is \emph{\Aone-homotopy invariant},
    sending \Aone-homotopy equivalences to homotopy equivalences.

    For a derived stack $\frX$ over $\bbC$,
    define the \emph{topological realization} of $\frX$
    to be the topological realization of its classical truncation,
    $|\frX| = |\frX_{\cl}|$.
\end{definition}

\begin{example}
    \label{eg-quot-st-top}
    For a quotient stack $[X/G]$,
    where $X$ is a $\bbC$-scheme and $G$ is an algebraic group over $\bbC$,
    we have 
    \begin{equation}
        |[X/G]| \simeq |X| / |G|,
    \end{equation}
    where we take the homotopy quotient of a topological space
    by a topological group.
    This is because $|{-}|$ commutes with homotopy colimits.
\end{example}

We also relate the algebraic $K$-theory of a $\bbC$-scheme
with the topological $K$-theory of its topological realization.

\begin{definition}
    \label{def-perf}
    Let $\Perf$ be the classifying stack of perfect complexes over~$\bbC$,
    which is a derived stack over~$\bbC$,
    as in Toën--Vaquié~\cite[p.\,440]{ToenVaquie2007}
    or Toën--Vezzosi~\cite[Definition~1.3.7.5]{ToenVezzosi2008}.
    As in Blanc~\cite[Theorems~4.5 and~4.21]{Blanc2016},
    we have a homotopy equivalence
    \begin{equation}
        \label{eq-perf-bu}
        |\Perf| \simeq \upB \upU \times \bbZ.
    \end{equation}
    This gives rise to a natural map
    \begin{equation}
        \label{eq-map-kalg-ktop}
        \Perf (\frX) \longrightarrow
        \cat{Map}_{\cat{Top}} (|\frX|, \BU \times \bbZ)
    \end{equation}
    for any derived $\bbC$-stack $\frX$,
    where the left-hand side denotes the mapping space from $\frX$ to $\Perf$
    in the $\infty$-category of derived stacks over $\bbC$,
    and the right-hand side denotes the mapping space
    in the $\infty$-category of topological spaces.
    In particular, we obtain a map of abelian groups
    \begin{equation}
        K_0 (\Perf (\frX)) \longrightarrow K (|\frX|),
    \end{equation}
    where the left-hand side is the Grothendieck group
    of the stable $\infty$-category of perfect complexes on $\frX$,
    and the right-hand side is the topological $K$-theory of $|\frX|$.
    This assigns to each perfect complex on $\frX$
    its corresponding topological $K$-theory class.

    This map also factors through the \emph{semi-topological $K$-theory},
    defined in Blanc~\cite[\S4]{Blanc2016}
    via the topological realization $|\cat{Perf} (\frX)|$.
\end{definition}

\section{Self-dual categories and moduli stacks}

\label{sect-cat}

\subsection{Self-dual categories}

We recall the definition of \emph{self-dual categories}
from the author~\cite[\S3]{Bu2023}.

\begin{definition}
    \label{def-sd-cat}
    A \emph{self-dual category} is a triple $(\calA, D, \eta)$, where

    \begin{itemize}
        \item
            $\calA$ is a category.

        \item
            $D \colon \calA \simto \calA^\op$
            is an equivalence,
            called the \emph{dual functor}.
            We often denote this functor by $(-)^\vee$.

        \item
            $\eta \colon D^\op \circ D \simTo \id_{\calA}$
            is a natural isomorphism,
            such that for any $E \in \calA$, we have
            \begin{equation}
                \eta_{E^\vee} = (\eta_E^\vee)^{-1} \colon
                E^{\vee\vee\vee} \longsimto E^\vee .
            \end{equation}
    \end{itemize}
    By abuse of language, we sometimes say that $(\calA, D)$,
    or simply $\calA$, is a self-dual category.
    
    If, in addition, $\calA$ is equipped with some additional structure,
    such as that of an additive category,
    a $\bbC$-linear category,
    or an exact category, etc.,
    and if such structure is respected by $D$ and $\eta$,
    meaning that $D$ and $\eta$ are $1$- and $2$-morphisms
    in the $2$-category of such categories with additional structure,
    then we say that $(\calA, D, \eta)$
    is, for example, a \emph{self-dual additive category}, 
    or a \emph{self-dual $\bbC$-linear category}, etc.
\end{definition}

\begin{definition}
    \label{def-sd-obj}
    Let $(\calA, D, \eta)$ be a self-dual category.
    A \emph{self-dual object} in $\calA$
    is a pair $(E, \phi)$, where
    \begin{itemize}
        \item 
            $E \in \calA$ is an object.
        \item
            $\phi \colon E \simto E^\vee$ is an isomorphism in $\calA$,
            such that $\phi^\vee = \phi \circ \eta_E$\,:
            \begin{equation} \begin{tikzcd}[row sep={2.2em,between origins}, column sep=3em]
                E \ar[dr, "\phi" {pos=.47},
                    "\textstyle \sim" {pos=.52, anchor=center, rotate=-20, shift={(0, -1ex)}}] \\
                & E^\vee \rlap{ .} \\
                E^{\vee\vee} \ar[uu, "\eta_E", 
                    "\textstyle \sim" {anchor=center, rotate=90, shift={(0, -1ex)}}]
                \ar[ur, "\phi^\vee"' {pos=.42, inner sep=.05em},
                    "\textstyle \sim" {pos=.47, anchor=center, rotate=20, shift={(0, .6ex)}}]
            \end{tikzcd} \end{equation}
    \end{itemize}
    
    A \emph{morphism of self-dual objects} from $(E, \phi)$ to $(F, \psi)$
    is an isomorphism $f \colon E \simto F$ in~$\calA$, such that the diagram
    \begin{equation} \begin{tikzcd}
        E \ar[d, "f"', "\textstyle \sim" {anchor=center, rotate=90, shift={(0, -1ex)}}]
        \ar[r, "\phi", "\textstyle \sim"']
        & E^\vee \\
        F \ar[r, "\psi"', "\textstyle \sim" {shift={(0, -.4ex)}}]
        & F^\vee \ar[u, "f^\vee"', "\textstyle \sim" {anchor=center, rotate=90, shift={(0, .6ex)}}]
    \end{tikzcd} \end{equation}
    commutes.
    
    Let $\calA^\sd$ denote the groupoid
    whose objects are self-dual objects in $\calA$,
    and morphisms are as above.
\end{definition}

\begin{definition}
    Let $(\calA, D, \eta)$ and $(\calA', D', \eta')$
    be self-dual categories.

    A \emph{self-dual functor} from $\calA$ to $\calA'$
    is a pair $(F, \kappa)$, where

    \begin{itemize}
        \item 
            $F \colon \calA \to \calA'$ is a functor.

        \item
            $\kappa \colon F^\op \circ D \simTo D' \circ F$
            is a natural isomorphism,
            such that if we compose the $2$-morphisms
            $\kappa, \kappa^\op, \eta, \eta'$
            as in the diagram
            \begin{equation} \begin{tikzcd}[row sep=2em]
                \calA \ar[r, "D"] \ar[d, "F"]
                \ar[rr, bend left=40, "\id"]
                \ar[rr, phantom, shift left=5.5, "\Uparrow"]
                \ar[rr, phantom, shift left=5.5, "\scriptstyle \mathllap{\eta} \hspace{.75em}"]
                & \calA^\op \ar[r, "D^\op"] \ar[d, "F^\op"]
                & \calA \ar[d, "F"]
                \\ \calA' \ar[r, "D'"']
                \ar[rr, bend right=40, "\id"']
                \ar[ur, phantom, "\Leftarrow" rotate=45]
                \ar[ur, phantom, shift left=2.5, "\scriptstyle \kappa"]
                \ar[rr, phantom, shift right=5, "\Downarrow"]
                \ar[rr, phantom, shift right=4.5, "\scriptstyle \mathllap{\eta'} \hspace{.5em}"]
                & \calA'^\op \ar[r, "D'^\op"']
                \ar[ur, phantom, shift right=1, "\Leftarrow" rotate=45]
                \ar[ur, phantom, shift left=1.5, "\scriptstyle \kappa^\op"]
                & \calA' \rlap{ ,}
            \end{tikzcd} \end{equation}
            then the resulting $2$-morphism is equal to $\id_F$\,.
    \end{itemize}
    If, in addition, $\calA$ and $\calA'$
    are self-dual categories equipped with extra structure,
    as in Definition~\ref{def-sd-cat},
    then we may define a self-dual functors respecting such structure
    as a pair $(F, \kappa)$ as above,
    such that the extra structure is respected by $F$ and $\kappa$.
\end{definition}

\begin{remark}
    \label{rmk-alt-def-sd-cat}
    Alternatively, and more concisely,
    one can define the \emph{$2$-category of self-dual categories}
    as the homotopy fixed points of the $\bbZ_2$-action
    on the $2$-category of categories, given by taking the opposite category.
    Moreover, the \emph{$2$-category of self-dual categories with extra structure},
    where the extra structure can refer to additive structure, etc.,
    can be defined as the $\bbZ_2$-fixed points in
    the $2$-category of categories with such extra structure.

    This approach is perhaps more natural,
    and indeed, we will take this approach when we define
    self-dual $\infty$-categories in \S\ref{sect-sd-infty-cat} below.
    However, we choose to spell out the explicit definition
    in the $1$-categorical case,
    to help the reader with understanding this type of structures.
\end{remark}

We introduce a construction that will become important later on.

\begin{example}
    \label{eg-oplus-sd}
    Let $\calA$ be a self-dual additive category,
    and $\calA^\sd$ its groupoid of self-dual objects.
    For each object $E \in \calA$, 
    define an object $(\bar{E}, \phi) \in \calA^\sd$ by
    \begin{equation}
        \bar{E} = E \oplus E^\vee, \qquad
        \phi = \begin{pmatrix}
            0 & \id_{E^\vee} \\
            \eta \sss{E}{-1} & 0
        \end{pmatrix},
    \end{equation}
    where $\eta_E \colon E^{\vee \vee} \simto E$
    is given by the self-dual structure of $\calA$.
    This defines a functor $\calA^\simeq \to \calA^\sd$,
    where $\calA^\simeq$ is the groupoid
    obtained from $\calA$ by discarding all non-invertible arrows.

    Moreover, there is a functor
    \begin{align*}
        \oplus^\sd \colon \calA^\simeq \times \calA^\sd
        & \longrightarrow \calA^\sd, \\
        (E, F) & \longmapsto \bar{E} \oplus F,
    \end{align*}
    where we omitted the structure maps for the self-dual objects.
    Later on, this will be used to establish
    moduli spaces of $\calA^\sd$ as modules over moduli spaces of $\calA$.
\end{example}

\subsection{Self-dual \texorpdfstring{$\infty$}{∞}-categories}
\label{sect-sd-infty-cat}

We introduce the notion of \emph{self-dual $\infty$-categories}.

To circumvent set-theoretic issues,
we assume the existence of sufficiently many \emph{Grothendieck universes},
as in Lurie~\cite[\S1.2.15]{LurieHTT}.
Let $\inftyCat$ be the $(\infty, 1)$-category of small $(\infty, 1)$-categories,
as in \cite[Chapter~3]{LurieHTT}.
Expanding the universe,
we may choose another $(\infty, 1)$-category $\inftyCAT$
of bigger $(\infty, 1)$-categories, so that
\[
    \inftyCat \in \inftyCAT.
\]

\begin{definition}
    \label{def-sd-infty-cat}
    Let $\bbZ_2$ act on $\inftyCat$ by sending an $\infty$-category $\calC$
    to its opposite category $\calC^\op$.
    The $\infty$-category of \emph{self-dual $\infty$-categories}
    is the homotopy $\bbZ_2$-fixed points of this action,
    \begin{equation}
        \inftyCat^\sd = (\inftyCat)^{\bbZ_2},
    \end{equation}
    as a homotopy limit in $\inftyCAT$.
\end{definition}

More explicitly, a self-dual $\infty$-category consists of
an $\infty$-category $\calC$,
an equivalence $D \colon \calC \simto \calC^\op$,
an equivalence $\eta \colon D^\op \circ D \simTo \id_{\calC}$\,,
and infinitely many higher homotopies satisfying coherence conditions.

\begin{definition}
    \label{def-sd-obj-infty-cat}
    Let $\calC$ be a self-dual $\infty$-category.
    The \emph{$\infty$-groupoid of self-dual objects} in $\calC$
    is the $\bbZ_2$-fixed points
    \begin{equation}
        \calC^\sd = (\calC^\simeq)^{\bbZ_2} ,
    \end{equation}
    as a homotopy limit in $\inftyGpd$,
    where $\calC^\simeq$ is the $\infty$-groupoid
    obtained from $\calC$ by removing all non-invertible arrows,
    and $\bbZ_2$ acts by composing the restriction of
    $D \colon \calC \simto \calC^\op$ to $\calC^\simeq$,
    which is an equivalence $\calC^\simeq \simto (\calC^\op)^\simeq$,
    with the natural identification
    $(\calC^\op)^\simeq \simeq \calC^\simeq$,
    sending each $1$-morphism to its inverse.
\end{definition}

\subsection{Linear categories and moduli stacks}

We define \emph{linear stacks},
which can be seen as moduli stacks of objects in a linear category.
These are similar to \emph{linear algebraic stacks}
in the sense of \cite[\S4]{Bu2023},
but we drop the algebraicity condition here,
as it is not needed for the present work.

Linear stacks are a type of \emph{stacks in categories},
as opposed to a usual stack,
which is a functor taking values in groupoids.
As discussed in \cite{Bu2023},
a stack in categories is,
in some sense, a more natural notion of a moduli stack of objects in a category,
as many operations,
such as taking the direct sum in an additive category,
can be obtained as intrinsic data of the stack,
rather than extra information.

\begin{definition}
    \label{def-lin-st}
    An \emph{additive stack} over $\bbC$,
    as in~\cite[\S4]{Bu2023},
    is a $2$-functor $\calX \colon \cat{Aff}_{\bbC}^\op \to \cat{AddCat}$,
    where $\cat{AddCat}$ is the $2$-category of additive categories,
    such that $\calX$ satisfies descent for the étale topology.
    We have the $2$-category of additive stacks over $\bbC$.

    An additive stack $\calX$ is equipped with canonical morphisms
    \begin{align*}
        0 \colon *
        & \longrightarrow \calX,
        \\
        \oplus \colon \calX \times \calX
        & \longrightarrow \calX,
    \end{align*}
    given by the zero object and the direct sum, respectively.

    A \emph{$\bbC$-linear stack}
    is an additive stack $\calX$ over $\bbC$, together with an action
    \[
        \odot \colon [*/\bbA^1] \times \calX
        \longrightarrow \calX,
    \]
    compatible with the additive structure on $\calX$,
    in the sense of~\cite[Definition~4.13]{Bu2023}.

    A \emph{morphism of $\bbC$-linear stacks}
    is a morphism of additive stacks, compatible with the $\bbA^1$-action.
    Similarly, a \emph{$2$-morphism} between such morphisms
    is a $2$-morphism in the $2$-category of additive stacks
    compatible with the $\bbA^1$-action.
    See~\cite[\S2.4]{Bu2023} for the precise definition
    of these compatibility conditions.

    These define a $2$-category of $\bbC$-linear stacks.
\end{definition}

We also define a notion of \emph{self-dual linear stacks},
which will serve as moduli stacks for self-dual $\bbC$-linear categories.

\begin{definition}
    \label{def-sd-lin-st}
    A \emph{self-dual $\bbC$-linear stack}
    is a triple $(\calX, D, \eta)$, where

    \begin{itemize}
        \item 
            $\calX$ is a $\bbC$-linear stack.

        \item
            $D \colon \calX \simto \calX^\op$
            is an isomorphism of linear stacks,
            where $\calX^\op$ is the opposite linear stack of $\calX$.

        \item
            $\eta \colon D^\op \circ D \simTo \id_{\calX}$
            is a $2$-morphism in the $2$-category of $\bbC$-linear stacks,
            such that $D \circ \eta = \eta^\op \circ D \colon
            D \circ D^\op \circ D \simTo D$.
    \end{itemize}
    Given this data, the \emph{stack of self-dual objects} in $\calX$,
    denoted by $\calX^\sd$,
    is an ordinary stack over $\bbC$, i.e.~a stack in groupoids,
    given by $\calX^\sd (S) = \calX (S)^\sd$
    for any affine $\bbC$-scheme $S$,
    where the right-hand side is as in Definition~\ref{def-sd-obj},
    and the self-dual structure on $\calX (S)$
    is induced by the data $D$ and $\eta$ above.
\end{definition}

\begin{remark}
    \label{rmk-alt-def-sd-lin-st}
    As in Remark~\ref{rmk-alt-def-sd-cat},
    one can define the $2$-category of self-dual $\bbC$-linear stacks
    as the homotopy $\bbZ_2$-fixed points
    in the $2$-category of $\bbC$-linear stacks,
    where the $\bbZ_2$-action is given by
    taking the opposite linear stack.
    Also, for a self-dual $\bbC$-linear stack $\calX$,
    the stack $\calX^\sd$ is the homotopy $\bbZ_2$-fixed points in $\calX$.
    But we choose to spell out the definition explicitly
    in the case of $1$-stacks.
\end{remark}

\subsection{Linear \texorpdfstring{$\infty$}{∞}-categories and moduli stacks}

We define \emph{additive $\infty$-categories},
following Gepner--Groth--Nikolaus~\cite{GepnerGrothNikolaus2015}.
The following characterization
can be found in \cite[Proposition~2-8]{GepnerGrothNikolaus2015}.

\begin{definition}
    An \emph{additive $\infty$-category}
    is an $\infty$-category $\calC$ with finite products and coproducts,
    such that its homotopy category $\cat{Ho} (\calC)$
    is an additive category.

    An \emph{additive functor} between additive $\infty$-categories
    is a functor preserving finite products and coproducts.

    The \emph{$\infty$-category of additive $\infty$-categories}
    $\infty \mathhyphen \cat{AddCat}$
    is the subcategory of $\inftyCat$ spanned by
    additive $\infty$-categories,
    additive functors between additive $\infty$-categories,
    and all higher morphisms between them.
\end{definition}

Note that we have used the fact that for an ordinary category,
being additive is a property instead of an extra structure.
From the definition,
we can see that the same is true for being an additive $\infty$-category.

\begin{definition}
    The \emph{$\infty$-category of self-dual additive $\infty$-categories}
    is the homotopy $\bbZ_2$-fixed points in
    the $\infty$-category of additive $\infty$-categories,
    with the $\bbZ_2$-action given by taking the opposite category.
\end{definition}

Next, we define a notion of \emph{linear $\infty$-stacks},
which will be \emph{stacks in $\infty$-categories}.
As in the $1$-categorical case,
we believe that stacks in $\infty$-categories
are more suitable for moduli problems,
compared to the usual notion of $\infty$-stacks,
which are functors taking values in $\infty$-groupoids.
The reason for this is that by encoding information of non-invertible morphisms,
many constructions for $\infty$-categories become
automatically available as intrinsic data of the stack,
such as taking direct sums inside the original $\infty$-category.

\begin{definition}
    \label{def-lin-infty-st}
    An \emph{additive $\infty$-stack} over $\bbC$ is a functor
    $\calX \colon \cat{Aff}_{\bbC}^\op \to \infty \mathhyphen \cat{AddCat}$,
    satisfying descent for the étale topology.

    An additive $\infty$-stack $\calX$ is equipped with canonical morphisms
    \begin{align*}
        0 \colon *
        & \longrightarrow \calX,
        \\
        \oplus \colon \calX \times \calX
        & \longrightarrow \calX,
    \end{align*}
    given by the zero object and the direct sum, respectively.

    A \emph{$\bbC$-linear $\infty$-stack}
    is an additive $\infty$-stack $\calX$,
    together with an action
    \[
        \odot \colon [*/\bbA^1] \times \calX
        \longrightarrow \calX,
    \]
    compatible with the additive structure on $\calX$.
    More precisely,
    this can be formulated as a morphism of monoid objects
    $[*/\bbA^1] \to \cat{End}^\oplus (\calX)$,
    where $\cat{End}^\oplus (\calX)$ is the $\infty$-stack of
    additive endomorphisms of $\calX$.

    One can form the $\infty$-category of $\bbC$-linear $\infty$-stacks,
    as the full subcategory of the
    $\infty$-category of $[*/\bbA^1]$-equivariant
    stacks in $\infty$-categories over $\bbC$,
    consisting of $\bbC$-linear $\infty$-stacks,
    additive morphisms, and all higher morphisms.

    We also have the notion of a \emph{$\bbC$-linear derived stack},
    where we use $\cat{dAff}_{\bbC}^\op$
    in place of $\cat{Aff}_{\bbC}^\op$
    as the domain of the functor.
\end{definition}

We also define self-dual $\bbC$-linear $\infty$-stacks.

\begin{definition}
    \label{def-sd-lin-infty-st}
    The $\infty$-category of \emph{self-dual $\bbC$-linear $\infty$-stacks}
    is the homotopy $\bbZ_2$-fixed points in the $\infty$-category
    of $\bbC$-linear $\infty$-stacks,
    where $\bbZ_2$ acts by taking the opposite stack.

    For a self-dual $\bbC$-linear $\infty$-stack $\calX$,
    its \emph{stack of self-dual objects}
    is an $\infty$-stack $\calX^\sd$ over $\bbC$,
    taking values in $\infty$-groupoids,
    which assigns to each affine $\bbC$-scheme $S$
    the $\infty$-groupoid
    $\calX^\sd (S) = \calX (S)^\sd$,
    as in Definition~\ref{def-sd-obj-infty-cat},
    where the self-dual structure on $\calX (S)$
    is induced from that of $\calX$.

    Likewise, we also have \emph{self-dual $\bbC$-linear derived stacks}.
    Their stacks of self-dual objects are also derived stacks.
\end{definition}

\begin{example}
    \label{def-operf}
    Let $\Perf$ be the classifying stack of perfect complexes over~$\bbC$,
    as in Definition~\ref{def-perf}.
    It extends to a $\bbC$-linear derived stack $\Perf^+$,
    associating each derived affine $\bbC$-scheme $S$
    with the $\bbC$-linear stable $\infty$-category of perfect complexes on $S$.

    We consider two self-dual structures on $\Perf^+$,
    depending on a parameter $\epsilon = \pm 1$.
    Namely, consider the involution $\Perf^+ \simto (\Perf^+)^\op$
    sending a perfect complex to its dual.
    Then, identify the double dual of a perfect complex with itself
    using the sign $\epsilon$.

    The stack of self-dual objects $\Perf^\sd$
    classifies perfect complexes $E$ equipped with an isomorphism
    $\phi \colon E \simto E^\vee$,
    with $\phi^\vee \simeq \epsilon \phi$.

    Define derived stacks~$\OPerf$ and $\SpPerf$
    to be the stack $\Perf^\sd$ when $\epsilon = +1$ or $-1$, respectively.
    They classify objects which one may call
    \emph{orthogonal} or \emph{symplectic perfect complexes}.
    They should also be, in some sense, representing objects
    for \emph{real algebraic $K$-theory},
    as in Hesselholt--Madsen~\cite{HesselholtMadsen}.
\end{example}

Recall that we have an equivalence $|\Perf| \simeq \upB \upU \times \bbZ$
in~\eqref{eq-perf-bu},
relating the \emph{semi-topological $K$-theory}
of a point to its topological $K$-theory,
We also hope that similar relations in the real case should hold.
We state this as the following conjecture.

\begin{conjecture}
    \label{conj-homology}
    We should have equivalences of spectra
    \begin{equation}
        |\OPerf|_{\mathbb{S}} \simeq \mathit{ko}, \qquad
        |\SpPerf|_{\mathbb{S}} \simeq \mathit{ksp},
    \end{equation}
    where $|{-}|_{\mathbb{S}}$ denotes spectral realization,
    as in Blanc~\textnormal{\cite[\S3.4]{Blanc2016},}
    and $\mathit{ko}$ and $\mathit{ksp}$
    are connective covers of the $K$-theory spectra
    $\mathit{KO}$ and $\mathit{KSp},$ representing
    orthogonal and symplectic topological $K$-theory, respectively.

    In particular, we should have homotopy equivalences
    \begin{equation}
        |\OPerf| \simeq \BO \times \bbZ, \qquad
        |\SpPerf| \simeq \BSp \times \bbZ.
    \end{equation}
\end{conjecture}

\subsection{Stability conditions on self-dual categories}

We discuss \emph{stability conditions} on self-dual exact categories,
as in~\cite[\S3.3]{Bu2023}.
Here, we also focus on describing \emph{stable self-dual objects},
which were not discussed there.
Some results appeared in Young~\cite[\S3.1]{Young2015}
in the special case of self-dual quiver representations.

We will use the notion of \emph{quasi-abelian categories},
as discussed in Bridgeland~\cite[\S4]{Bridgeland2002},
as it interacts well with stability conditions.

We first recall the notion of a \emph{self-dual stability condition}
on a self-dual exact category, as in~\cite[Definition~3.12]{Bu2023}.

\begin{definition}
    \label{def-sd-stab}
    Let $\calA$ be a self-dual $\bbC$-linear exact category.
    Choose a quotient group $K (\calA)$
    of the Grothendieck group of $\calA$,
    compatible with the self-dual structure,
    and write $C (\calA) \subset K (\calA)$ for the classes
    realized by non-zero objects. We require that
    the class $0 \in K (\calA)$ only contains the zero object.

    A \emph{self-dual stability condition} on $\calA$
    consists of the following data.
    \begin{itemize}
        \item 
            A totally ordered set $T$,
            equipped with an order-reversing involution $t \mapsto -t$,
            fixing a unique element $0 \in T$.
        \item
            A map $\tau \colon C (\calA) \to T$,
            satisfying $\tau (\alpha) = -\tau (\alpha^\vee)$
            for all $\alpha \in C (\calA)$,
            such that for any $\alpha, \beta, \gamma \in C (\calA)$
            with $\beta = \alpha + \gamma$, we have either
            \[
                \tau (\alpha) < \tau (\beta) < \tau (\gamma)
                \quad \text{or} \quad
                \tau (\alpha) = \tau (\beta) = \tau (\gamma)
                \quad \text{or} \quad
                \tau (\alpha) > \tau (\beta) > \tau (\gamma).
            \]
    \end{itemize}
\end{definition}

\begin{situation}
    \label{sit-quasi-ab-sd}
    Let $\calA$ be a self-dual $\bbC$-linear quasi-abelian category,
    with a choice of $K (\calA)$, as in Definition~\ref{def-sd-stab}.
    Let $\tau$ be a self-dual stability condition on $\calA$.

    We assume that $\calA$ is \emph{noetherian}, that is,
    for any object $E \in \calA$ and an infinite chain
    $E_0 \subset E_1 \subset \cdots \subset E$ of subobjects of $E$,
    there exists an integer $n$ such that
    the inclusion $E_i \hookrightarrow E_{i+1}$ is an isomorphism when $i > n$.

    As $\calA$ is self-dual, this also implies that $\calA$ is \emph{artinian},
    that is, for an infinite chain
    $E_0 \supset E_1 \supset \cdots$ in $\calA$,
    there exists an integer $n$ such that
    the inclusion $E_{i+1} \hookrightarrow E_i$ is an isomorphism when $i > n$.
\end{situation}

\begin{remark}
    One could also consider the $\infty$-categorical setting,
    where we have a self-dual $\bbC$-linear stable $\infty$-category,
    and a Bridgeland stability condition \cite{Bridgeland2002} on it,
    which is self-dual.
    However, by choosing a slice in the stable $\infty$-category,
    we obtain a quasi-abelian category, as in \cite[\S4]{Bridgeland2002},
    which contains enough data to perform the constructions in this section.
    In other words, considering the full stable $\infty$-category
    does not give us extra generality.
\end{remark}

Recall the notion of \emph{isotropic subobjects}
of a self-dual object, as in~\cite[Definition~3.7]{Bu2023}.
Namely, for a self-dual object $(E, \phi) \in \calA^\sd$,
a subobject $F \subset E$ is \emph{isotropic}
if the composition $F \hookrightarrow E \simto E^\vee \to F^\vee$
is zero, where the last map is the dual of the inclusion $F \hookrightarrow E$.

\begin{definition}
    In Situation~\ref{sit-quasi-ab-sd},
    a self-dual object $(E, \phi) \in \calA^\sd$ is said to be
    \begin{itemize}
        \item 
            \emph{$\tau$-semistable},
            if for any isotropic subobject $F \subset E$, we have $\tau (F) \leq 0$.
        \item
            \emph{$\tau$-stable},
            if for any non-zero isotropic subobject $F \subset E$,
            we have $\tau (F) < 0$.
    \end{itemize}
\end{definition}

Note that using our definition,
the object $(0, \id_0) \in \calA^\sd$ is stable as a self-dual object,
but $0 \in \calA$ is not stable as an ordinary object.
This is a shadow of a strange behaviour of
stable self-dual objects, in that they are in general decomposable,
unlike ordinary stable objects.
The precise statement will be given in Theorem~\ref{thm-st-sd-decomp} below.

We state an important characterization of semistable self-dual objects.

\begin{theorem}
    \label{thm-ss-sd}
    In Situation~\textnormal{\ref{sit-quasi-ab-sd},}
    a self-dual object $(E, \phi) \in \calA^\sd$ is semistable,
    if and only if the object $E \in \calA$ is semistable.
\end{theorem}

\begin{proof}
    See~\cite[Corollary~3.16]{Bu2023}.
\end{proof}

\begin{definition}
    Let $E \in \calA$ be a $\tau$-semistable object.
    A \emph{$\tau$-Jordan--Hölder filtration} for $E$
    is a sequence of inclusions
    \[
        0 = E_0 \hookrightarrow E_1 \hookrightarrow \cdots \hookrightarrow E_m = E,
    \]
    such that each quotient $F_i = E_i / E_{i-1}$ is $\tau$-stable, and
    \[
        \tau (F_1) = \cdots = \tau (F_m).
    \]
\end{definition}

The following result is essentially due to
Rudakov~\cite[Theorem~3]{Rudakov1997}, also mentioned in
Joyce~\cite[Theorem~4.5]{Joyce2007III}.
Here, we only assume that $\calA$ is a quasi-abelian category,
instead of an abelian category.

\begin{lemma}
    \label{lem-tau-jh}
    In Situation~\textnormal{\ref{sit-quasi-ab-sd},}
    every $\tau$-semistable object $E \in \calA$
    has a $\tau$-Jordan--Hölder filtration.
    Moreover, for any two such filtrations,
    the quotients for one filtration is always
    a permutation of the quotients for the other one.
\end{lemma}

\begin{proof}
    We may assume that $E \nsimeq 0$.
    Consider the partially ordered set $P$ of subobjects $F \subset E$
    with $\tau (F) = \tau (E)$.
    Then, by the noetherian and artinian assumptions in
    Situation~\ref{sit-quasi-ab-sd},
    $P$ satisfies both the ascending and the descending chain conditions.
    Therefore, there exists a maximal chain,
    which verifies the existence part.

    The uniqueness part follows from standard arguments.
\end{proof}

\begin{theorem}
    \label{thm-st-sd-decomp}
    In Situation~\textnormal{\ref{sit-quasi-ab-sd},}
    if $(E, \phi) \in \calA^\sd$ is a $\tau$-stable self-dual object,
    then there is a decomposition
    \begin{equation}
        (E, \phi) \simeq
        (E_1, \phi_1) \oplus \cdots \oplus (E_m, \phi_m),
    \end{equation}
    where the self-dual objects $(E_i, \phi_i) \in \calA^\sd$ are pairwise non-isomorphic,
    and each $E_i$ is a $\tau$-stable object in $\calA$.
    Such a decomposition is unique up to a permutation.

    In particular, we have
    \begin{equation}
        \Aut (E, \phi) \simeq \bbZ_2^m .
    \end{equation}
\end{theorem}

\begin{proof}
    Consider a $\tau$-Jordan--Hölder filtration of $E$,
    \[
        0 = F_0 \hookrightarrow F_1 \hookrightarrow \cdots \hookrightarrow F_m = E.
    \]
    If $m \leq 1$, then we are already done.
    Otherwise, let $G_1, \dotsc, G_m$ denote the stepwise quotients.
    There is another $\tau$-Jordan--Hölder filtration
    \[
        0 = F_m^\perp \hookrightarrow F_{m-1}^\perp \hookrightarrow \cdots
        \hookrightarrow F_0^\perp = E,
    \]
    where $F_i^\perp \simeq (E/F_i)^\vee$, and we have used the self-dual structure
    $\phi \colon E \simto E^\vee$
    to identify $F_0^\perp \simeq E$.
    The stepwise quotients of this filtration are $G_m^\vee, \dotsc, G_1^\vee$.

    Since $(E, \phi)$ is $\tau$-stable,
    the map $F_1 \hookrightarrow E$ cannot factor through $F_1^\perp$.
    Thus, the composition $F_1 \hookrightarrow E \to E/F_1^\perp \simeq F_1^\vee$
    is non-zero, and hence is an isomorphism.
    One can verify that this defines a self-dual structure of $F_1$\,.
    
    Now, the inclusion $F_1 \hookrightarrow E$
    splits the short exact sequence
    $F_1^\perp \hookrightarrow E \twoheadrightarrow F_1^\vee$, so that
    $E \simeq F_1 \oplus F_1^\perp$.
    We can then identify $E^\vee \simeq F_1^\vee \oplus F_1^{\perp \vee}$.
    Since $\phi \colon E \simto E^\vee$ sends $F_1$ into $F_1^\vee$,
    and $F_1^\perp$ into $F_1^{\perp \vee}$, it follows that
    $\smash{\phi|_{F_1}}$ and $\smash{\phi|_{F_1^\perp}}$ must be isomorphisms onto
    $F_1^\vee$ and $F_1^{\perp \vee}$, respectively.
    They define self-dual structures on $F_1$ and $F_1^\perp$, and we have
    \[
        (E, \phi) \simeq (F_1, \phi|_{F_1}) \oplus (F_1^\perp, \phi|_{F_1^\perp}),
    \]
    with $F_1$ stable as an object of $\calA$.
    By the $\tau$-artinian assumption, 
    repeating the process with $F_1^\perp$ in place of $E$,
    we will eventually obtain a desired decomposition.

    To see that the $(E_i, \phi_i)$ are pairwise non-isomorphic, 
    suppose the contrary.
    We may then assume that there is an isomorphism
    $f \colon E_1 \simto E_2$\,, compatible with $\phi_1$ and $\phi_2$\,.
    The map $\id_{E_1} + \upi f \colon E_1 \to E$
    then gives an isotropic subobject,
    contradicting with the assumption that $(E, \phi)$ is $\tau$-stable.

    The uniqueness follows from the uniqueness part of
    Lemma~\ref{lem-tau-jh}, as any such decomposition
    gives rise to a $\tau$-Jordan--Hölder filtration of $E$.

    Finally, for the last statement,
    suppose that $f \colon E \simto E$ is an isomorphism.
    Then we must have $f (E_i) = E_i$ for all $i$.
    Since $\Aut (E_i) \simeq \Gm$\,,
    we have $\Aut (E_i, \phi_i) \simeq \bbZ_2$\,,
    as a scalar automorphism of $E_i$ preserves $\phi_i$
    if and only if it squares to the identity.
    This completes the proof.
\end{proof}

\section{Vertex algebras and modules}

\label{sect-vertex}

In this section, we discuss \emph{vertex algebras}
and related topics, including modules for vertex algebras.
A standard reference for this is Frenkel--Ben-Zvi~\cite{FrenkelBenZvi2004}.
We also introduce new notions including
\emph{involutive vertex algebras} and \emph{twisted modules} for vertex algebras,
which will become important algebraic structures
for studying enumerative invariants in self-dual categories.

Throughout this section,
we fix a field $\bbK$ of characteristic zero.

\subsection{Vertex algebras}

\begin{definition}
    \label{def-va}
    A \emph{graded vertex algebra} over $\bbK$ consists of the following data:

    \begin{itemize}
        \item 
            A $\bbZ$-graded $\bbK$-vector space
            \[
                V = \bigoplus_{n \in \bbZ} V_{(n)} \, .
            \]
            
        \item 
            An element $1 \in V_{(0)}$\,, called the \emph{unit}, or the \emph{vacuum}.
            
        \item 
            A linear map $D \colon V \to V$ of degree $2$ with respect to the grading,
            called the \emph{translation operator}.
            
        \item 
            A linear map
            \[
                Y (-, z) (-) \colon
                V \otimes V \longrightarrow
                V \llparen z \rrparen,
            \]
            respecting the grading, with $\deg z = -2$.
            This is called the \emph{multiplication},
            or the \emph{state--field correspondence}.
    \end{itemize}
    They should satisfy the following axioms:

    \begin{enumerate}
        \item 
            (\emph{Identity})
            For any $A \in V$, we have 
            \begin{align}
                Y (1, z) \, A & = A, \\
                Y (A, z) \, 1 & \in A + z \, V \llbr z \rrbr.
            \end{align}
            
        \item 
            (\emph{Translation})
            For any $A, B \in V$, we have 
            \begin{equation}
                \label{eq-def-va-transl}
                [D, Y (A, z)] \, B = 
                \frac{\partial}{\partial z} \, Y (A, z) \, B,
            \end{equation}
            where $[-,-]$ denotes the commutator,
            and we have $D (1) = 0$.
            
        \item 
            \label{itm-va-locality}
            (\emph{Locality})
            For any $A, B, C \in V$,
            there exists an element
            \begin{equation}
                X (A, B, C; y, z) \in
                V \llbr y, z \rrbr [y^{-1}, z^{-1}, (y - z)^{-1}],
            \end{equation}
            such that
            \begin{equation}
                Y (A, y) \, Y (B, z) \, C
                = \iota_{y, \> z} \bigl( X (A, B, C; y, z) \bigr)
                \ \in V \llparen y \rrparen \llparen z \rrparen,
            \end{equation}
            where $\iota_{y, \> z}$ is as in Definition~\ref{def-iota}.
            Moreover, if $A \in V_{(a)}$ and $B \in V_{(b)}$\,, then
            \begin{equation}
                X (A, B, C; y, z) = (-1)^{a b} \cdot X (B, A, C; z, y).
            \end{equation}
    \end{enumerate}
    It is common to introduce the following notation.
    For $A, B \in V$ and $n \in \bbZ$, write
    \begin{equation}
        A_{(n)} B = \res_{z=0} \bigl( z^n \, Y (A, z) \, B \bigr),
    \end{equation}
    so that
    \begin{equation}
        Y (A, z) \, B = \sum_{n \in \bbZ} A_{(n)} B \cdot z^{-n-1}.
    \end{equation}
    Note that if $A \in V_{(a)}$ and $B \in V_{(b)}$\,, 
    then $A_{(n)} B \in V_{(a+b-2n-2)}$\,.

    A \emph{morphism of vertex algebras}
    is a morphism of the underlying graded vector spaces,
    compatible with the unit elements,
    the translation operators,
    and the multiplication maps.
\end{definition}

\begin{remark}
    \label{rmk-va-loc}
    The locality axiom, Definition~\ref{def-va}~\ref{itm-va-locality},
    is sometimes stated as follows:

    \begin{itemize}
        \item 
            For any $A \in V_{(a)}$\,, $B \in V_{(b)}$\,, and $C \in V$,
            there exists an integer $N > 0$, such that
            \begin{equation}
                (y - z)^N \cdot \Bigl(
                    Y (A, y) \, Y (B, z) \, C -
                    (-1)^{ab} \cdot Y (B, z) \, Y (A, y) \, C
                \Bigr) = 0
            \end{equation}
            in $V \llbr y^{\pm 1}, z^{\pm 1} \rrbr$.
    \end{itemize}
    To see that these are equivalent, suppose that this condition is satisfied.
    Note that
    \begin{equation}
        V \llparen y \rrparen \llparen z \rrparen \cap
        V \llparen z \rrparen \llparen y \rrparen =
        V \llbr y, z \rrbr [y^{-1}, z^{-1}]
    \end{equation}
    in $V \llbr y^{\pm 1}, z^{\pm 1} \rrbr$,
    so both $(y - z)^N \cdot Y (A, y) \, Y (B, z) \, C$
    and $(y - z)^N \cdot Y (B, z) \, Y (A, y) \, C$
    must lie in this subspace, giving the required element $X (A, B, C; y, z)$.
    Conversely, assuming Definition~\ref{def-va}~\ref{itm-va-locality},
    we may choose $N$ such that
    $(y - z)^N \cdot X (A, B, C; y, z) \in V \llbr y, z \rrbr [y^{-1}, z^{-1}]$.
    Then, we have 
    \begin{equation}
        \iota_{y, \> z} \bigl( (y - z)^N \cdot X (A, B, C; y, z) \bigr) =
        \iota_{z, \> y} \bigl( (y - z)^N \cdot X (A, B, C; y, z) \bigr)
    \end{equation}
    in $V \llbr y^{\pm 1}, z^{\pm 1} \rrbr$,
    which implies the condition given here.
\end{remark}

We state a few standard results about vertex algebras.

\begin{lemma}
    \label{lem-va-facts}
    Let $V$ be a graded vertex algebra over $\bbK$.
    \begin{enumerate}
        \item
            \label{itm-va-az1}
            For any $A \in V$, we have
            \begin{equation}
                Y (A, z) \, 1 =
                \exp (z D) \, A.
            \end{equation}
            In particular, $D (A) = A_{(-2)} 1,$ and $D$ is determined by $Y$.
        \item
            \label{itm-va-ydaz}
            For any $A, B \in V$, we have
            \begin{equation}
                Y (D (A), z) \, B =
                \frac{\partial}{\partial z} \, Y (A, z) \, B .
            \end{equation}
        \item
            \label{itm-va-assoc}
            For any $A, B, C \in V$, we have
            \begin{equation}
                Y \bigl( Y (A, y - z) \, B, \, z \bigr) \, C
                = \iota_{z, \> y-z} \bigl( X (A, B, C; y, z) \bigr)
            \end{equation}
            in $V \llparen z \rrparen \llparen y-z \rrparen$.
    \end{enumerate}
\end{lemma}

\begin{proof}
    These can be found in Frenkel--Ben-Zvi
    \cite[Lemma~3.1.3, Corollary~3.1.6, Theorem~3.2.1]{FrenkelBenZvi2004}.
\end{proof}

Next, we recall a well-known construction,
associating to each graded vertex algebra a graded Lie algebra.

\begin{definition}
    \label{def-lie-alg}
    A \emph{graded Lie algebra} over $\bbK$ consists of the following data:

    \begin{itemize}
        \item 
            A $\bbZ$-graded $\bbK$-vector space
            \[
                L = \bigoplus_{n \in \bbZ} L_{(n)} \, .
            \]
            
        \item 
            A linear map
            \[
                [-, -] \colon L \otimes L \longrightarrow L,
            \]
            preserving the grading, called the \emph{Lie bracket}.
    \end{itemize}
    They should satisfy the following axioms:

    \begin{enumerate}
        \item 
            (\emph{Antisymmetry})
            For any $A \in L_{(a)}$ and $B \in L_{(b)}$\,, we have
            \begin{equation}
                [A, B] + (-1)^{a b} \cdot [B, A] = 0.
            \end{equation}
        \item 
            (\emph{Jacobi identity})
            For any $A, B, C \in L$, we have
            \begin{equation}
                [[A, B], C] + [[B, C], A] + [[C, A], B] = 0.
            \end{equation}
    \end{enumerate}
\end{definition}

\begin{theorem}
    \label{thm-lie-alg}
    Let $V$ be a graded vertex algebra over $\bbK$.
    Define a graded $\bbK$-vector space
    \[
        L = V / D (V),
    \]
    with grading given by $L_{(n)} = V_{(n+2)} / D (V_{(n)})$
    for $n \in \bbZ$.
    Define a Lie bracket on $L$ by setting
    \begin{equation}
        [A + D (V), B + D (V)] = A_{(0)} B + D (V)
    \end{equation}
    for all $A, B \in V$.
    Then the Lie bracket is well-defined,
    and makes $L$ into a graded Lie algebra over $\bbK$.

    Moreover, a morphism of graded vertex algebras
    induces a morphism of the associated graded Lie algebras.
\end{theorem}

See, for example, Frenkel--Ben-Zvi \cite[\S4.1]{FrenkelBenZvi2004}.

\subsection{Modules and twisted modules}
\label{sect-va-mod}

We introduce \emph{twisted modules} for a vertex algebra,
in a similar spirit to
\emph{twisted vertex algebras} in the sense of Anguelova~\cite{Anguelova2013},
based on ideas of Frenkel and Reshetikhin~\cite{FrenkelReshetikhin1997}.

We note that there is another notion of twisted modules for vertex algebras,
as in \cite[\S5.6]{FrenkelBenZvi2004},
and it is different from ours.
Namely, the existing notion of twisted modules
interacts with an automorphism of the vertex algebra,
while our notion of twisted modules interacts with an \emph{involution}
of the vertex algebra, in the sense of Definition~\ref{def-va-invol} below,
which will not be an automorphism in general.

Before we discuss twisted modules,
we give a definition of ordinary modules for a vertex algebra,
following~\cite[\S5.1]{FrenkelBenZvi2004}.

\begin{definition}
    \label{def-va-module}
    Let $V$ be a graded vertex algebra over $\bbK$.
    A \emph{graded $V \<$-module} consists of the following data:

    \begin{itemize}
        \item 
            A $\bbZ$-graded $\bbK$-vector space
            \[
                W = \bigoplus_{n \in \bbZ} W_{(n)} \, .
            \]
        \item 
            A linear map
            \[
                Y^W (-, z) (-) \colon
                V \otimes W \longrightarrow
                W \llparen z \rrparen,
            \]
            respecting the gradings on $V$ and $W$, with $\deg z = -2$.
    \end{itemize}
    They should satisfy the following conditions:

    \begin{enumerate}
        \item 
            \label{itm-vm-unit}
            For any $M \in W$, we have 
            \begin{equation}
                \label{eq-vm-unit}
                Y^W (1, z) \, M = M.
            \end{equation}
        \item 
            \label{itm-vm-assoc}
            For any $A, B \in V$ and $M \in W$,
            there exists an element
            \begin{equation}
                \label{eq-vm-xw}
                X^W (A, B, M; y, z) \in
                W \llbr y, z \rrbr [y^{-1}, z^{-1}, (y - z)^{-1}],
            \end{equation}
            such that
            \begin{alignat}{2}
                \label{eq-tw-mod-assoc-1}
                Y^W (A, y) \, Y^W (B, z) \, M
                & = \iota_{y, \> z} \bigl( X^W (A, B, M; y, z) \bigr)
                && \ \in W \llparen y \rrparen \llparen z \rrparen, \\
                \label{eq-tw-mod-assoc-2}
                Y^W \bigl( Y (A, y - z) \, B, \, z \bigr) \, M
                & = \iota_{z, \> y-z} \bigl( X^W (A, B, M; y, z) \bigr)
                && \ \in W \llparen z \rrparen \llparen y-z \rrparen.
            \end{alignat}
    \end{enumerate}
\end{definition}

Modules for a vertex algebra give rise to
modules for the associated Lie algebra,
introduced in Theorem~\ref{thm-lie-alg}.

\begin{definition}
    \label{def-lie-mod}
    Let $L$ be a graded Lie algebra over $\bbK$.
    A \emph{graded module} over $L$ 
    is a $\bbZ$-graded $\bbK$-vector space
    \[
        W = \bigoplus_{n \in \bbZ} W_{(n)} \, ,
    \]
    together with a linear map
    $\cdot \ \colon L \otimes W \to W$
    preserving the grading,
    such that for any $A \in L_{(a)}$\,, $B \in L_{(b)}$\,,
    and $M \in W$, we have
    \begin{equation}
        A \cdot (B \cdot M) - (-1)^{a b} \, B \cdot (A \cdot M)
        = [A, B] \cdot M.
    \end{equation}
\end{definition}

\begin{lemma}
    Let $V$ be a graded vertex algebra over $\bbK,$
    and $W$ a graded $V \<$-module.
    Define a linear map
    \begin{align*}
        \cdot \ \colon (V / D (V)) \otimes W
        & \longrightarrow W, \\
        (A + D (V), M) & \longmapsto A \cdot M 
    \end{align*}
    by
    \begin{equation}
        A \cdot M = \res_{z = 0} Y^W (A, z) \, M.
    \end{equation}
    Then this map is well-defined, 
    and equips $W$ with the structure of a graded module
    over the Lie algebra $V / D (V)$.
\end{lemma}

\begin{proof}
    This is a consequence of~\cite[Theorem~5.1.6]{FrenkelBenZvi2004}.
\end{proof}

Next, we introduce \emph{twisted modules} for a vertex algebra.
Roughly speaking, a twisted module
is similar to an ordinary module,
as in Definition~\ref{def-va-module},
but we allow the product $Y^W (A, y) \, Y^W (B, z) \, M$
to be singular not only along the loci $y = 0$, $z = 0$, and $y - z = 0$,
but also along the locus $y + z = 0$.

\begin{definition}
    \label{def-tw-module}
    Let $V$ be a graded vertex algebra over $\bbK$.
    A \emph{graded twisted $V \<$-module} consists of the data
    as in Definition~\ref{def-va-module},
    except that we replace~\eqref{eq-vm-xw} by
    \begin{equation}
        X^W (A, B, M; y, z) \in
        W \llbr y, z \rrbr [y^{-1}, z^{-1}, (y \pm z)^{-1}].
    \end{equation}
\end{definition}

Note that every graded $V \<$-module is automatically a graded twisted $V \<$-module.

\begin{lemma}
    \label{lem-tw-mod}
    Let $V$ be a graded vertex algebra over $\bbK,$
    and let $W$ be a graded twisted $V \<$-module.
    \begin{enumerate}
        \item 
            \label{itm-tw-mod-translation}
            For any $A \in V$ and $M \in W,$ we have
            \begin{equation}
                Y^W (D (A), z) \, M = \frac{\partial}{\partial z} \, Y^W (A, z) \, M.
            \end{equation}
        \item
            \label{itm-tw-mod-comm}
            For any $A \in V_{(a)}\,,$ $B \in V_{(b)}\,,$ and $M \in W,$ we have
            \begin{equation}
                X^W (A, B, M; y, z) = (-1)^{a b} \cdot X^W (B, A, M; z, y).
            \end{equation}
    \end{enumerate}
\end{lemma}

\begin{proof}
    For~\ref{itm-tw-mod-translation}, for any $M \in W$,
    consider the element $N (y, z) = X^W (A, 1, M; y, z)$.
    Applying~\eqref{eq-tw-mod-assoc-1}, we obtain
    \begin{equation}
        Y^W (A, y) \, M = \iota_{y, \> z} ( N (y, z) ).
    \end{equation}
    In particular, $N (y, z)$ does not depend on $z$, 
    and we can identify it with an element $N (y) \in W \llparen y \rrparen$.
    Applying~\eqref{eq-tw-mod-assoc-2},
    using the identity $A_{(-2)} 1 = D (A)$
    from Lemma~\ref{lem-va-facts}~\ref{itm-va-az1}, we see that
    \begin{equation}
        Y^W (D (A), z) \, M
        = \res_{y-z=0} \bigl( (y-z)^{-2} \, N (y) \bigr)
        = \frac{\partial}{\partial z} \, N (z),
    \end{equation}
    as desired.

    For~\ref{itm-tw-mod-comm}, the proof in
    \cite[Remark~5.1.4]{FrenkelBenZvi2004}
    still applies in our case.
\end{proof}

\begin{remark}
    \label{rmk-tw-mod-loc}
    The locality property in a twisted module,
    Lemma~\ref{lem-tw-mod}~\ref{itm-tw-mod-comm}, 
    can also be written as follows:
    \begin{itemize}
        \item 
            For any $A \in V_{(a)}$\,, $B \in V_{(b)}$\,, and $M \in W$,
            there exists $N > 0$, such that
            \begin{equation}
                (y^2 - z^2)^N \cdot \Bigl(
                    Y^W (A, y) \, Y^W (B, z) \, M -
                    (-1)^{ab} \cdot Y^W (B, z) \, Y^W (A, y) \, M
                \Bigr) = 0
            \end{equation}
            in $M \llbr y^{\pm1}, z^{\pm1} \rrbr$.
    \end{itemize}
    The reason for this is analogous to Remark~\ref{rmk-va-loc}.
\end{remark}

\subsection{Involutive vertex algebras and modules}

\begin{definition}
    \label{def-va-invol}
    An \emph{involutive graded vertex algebra} consists of the following data:
    \begin{itemize}
        \item 
            A graded vertex algebra $V$.
        \item
            A linear map $I \colon V \to V$,
            preserving the grading,
            such that $I^2 = \id_{V}$\,.
            This is called an \emph{involution} of $V$, and is also denoted by $(-)^\vee$.
    \end{itemize}
    They should satisfy the following properties:
    \begin{itemize}
        \item 
            We have $1^\vee = 1$.
        \item
            For any $A \in V$, we have $D (A^\vee) = - (D A)^\vee$.
        \item
            For any $A, B \in V$, we have
            $Y (A^\vee, z) \, B^\vee = \bigl( Y (A, -z) \, B \bigr)^\vee$.
    \end{itemize}
\end{definition}

\begin{definition}
    \label{def-tw-mod-invol}
    Let $V$ be an involutive graded vertex algebra,
    and let $W$ be a graded twisted $V \<$-module.
    We say that $W$ is \emph{involutive}, if it satisfies the following property:
    \begin{itemize}
        \item 
            For any $A \in V$ and $M \in W$, we have
            \begin{equation}
                \label{eq-def-tw-mod-antisym}
                Y^W (A^\vee, z) \, M = Y^W (A, -z) \, M.
            \end{equation}
    \end{itemize}
\end{definition}

We will show that given an involutive graded twisted module over a vertex algebra,
one can construct a twisted module over the associated Lie algebra.
Let us first introduce the involutive versions of
Lie algebras and twisted modules.

\begin{definition}
    An \emph{involutive graded Lie algebra} consists of the following data:
    \begin{itemize}
        \item 
            A graded Lie algebra $L$.
        \item
            A linear map $I \colon L \to L$,
            preserving the grading,
            such that $I^2 = \id_{L}$\,.
            This is called an \emph{involution} of $V$, and is also denoted by $(-)^\vee$.
    \end{itemize}
    They should satisfy the following property:
    \begin{itemize}
        \item
            For any $A, B \in L,$ we have
            $[A^\vee, B^\vee] = -[A, B]^\vee$.
    \end{itemize}
\end{definition}

The following is a graded version of~\cite[Definition~7.10]{Bu2023}.

\begin{definition}
    \label{def-tw-lie-mod}
    Let $L$ be an involutive graded Lie algebra over $\bbK$.
    A \emph{graded twisted module} over $L$ 
    is a $\bbZ$-graded $\bbK$-vector space
    \[
        W = \bigoplus_{n \in \bbZ} W_{(n)} \, ,
    \]
    together with a linear map
    $\heart \colon L \otimes W \to W$
    preserving the grading,
    such that for any $A \in L_{(a)}$\,, $B \in L_{(b)}$\,,
    and $M \in W$, we have
    \begin{align}
        \label{eq-tw-lie-mod-antisym}
        A \heart M & = -A^\vee \heart M, \\
        \label{eq-tw-lie-mod-jacobi}
        A \heart (B \heart M) - (-1)^{a b} \cdot B \heart (A \heart M) & =
        [A, B] \heart M - [A, B^\vee] \heart M.
    \end{align}
\end{definition}

\begin{remark}
    As discussed in~\cite[\S7.3]{Bu2023},
    for an involutive graded Lie algebra $L$,
    giving a graded twisted $L$-module $W$
    is equivalent to giving a graded module over the subalgebra
    $L^+ \subset L$ defined by
    \begin{equation}
        L^+ = \{ A \in L \mid A^\vee = -A \},
    \end{equation}
    with the action $\cdot \ \colon L^+ \otimes W \to W$ given by
    \begin{equation}
        A \cdot M = \frac{1}{2} A \heart M
    \end{equation}
    for $A \in L^+$ and $M \in W$.
    
    For this reason, Definition~\ref{def-tw-lie-mod}
    might seem unnecessary,
    as one could have defined a twisted $L$-module
    simply as an $L^+$-module.
    However, we believe that the relations
    \eqref{eq-tw-lie-mod-antisym}--\eqref{eq-tw-lie-mod-jacobi}
    are more natural than the usual Jacobi identity for an $L^+$-action,
    since similar relations appear in the vertex algebra picture,
    as will be seen in the proof of Theorem~\ref{thm-tw-lie-mod} below.
\end{remark}

\begin{theorem}
    \label{thm-tw-lie-mod}
    Let $V$ be an involutive graded vertex algebra,
    and let $W$ be an involutive graded twisted $V \<$-module.
    Define a linear map
    \begin{align*}
        {\heart} \colon (V / D (V)) \otimes W
        & \longrightarrow W, \\
        (A + D (V), M) & \longmapsto A \heart M 
    \end{align*}
    by
    \begin{equation}
        A \heart M = \res_{z = 0} Y^W (A, z) \, M.
    \end{equation}
    Then $\heart$ is well-defined,
    and equips $W$ with the structure of a graded twisted module 
    over the involutive graded Lie algebra $V / D (V)$.
\end{theorem}

\begin{proof}
    To see that $\heart$ is well-defined, it is enough to prove that
    \[
        \res_{z = 0} Y^W (D A, z) \, M = 0
    \]
    for any $A \in V$ and $M \in W$.
    But this follows from Lemma~\ref{lem-tw-mod}~\ref{itm-tw-mod-translation},
    as the residue of a derivative in $z$ has to be zero.

    Equation~\eqref{eq-tw-lie-mod-antisym}
    follows from taking the $z = 0$ residue of~\eqref{eq-def-tw-mod-antisym}.
    
    To prove~\eqref{eq-tw-lie-mod-jacobi}, write
    \[
        F (y, z) =
        X^W (A, B, M; y, z) \in
        W \llbr y, z \rrbr [y^{-1}, z^{-1}, (y \pm z)^{-1}].
    \]
    By the residue theorem, we have
    \begin{equation}
        \label{eq-pf-lie-mod}
        \res_{y=0} \res_{z=0} F (y, z) =
        \res_{z=0} \Bigl(
            \res_{y=0} F (y, z) +
            \res_{y=z} F (y, z) +
            \res_{y=-z} F (y, z)
        \Bigr).
    \end{equation}
    Analysing each term of~\eqref{eq-pf-lie-mod},
    and using Lemma~\ref{lem-tw-mod}~\ref{itm-tw-mod-comm},
    we have
    \begin{align}
        \res_{y=0} \res_{z=0} F (y, z)
        & = A \heart (B \heart M), \\
        \res_{z=0} \res_{y=0} F (y, z)
        & = (-1)^{a b} \cdot B \heart (A \heart M), \\
        \res_{z=0} \res_{y=z} F (y, z)
        & = [A, B] \heart M, \\
        \res_{z=0} \res_{y=-z} F (y, z)
        & = -[A, B^\vee] \heart M.
    \end{align}
    Here, the last equality is because
    \begin{align*}
        & \phantom{{} = {}}
        {\res_{z=0} \res_{y=-z} F (y, z)} \\
        & =
        \res_{z=0} \res_{y=-z} Y^W (A, y) \, Y^W (B, z) \, M \\
        & =
        \res_{z=0} \res_{y=-z} Y^W (A, y) \, Y^W (B^\vee, -z) \, M \\
        & =
        \res_{z=0} \res_{y=-z} Y^W ( Y (A, y+z) \, B^\vee, -z ) \, M \\
        & =
        -[A, B^\vee] \heart M,
        \numberthis
    \end{align*}
    where we used~\eqref{eq-def-tw-mod-antisym}
    in the second step.
    Substituting these into~\eqref{eq-pf-lie-mod}
    proves~\eqref{eq-tw-lie-mod-jacobi}.
\end{proof}

\subsection{Equivalent formulations}

We discuss an alternative way to define vertex algebras
and twisted modules.
For vertex algebras, this is essentially due to Kim~\cite{Kim2011},
and we state this alternative formulation below.

\begin{theorem}
    \label{thm-va-alt-def}
    A graded vertex algebra over $\bbK$
    can be equivalently defined as the following data:

    \begin{itemize}
        \item 
            A $\bbZ$-graded $\bbK$-vector space $V$.
        \item
            For each $n \geq 0$, a linear map
            \[
                X_n (-, \dotsc, -; z_1, \dotsc, z_n) \colon
                V^{\otimes n} \longrightarrow
                V \llbr z_1, \dotsc, z_n \rrbr [z_1^{-1}, (z_i - z_j)^{-1}],
            \]
            where we invert $z_i$ for all $1 \leq i \leq n,$
            and $z_i - z_j$ for all $1 \leq i < j \leq n$.
            This map is required to preserve grading,
            with $\deg z_i = -2$ for all $i$.
    \end{itemize}
    We require that they satisfy the following properties:

    \begin{enumerate}
        \item
            \label{itm-va-multi-unit}
            For any $A \in V,$ we have
            \begin{equation}
                X_1 (A; z) \in A + z \, V \llbr z \rrbr.
            \end{equation}
            
        \item
            \label{itm-va-multi-equivar}
            Each $X_n$ is equivariant under the symmetric group $\frS_n$.
            That is, for any $\sigma \in \frS_n$\,, 
            given $A_i \in V_{(a_i)}$ for $i = 1, \dotsc, n,$ we have
            \begin{equation}
                X_n (A_{\sigma (1)} \, , \dotsc, A_{\sigma (n)};
                z_{\sigma (1)} \, , \dotsc, z_{\sigma (n)}) =
                \pm X_n (A_1, \dotsc, A_n; z_1, \dotsc, z_n),
            \end{equation}
            where the sign is \textnormal{`$-$'} if and only if $\sigma$
            restricts to an odd permutation on the odd elements.
            
        \item 
            \label{itm-va-multi-assoc}
            For integers $1 \leq \ell \leq m$ and $n \geq 0,$ 
            and $A_1, \dotsc, A_{\ell-1}, A_{\ell+1}, \dotsc, A_m, B_1, \dotsc, B_n \in V,$
            we have
            \begin{multline}
                \label{eq-va-multi-assoc}
                X_m \Bigl( 
                    A_1, \dotsc, A_{\ell-1}, \,
                    X_n (B_1, \dotsc, B_n; z_1, \dotsc, z_n), \,
                    A_{i+1}, \dotsc, A_m \, ; \ 
                    y_1, \dotsc, y_m
                \Bigr)
                \\[-1ex] = 
                \iota_{ \{ y_1, \> \dotsc \> , \> y_m \}, \{ z_1, \> \dotsc \> , \> z_n \} }
                \Bigl( 
                    X_{m+n-1} \bigl(
                        A_1, \dotsc, A_{\ell-1}, B_1, \dotsc, B_n,
                        A_{\ell+1}, \dotsc, A_m \, ;
                        \\[-1ex] 
                        y_1, \dotsc, y_{\ell-1}, \,
                        y_\ell + z_1, \dotsc, y_\ell + z_n, \,
                        y_{\ell+1}, \dotsc, y_m
                    \bigr)
                \Bigr), 
            \end{multline}
            where 
            $\iota_{ \{ y_1, \> \dotsc \> , \> y_m \}, \{ z_1, \> \dotsc \> , \> z_n \} }$
            is the map defined in \textnormal{\S\ref{sect-pow-ser}}.
    \end{enumerate}
\end{theorem}

\begin{proof}
    Given a graded vertex algebra $V$, define the maps $X_n$ by
    \begin{equation}
        X_n (A_1, \dotsc, A_n; z_1, \dotsc, z_n) =
        \iota \sss{z_1, \> \dotsc, \> z_n}{-1}
        \bigl( Y (A_1, z_1) \cdots Y (A_n, z_n) \, 1 \bigr) ,
    \end{equation}
    where we take the unique preimage under the embedding
    \[
        \iota_{z_1, \> \dotsc, \> z_n} \colon
        V \llbr z_1, \dotsc, z_n \rrbr [z_i^{-1}, (z_i - z_j)^{-1}]
        \longhookrightarrow
        V \llparen z_1 \rrparen \cdots \llparen z_n \rrparen .
    \]
    It follows from Kim~\cite[Theorem~3.14]{Kim2011}
    that this is well-defined,
    and that it satisfies~\ref{itm-va-multi-assoc}.
    The properties~\ref{itm-va-multi-unit}
    and~\ref{itm-va-multi-equivar}
    follow from the identity and locality axioms for a vertex algebra.

    Conversely, given this new data,
    we may define a graded vertex algebra structure on $V$ by taking
    \begin{align}
        \label{eq-va-alt-unit}
        1 & = X_0 ( \, ; \, ) , \\
        \label{eq-va-alt-transl}
        D (A) & = \res_{z=0} z^{-2} \, X_1 (A; z), \\
        \label{eq-va-alt-mult}
        Y (A, z) \, B & = X_2 (A, B; z, 0)
    \end{align}
    for all $A, B \in V$.
    The axioms for a graded vertex algebra are easily verified,
    except for~\eqref{eq-def-va-transl}, which we prove as follows.
    Applying~\eqref{eq-va-multi-assoc} in two different ways, we have
    \begin{align}
        X_1 (Y (A, z) \, B; y) & =
        \iota_{y, z} \bigl(
            X_2 (A, B; y + z, y)
        \bigr), \\
        Y (A, z) \, X_1 (B; y) & =
        \iota_{z, y} \bigl(
            X_2 (A, B; z, y)
        \bigr).
    \end{align}
    Multiplying by $y^{-2}$ and taking the residue, we obtain
    \begin{align*}
        [D, Y (A, z)] \, B & =
        \res_{y=0} y^{-2} \bigl(
            X_2 (A, B; y + z, y) - X_2 (A, B; z, y)
        \bigr) \\
        & = \frac{\partial}{\partial z} \, Y (A, z) \, B.
    \end{align*}
\end{proof}

This can be naturally generalized to twisted modules.

\begin{theorem}
    \label{thm-vm-alt-def}
    \allowdisplaybreaks
    Let $V$ be graded vertex algebra over $\bbK,$
    defined as in Theorem~\textnormal{\ref{thm-va-alt-def}}.
    Then a graded twisted $V \<$-module
    can be equivalently defined as the following data:

    \begin{itemize}
        \item 
            A $\bbZ$-graded $\bbK$-vector space $W$.
        \item
            For each $n \geq 0$, a linear map
            \[
                X^W_n (-, \dotsc, -; z_1, \dotsc, z_n) \colon
                V^{\otimes n} \otimes W \longrightarrow
                W \llbr z_1, \dotsc, z_n \rrbr [z_i^{-1}, (z_i \pm z_j)^{-1}],
            \]
            where we invert $z_i$ for all $1 \leq i \leq n,$
            and $z_i \pm z_j$ for all $1 \leq i < j \leq n$.
    \end{itemize}
    We require that they satisfy the following properties:

    \begin{enumerate}
        \item 
            For any $1 \leq \ell \leq m$ and $n \geq 0$,
            $A_1, \dotsc, A_{\ell-1}, A_{\ell+1}, \dotsc, A_m, B_1, \dotsc, B_n \in V,$
            and $M \in W,$ we have
            \begin{multline}
                \label{eq-tw-mod-assoc-x1}
                X^W_m \Bigl( 
                    A_1, \dotsc, A_{\ell-1}, \,
                    X_n (B_1, \dotsc, B_n; z_1, \dotsc, z_n), \,
                    A_{\ell+1}, \dotsc, A_m, M; \ 
                    y_1, \dotsc, y_m
                \Bigr)
                \\*[-1ex] 
                \shoveleft{
                = \iota_{ \{ y_1, \> \dotsc \> , \> y_m \}, \{ z_1, \> \dotsc \> , \> z_n \} } \Bigl( 
                    X^W_{m+n-1}
                    \bigl(
                        A_1, \dotsc, A_{\ell-1}, B_1, \dotsc, B_n, 
                        A_{\ell+1}, \dotsc, A_m, M;
                } \\*[-1ex]
                        y_1, \dotsc, y_{\ell-1}, \,
                        y_\ell + z_1, \dotsc, y_\ell + z_n, \,
                        y_{\ell+1}, \dotsc, y_m
                    \bigr)
                \Bigr) ,
            \end{multline}
            where 
            $\iota_{ \{ y_1, \> \dotsc \> , \> y_m \}, \{ z_1, \> \dotsc \> , \> z_n \} }$
            is the map defined in \textnormal{\S\ref{sect-pow-ser}}.
            
        \item 
            For any $m, n \geq 0$,
            $A_1, \dotsc, A_m, B_1, \dotsc, B_n \in V,$
            and $M \in W,$ we have
            \begin{multline}
                \label{eq-tw-mod-assoc-x2}
                X^W_m \Bigl( 
                    A_1, \dotsc, A_m, \,
                    X^W_n (B_1, \dotsc, B_n, M; z_1, \dotsc, z_n); \ 
                    y_1, \dotsc, y_m
                \Bigr)
                \\*[-1ex] 
                \shoveleft{
                = \iota_{ \{ y_1, \> \dotsc \> , \> y_m \}, \{ z_1, \> \dotsc \> , \> z_n \} } \Bigl( 
                    X^W_{m+n} (
                        A_1, \dotsc, A_m, B_1, \dotsc, B_n, M;
                } \\*[-1ex]
                        y_1, \dotsc, y_m,
                        z_1, \dotsc, z_n
                    )
                \Bigr) .
            \end{multline}
    \end{enumerate}
\end{theorem}

\begin{proof}
    We essentially follow the argument of Kim~\cite{Kim2011},
    making modifications as necessary.
    For simplicity, we omit the signs and assume that all elements are even,
    but the argument can be adapted to the general case
    by adding the appropriate signs.

    Given a graded twisted $V$-module $W$, define the maps $X^W_n$ by
    \begin{equation}
        X^W_n (A_1, \dotsc, A_n, M; z_1, \dotsc, z_n) =
        \iota \sss{z_1, \> \dotsc, \> z_n}{-1}
        \bigl( Y^W (A_1, z_1) \cdots Y^W (A_n, z_n) \, M \bigr),
    \end{equation}
    where we take the unique preimage under the embedding
    \begin{equation}
        \iota_{z_1, \> \dotsc, \> z_n} \colon
        W \llbr z_1, \dotsc, z_n \rrbr [z_i^{-1}, (z_i \pm z_j)^{-1}]
        \longhookrightarrow
        W \llparen z_1 \rrparen \cdots \llparen z_n \rrparen .
    \end{equation}
    To see that such a preimage exists,
    note that by Remark~\ref{rmk-tw-mod-loc},
    for fixed elements $A_1, \dotsc, A_n, M$ as above,
    there exists $N > 0$ such that in the expression
    \begin{equation}
        \biggl( \prod_{1 \leq i < j \leq n} {} (z_i^2 - z_j^2) \biggr)^N \cdot
        Y^W (A_1, z_1) \cdots Y^W (A_n, z_n) \, M,
    \end{equation}
    the order of the operators $Y^W (A_i, z_i)$ can be permuted freely.
    This means that this expression lies in the subspace 
    $W \llparen z_{\smash{\sigma (1)}} \rrparen
    \cdots \llparen z_{\smash{\sigma (n)}} \rrparen \subset
    W \llbr z_1^{\pm 1}, \dotsc, z_n^{\pm1} \rrbr$
    for all $\sigma \in \frS_n$\,.
    Consequently, it belongs to their intersection,
    $W \llparen z_1, \dotsc, z_n \rrparen$.

    From here, verifying
    \eqref{eq-tw-mod-assoc-x1}--\eqref{eq-tw-mod-assoc-x2}
    is straightforward.

    Conversely, given this new data,
    if we define
    \begin{equation}
        Y^W (A, z) \, M = X^W_1 (A, M; z)
    \end{equation}
    for all $A \in V$ and $M \in W$,
    then this defines a twisted $V$-module structure on $W$.
    Note that \eqref{eq-vm-unit} holds since
    \[
        Y^W (1, z) \, M = X^W_1 (X_0 ( \, ; \, ), M; z)
        = X^W_0 (M; \, ) = M.
    \]
\end{proof}

\section{Vertex algebras and modules from moduli stacks}

\label{sect-vertex-moduli}

\subsection{\texorpdfstring{$H$}{H}-spaces and modules}

We briefly discuss several notions related to $H$-spaces,
which will be a key ingredient in constructing
vertex algebras and modules from moduli stacks.
For a general account on $H$-spaces, see, for example,~\cite[\S3.C]{Hatcher2002}.

In the following,
write $\cat{Ho} (\cat{Top})$ for the homotopy category of topological spaces,
whose objects are topological spaces,
which we assume to have homotopy types of CW complexes,
and morphisms are homotopy classes of continuous maps.

\begin{definition}
    \label{def-h-sp}
    An \emph{$H$-space} is a unital magma object in $\cat{Ho} (\cat{Top})$,
    that is, a triple $(X, m, e)$, where
    \begin{itemize}
        \item
            $X$ is a topological space,
        \item
            $m \colon X \times X \to X$ is a continuous map, called \emph{multiplication}, and
        \item
            $e \colon {*} \to X$ is a map from the singleton space $*$,
            called the \emph{unit},
    \end{itemize}
    such that multiplication is left and right unital up to homotopies.
    Similarly,
    \begin{itemize}
        \item 
            An \emph{associative $H$-space}
            is a monoid object in $\cat{Ho} (\cat{Top})$,
            i.e.\ an $H$-space satisfying associativity up to a homotopy.
        \item 
            A \emph{commutative $H$-space}
            is a unital commutative monoid object in $\cat{Ho} (\cat{Top})$,
            i.e.\ an associative $H$-space
            satisfying commutativity up to a homotopy.
    \end{itemize}
\end{definition}

\begin{definition}
    \label{def-pi0-h-sp}
    Let $X$ be a commutative $H$-space.
    Then the set $\pi_0 (X)$ of connected components of $X$
    has the structure of a commutative monoid.
    
    For $\alpha, \beta \in \pi_0 (X)$, write
    \[
        m_{\alpha, \> \beta} \colon X_\alpha \times X_\beta
        \longrightarrow X_{\alpha + \beta}
    \]
    for the restriction of the multiplication map $m$.
\end{definition}

In our applications,
the $H$-space $X$ will be taken to be the topological realization
of a moduli stack of objects in an additive category
or an additive $\infty$-category,
with the $H$-space structure given by the direct sum.

\begin{definition}
    \label{def-h-sp-mod}
    Let $X$ be an associative $H$-space.
    A \emph{module} over $X$ is a module object for $X$ in $\cat{Ho} (\cat{Top})$,
    which consists of a topological space $Y$, and a map
    \[
        m^Y \colon X \times Y \longrightarrow Y
    \]
    satisfying the unit law and the associativity law.

    If, moreover, $Y$ is also an $H$-space,
    and the $H$-space structure is respected by $m^Y$,
    then we say that $m^Y$ is an \emph{$H$-action} of $X$ on $Y$.
\end{definition}

We also define involutive $H$-spaces and modules,
which capture the key properties of 
the moduli stack of objects in a self-dual category,
and the moduli stack of self-dual objects,
respectively.
    
\begin{definition}
    \label{def-h-sp-invol}
    An \emph{involutive commutative $H$-space}
    is a commutative $H$-space $X$,
    together with an $H$-action by $\bbZ_2$\,,
    that is, an isomorphism of $H$-spaces 
    \[
        I \colon X \longrightarrow X,
    \]
    such that $I^2$ is homotopic to $\id_X$\,.

    An \emph{involutive module} over an involutive commutative $H$-space $X$
    is a module $Y$ over the $X$, such that the diagram
    \[ \begin{tikzcd}[row sep=tiny]
        X \times Y
        \ar[dr, pos=.3, "m^Y"]
        \ar[dd, "I \> \times \> \id_Y"']
        \\ & Y
        \\ X \times Y
        \ar[ur, pos=.3, "m^Y"']
    \end{tikzcd} \]
    is commutative in $\cat{Ho} (\cat{Top})$.

    In particular, when $X, Y$ are discrete spaces,
    we say that $X$ is an \emph{involutive abelian monoid},
    and that $Y$ is an \emph{involutive $X$-module}.
\end{definition}

\subsection{Vertex algebras from moduli stacks}
\label{sect-cons-va}

We discuss the vertex algebra construction of
Joyce~\cite[Theorem~3.14]{JoyceHall},
also discussed in~\cite{GrossJoyceTanaka, Joyce2021, Latyntsev2021}.
We use a slightly different setting compared to Joyce's,
in that Joyce works with a notion of
\emph{homology theories for algebraic stacks},
while we work with $H$-spaces instead of stacks.

To construct this vertex algebra,
we will need a complicated set of data and conditions stated below.
However, in the context of moduli theory,
there are often natural choices for this data,
as will be discussed in Example~\ref{eg-lin-cat} below,
which the reader might find helpful.

\begin{assumption}
    \label{asn-h-sp}
    Suppose that we have the following.

    \begin{enumerate}
        \item 
            \label{itm-h-sp}
            Let $X$ be a commutative $H$-space,
            as in Definition~\ref{def-h-sp}, with structure map
            \[
                {\oplus} \colon X \times X \longrightarrow X,
            \]
            and unit element $0 \in X$.
            Suppose that there is an $H$-action of $\BU (1)$ on $X$, denoted by
            \[
                {\odot} \colon \BU (1) \times X \longrightarrow X,
            \]
            that respects the $H$-space structure on $X$.
    \end{enumerate}
    A typical example of this is when
    $X$ is the topological realization of a moduli stack,
    as will be discussed in Example~\ref{eg-lin-cat} below.

    We assume the following conditions and extra data.

    \begin{enumerate}[resume]
        \item 
            \label{itm-theta}
            We have a class in topological $K$-theory,
            \[
                \Theta \in K (X \times X),
            \]
            of $\upU (1)$-weight $(1, -1)$, as in Definition~\ref{def-bu1-wt}.
            Write
            \[
                \Theta_{\alpha, \> \beta} =
                \Theta |_{ X_\alpha \times X_\beta }
                \in K ( X_\alpha \times X_\beta )
            \]
            for $\alpha, \beta \in \pi_0 (X)$.
            Define a map $\chi \colon \pi_0 (X) \times \pi_0 (X) \to \bbZ$ by
            \[
                \chi (\alpha, \beta) = \rank \Theta_{\alpha, \> \beta}\,.
            \]

        \item
            \label{itm-theta-linear}
            For all $\alpha, \beta, \gamma \in \pi_0 (X)$, we have
            \begin{align}
                \label{eq-theta-linear-first}
                ({\oplus}_{\alpha, \> \beta} \times \id_{X_\gamma})^*
                (\Theta_{\alpha + \beta, \> \gamma})
                & =
                \pr_{1,3}^* (\Theta_{\alpha, \> \gamma}) +
                \pr_{2,3}^* (\Theta_{\beta, \> \gamma}), \\
                \label{eq-theta-linear-second}
                (\id_{X_\alpha} \times {\oplus}_{\beta, \> \gamma})^*
                (\Theta_{\alpha, \> \beta + \gamma})
                & =
                \pr_{1,2}^* (\Theta_{\alpha, \> \beta}) +
                \pr_{1,3}^* (\Theta_{\alpha, \> \gamma})
            \end{align}
            in $K (X_\alpha \times X_\beta \times X_\gamma)$,
            where the $\pr_{i, \> j}$ are projections from $X_\alpha \times X_\beta \times X_\gamma$
            to its factors.
            In particular, we have
            \begin{align}
                \chi (\alpha + \beta, \gamma)
                & = \chi (\alpha, \gamma) + \chi (\beta, \gamma), \\
                \chi (\alpha, \beta + \gamma)
                & = \chi (\alpha, \beta) + \chi (\alpha, \gamma)
            \end{align}
            for all $\alpha, \beta, \gamma \in \pi_0 (X)$.

        \item
            \label{itm-theta-sym}
            For all $\alpha, \beta \in \pi_0 (X)$, we have
            \begin{equation}
                \label{eq-theta-sym}
                \sigma \sss{\alpha, \> \beta}{*} ( \Theta_{\beta, \> \alpha} ) =
                \Theta \sss{\alpha, \> \beta}{\vee}
            \end{equation}
            in $K (X_\alpha \times X_\beta)$,
            where $\sigma_{\alpha, \> \beta} \colon X_\alpha \times X_\beta \to
            X_\beta \times X_\alpha$
            is the map exchanging the two factors.
            In particular, we have
            \begin{equation}
                \chi (\alpha, \beta) = \chi (\beta, \alpha)
            \end{equation}
            for all $\alpha, \beta \in \pi_0 (X)$.

        \item
            \label{itm-chi-check}
            We have a map
            \[
                \hat{\chi} \colon \pi_0 (X) \longrightarrow \bbZ,
            \]
            satisfying
            \begin{equation}
                \hat{\chi} (\alpha + \beta) =
                \hat{\chi} (\alpha) + \hat{\chi} (\beta) + 2 \chi (\alpha, \beta)
            \end{equation}
            for all $\alpha, \beta \in \pi_0 (X)$.
            
        \item 
            \label{itm-epsilon}
            We have a map
            \[
                \varepsilon \colon \pi_0 (X) \times \pi_0 (X)
                \longrightarrow \{ 1, -1 \},
            \]
            satisfying
            \begin{align}
                \varepsilon (0, 0) & = 1,
                \\
                \varepsilon (\alpha, \beta) \cdot \varepsilon (\beta, \alpha)
                & = (-1)^{\chi (\alpha, \> \beta) + \hat{\chi} (\alpha) \hat{\chi} (\beta)},
                \\
                \varepsilon (\alpha, \beta) \cdot \varepsilon (\alpha + \beta, \gamma)
                & = \varepsilon (\alpha, \beta + \gamma) \cdot \varepsilon (\beta, \gamma)
                \numberthis \label{eq-epsilon-assoc}
            \end{align}
            for all $\alpha, \beta, \gamma \in \pi_0 (X)$.
    \end{enumerate}
    We introduce the following notations.
    
    \begin{itemize}
        \item 
            Define a shift in grading
            \begin{equation}
                \hat{H}_* (X_\alpha; \bbQ) =
                H_{* - \hat{\chi} (\alpha)} (X_\alpha; \bbQ),
            \end{equation}
            and write
            \begin{equation}
                \hat{H}_* (X; \bbQ) =
                \bigoplus_{\alpha \in \pi_0 (X)} \hat{H}_* (X_\alpha; \bbQ).
            \end{equation}
            The reason for doing this shift is that in the following,
            the space $\hat{H}_* (X; \bbQ)$ will acquire
            the structure of a graded vertex algebra,
            and only the shifted grading is compatible
            with the vertex algebra structure.
        \item 
            Write
            \begin{equation}
                \varepsilon (\alpha_1, \dotsc, \alpha_n) =
                \varepsilon (\alpha_1, \alpha_2) \cdot
                \varepsilon (\alpha_1 + \alpha_2, \alpha_3) \cdots
                \varepsilon (\alpha_1 + \cdots + \alpha_{n-1}, \alpha_n)
            \end{equation}
            for any $n > 0$ and $\alpha_1, \dotsc, \alpha_n \in \pi_0 (X)$.
    \end{itemize}
\end{assumption}

\begin{remark}
    Note that we could have written~%
    \eqref{eq-theta-linear-first}--\eqref{eq-theta-sym}
    globally on $X \times X \times X$ or $X \times X$,
    without the subscripts $\alpha, \beta$, etc.,
    which is perhaps the more natural way to write them.
    However, we choose to include the subscripts to help with understanding.
    The same remark applies to many axioms below as well.
\end{remark}

\begin{remark}
    In Assumption~\ref{asn-h-sp}~\ref{itm-chi-check},
    an elementary argument shows that the map $\hat{\chi}$ must have the form
    \begin{equation}
        \hat{\chi} (\alpha) = \chi (\alpha, \alpha) + f (\alpha),
    \end{equation}
    where $f \colon \pi_0 (X) \to \bbZ$ is a monoid homomorphism.
    In fact, we will take $f = 0$ in all our examples,
    so that $\hat{\chi} (\alpha) = \chi (\alpha, \alpha)$ for all $\alpha \in \pi_0 (X)$.
\end{remark}

Before we state the vertex algebra construction,
we set up some notations regarding the homology of $\BU (1)$.

\begin{definition}
    \label{def-h-bu-1}
    We define isomorphisms
    \[
        H^* (\BU (1); \bbQ) \simeq \bbQ [T], \qquad
        H_* (\BU (1); \bbQ) \simeq \bbQ [t],
    \]
    where the first one is an isomorphism of graded $\bbQ$-algebras,
    and the second one is an isomorphism of graded $\bbQ$-vector spaces,
    with $\deg T = \deg t = 2$.
    Namely, let $T \in H^2 (\BU (1); \bbQ)$ be the universal first Chern class,
    i.e.~the first Chern class of the universal line bundle.
    This extends to a $\bbQ$-algebra isomorphism
    $H^* (\BU (1); \bbQ) \simeq \bbQ [T]$ as required.
    For each $k \geq 0$, define $t^k \in H_{2k} (\BU (1); \bbQ)$
    to be the unique element satisfying
    \begin{equation}
        T^k \cdot t^k = k!.
    \end{equation}
    This choice has the property that,
    for any polynomials $p (T) \in H^* (\BU (1); \bbQ)$
    and $q (t) \in H_* (\BU (1); \bbQ)$, we have
    \begin{equation}
        q (t) \cap P (T) =
        P \Bigl( \frac{\partial}{\partial t} \Bigr) \, q (t).
    \end{equation}
\end{definition}

We can now state Joyce's vertex algebra construction.
Although our setting is different,
Joyce's proof applies to our case without any substantial modification.
We will present a proof of this in \S\ref{sect-proof-va}.

\begin{theorem}
    \label{thm-va}
    Under Assumption~\textnormal{\ref{asn-h-sp},}
    there is a graded vertex algebra structure
    on the graded vector space $\hat{H}_* (X; \bbQ)$,
    defined as follows.

    \begin{itemize}
        \item 
            The unit element $1 \in \hat{H}_0 (X; \bbQ)$
            is given by the class 
            \begin{equation}
                \label{eq-constr-unit}
                1 = [0] \in H_0 (X_0; \bbQ),
            \end{equation}
            where $0 \in X$ is the unit element of $X$.
        
        \item 
            For $A \in H_a (X_\alpha; \bbQ)$, we have
            \begin{equation}
                \label{eq-constr-transl}
                D (A) = \odot_* (t \boxtimes A) \in H_{a+2} (X_\alpha; \bbQ),
            \end{equation}
            where $t \in H_2 (\BU (1); \bbQ)$ is as in
            Definition~\textnormal{\ref{def-h-bu-1}}.
            
        \item 
            For $A \in H_a (X_\alpha; \bbQ)$ and $B \in H_b (X_\beta; \bbQ)$, we have
            \begin{multline}
                \label{eq-constr-va}
                Y (A, z) \, B =
                (-1)^{a \hat{\chi} (\beta)} 
                \varepsilon (\alpha, \beta) \cdot
                z^{\chi (\alpha, \> \beta)} \cdot {} \\
                (\oplus_{\alpha, \> \beta})_* \circ
                (\upe^{z D} \otimes \id) \bigl(
                (A \boxtimes B) \cap c_{1/z} (\Theta_{\alpha, \> \beta}) \bigr).
            \end{multline}
    \end{itemize}
\end{theorem}

Note that in~\eqref{eq-constr-va},
only finitely many terms of $c_{1/z} (\Theta_{\alpha, \> \beta})$
are involved, as higher terms will vanish after intersecting with $A \boxtimes B$.
This ensures that the expression is well-defined
and lies in $\hat{H}_* (X; \bbQ) \llparen z \rrparen$.

Next, we formulate a situation where the Lie algebra
associated to the vertex algebra $\hat{H}_* (X; \bbQ)$,
as in Theorem~\ref{thm-lie-alg},
can be endowed with a geometric meaning.

\begin{assumption}
    \label{asn-rat-triv}
    Under Assumption~\ref{asn-h-sp}, we further assume the following.

    \begin{enumerate}
        \item 
            We are given a map
            \[
                \pi \colon X \longrightarrow X^\pl ,
            \]
            inducing a bijection on connected components.
            
        \item
            \label{itm-rat-triv}
            For each $\alpha \in \pi_0 (X)$, one of the following is true.
            \begin{enumerate}[label={(\alph*)}]
                \item
                    \label{itm-rat-triv-a}
                    The $\BU (1)$-action on $X_\alpha$ is trivial,
                    and the restriction $\pi_\alpha \colon X_\alpha \to X^\pl_\alpha$
                    is an equivalence.
                \item
                    \label{itm-rat-triv-b}
                    The $\BU (1)$-action on $X_\alpha$
                    establishes the restriction
                    $\pi_\alpha \colon X_\alpha \to X^\pl_\alpha$
                    as a principal $\BU (1)$-bundle.
                    Moreover, there exists a $K$-theory class $E_\alpha \in K (X_\alpha)$
                    of non-zero rank and non-zero $\upU (1)$-weight.
            \end{enumerate}
            This is rather technical,
            and see below for more explanations.
    \end{enumerate}
    We introduce the following notations.
    
    \begin{itemize}
        \item 
            Define a shift in grading
            \begin{equation}
                \check{H}_* (X^\pl_\alpha; \bbQ) =
                H_{* + 2 - \hat{\chi} (\alpha)} (X^\pl_\alpha; \bbQ),
            \end{equation}
            and write
            \begin{equation}
                \check{H}_* (X^\pl; \bbQ) =
                \bigoplus_{\alpha \in \pi_0 (X)} \check{H}_* (X^\pl_\alpha; \bbQ).
            \end{equation}
    \end{itemize}
\end{assumption}

Here, we think of the space $X^\pl$ roughly as a homotopy quotient $X / \BU (1)$,
or the space obtained from a moduli stack
by removing scalar automorphisms of objects,
i.e.~the $\Gm$-\emph{rigidification} of a moduli stack.
See Example~\ref{eg-lin-cat} below for more details.

The technical condition in Assumption~\ref{asn-rat-triv}~\ref{itm-rat-triv}
comes from the requirement that $\pi \colon X \to X^\pl$
should be \emph{rationally trivial},
in the sense of~\cite[Definition~2.26]{JoyceHall},
which makes Theorem~\ref{thm-pl-lie-alg} below possible.
This assumption ensures that,
from the viewpoint of rational homology $H_* (-; \bbQ)$,
the map $\pi$ is an isomorphism on some components,
and a trivial $\BU (1)$-bundle on the other components.

\begin{theorem}
    \label{thm-pl-lie-alg}
    Under Assumptions~\textnormal{\ref{asn-h-sp}}
    and~\textnormal{\ref{asn-rat-triv},} we have an isomorphism
    \begin{equation}
        H_* (X; \bbQ) / D (H_{*-2} (X; \bbQ)) \simeq
        H_* (X^\pl; \bbQ).
    \end{equation}
    Therefore, the graded vertex algebra structure on $\hat{H}_* (X; \bbQ)$
    induces a graded Lie algebra structure on $\check{H}_* (X^\pl; \bbQ)$.
\end{theorem}

\begin{proof}
    We follow ideas from~\cite{JoyceHall}.
    Fix $\alpha \in \pi_0 (X)$, and write $V_\alpha = H_* (X_\alpha; \bbQ)$.
    We need to show that $V_\alpha / D (V_\alpha) \simeq
    H_* (X_\alpha^\pl; \bbQ)$.
    We split into two cases.

    If Assumption~\ref{asn-rat-triv}~\ref{itm-rat-triv}~\ref{itm-rat-triv-a}
    holds, then using the notations of Theorem~\ref{thm-va},
    we have $\odot_* (t \boxtimes A) = 0$ for all $A \in H_* (X_\alpha; \bbQ)$.
    This means that $D = 0$, and
    $V_\alpha / D (V_\alpha) \simeq V_\alpha \simeq
    H_* (X_\alpha^\pl; \bbQ)$.

    If Assumption~\ref{asn-rat-triv}~\ref{itm-rat-triv}~\ref{itm-rat-triv-b}
    holds, let $L = \det (E_\alpha)$, where $E_\alpha$ is the class
    given by the assumption.
    Then $L$ has $\upU (1)$-weight $d \neq 0$,
    where $d$ is the product
    of the rank and the $\upU (1)$-weight of $E_\alpha$\,.
    We view $L$ as a line bundle over $X_\alpha$\,,
    and let $L^\circ = L \setminus X_\alpha$
    be the complement of the zero section,
    and let $p \colon L^\circ \to X_\alpha$ denote the projection.
    Consider the diagram
    \[ \begin{tikzcd}
        L^\circ \times \BU (1)
        \ar[r, "p \times \id"]
        \ar[d, "\pr_1"]
        & X_\alpha \times \BU (1)
        \ar[r, "\odot"]
        \ar[d, "\pr_1"]
        & X_\alpha \ar[d, "\pi_\alpha"]
        \\ L^\circ \ar[r, "p"]
        & X_\alpha \ar[r, "\pi_\alpha"]
        & X^\pl_\alpha \rlap{ ,}
    \end{tikzcd} \]
    where $\pi_\alpha$ is the restriction of the map $\pi$
    in Assumption~\ref{asn-rat-triv} to $X_\alpha$\,.
    We see that both squares are homotopy pullbacks.
    The composition of the bottom maps is a principal $\upB \bbZ_{|d|}$-bundle,
    and induces an isomorphism on rational homology.
    Moreover, the pullback of $\pi_\alpha$ along this map
    is a trivial $\BU (1)$-bundle.
    Therefore, we have an isomorphism
    \begin{align*}
        V_\alpha & \simeq
        H_* (X_\alpha^\pl; \bbQ) \otimes H_* (\BU (1); \bbQ) \\
        & \simeq
        H_* (X_\alpha^\pl; \bbQ) [t],
        \numberthis
    \end{align*}
    and using this isomorphism,
    the operator $D$ sends is multiplication by $t$,
    so that $V_\alpha / D (V_\alpha) \simeq
    H_* (X_\alpha^\pl; \bbQ)$.
\end{proof}

\begin{example}
    \label{eg-va-triv}
    A somewhat trivial example satisfying Assumption~\ref{asn-h-sp} is given as follows.
    Let $X$ be any commutative $H$-space with an $H$-action by $\BU (1)$.
    Set $\Theta = 0$, $\hat{\xi} (\alpha) = 0$,
    and $\varepsilon (\alpha, \beta) = 1$ for all $\alpha, \beta \in \pi_0 (X)$.
    Then, Theorem~\ref{thm-va} produces a \emph{commutative vertex algebra},
    in the sense of~\cite[\S1.4]{FrenkelBenZvi2004},
    which is the same as a commutative algebra equipped with a derivation.
\end{example}

\begin{example}
    \label{eg-lin-cat}
    We describe an important class of examples satisfying Assumption~\ref{asn-h-sp},
    arising from moduli stacks of objects in $\bbC$-linear categories.
    
    Let $\calM^+$ be a $\bbC$-linear stack, as in Definition~\ref{def-lin-st},
    or a $\bbC$-linear $\infty$-stack, as in Definition~\ref{def-lin-infty-st}.
    Let $\calM = (\calM^+)^\simeq$ be the underlying stack of $\calM$,
    as a functor taking values in groupoids or $\infty$-groupoids,
    obtained from $\calM^+$ by discarding the non-invertible arrows.
    We view $\calM$ as a \emph{moduli stack}
    of objects of the $\bbC$-linear category
    or $\bbC$-linear $\infty$-category
    $\calA = \calM^+ (\Spec \bbC)$.

    For Assumption~\ref{asn-h-sp}~\ref{itm-h-sp}, take 
    \[
        X = |\calM|,
    \]
    the topological realization of $\calM$.
    The $H$-space structure map, $\oplus \colon X \times X \to X$,
    is given by the direct sum morphism 
    \[
        \oplus \colon \calM \times \calM \to \calM.
    \]
    The $\BU (1)$-action on $X$ is given by the $[*/\Gm]$-action on $\calM$
    via scalar automorphisms,
    \[
        \odot \colon [*/\Gm] \times \calM \to \calM,
    \]
    as we have an equivalence $|[*/\Gm]| \simeq \BU (1)$
    respecting the group structures.

    For Assumption~\ref{asn-h-sp}~\ref{itm-theta},
    we assume that there is a perfect complex
    \[
        \Thetaul \in \cat{Perf} (\calM \times \calM),
    \]
    of $\Gm$-weight $(1, -1)$,
    giving a class in $K$-theory, $\Theta \in K (X \times X)$,
    as in Definition~\ref{def-perf}.
    Here, $\Gm$-weights are defined analogously to
    $\upU (1)$-weights in Definition~\ref{def-bu1-wt}.
    Informally, the complex $\Thetaul$ assigns to each pair of objects
    $E, F \in \calA$ a chain complex of vector spaces $\Thetaul (E, F)$,
    satisfying the conditions below.

    For Assumption~\ref{asn-h-sp}~\ref{itm-theta-linear},
    we require that $\Thetaul$ satisfies the relations~%
    \eqref{eq-theta-linear-first}--\eqref{eq-theta-linear-second}
    in $K$-theory.
    In fact, in all our examples, these relations will be satisfied
    as isomorphisms of perfect complexes.
    In the latter case, \eqref{eq-theta-linear-first}--\eqref{eq-theta-linear-second}
    translate informally to the relations
    \begin{align}
        \Thetaul (E \oplus F, G) & \simeq \Thetaul (E, G) \oplus \Thetaul (F, G), \\
        \Thetaul (E, F \oplus G) & \simeq \Thetaul (E, F) \oplus \Thetaul (E, G)
    \end{align}
    for objects $E, F, G \in \calA$.

    For Assumption~\ref{asn-h-sp}~\ref{itm-theta-sym},
    we require that $\Thetaul$ satisfies~\eqref{eq-theta-sym}
    in $K$-theory.
    In examples, this is satisfied either as an isomorphism of perfect complexes,
    or as an isomorphism after shifting by an even number,
    \begin{equation}
        \sigma \sss{\alpha, \> \beta}{*} ( \Thetaul_{\beta, \> \alpha} ) \simeq
        \Thetaul \sss{\alpha, \> \beta}{\vee} [2n],
    \end{equation}
    where $n \in \bbZ$. More remarks on this can be found in~\cite[Remark~3.3]{JoyceHall}.

    For Assumption~\ref{asn-h-sp}~\ref{itm-chi-check}, we often take
    \begin{equation}
        \hat{\chi} (\alpha) = \chi (\alpha, \alpha),
    \end{equation}
    and for Assumption~\ref{asn-h-sp}~\ref{itm-epsilon}, we often take
    \begin{equation}
        \varepsilon (\alpha, \beta) = (-1)^{\chi_\calA (\alpha, \> \beta)},
    \end{equation}
    for a system of integers $\chi_\calA (\alpha, \beta)$,
    bilinear in $\alpha, \beta$, and satisfying
    $\chi (\alpha, \beta) = \chi_\calA (\alpha, \beta) + \chi_\calA (\beta, \alpha)$.
    See also~\cite[Remark~3.4]{JoyceHall}.

    A typical choice of this set of data is to take
    \begin{equation}
        \label{eq-choice-theta-ul}
        \Thetaul = \calExt^\vee \oplus \sigma^* (\calExt)
    \end{equation}
    when there exists an Ext-complex, $\calExt \in \cat{Perf} (\calM \times \calM)$,
    assigning to $E, F \in \calA$ the complex $\calExt (E, F)$
    consisting of the Ext groups $\Ext^i (E, F)$.
    In this case, \eqref{eq-theta-linear-first}--\eqref{eq-theta-sym}
    are satisfied as isomorphisms of perfect complexes,
    and there is a natural choice of
    $\varepsilon (\alpha, \beta) = (-1)^{\chi_\calA (\alpha, \> \beta)}$,
    by taking
    \begin{equation}
        \chi_\calA (\alpha, \beta) = \rank \calExt_{\alpha, \> \beta} \, .
    \end{equation}

    Now, one can apply Theorem~\ref{thm-va} 
    to obtain a graded vertex algebra structure on the homology
    $\hat{H}_* (\calM; \bbQ)$.

    In addition, if Assumption~\ref{asn-rat-triv} is also satisfied,
    then one can obtain a graded Lie algebra structure on the homology
    $\check{H}_* (\calM; \bbQ)$
\end{example}

\begin{remark}
    \label{rmk-rat-triv}
    In Example~\ref{eg-lin-cat},
    Assumption~\ref{asn-rat-triv} is often satisfied
    for moduli stacks for a $\bbC$-linear quasi-abelian or abelian category,
    but not for a derived category.
    The reason is often that the class $\alpha = 0$ does not satisfy
    Assumption~\ref{asn-rat-triv}~\ref{itm-rat-triv} in the derived case.
    However, one can often find a subspace
    $\calM_{>0} \subset \calM$ consisting of components
    that contain non-zero objects of a quasi-abelian or abelian subcategory.
    We are often able to show that these components satisfy
    Assumption~\ref{asn-rat-triv}~\ref{itm-rat-triv}~\ref{itm-rat-triv-b},
    so we obtain a graded Lie subalgebra
    \begin{equation}
        \label{eq-mpl-gt0}
        \check{H}_* (\calM^\pl_{>0}; \bbQ) \subset
        \hat{H}_{*+2} (\calM; \bbQ) / D (\hat{H}_* (\calM; \bbQ))
    \end{equation}
    and this subalgebra is sufficient for working with enumerative invariants.

    Another workaround for this problem is given by Upmeier~\cite[Theorem~1.4]{Upmeier2021},
    who constructed a graded Lie algebra structure on
    $\check{H}_* (\calM^\pl; \bbQ)$ directly,
    without using Assumption~\ref{asn-rat-triv}.
    Thus, when this assumption is not satisfied,
    $\check{H}_* (\calM^\pl; \bbQ)$ and
    $\hat{H}_{*+2} (\calM; \bbQ) / D (\hat{H}_* (\calM; \bbQ))$
    are two different graded Lie algebras.
    We do not take this approach in the present work,
    but we expect that certain constructions below
    involving this Lie algebra
    should be able to be done using Upmeier's construction as well.
\end{remark}

\subsection{Twisted modules from moduli stacks}
\label{sect-cons-tw-mod}

In the following, we construct twisted modules,
introduced in Definition~\ref{def-tw-module},
for the vertex algebras constructed in \S\ref{sect-cons-va}.
This will need extra data that may look rather complicated at the first glance,
but we explain how a choice of this data
can be naturally extracted from the information of a moduli stack.

Before stating the main construction,
we introduce a piece of notation.

\begin{definition}
    \label{def-z2-bu-1}
    Let $\bbZ_2$ act on the topological group $\BU (1)$
    by assigning to the non-trivial element of $\bbZ_2$
    the map $\BU (1) \to \BU (1)$ induced by the map
    $\upU (1) \to \upU (1)$, $\upe^{\upi t} \mapsto \upe^{-\upi t}$.
    In other words, $\bbZ_2$ acts by taking the dual line bundle.
    Let $\bbZ_2 \rtimes \BU (1)$ be the semi-direct product
    associated to this action, as a topological group.
\end{definition}

\begin{assumption}
    \label{asn-h-sp-mod}
    Under Assumption~\ref{asn-h-sp},
    we further assume the following.
    
    \begin{enumerate}
        \item
            \label{itm-h-sp-invol}
            We have an involution
            \[
                I \colon X \longrightarrow X
            \]
            of $X$, giving an involutive commutative $H$-space,
            in the sense of Definition~\ref{def-h-sp-invol}.
            Denote the involution on $\pi_0 (X)$ by
            $(-)^\vee \colon \pi_0 (X) \to \pi_0 (X)$.
            Write
            \[
                I_\alpha = I|_{X_\alpha} \colon
                X_\alpha \longrightarrow X_{\smash{\alpha^\vee}}
            \]
            for $\alpha \in \pi_0 (X)$.
        
            Suppose that the involution of $X$ and the $H$-action by $\BU (1)$
            together give an $H$-action by the group
            $\bbZ_2 \rtimes \BU (1)$ in Definition~\ref{def-z2-bu-1}.
            
        \item 
            \label{itm-h-sp-mod}
            We have an involutive module $X^\sd$ over $X$,
            with structure map 
            \[
                {\oplus}^\sd \colon X \times X^\sd \to X^\sd.
            \]
            Denote the action of $\pi_0 (X)$ on $\pi_0 (X^\sd)$ by
            $(\alpha, \theta) \mapsto \bar{\alpha} + \theta$.
            This gives $\pi_0 (X^\sd)$ the structure of
            an involutive $\pi_0 (X)$-module,
            in the sense of Definition~\ref{def-h-sp-invol}.
            
        \item
            \label{itm-theta-equivar}
            Consider the map
            \[
                \delta = (I \times I) \circ \sigma \colon
                X \times X \longrightarrow
                X \times X,
            \]
            where $\sigma$ exchanges the two factors,
            and $I \colon X \to X$ is the involution.
            Then $\delta$ defines an $H$-action by $\bbZ_2$ on $X \times X$.
            We assume that the class $\Theta \in K (X \times X)$
            in Assumption~\ref{asn-h-sp}~\ref{itm-theta} is $\bbZ_2$-invariant, that is, we have
            \begin{equation}
                \delta^* (\Theta) = \Theta.
            \end{equation}
            We require that
            \begin{align}
                \chi (\alpha, \beta) & = \chi (\beta^\vee, \alpha^\vee), \\
                \hat{\chi} (\alpha) & = \hat{\chi} (\alpha^\vee)
            \end{align}
            for all $\alpha, \beta \in \pi_0 (X)$,
            where the first equation is implied by the condition on $\Theta$.
            
        \item 
            \label{itm-theta-dots}
            We have classes
            \begin{align*}
                \dot{\Theta} & \in K (X \times X^\sd), \\
                \ddot{\Theta} & \in K (X),
            \end{align*}
            of $\upU (1)$-weights $1$ and $2$, respectively.
            For $\alpha \in \pi_0 (X)$ and $\theta \in \pi_0 (X^\sd)$, write
            \begin{alignat*}{2}
                \dot{\Theta}_{\alpha, \> \theta} & =
                \dot{\Theta} |_{X_\alpha \times X^\sd_\theta}
                && \in K (X_\alpha \times X^\sd_\theta), \\
                \ddot{\Theta}_{\alpha} & =
                \ddot{\Theta} |_{X_\alpha}
                && \in K (X_\alpha).
            \end{alignat*}
            Write
            \begin{align}
                \dot{\chi} (\alpha, \theta) & =
                \rank \dot{\Theta}_{\alpha, \> \theta} \, , \\
                \ddot{\chi} (\alpha) & =
                \rank \ddot{\Theta}_{\alpha} \, , \\
                \chi^\sd (\alpha, \theta) & =
                \dot{\chi} (\alpha, \theta) +
                \ddot{\chi} (\alpha)
            \end{align}
            for all $\alpha \in \pi_0 (X)$ and $\theta \in \pi_0 (X^\sd)$.
            
        \item
            \label{itm-theta-dots-linear}
            For all $\alpha, \beta \in \pi_0 (X)$ and $\theta \in \pi_0 (X^\sd)$, we have
            \begin{align}
                \label{eq-theta-dot-linear-first}
                ({\oplus}_{\alpha, \> \beta} \times \id_{X\sss{\theta}{\sd}})^*
                (\dot{\Theta}_{\alpha + \beta, \> \theta})
                & =
                \pr_{1, 3}^* (\dot{\Theta}_{\alpha, \> \theta}) 
                + \pr_{2, 3}^* (\dot{\Theta}_{\beta, \> \theta}),
                \\
                (\id_{X_\alpha} \times {\oplus}\sss{\beta, \> \theta}{\sd})^*
                (\dot{\Theta}_{\alpha, \> \bar{\beta} + \theta})
                & =
                \pr_{1, 2}^* (\Theta_{\alpha, \> \beta}) + {} \notag \\*[-.5ex]
                \label{eq-theta-dot-linear-second}
                & \hspace{2em}
                {\pr_{1, 2}^*} \circ (\id_{X_\alpha} \times I_\beta)^*
                (\Theta_{\smash{\alpha, \> \beta^\vee}})
                + \pr_{1, 3}^* (\dot{\Theta}_{\alpha, \> \theta}),
                \\
                \label{eq-theta-ddot-linear}
                ({\oplus}_{\alpha, \> \beta})^* (\ddot{\Theta}_{\alpha + \beta})
                & =
                \pr_1^* (\ddot{\Theta}_\alpha) + \pr_2^* (\ddot{\Theta}_\beta)
                + (\id_{X_\alpha} \times I_\beta)^* (\Theta_{\smash{\alpha, \> \beta^\vee}}),
            \end{align}
            where~\eqref{eq-theta-dot-linear-first}--\eqref{eq-theta-dot-linear-second}
            are in $K (X_\alpha \times X_\beta \times X\sss{\theta}{\sd})$,
            and~\eqref{eq-theta-ddot-linear} in
            $K (X_\alpha \times X_\beta)$.
            The maps $\pr_i$ and $\pr_{i, \> j}$ are projections onto factors.
            In particular, we have
            \begin{align}
                \dot{\chi} (\alpha + \beta, \theta)
                & = \dot{\chi} (\alpha, \theta) + \dot{\chi} (\beta, \theta), \\
                \dot{\chi} (\alpha, \bar{\beta} + \theta)
                & = \chi (\alpha, \beta) + \chi (\alpha, \beta^\vee)
                + \dot{\chi} (\alpha, \theta), \\
                \ddot{\chi} (\alpha + \beta)
                & = \ddot{\chi} (\alpha) + \ddot{\chi} (\beta) + \chi (\alpha, \beta^\vee)
            \end{align}
            for all $\alpha, \beta \in \pi_0 (X)$ and $\theta \in \pi_0 (X^\sd)$.

        \item
            \label{itm-theta-dots-sym}
            For all $\alpha \in \pi_0 (X)$ and $\theta \in \pi_0 (X^\sd)$, we have
            \begin{align}
                (I_\alpha \times \id_{X_\theta^\sd})^* (\dot{\Theta}_{\smash{\alpha^\vee, \> \theta}})
                & = \dot{\Theta} \sss{\alpha, \> \theta}{\vee}, \\
                I_\alpha^* (\ddot{\Theta}_{\smash{\alpha^\vee}})
                & = \ddot{\Theta}_\alpha^\vee.
            \end{align}
            In particular, we have
            \begin{align}
                \dot{\chi} (\alpha, \theta)
                & = \dot{\chi} (\alpha^\vee, \theta), \\
                \ddot{\chi} (\alpha)
                & = \ddot{\chi} (\alpha^\vee)
            \end{align}
            for all $\alpha \in \pi_0 (X)$ and $\theta \in \pi_0 (X^\sd)$.

        \item
            \label{itm-chi-ring}
            We have a map
            \[
                \ring{\chi} \colon \pi_0 (X^\sd)
                \longrightarrow \bbZ,
            \]
            satisfying
            \begin{equation}
                \ring{\chi} (\bar{\alpha} + \theta) =
                \ring{\chi} (\theta) + \hat{\chi} (\alpha) + 
                2 \chi^\sd (\alpha, \theta)
            \end{equation}
            for all $\alpha \in \pi_0 (X)$ and
            $\theta \in \pi_0 (X^\sd)$.
            
        \item 
            \label{itm-epsilon-sd}
            We have a map
            \[
                \varepsilon^\sd \colon \pi_0 (X) \times \pi_0 (X^\sd)
                \longrightarrow \{ 1, -1 \},
            \]
            satisfying
            \begin{align}
                \varepsilon^\sd (\alpha, \theta) \cdot
                \varepsilon^\sd (\alpha^\vee , \theta)
                & = (-1)^{\chi^\sd (\alpha, \theta) +
                    \hat{\chi} (\alpha) \ring{\chi} (\theta)},
                \\
                \label{eq-epsilon-sd-assoc}
                \varepsilon (\alpha, \beta) \cdot
                \varepsilon^\sd (\alpha + \beta, \theta)
                & = \varepsilon^\sd (\beta, \theta) \cdot
                \varepsilon^\sd (\alpha, \bar{\beta} + \theta)
            \end{align}
            for all $\alpha, \beta \in \pi_0 (X)$ and $\theta \in \pi_0 (X^\sd)$.
    \end{enumerate}
    We introduce the following notations.
    
    \begin{itemize}
        \item 
            Define a shift in grading
            \begin{equation}
                \ring{H}_* (X^\sd_\theta; \bbQ) =
                H_{* - \ring{\chi} (\theta)} (X^\sd_\theta; \bbQ),
            \end{equation}
            and write
            \begin{equation}
                \ring{H}_* (X^\sd; \bbQ) =
                \bigoplus_{\theta \in \pi_0 (X^\sd)} \ring{H}_* (X^\sd_\theta; \bbQ).
            \end{equation}
        \item 
            Write
            \begin{equation}
                \varepsilon^\sd (\alpha_1, \> \dotsc, \> \alpha_n, \> \theta) =
                \varepsilon (\alpha_1, \dotsc, \alpha_n) \cdot
                \varepsilon^\sd (\alpha_1 + \cdots + \alpha_n, \> \theta)
            \end{equation}
            for any $n \geq 0$, $\alpha_1, \dotsc, \alpha_n \in \pi_0 (X)$,
            and $\theta \in \pi_0 (X^\sd)$.
    \end{itemize}
\end{assumption}

\begin{remark}
    From Assumption~\ref{asn-h-sp-mod}~\ref{itm-theta-dots-linear}--\ref{itm-theta-dots-sym},
    one can deduce the relations
    \begin{align}
        \chi^\sd (\alpha, \theta) 
        & = \chi^\sd (\alpha^\vee, \theta), \\
        \chi (\alpha, \beta) + \chi^\sd (\alpha + \beta, \theta) 
        & = \chi^\sd (\alpha, \bar{\beta} + \theta) +
        \chi^\sd (\beta, \theta)
    \end{align}
    for $\alpha, \beta \in \pi_0 (X)$ and $\theta \in \pi_0 (X^\sd)$.
\end{remark}

\begin{theorem}
    \label{thm-tw-mod}
    Under Assumption~\textnormal{\ref{asn-h-sp-mod},}
    the graded vertex algebra $\hat{H}_* (X; \bbQ)$
    is equipped with an involution, 
    in the sense of Definition~\textnormal{\ref{def-va-invol},} given by
    \begin{equation}
        A^\vee = I_* (A)
    \end{equation}
    for $A \in \hat{H}_* (X; \bbQ)$,
    where $I \colon X \to X$ is the involution.
    
    The graded vector space
    $\ring{H}_* (X^\sd; \bbQ)$ has the structure of
    an involutive graded twisted module
    for the involutive graded vertex algebra $\hat{H}_* (X; \bbQ),$
    in the sense of Definition~\textnormal{\ref{def-tw-mod-invol},}
    with structure map given by
    \begin{multline}
        \label{eq-def-ysd}
        Y^\sd (A, z) \, M =
        (-1)^{a \ring{\chi} (\theta)}
        \varepsilon^\sd (\alpha, \theta) \cdot
        z^{\dot{\chi} (\alpha, \> \theta)} \cdot
        (2z)^{\ddot{\chi} (\alpha)} \cdot {} \\
        (\oplus_{\alpha, \> \theta}^\sd)_* \circ
        (\upe^{z D} \otimes \id) \bigl(
        (A \boxtimes M) \cap (
            c_{1/z} (\dot{\Theta}_{\alpha, \> \theta}) \cdot c_{1/(2z)} (\ddot{\Theta}_\alpha)
        ) \bigr)
    \end{multline}
    for $A \in H_a (X_\alpha; \bbQ)$ and $M \in H_m (X^\sd_\theta; \bbQ)$.

    In particular, if Assumption~\textnormal{\ref{asn-rat-triv}}
    is satisfied, then $\ring{H}_* (X^\sd; \bbQ)$
    is also a graded twisted module for the involutive graded Lie algebra
    $\check{H}_* (X^\pl; \bbQ)$.
\end{theorem}

The proof will be given in \S\ref{sect-proof-tw-mod}.

\begin{example}
    \label{eg-sd-lin-cat}
    We describe a class of examples satisfying Assumption~\ref{asn-h-sp-mod},
    which comes from moduli stacks of objects in a self-dual category.
    
    Let $\calM^+$ be a self-dual $\bbC$-linear stack,
    as in Definition~\ref{def-sd-lin-st},
    or a self-dual $\bbC$-linear $\infty$-stack,
    as in Definition~\ref{def-sd-lin-infty-st}.
    Let $\calM = (\calM^+)^\simeq$ be the underlying stack of $\calM$,
    as in Example~\ref{eg-lin-cat}.
    We view $\calM$ as a \emph{moduli stack}
    of objects of $\calA = \calM^+ (\Spec \bbC)$,
    and we view the stack $\calM^\sd$,
    defined in Definitions~\ref{def-sd-lin-st} and~\ref{def-sd-lin-infty-st},
    as a \emph{moduli stack of self-dual objects} in $\calA$.

    We assume that we are given data as in Example~\ref{eg-lin-cat},
    so that Assumption~\ref{asn-h-sp} is satisfied,
    where $X = |\calM|$ is the topological realization of $\calM$.

    Assumption~\ref{asn-h-sp-mod}~\ref{itm-h-sp-invol}
    is satisfied,
    with the involution $I \colon X \to X$ induced by the dual morphism
    \[
        (-)^\vee \colon \calM \longrightarrow \calM.
    \]

    For Assumption~\ref{asn-h-sp-mod}~\ref{itm-h-sp-mod}, we take
    \[
        X^\sd = |\calM^\sd|,
    \]
    the topological realization of $\calM^\sd$.
    The $H$-space action
    $\oplus^\sd \colon X \times X^\sd \to X^\sd$ is induced by the morphism
    \begin{align*}
        \oplus^\sd \colon \calM \times \calM^\sd
        & \longrightarrow \calM^\sd, \\
        (E, F) & \longmapsto \bar{E} \oplus F,
    \end{align*}
    as in Example~\ref{eg-oplus-sd}.

    For Assumption~\ref{asn-h-sp-mod}~\ref{itm-theta-equivar},
    we require that the perfect complex $\Thetaul$ in Example~\ref{eg-lin-cat}
    is equivariant with respect to the $\bbZ_2$-action on $\calM \times \calM$
    sending $(E, F)$ to $(F^\vee, E^\vee)$. Informally, this means that
    \begin{equation}
        \Thetaul (E, F) \simeq \Thetaul (F^\vee, E^\vee)
    \end{equation}
    for all $E, F \in \calA$. In particular,
    if $\Thetaul$ is given by an Ext complex, as in Example~\ref{eg-lin-cat},
    then this is automatically satisfied.

    For Assumption~\ref{asn-h-sp-mod}~\ref{itm-theta-dots}, we choose perfect complexes
    \begin{align*}
        \Thetaul[\dot] & \in \cat{Perf} (\calM \times \calM^\sd), \\
        \Thetaul[\ddot] & \in \cat{Perf} (\calM).
    \end{align*}
    A typical choice is to take
    \begin{align}
        \label{eq-choice-theta-dot}
        \Thetaul[\dot]_{\alpha, \> \theta} & =
        (\id_{\calM_\alpha} \times j_\theta)^* (\Theta_{\alpha, \> j (\theta)}), \\
        \label{eq-choice-theta-ddot}
        \Thetaul[\ddot]_\alpha & =
        \bigl( (\id_{\calM_\alpha} , I_\alpha)^* 
        (\Theta_{\smash{\alpha, \> \alpha^\vee}}) \bigr)^{\bbZ_2},
    \end{align}
    where $j \colon \calM^\sd \to \calM$ is the projection,
    sending a self-dual object to its underlying object,
    and $j_\theta = j|_{\smash{\calM^\sd_\theta}} \colon
    \calM^\sd_\theta \to \calM_{\smash{j (\theta)}}$\,,
    where we have assumed the existence of a class $j (\theta)$,
    such that the image of $j_\theta$ lies in $\calM_{\smash{j (\theta)}}$\,.
    The notation $(-)^{\bbZ_2}$ means taking the homotopy $\bbZ_2$-fixed points.
    Informally, we are taking
    $\Thetaul[\ddot] (E) = \Thetaul (E, E^\vee)^{\bbZ_2}$ for $E \in \calA$,
    where the $\bbZ_2$-action comes from
    the $\bbZ_2$-equivariance of the complex $\Thetaul$ discussed above.

    For these specific choices of $\Thetaul[\dot]$ and $\Thetaul[\ddot]$,
    Assumption~\ref{asn-h-sp-mod}~\ref{itm-theta-dots-linear}--\ref{itm-theta-dots-sym}
    follow directly from the assumptions on $\Thetaul$.

    For Assumption~\ref{asn-h-sp-mod}~\ref{itm-chi-ring},
    if we have chosen $\Thetaul[\dot]$ and $\Thetaul[\ddot]$ as above,
    one natural choice of $\ring{\chi}$ is to take
    $\ring{\chi} (\theta) = \ddot{\chi} (j (\theta))$ for all $\theta \in C^\sd (\calA)$.
    The assumptions on $\ddot{\chi}$ imply the required property of $\ring{\chi}$,
    using that $j (\bar{\alpha}) = \alpha + \alpha^\vee$ for all $\alpha \in C^\circ (\calA)$.

    Finally, for Assumption~\ref{asn-h-sp-mod}~\ref{itm-epsilon-sd},
    as in Example~\ref{eg-lin-cat},
    if we have access to the extra data
    $\chi_{\calA} (\alpha, \beta)$,
    we may define
    $\dot{\chi}_{\calA} (\alpha, \theta) = \chi_{\calA} (\alpha, j (\theta))$.
    We also need an extra piece of data
    $\ddot{\chi}_{\calA} (\alpha) \in \bbZ$ for $\alpha \in C^\circ (\calA)$,
    satisfying $\ddot{\chi} (\alpha) =
    \ddot{\chi}_{\calA} (\alpha) + \ddot{\chi}_{\calA} (\alpha^\vee)$
    for all $\alpha \in C^\circ (\calA)$.
    We also require that
    \begin{equation}
        \ddot{\chi}_{\calA} (\alpha + \beta)
        = \ddot{\chi}_{\calA} (\alpha) + 
        \ddot{\chi}_{\calA} (\beta) + 
        \chi_{\calA} (\alpha, \beta^\vee)
    \end{equation}
    for all $\alpha, \beta \in C^\circ (\calA)$.
    We can then define $\chi \sss{\calA}{\sd} (\alpha, \theta) =
    \dot{\chi}_{\calA} (\alpha, \theta) + \ddot{\chi}_{\calA} (\alpha)$,
    and $\varepsilon^\sd (\alpha, \theta) =
    (-1)^{\chi \sss{\calA}{\sd} (\alpha, \theta)}$
    for all $\alpha \in C^\circ (\calA)$ and $\theta \in C^\sd (\calA)$.

    A typical choice of the optional extra data $\ddot{\chi}_{\calA} (\alpha)$
    is described as follows.
    If we choose $\Thetaul$ as in~\eqref{eq-choice-theta-ul},
    we may take
    \begin{equation}
        \ddot{\chi}_{\calA} (\alpha) =
        \rank (\calExt_{\smash{\alpha, \> \alpha^\vee}})^{\bbZ_2},
    \end{equation}
    where we take the $\bbZ_2$-fixed part in the Ext-complex,
    and the $\bbZ_2$-action comes from the natural identification
    \begin{equation}
        \calExt (E, F) \simeq \calExt (F^\vee, E^\vee)
    \end{equation}
    for $E, F \in \calA$.

    Now, one can apply Theorem~\ref{thm-tw-mod}
    to equip the graded vector space $\ring{H}_* (\calM^\sd; \bbQ)$
    with the structure of an involutive graded twisted module
    for the involutive graded vertex algebra $\hat{H}_* (\calM; \bbQ)$
    defined in Example~\ref{eg-lin-cat}.

    In addition, if Assumption~\ref{asn-rat-triv} is also satisfied,
    then $\ring{H}_* (\calM^\sd; \bbQ)$
    is also a graded twisted module
    for the involutive graded Lie algebra
    $\check{H}_* (\calM^\pl; \bbQ)$.

    As in Remark~\ref{rmk-rat-triv},
    even when Assumption~\ref{asn-rat-triv}~\ref{itm-rat-triv}
    is not satisfied for all classes $\alpha \in \pi_0 (\calM)$,
    it might still be possible to obtain a graded Lie algebra
    $\check{H}_* (\calM^\pl_{>0}; \bbQ)$.
    In this case, if $\calM_{>0}$ is chosen to be invariant under the involution,
    then this graded Lie algebra is also involutive,
    and $\ring{H}_* (\calM^\sd; \bbQ)$ is a graded twisted module for it.
\end{example}

\begin{example}
    \label{eg-sd-lin-cat-semitop}
    In Example~\ref{eg-sd-lin-cat},
    one could also consider the space
    \[
        X^\sd = |\calM|^{\bbZ_2} ,
    \]
    the homotopy fixed points of the $\bbZ_2$-action on $|\calM|$.

    We note that in order to take the homotopy fixed points,
    one has to have a \emph{full} $\bbZ_2$-action on $|\calM|$,
    i.e.\ an action up to coherent homotopy,
    rather than an $H$-action.
    Here, the full $\bbZ_2$-action is given by the $\bbZ_2$-action on $\calM$
    coming from the self-dual structure of $\calM^+$.

    Note also that there is a map
    \[
        \pi \colon |\calM^\sd| \longrightarrow |\calM|^{\bbZ_2},
    \]
    induced by the $\bbZ_2$-invariant map $|\calM^\sd| \to |\calM|$.
    The map $\pi$ is sometimes a homotopy equivalence,
    in which case we this example would be the same as
    Example~\ref{eg-sd-lin-cat},
    but in other cases, these two constructions can be different.

    Let $\bbZ_2$ act on the space $|\calM| \times |\calM| \times |\calM|$
    by swapping the first and third factors, 
    and applying the involution $I \colon |\calM| \to |\calM|$
    to the second factor.
    Then, the map
    \[
        (-) \oplus (-) \oplus (-) \colon
        |\calM| \times |\calM| \times |\calM| \longrightarrow |\calM|
    \]
    is $\bbZ_2$-equivariant, where $\bbZ_2$ acts on the
    right-hand-side $|\calM|$ by $I$.
    Taking the homotopy fixed points, we obtain a map
    \[
        \oplus^\sd \colon |\calM| \times |\calM|^{\bbZ_2}
        \longrightarrow |\calM|^{\bbZ_2},
    \]
    which we use for Assumption~\ref{asn-h-sp-mod}~\ref{itm-h-sp-mod}.

    More generally, we may take $X^\sd$ to be any space $|\calM|^\sd$,
    equipped with maps
    $\oplus^\sd \colon |\calM| \times |\calM|^\sd \to |\calM|^\sd$
    and $i \colon |\calM|^\sd \longrightarrow |\calM|^{\bbZ_2}$,
    such that the diagram
    \[ \begin{tikzcd}
        {} |\calM| \times |\calM|^\sd
        \ar[r, "\oplus^\sd"] \ar[d, "\id \times i"']
        &{} |\calM|^\sd \ar[d, "i"]
        \\{} |\calM| \times |\calM|^{\bbZ_2} \ar[r, "\oplus^\sd"]
        &{} |\calM|^{\bbZ_2}
    \end{tikzcd} \]
    commutes up to a homotopy.

    Following a similar procedure as in Example~\ref{eg-lin-cat},
    we may choose $K$-theory classes $\dot{\Theta}, \ddot{\Theta}$
    based on the choice of $\Theta$,
    and then choose $\ring{\chi}$ and $\varepsilon^\sd$ based on this.
    
    Now, Theorem~\ref{thm-tw-mod} equips
    the graded vector space $\ring{H}_* (X^\sd; \bbQ)$
    with the structure of an involutive graded twisted module
    for the involutive graded vertex algebra $\hat{H}_* (\calM; \bbQ)$
    defined in Example~\ref{eg-lin-cat}.
\end{example}

\subsection{Morphisms of vertex algebras}

\begin{assumption}
    \label{asn-h-sp-mor}
    Let $(X, \Theta, \varepsilon)$ and
    $(X', \Theta', \varepsilon')$
    be two sets of data, satisfying Assumption~\ref{asn-h-sp}.
    Suppose that we have the following.
    
    \begin{enumerate}
        \item
            \label{itm-h-sp-mor}
            We have a map of $H$-spaces
            \[
                F \colon X \longrightarrow X'.
            \]
            
        \item
            \label{itm-xi}
            We have classes
            \begin{align*}
                \Xi & \in K (X \times X), \\
                \hat{\Xi} & \in K (X),
            \end{align*}
            of $\upU (1)$-weights $(1, -1)$ and $0$, respectively.
            We assume that $\hat{\Xi}$ can be realized as
            the class of a vector bundle,
            possibly of mixed rank. Write
            \begin{alignat*}{2}
                \Xi_{\alpha, \> \beta} & =
                \Xi |_{ X_\alpha \times X_\beta }
                && \in K ( X_\alpha \times X_\beta ), \\
                \hat{\Xi}_{\alpha} & =
                \hat{\Xi} |_{ X_\alpha }
                && \in K ( X_\alpha )
            \end{alignat*}
            for $\alpha, \beta \in \pi_0 (X)$.
            Define
            \begin{align}
                \xi (\alpha, \beta) & = \rank \Xi_{\alpha, \> \beta} \, , \\
                \hat{\xi} (\alpha) & = \rank \hat{\Xi}_\alpha \, .
            \end{align}

        \item
            \label{itm-xi-check}
            For all $\alpha, \beta \in \pi_0 (X)$, we have
            \begin{equation}
                \label{eq-xi-check}
                ({\oplus}_{\alpha, \> \beta})^*
                (\hat{\Xi}_{\alpha + \beta}) =
                \pr_1^* (\hat{\Xi}_\alpha) +
                \pr_2^* (\hat{\Xi}_\beta) +
                \Xi_{\alpha, \> \beta} +
                \sigma \sss{\alpha, \> \beta}{*} (\Xi_{\beta, \> \alpha})
            \end{equation}
            in $K (X_\alpha \times \smash{X_\beta})$,
            where $\smash{\sigma_{\alpha, \> \beta}} \colon X_\alpha \times \smash{X_\beta}
            \to \smash{X_\beta} \times X_\alpha$ is the map
            swapping the two factors.
            In particular, we have
            \begin{equation}
                \label{eq-xi-check-rk}
                \hat{\xi} (\alpha + \beta) =
                \hat{\xi} (\alpha) + \hat{\xi} (\beta) +
                \xi (\alpha, \beta) + \xi (\beta, \alpha)
            \end{equation}
            for all $\alpha, \beta \in \pi_0 (X)$.

        \item
            \label{itm-xi-compat-theta}
            For all $\alpha, \beta \in \pi_0 (X)$, we have
            \begin{equation}
                (F_\alpha \times F_\beta)^*
                (\Theta \sss{F (\alpha), \> F (\beta)}{\prime}) =
                \Theta_{\alpha, \> \beta} +
                \Xi_{\alpha, \> \beta} +
                \sigma \sss{\alpha, \> \beta}{*} (\Xi \sss{\beta, \> \alpha}{\vee})
            \end{equation}
            in $K (X_\alpha \times X_\beta)$. 
            In particular, we have
            \begin{equation}
                \chi' (F (\alpha), F (\beta)) =
                \chi (\alpha, \beta) + \xi (\alpha, \beta) + \xi (\beta, \alpha)
            \end{equation}
            for all $\alpha, \beta \in \pi_0 (X)$.
            We also require that
            \begin{equation}
                \label{eq-chi-check-diff}
                \hat{\chi}' (F (\alpha)) =
                \hat{\chi} (\alpha) + 2 \hat{\xi} (\alpha)
            \end{equation}
            for all $\alpha \in \pi_0 (X)$.

        \item
            \label{itm-xi-compat-epsilon}
            For all $\alpha, \beta \in \pi_0 (X)$, we have
            \begin{equation}
                \varepsilon' (F (\alpha), F (\beta)) =
                (-1)^{\xi (\alpha, \> \beta)} \cdot \varepsilon (\alpha, \beta).
            \end{equation}
    \end{enumerate}
\end{assumption}

\begin{remark}
    Assumption~\ref{asn-h-sp-mor}~\ref{itm-xi-check}
    will be automatically satisfied
    when $\hat{\Xi}_\alpha = \Delta_\alpha^* (\Xi_{\alpha, \> \alpha})$,
    where $\Delta_\alpha \colon X_\alpha \to X_\alpha \times X_\alpha$
    is the diagonal map.
    This will be the case in all our examples.
\end{remark}

\begin{theorem}
    \label{thm-mor-va}
    Under Assumption~\textnormal{\ref{asn-h-sp-mor},}
    consider the graded vertex algebras
    $\hat{H}_* (X; \bbQ)$ and
    $\hat{H}_* (X'; \bbQ),$
    defined in Theorem~\textnormal{\ref{thm-va}.}
    Then there is a morphism of graded vertex algebras
    \[
        \Omega \colon \hat{H}_* (X; \bbQ) \longrightarrow
        \hat{H}_* (X'; \bbQ),
    \]
    given by
    \begin{equation}
        \Omega (A) = F_* (A \cap c_{\hat{\xi} (\alpha)} (\hat{\Xi} \sss{\alpha}{\vee}))
    \end{equation}
    for $A \in \hat{H}_* (X_\alpha; \bbQ)$.

    In particular, if Assumption~\textnormal{\ref{asn-rat-triv}}
    is satisfied for both $X$ and $X',$
    then there is an induced morphism of graded Lie algebras,
    \[
        \Omega^\pl \colon \check{H}_* (X^\pl; \bbQ) \longrightarrow
        \check{H}_* (X'^\pl; \bbQ).
    \]
\end{theorem}

The proof will be given in \S\ref{sect-proof-mor-va}.

\begin{theorem}
    \label{thm-comp-mor-va}
    Suppose we are given morphisms of $H$-spaces
    \[
        X \overset{F}{\longrightarrow}
        X' \overset{F'}{\longrightarrow}
        X \second,
    \]
    where $X, X', X \second$ are equipped with extra data
    satisfying Assumption~\textnormal{\ref{asn-h-sp},}
    and $F, F'$ are equipped with extra data
    satisfying Assumption~\textnormal{\ref{asn-h-sp-mor}}.

    Let $F \second = F' \circ F,$ and set
    \begin{alignat}{2}
        \Xi \sss{\alpha, \> \beta}{\dprime} & =
        \Xi_{\alpha, \> \beta} +
        (F_\alpha \times F_\beta)^* 
        (\Xi \sss{F (\alpha), \> F (\beta)}{\prime})
        && \ \in K (X_\alpha \times X_\beta),
        \\
        \hat{\Xi} \sss{\alpha}{\dprime} & =
        \hat{\Xi}_\alpha + F_\alpha^* (\hat{\Xi} \sss{F (\alpha)}{\prime})
        && \ \in K (X_\alpha)
    \end{alignat}
    for all $\alpha, \beta \in \pi_0 (X)$.
    Then $F \second$ equipped with this data
    satisfies Assumption~\textnormal{\ref{asn-h-sp-mor}}.

    Moreover, if $\Omega, \Omega', \Omega \second$
    are maps induced by $F, F', F \second,$
    as in Theorem~\textnormal{\ref{thm-mor-va},} then
    \begin{equation}
        \Omega \second = \Omega' \circ \Omega.
    \end{equation} 
\end{theorem}

The proof will be given in \S\ref{sect-proof-comp-mor-va}.

\subsection{Morphisms of twisted modules}

\begin{assumption}
    \label{asn-h-sp-mod-mor}
    Let
    $(X, X^\sd, \Theta, \dot{\Theta}, \ddot{\Theta}, \ring{\chi}, \varepsilon^\sd)$ and
    $(X', X'^\sd, \Theta', \dot{\Theta}', \ddot{\Theta}', \ring{\chi}', \varepsilon'^\sd)$
    be two sets of data, satisfying Assumption~\ref{asn-h-sp-mod}.
    Suppose that we have the following.

    We assume the following extra data and conditions.
    
    \begin{enumerate}
        \item
            \label{itm-f-invol}
            We have a map of involutive commutative $H$-spaces
            \begin{align*}
                F \colon X & \longrightarrow X' , 
            \end{align*}
            equipped with extra data
            satisfying Assumption~\ref{asn-h-sp-mor}.

        \item
            \label{itm-fsd}
            We have a map
            \begin{align*}
                F^\sd \colon X^\sd & \longrightarrow X'^\sd ,
            \end{align*}
            such that the diagram
            \[ \begin{tikzcd}
                X \times X^\sd
                \ar[r, "F \times F^\sd"]
                \ar[d, "\oplus^\sd"']
                & X' \times X'^\sd
                \ar[d, "\oplus^\sd"]
                \\ X^\sd \ar[r, "F^\sd"]
                & X'^\sd
            \end{tikzcd} \]
            commutes up to a homotopy.
            
        \item
            \label{itm-xi-dots}
            We have classes
            \begin{align*}
                \dot{\Xi} & \in K (X \times X^\sd), \\
                \ddot{\Xi} & \in K (X), \\
                \ring{\Xi} & \in K (X^\sd),
            \end{align*}
            where $\dot{\Xi}$ and $\ddot{\Xi}$ have weights~$1$ and~$2$,
            respectively.
            We further assume that $\ring{\Xi}$ is realized by a vector bundle.
            Write
            \begin{alignat*}{2}
                \dot{\Xi}_{\alpha, \> \theta} & =
                \dot{\Xi} |_{ X_\alpha \times X^\sd_\theta }
                && \in K ( X_\alpha \times X^\sd_\theta ), \\
                \ddot{\Xi}_{\alpha} & =
                \ddot{\Xi} |_{ X_\alpha }
                && \in K ( X_\alpha ), \\
                \ring{\Xi}_{\theta} & =
                \ring{\Xi} |_{ X^\sd_\theta }
                && \in K ( X^\sd_\theta )
            \end{alignat*}
            for $\alpha \in \pi_0 (X)$ and $\theta \in \pi_0 (X^\sd)$.
            Define
            \begin{align}
                \dot{\xi} (\alpha, \theta) & = \rank \dot{\Xi}_{\alpha, \> \theta} \, , \\
                \ddot{\xi} (\alpha) & = \rank \ddot{\Xi}_{\alpha} \, , \\
                \ring{\xi} (\theta) & = \rank \ring{\Xi}_\theta \, .
            \end{align}

        \item
            \label{itm-xi-ring}
            For all $\alpha \in \pi_0 (X)$ and $\theta \in \pi_0 (X^\sd)$,
            we have
            \begin{multline}
                \label{eq-xi-ring}
                ({\oplus}\sss{\alpha, \> \theta}{\sd})^* (\ring{\Xi}_{\bar{\alpha} + \theta})
                =
                \pr_1^* (\hat{\Xi}_{\alpha}) +
                \pr_2^* (\ring{\Xi}_\theta) + {}
                \\*
                \dot{\Xi}_{\alpha, \> \theta} +
                (I_\alpha \times \id_{X^\sd_\theta})^* (\dot{\Xi}_{\alpha^\vee, \> \theta}) +
                \pr_1^* (\ddot{\Xi}_\alpha) +
                I_\alpha^* \pr_1^* (\ddot{\Xi}_{\alpha^\vee})
            \end{multline}
            in $K (X_\alpha \times X^\sd_\theta)$. In particular, we have
            \begin{equation}
                \label{eq-xi-ring-rk}
                \ring{\xi} (\bar{\alpha} + \theta) =
                \hat{\xi} (\alpha) + \ring{\xi} (\theta) +
                \dot{\xi} (\alpha, \theta) + \dot{\xi} (\alpha^\vee, \theta) +
                \ddot{\xi} (\alpha) + \ddot{\xi} (\alpha^\vee)
            \end{equation}
            for all $\alpha \in \pi_0 (X)$ and $\theta \in \pi_0 (X^\sd)$.

        \item
            \label{itm-xi-dots-compat-theta}
            For all $\alpha \in \pi_0 (X)$ and $\theta \in \pi_0 (X^\sd)$,
            we have
            \begin{align}
                (F_\alpha \times F^\sd_\theta)^* (\dot{\Theta}'_{F (\alpha), \> F^\sd (\theta)}) & =
                \dot{\Theta}_{\alpha, \> \theta} +
                \dot{\Xi}_{\alpha, \> \theta} +
                (I_\alpha \times \id_{X^\sd_\theta})^* (\dot{\Xi} \sss{\alpha^\vee, \> \theta}{\vee}),
                \\
                F_\alpha^* (\ddot{\Theta}'_{F (\alpha)}) & =
                \ddot{\Theta}_{\alpha} +
                \ddot{\Xi}_{\alpha} +
                I_\alpha^* (\ddot{\Xi} \sss{\alpha^\vee}{\vee})
            \end{align}
            in $K (X_\alpha \times X^\sd_\theta)$
            and $K (X_\alpha)$, respectively. 
            In particular, we have
            \begin{align}
                \dot{\chi}' ( F (\alpha), F^\sd (\theta) ) & =
                \dot{\chi} (\alpha, \theta) +
                \dot{\xi} (\alpha, \theta) +
                \dot{\xi} (\alpha^\vee, \theta), \\
                \ddot{\chi}' ( F (\alpha) ) & =
                \ddot{\chi} (\alpha) + \ddot{\xi} (\alpha) + \ddot{\xi} (\alpha^\vee) 
            \end{align}
            for all $\alpha \in \pi_0 (X)$ and $\theta \in \pi_0 (X^\sd)$.
            We also require that
            \begin{equation}
                \label{eq-chi-ring-diff}
                \ring{\chi}' ( F^\sd (\theta) ) =
                \ring{\chi} (\theta) + 2 \ring{\xi} (\theta)
            \end{equation}
            for all $\theta \in \pi_0 (X^\sd)$.

        \item
            \label{itm-xi-dots-compat-epsilon}
            For all $\alpha \in \pi_0 (X)$ and $\theta \in \pi_0 (X^\sd)$, we have
            \begin{equation}
                \varepsilon'^\sd (F (\alpha), F^\sd (\theta)) =
                (-1)^{\dot{\xi} (\alpha, \> \theta) + \ddot{\xi} (\alpha)} \cdot
                \varepsilon^\sd (\alpha, \theta).
            \end{equation}
    \end{enumerate}
\end{assumption}

\begin{theorem}
    \label{thm-mor-tw-mod}
    Under Assumption~\textnormal{\ref{asn-h-sp-mod-mor},}
    we have a morphism of graded vector spaces
    \[
        \Omega^\sd \colon \ring{H}_* (X^\sd; \bbQ) \longrightarrow
        \ring{H}_* (X'^\sd; \bbQ),
    \]
    given by
    \begin{equation}
        \Omega^\sd (M) =
        F^\sd_* (M \cap c_{\ring{\xi} (\theta)} (\ring{\Xi} \sss{\theta}{\vee})),
    \end{equation}
    such that we have a commutative diagram
    \begin{equation}
    \label{eq-cd-omega-sd}
    \begin{tikzcd}
        \hat{H}_* (X; \bbQ) \otimes \ring{H}_* (X^\sd; \bbQ)
        \ar[r, "\Omega \otimes \Omega^\sd"]
        \ar[d, "{Y^\sd (-, \> z) (-)}"']
        & \hat{H}_* (X'; \bbQ) \otimes \ring{H}_* (X'^\sd; \bbQ)
        \ar[d, "{Y'^\sd (-, \> z) (-)}"]
        \\
        \ring{H}_* (X^\sd; \bbQ) \llparen z \rrparen
        \ar[r, "\Omega^\sd"]
        & \ring{H}_* (X'^\sd; \bbQ) \llparen z \rrparen \rlap{ ,}
    \end{tikzcd}
    \end{equation}
    where $\Omega$ is the morphism given by
    Theorem~\textnormal{\ref{thm-mor-va}.}

    In particular, if Assumption~\textnormal{\ref{asn-rat-triv}}
    is satisfied for both $X$ and $X',$
    we have a commutative diagram
    \begin{equation}
    \label{eq-cd-omega-pl-sd}
    \begin{tikzcd}
        \check{H}_* (X^\pl; \bbQ) \otimes \ring{H}_* (X^\sd; \bbQ)
        \ar[r, "\Omega^\pl \otimes \Omega^\sd"]
        \ar[d, "{\heart}"']
        & \check{H}_* (X'^\pl; \bbQ) \otimes \ring{H}_* (X'^\sd; \bbQ)
        \ar[d, "{\heart}"]
        \\
        \ring{H}_* (X^\sd; \bbQ) \llparen z \rrparen
        \ar[r, "\Omega^\sd"]
        & \ring{H}_* (X'^\sd; \bbQ) \llparen z \rrparen \rlap{ ,}
    \end{tikzcd}
    \end{equation}
    where $\Omega^\pl$ is the morphism of graded Lie algebras
    in Theorem~\textnormal{\ref{thm-mor-va}.}
\end{theorem}

The proof will be given in \S\ref{sect-proof-mor-tw-mod}.

\begin{theorem}
    \label{thm-comp-mor-tw-mod}
    Suppose that we are in the situation of
    Theorem~\textnormal{\ref{thm-comp-mor-va},}
    where we are given morphisms of commutative $H$-spaces
    \[
        X \overset{F}{\longrightarrow}
        X' \overset{F'}{\longrightarrow}
        X \second.
    \]
    Suppose that $X, X', X \second$ are involutive,
    and are equipped with extra data
    satisfying Assumption~\textnormal{\ref{asn-h-sp-mod},}
    including modules $X^\sd, X'^\sd, X^{\dprime \sd}$.
    Suppose also that $F, F'$ 
    are equipped with extra data satisfying
    Assumption~\textnormal{\ref{asn-h-sp-mod-mor},} including maps
    \[
        X^\sd \overset{F^\sd}{\longrightarrow}
        X'^\sd \overset{F'^\sd}{\longrightarrow}
        X^{\dprime \sd},
    \]
    compatible with $F$ and $F'$.

    Let $F \second = F' \circ F,$
    $F^{\dprime \sd} = F'^\sd \circ F^\sd,$ and set
    \begin{alignat}{2}
        \dot{\Xi} \sss{\alpha, \> \theta}{\dprime} & =
        \dot{\Xi}_{\alpha, \> \theta} +
        (F_\alpha \times F^\sd_\theta)^* 
        (\dot{\Xi} \sss{F (\alpha), \> F^\sd (\theta)}{\prime})
        && \ \in K (X_\alpha \times X^\sd_\theta),
        \\
        \ddot{\Xi} \sss{\alpha}{\dprime} & =
        \ddot{\Xi}_\alpha + F_\alpha^* (\ddot{\Xi} \sss{F (\alpha)}{\prime})
        && \ \in K (X_\alpha),
        \\
        \ring{\Xi} \sss{\theta}{\dprime} & =
        \ring{\Xi}_\theta +
        (F^\sd_\theta)^* (\ring{\Xi} \sss{F^\sd (\theta)}{\prime})
        && \ \in K (X^\sd_\theta)
    \end{alignat}
    for all $\alpha \in \pi_0 (X)$ and $\theta \in \pi_0 (X^\sd)$.
    Then $F \second$ and $F^{\dprime \sd}$ equipped with this data
    satisfies Assumption~\textnormal{\ref{asn-h-sp-mod-mor}}.

    Moreover, if $\Omega^\sd, \Omega'^\sd, \Omega^{\dprime \sd}$
    are maps induced by $F^\sd, F'^\sd, F^{\dprime \sd},$
    as in Theorem~\textnormal{\ref{thm-mor-tw-mod},} then
    \begin{equation}
        \Omega^{\dprime \sd} = \Omega'^\sd \circ \Omega^\sd.
    \end{equation} 
\end{theorem}

The proof will be given in \S\ref{sect-proof-comp-mor-tw-mod}.

\section{Enumerative invariants}

\label{sect-inv}

\subsection{The main conjecture}

We state the main conjecture,
Conjecture~\ref{conj-main},
on the structure of enumerative invariants
counting self-dual objects in a self-dual $\bbC$-linear category.
This will be an analogue of the conjectures in
Gross--Joyce--Tanaka~\cite[\S4]{GrossJoyceTanaka}
involving the ordinary enumerative invariants,
partially proved in Joyce~\cite{Joyce2021}.

We state the main result of Joyce~\cite{Joyce2021} below,
but in an imprecise manner,
as the precise statements are very long and technical.
They can be found in~\cite[Theorems~5.7--5.9]{Joyce2021}.

\begin{theorem}
    \label{thm-joyce-main}
    Let $\calA$ be a $\bbC$-linear quasi-abelian category,
    equipped with a moduli stack $\calM,$
    and extra data as in Example~\textnormal{\ref{eg-lin-cat},}
    satisfying Assumptions~\textnormal{\ref{asn-h-sp}}
    and~\textnormal{\ref{asn-rat-triv},} so that
    \[
        \hat{H}_* (\calM; \bbQ) \quad \text{and} \quad
        \check{H}_* (\calM^\pl; \bbQ)
    \]
    are equipped with the structure of
    a graded vertex algebra and
    a graded Lie algebra, respectively,
    by Theorems~\textnormal{\ref{thm-va}}
    and~\textnormal{\ref{thm-pl-lie-alg}}.
    Write $C^\circ (\calA) = \pi_0 (\calM),$
    which is an abelian monoid,
    and write $C (\calA) = C^\circ (\calA) \setminus \{ 0 \}$.

    We assume that we are given a derived enhancement $\calMul$ of $\calM$
    that is quasi-smooth.

    For stability conditions
    $\tau$ on $\calA$ satisfying certain conditions,
    one can define homology classes
    \[
        \inv{ \calM^\plss_{\alpha} (\tau) }
        \in H_{2 - \hat{\chi} (\alpha)} (\calM^\pl_{\alpha}; \bbQ)
    \]
    for $\alpha \in C (\calA)$ satisfying certain conditions,
    which we call the \emph{Joyce invariants.}
    Here, we assume the existence of an open substack
    $\calM^\plss_{\alpha} (\tau) \subset \calM^\pl_\alpha$
    consisting of the $\tau$-semistable objects.
    Moreover, we have the following:

    \begin{enumerate}
        \item 
            For any class $\alpha \in C (\calA),$
            such that all $\tau$-semistable objects of class~$\alpha$
            are $\tau$-stable, the stack
            $\calM^\plss_{\alpha} (\tau)$
            is a proper Deligne--Mumford stack, and we have
            \begin{equation}
                \inv{ \calM^\plss_{\alpha} (\tau) } =
                \iota_* \virt{ \calM^\plss_{\alpha} (\tau) } \, ,
            \end{equation}
            where $\iota \colon \calM^\plss_{\alpha} (\tau)
            \hookrightarrow \calM^\pl_\alpha$
            denotes the inclusion,
            and the right-hand side is the virtual fundamental class
            obtained from the derived enhancement,
            as in Definition~\textnormal{\ref{def-virt}}.

        \item
            For stability conditions $\tau, \tilde{\tau}$ on $\calA$
            satisfying certain conditions,
            and any $\alpha \in C (\calA),$
            we have the \emph{wall-crossing formula}
            \begin{multline}
                \label{eq-wcf-gl}
                \inv{ \calM^\plss_\alpha (\tilde{\tau}) } =
                \sum_{ \leftsubstack[8em]{
                    \\[-3ex]
                    & n > 0;
                    \alpha_1, \dotsc, \alpha_n \in C (\calA) \colon \\[-.5ex]
                    & \alpha = \alpha_1 + \cdots + \alpha_n
                } } {} 
                \tilde{U} (\alpha_1, \dotsc, \alpha_n; \tau, \tilde{\tau}) \cdot {} 
                \\[.5ex] 
                \bigl[ \bigl[ \bigl[ 
                \inv{ \calM^\plss_{\alpha_1} (\tau) } \, , 
                \inv{ \calM^\plss_{\alpha_2} (\tau) } \bigr] , \dotsc \bigr] ,
                \inv{ \calM^\plss_{\alpha_n} (\tau) } \bigr] 
            \end{multline}
            in the Lie algebra $\check{H}_* (\calM^\pl; \bbQ),$
            where $\tilde{U} ({\cdots})$ are certain combinatorial coefficients,
            and there are only finitely many non-zero terms.
    \end{enumerate}
\end{theorem}

\begin{remark}
    \label{rmk-dg}
    In Theorem~\ref{thm-joyce-main},
    one can also consider another setting,
    also discussed in~\cite[\S5.5]{Joyce2021},
    where $\calM$ is replaced by a moduli stack $\frM$ of objects
    in a $\bbC$-linear stable $\infty$-category $\calC$,
    equipped with a Bridgeland stability condition~\cite{Bridgeland2002}.
    We then take $\calA = \calP (I) \subset \calC$
    for an interval $I \subset \bbR$,
    as in Bridgeland~\cite[\S4]{Bridgeland2002},
    which is a quasi-abelian category
    provided that $I \cap (I + 1) = \varnothing$.
    The Bridgeland stability condition on $\calC$
    restricts to a stability condition $\tau$ on $\calA$.

    We may reformulate Theorem~\ref{thm-joyce-main} under this setting,
    with $\calA$ and $\frM$ in place of $\calA$ and $\calM$.
    The only essential difference is that
    the invariants now lie in the homology of $\frM$ instead of $\calM$.
    This does not provide any new information,
    as the new invariants can be obtained by
    pushing forward the old invariants
    along the inclusion $\calM \hookrightarrow \frM$.
\end{remark}

\begin{remark}
    \label{rmk-choice-ca}
    In Theorem~\ref{thm-joyce-main},
    we can sometimes take $C^\circ (\calA)$ to be
    a non-trivial quotient of $\pi_0 (\calM)$,
    as long as the data $\chi, \hat{\chi}, \varepsilon$, etc.,
    descend to $C^\circ (\calA)$,
    and the class $0 \in C^\circ (\calA)$ only contains the zero object.
    The reason for doing this is that sometimes,
    it can be tricky to determine the connected components of $\calM$,
    but one can simply work with a choice of $C^\circ (\calA)$
    without worrying about this problem.
\end{remark}

We make the following conjecture on the
structure of enumerative invariants
counting self-dual objects in a self-dual $\bbC$-linear category,
which is a counterpart of Gross--Joyce--Tanaka~\cite[\S4.3]{GrossJoyceTanaka}
and Theorem~\ref{thm-joyce-main}
for self-dual $\bbC$-linear categories.
We do not attempt to make a precise conjecture,
and instead, we expect that in different situations,
one would need extra conditions that depend on the specific case.

\begin{conjecture}
    \label{conj-main}
    Let $\calA$ be a self-dual $\bbC$-linear quasi-abelian category,
    equipped with extra data such that
    Theorem~\textnormal{\ref{thm-joyce-main}} holds.
    
    Suppose that the moduli stack $\calM$
    is equipped with extra data
    as in Example~\textnormal{\ref{eg-sd-lin-cat},}
    satisfying Assumption~\textnormal{\ref{asn-h-sp-mod}}.
    Write $C^\sd (\calA) = \pi_0 (\calM^\sd)$.

    We further assume that the self-dual structure of $\calM$
    extends to its derived enhancement $\calMul,$
    giving a derived enhancement $\calMul^\sd$ of $\calM^\sd,$
    as the $\bbZ_2$-fixed points in $\calMul$.
    We assume that $\calMul^\sd$ is also quasi-smooth.
    
    Then, for self-dual stability conditions $\tau$ on $\calA$
    satisfying certain conditions,
    one can define homology classes
    \[
        \inv{ \calM^\sdss_{\theta} (\tau) }
        \in H_{-\ring{\chi} (\theta)} (\calM^\sd_{\theta}; \bbQ)
    \]
    for $\theta \in C^\sd (\calA)$ satisfying certain conditions.
    Here, we assume the existence of an open substack
    $\calM^\sdss_{\theta} (\tau) \subset \calM^\sd_\theta$
    consisting of the $\tau$-semistable self-dual objects.
    Moreover, we should have the following:

    \begin{enumerate}
        \item 
            \label{itm-conj-fund}
            For any class $\theta \in C^\sd (\calA),$
            such that all $\tau$-semistable self-dual objects
            of class~$\theta$ are $\tau$-stable,
            the stack $\calM^\sdss_{\theta} (\tau)$ is a proper
            Deligne--Mumford stack, and we have
            \begin{equation}
                \inv{ \calM^\sdss_{\theta} (\tau) } =
                \iota_* \virt{ \calM^\sdss_{\theta} (\tau) } \, ,
            \end{equation}
            where $\iota \colon \calM^\sdss_{\theta} (\tau)
            \hookrightarrow \calM^\sd_\theta$
            denotes the inclusion,
            and the right-hand side is the virtual fundamental class
            obtained from the derived enhancement,
            as in Definition~\textnormal{\ref{def-virt}}.
            
        \item 
            \label{itm-conj-wcf}
            For self-dual stability conditions $\tau, \tilde{\tau}$ on $\calA$
            satisfying certain conditions,
            and any $\theta \in C^\sd (\calA),$
            we have the \emph{wall-crossing formula}
            \begin{align*}
                \numberthis
                \label{eq-wcf-sd}
                \hspace{2em} & \hspace{-2em}
                \inv{ \calM^\sdss_\theta (\tilde{\tau}) } = {}
                \\[.5ex] &
                \sum_{ \leftsubstack[8em]{
                    \\[-3ex]
                    & n \geq 0; \, m_1, \dotsc, m_n > 0; \\[-.5ex]
                    & \alpha_{1,1}, \dotsc, \alpha_{1,m_1}; \dotsc;
                    \alpha_{n,1}, \dotsc, \alpha_{n,m_n} \in C (\calA); \,
                    \rho \in \smash{C^\sd (\calA)} \colon \\[-.5ex]
                    & \theta = (\bar{\alpha}_{1,1} + \cdots + \bar{\alpha}_{1,m_1})
                    + \cdots + (\bar{\alpha}_{n,1} + \cdots + \bar{\alpha}_{n,m_n}) + \rho
                } } {} 
                \tilde{U}^\sd (\alpha_{1,1}, \dotsc, \alpha_{1,m_1}; \dotsc;
                \alpha_{n,1}, \dotsc, \alpha_{n,m_n}; \tau, \tilde{\tau}) \cdot {} 
                \\[.5ex] &
                \bigl[ \bigl[ \bigl[ 
                \inv{ \calM^\plss_{\alpha_{1,1}} (\tau) } \, , 
                \inv{ \calM^\plss_{\alpha_{1,2}} (\tau) } \bigr] , \dotsc \bigr] ,
                \inv{ \calM^\plss_{\alpha_{1,m_1}} (\tau) } \bigr] \heart \cdots \heart {}
                \\ &
                \bigl[ \bigl[ \bigl[ 
                \inv{ \calM^\plss_{\alpha_{n,1}} (\tau) } \, , 
                \inv{ \calM^\plss_{\alpha_{n,2}} (\tau) } \bigr] , \dotsc \bigr] ,
                \inv{ \calM^\plss_{\alpha_{n,m_n}} (\tau) } \bigr] \heart 
                \inv{ \calM^\sdss_{\rho} (\tau) } \ ,
            \end{align*}
            written in terms of
            the Lie algebra $\check{H}_* (\calM^\pl; \bbQ)$
            and its twisted module $\ring{H}_* (\calM^\sd; \bbQ),$
            where $\tilde{U}^\sd ({\cdots})$ are combinatorial coefficients
            defined in~\textnormal{\cite[Theorem~8.16]{Bu2023},}
            and there are only finitely many non-zero terms.
    \end{enumerate}
\end{conjecture}

We will give a partial proof of this conjecture
in the special case of the self-dual category of representations
of a self-dual quiver, in \S\ref{sect-quiver} below.

\begin{remark}
    Regarding Conjecture~\ref{conj-main}~\ref{itm-conj-fund},
    by Theorem~\ref{thm-st-sd-decomp},
    if all $\tau$-semistable self-dual objects of class~$\theta$
    are $\tau$-stable, then the stack
    $\calM^\sdss_\theta (\tau)$ has finite stabilizers,
    and hence, is Deligne--Mumford.
    Therefore, for the virtual fundamental class to be defined,
    we only need this stack to be proper.
    We note that this properness condition is far from automatic,
    and often needs to be checked on a case-by-case basis.
\end{remark}

Here are a few remarks on
variations of the statement of Conjecture~\ref{conj-main}.

\begin{remark}
    \label{rmk-bridgeland}
    In Conjecture~\ref{conj-main},
    as in Remark~\ref{rmk-dg},
    we can also consider a self-dual $\bbC$-linear stable $\infty$-category $\calC$,
    equipped with a moduli stack $\frM$,
    which is assumed to come from a self-dual $\bbC$-linear $\infty$-stack
    (Definition~\ref{def-sd-lin-infty-st}),
    and a self-dual Bridgeland stability condition,
    meaning that $\calP (t)^\vee = \calP (-t)$ for all $t \in \bbR$.
    Then, we may take $\calA = \calP (I) \subset \calC$
    for an interval $I \subset \bbR$ satisfying $I = -I$
    and $\pm \frac{1}{2} \notin I$.
    Then, $\calA$ becomes a self-dual quasi-abelian category,
    and the Bridgeland stability condition
    restricts to a self-dual stability condition on $\calA$.

    As in Remark~\ref{rmk-dg},
    we may reformulate Conjecture~\ref{conj-main},
    so that the invariants live in the homology of $\frM^\sd$
    instead of $\calM^\sd$.
    Again, this does not provide new information,
    but it sometimes facilitates computation.
\end{remark}

\begin{remark}
    \label{rmk-choice-csda}
    In Conjecture~\ref{conj-main},
    as in Remark~\ref{rmk-choice-ca},
    we can sometimes take $C^\circ (\calA)$
    to be a non-trivial quotient of $\pi_0 (\calM)$,
    and take $C^\sd (\calA)$
    to be a non-trivial quotient set of $\pi_0 (\calM^\sd)$,
    provided that $C^\circ (\calA)$ inherits
    an involutive abelian monoid structure,
    that $C^\sd (\calA)$ inherits
    an involutive $C^\circ (\calA)$-module structure,
    as in Definition~\ref{def-h-sp-invol},
    and that all the relevant data descends to $C^\sd (\calA)$.
    This can be helpful in cases where
    $\pi_0 (\calM)$ and $\pi_0 (\calM^\sd)$
    are difficult to determine.
\end{remark}

Finally, we make a remark on another situation
where a similar wall-crossing structure may be found,
which is a self-dual analogue of the conjecture in
Gross--Joyce--Tanaka~\cite[\S4.4]{GrossJoyceTanaka}.

\begin{remark}
    An interesting variant of Conjecture~\ref{conj-main}
    is given by considering the counting problem
    of coherent sheaves on a Calabi--Yau $4$-fold~$X$,
    also known as the theory of \emph{DT4 invariants},
    studied by Cao--Leung \cite{CaoLeung},
    Borisov--Joyce \cite{BorisovJoyce2017},
    Oh--Thomas \cite{OhThomas2023}, and others.

    In this case, we construct the self-dual category
    $\calA \subset \Db \cat{Coh} (X)$
    via a self-dual stability condition on $\Db \cat{Coh} (X)$,
    as in Remark~\ref{rmk-bridgeland}.
    For example, Bayer~\cite{Bayer2009} defined
    \emph{polynomial Bridgeland stability conditions}
    on $\Db \cat{Coh} (X)$.
    Some of these stability conditions are \emph{self-dual},
    in that the slice $\calP (\phi)$ is the dual of $\calP (-\phi)$ for all~$\phi$,
    which can be deduced from the discussion in~\cite[\S3.3]{Bayer2009},
    except that we may need to introduce a shift in~$\phi$ by a constant.
    One can then take
    \[
        \calA = \calP ((-1/2, 1/2)),
    \]
    the quasi-abelian category spanned by $\calP (\phi)$ with
    $-1/2 < \phi < 1/2$ under extensions.

    We hope that Theorem~\ref{thm-joyce-main} and Conjecture~\ref{conj-main}
    still hold in this case, but with the following modifications:

    \begin{enumerate}
        \item
            When defining the graded vertex algebra structure
            on $\hat{H}_* (\calM; \bbQ)$,
            we do not choose the data~$\Theta$
            as in~\eqref{eq-choice-theta-ul} in Example~\ref{eg-lin-cat}.
            Instead, we choose
            \[
                \Thetaul = \calExt^\vee,
            \]
            without symmetrizing it.
            Note that this only works when $X$ is even-dimensional and Calabi--Yau,
            since otherwise, \eqref{eq-theta-sym} would fail.
            The data $\dot{\Theta}, \ddot{\Theta}, \hat{\chi}, \ring{\chi}$, etc.,
            are to be chosen according to the new choice of~$\Theta$.

        \item 
            The derived enhancement $\calMul$ of the moduli stack $\calM$
            is no longer quasi-smooth,
            but it has a \emph{$(-2)$-shifted symplectic structure},
            in the sense of Pantev--Toën--Vaquié--Vezzosi~\cite{PTVV2013}.
            In the self-dual case, we do not require $\calMul^\sd$
            to be quasi-smooth, but we require that it has
            a $(-2)$-shifted symplectic structure 
            induced from that of $\calMul$.

        \item 
            \label{itm-cy4-orientation}
            Cao--Gross--Joyce~\cite{CaoGrossJoyce2020}
            showed that $\calMul$ and its rigidification
            $\calMul^\pl = \calMul / [*/\Gm]$ are \emph{orientable},
            in the sense of Borisov--Joyce~\cite[\S2.4]{BorisovJoyce2017}.
            We require that the self-dual moduli stacks $\calMul^\sd$ 
            are also orientable,
            and we will need to choose orientations on
            $\calMul^\pl$ and $\calMul^\sd$
            satisfying compatibility conditions.

        \item 
            The Behrend--Fantechi virtual fundamental class $\virt{-}$
            should be replaced by the virtual classes
            constructed by Borisov--Joyce~\cite{BorisovJoyce2017}
            or Oh--Thomas~\cite{OhThomas2023}
            for oriented $(-2)$-shifted symplectic
            derived manifolds or derived schemes,
            and suitable generalizations to
            derived orbifolds or derived Deligne--Mumford stacks.
    \end{enumerate}
    On the other hand, one might face the following problem:

    \begin{enumerate}[resume]
        \item 
            An orientation of the self-dual moduli stack $\calMul^\sd$,
            as required by~\ref{itm-cy4-orientation} above,
            might not exist, especially in the Type~B/C case,
            i.e.~the symplectic or odd rank orthogonal case,
            as suggested by Joyce--Upmeier~\cite[\S3.6]{JoyceUpmeier},
            who showed that certain related moduli spaces are unorientable,
            although they do not directly come from Calabi--Yau $4$-folds.
            However, it is possible that our theory can still work
            in the Type D case, or the even rank orthogonal case.
    \end{enumerate}
\end{remark}

\subsection{An alternative form of the wall-crossing formula}

We explain how to rewrite our wall-crossing formula~\eqref{eq-wcf-sd}
in the style of Kontsevich--Soibelman~\cite[\S1.4]{KontsevichSoibelman2008},
as the self-dual counterpart of
Gross--Joyce--Tanaka~\cite[\S4.2]{GrossJoyceTanaka}.
The main result is Theorem~\ref{thm-wcf-ks} below.
Compare also Young~\cite[Theorem~4.5]{Young2015}.

\begin{definition}
    \label{def-univ-env-tw-mod}
    Let $L$ be an involutive graded Lie algebra,
    and let $M$ be a twisted module for $L$.
    Define a $U (L)$-module
    \begin{equation}
        U^\tw (L; M) = (U (L) \otimes M) / I,
    \end{equation}
    where $U (L)$ is the universal enveloping algebra of $L$,
    and $I$ is the $U (L)$-submodule generated by elements of the form
    \[
        x \otimes m - x^\vee \otimes m - 1 \otimes (x \heart m),
    \]
    where $x \in L$ and $m \in M$.
    Here, $\heart$ denotes the twisted module structure on $M$.
    We denote by
    \[
        {\diamond} \colon U (L) \otimes U^\tw (L; M)
        \longrightarrow U^\tw (L; M)
    \]
    the $U (L)$-action on $U^\tw (L; M)$.
\end{definition}

\begin{definition}
    \label{def-moduli-utw}
    Under Assumption~\ref{asn-h-sp-mod},
    consider the Lie algebra $\check{H}_* (\calM^\pl; \bbQ)$
    and its module
    \begin{equation}
        U^\tw =
        U^\tw \bigl( \check{H}_* (\calM^\pl; \bbQ); \,
        \ring{H}_* (\calM^\sd; \bbQ) \bigr).
    \end{equation}
    We define a decomposition
    \begin{equation}
        U^\tw = \bigoplus_{\theta \in \pi_0 (X^\sd)} U^\tw_\theta
    \end{equation}
    by setting $U^\tw_\theta$ to be the sum of the images of
    \[
        \check{H}_* (\calM_{\alpha_1}^\pl; \bbQ) \otimes \cdots \otimes
        \check{H}_* (\calM_{\alpha_n}^\pl; \bbQ) \otimes
        \ring{H}_* (\calM^\sd_{\rho}; \bbQ),
    \]
    for all choices of $n \geq 0$ and
    $\alpha_1, \dotsc, \alpha_n \in \pi_0 (X)$ and $\rho \in \pi_0 (X^\sd)$
    such that
    $\theta = \bar{\alpha}_1 + \cdots + \bar{\alpha}_n + \rho$.

    We also define a completed version,
    \begin{equation}
        \hat{U}^\tw = \prod_{\theta \in \pi_0 (X^\sd)} U^\tw_\theta .
    \end{equation}
    The Lie algebra~$\check{H}_* (\calM^\pl; \bbQ)$
    also acts on the space $\hat{U}^\tw$,
    and we also denote this action by $\diamond$.
\end{definition}

\begin{theorem}
    One can rewrite the wall-crossing formulae~%
    \eqref{eq-wcf-gl} and~\eqref{eq-wcf-sd} as follows.

    \begin{enumerate}
        \item 
            \label{itm-wcf-star}
            The formula~\eqref{eq-wcf-gl} is equivalent to the formula
            \begin{multline}
                \label{eq-wcf-star}
                \inv{ \calM^\plss_\theta (\tilde{\tau}) } =
                \sum_{ \leftsubstack[8em]{
                    \\[-3ex]
                    & n \geq 1; \, 
                    \alpha_1, \dotsc, \alpha_n \in C (\calA) \colon \\[-.5ex]
                    & \alpha = \alpha_1 + \cdots + \alpha_n
                } } {} 
                U (\alpha_1, \dotsc, \alpha_n; \tau, \tilde{\tau}) \cdot {} 
                \\[.5ex]
                \inv{ \calM^\plss_{\alpha_1} (\tau) } * \cdots *
                \inv{ \calM^\plss_{\alpha_n} (\tau) }
            \end{multline}
            in~$U (\check{H}_* (\calM^\pl; \bbQ)),$
            where $U ({\cdots})$ are coefficients
            defined in~\textnormal{\cite[Definition~8.8]{Bu2023},}
            from Joyce~\textnormal{\cite[Definition~4.4]{Joyce2008IV}}.
        \item 
            \label{itm-wcf-diamond}
            The formula~\eqref{eq-wcf-sd} is equivalent to the formula
            \begin{multline}
                \label{eq-wcf-diamond}
                \inv{ \calM^\sdss_\theta (\tilde{\tau}) } =
                \sum_{ \leftsubstack[8em]{
                    \\[-3ex]
                    & n \geq 0; \, 
                    \alpha_1, \dotsc, \alpha_n \in C (\calA); \,
                    \rho \in \smash{C^\sd (\calA)} \colon \\[-.5ex]
                    & \theta = \bar{\alpha}_1 + \cdots + \bar{\alpha}_n + \rho
                } } {} 
                U^\sd (\alpha_1, \dotsc, \alpha_n; \tau, \tilde{\tau}) \cdot {} 
                \\[.5ex]
                \inv{ \calM^\plss_{\alpha_1} (\tau) } \diamond \cdots \diamond
                \inv{ \calM^\plss_{\alpha_n} (\tau) } \diamond 
                \inv{ \calM^\sdss_{\rho} (\tau) } 
            \end{multline}
            in~$U^\tw,$
            where $U^\sd ({\cdots})$ are coefficients
            defined in~\textnormal{\cite[Definition~8.8]{Bu2023}}.
    \end{enumerate}
\end{theorem}

\begin{proof}
    This conversion between the two versions
    is entirely formal and does not involve
    any properties of the invariants.
    Part~\ref{itm-wcf-star} follows from
    the definition of the coefficients~$\tilde{U} ({\cdots})$
    in~\cites[Theorem~3.12]{Joyce2021}[Theorem~8.14]{Bu2023},
    using Joyce~\cite[Theorem~5.4]{Joyce2008IV}.
    Part~\ref{itm-wcf-diamond}
    follows from~\cite[Theorem~8.16]{Bu2023}.
\end{proof}

\begin{definition}
    Assume that we are in a situation where Conjecture~\ref{conj-main} holds.
    Define elements
    \begin{align}
        \label{eq-def-delta}
        \delta_\alpha (\tau) & = 
        \sum_{ \leftsubstack[6em]{
            & \alpha_1, \dotsc, \alpha_n \in C (\calA) \colon \\[-.5ex]
            & \alpha = \alpha_1 + \cdots + \alpha_n \, , \\[-.5ex]
            & \tau (\alpha_1) = \cdots = \tau (\alpha_n)
        } } \hspace{.5em}
        \frac{1}{n!} \cdot
        \inv{\calM^\plss_{\alpha_1} (\tau)} * \cdots *
        \inv{\calM^\plss_{\alpha_n} (\tau)} \, ,
        \\[1ex]
        \label{eq-def-delta-sd}
        \delta^\sd_\theta (\tau) & =
        \sum_{ \leftsubstack[6em]{
            \\[-1.5ex]
            & \alpha_1, \dotsc, \alpha_n \in C (\calA), \,
            \rho \in C^\sd (\calA) \colon \\[-.5ex]
            & \theta = \bar{\alpha}_1 + \cdots + \bar{\alpha}_n + \rho \\[-.5ex]
            & \tau (\alpha_1) = \cdots = \tau (\alpha_n) = 0
        } } {}
        \frac{1}{2^n \, n!} \cdot
        \inv{\calM^\plss_{\alpha_1} (\tau)} \diamond \cdots \diamond
        \inv{\calM^\plss_{\alpha_n} (\tau)} \diamond
        \inv{\calM^\sdss_{\rho} (\tau)} \, .
        \raisetag{2.5ex}
    \end{align}
    Here, Equation~\eqref{eq-def-delta} is in
    $U (\check{H}_* (\calM^\pl; \bbQ))$,
    where the multiplication is denoted by $*$.
    Equation~\eqref{eq-def-delta-sd} is in
    $U^\tw (\check{H}_* (\calM^\pl; \bbQ); \ring{H}_* (\calM^\sd; \bbQ))$,
    as in Definition~\ref{def-univ-env-tw-mod}.

    The elements $\delta_\alpha (\tau)$ in~\eqref{eq-def-delta}
    were defined in Gross--Joyce--Tanaka~\cite[\S4.2]{GrossJoyceTanaka},
    in analogy to the identity in Joyce's work~\cite[Theorem~7.7]{Joyce2007III}
    on motivic invariants.
    Similarly, the elements $\delta^\sd_\theta (\tau)$
    in~\eqref{eq-def-delta-sd} are defined
    in analogy to the author's constructions
    of motivic invariants in~\cite[\S8.1]{Bu2023}.
\end{definition}

\begin{theorem}
    \label{thm-wcf-ks}
    Suppose we are in a situation
    where~\eqref{eq-wcf-gl} and~\eqref{eq-wcf-sd} hold.
    Then we have an identity
    \begin{multline}
        \label{eq-ks-wcf-sd}
        \sum_{ \leftsubstack[6em]{
            \\[-2ex]
            & n \geq 0, \, \alpha_1, \dotsc, \alpha_n \in C (\calA), \,
            \rho \in C^\sd (\calA) \colon \\[-.5ex]
            & \tau (\alpha_1) > \cdots > \tau (\alpha_n) > 0
        } } {}
        \delta_{\alpha_1} (\tau) \diamond \cdots \diamond
        \delta_{\alpha_n} (\tau) \diamond
        \delta^\sd_{\rho} (\tau)
        \\
        =
        \sum_{ \leftsubstack[6em]{
            \\[-2ex]
            & n \geq 0, \, \alpha_1, \dotsc, \alpha_n \in C (\calA), \,
            \rho \in C^\sd (\calA) \colon \\[-.5ex]
            & \tilde{\tau} (\alpha_1) > \cdots > \tilde{\tau} (\alpha_n) > 0
        } } {}
        \delta_{\alpha_1} (\tilde{\tau}) \diamond \cdots \diamond
        \delta_{\alpha_n} (\tilde{\tau}) \diamond
        \delta^\sd_{\rho} (\tilde{\tau})
    \end{multline}
    in the space $\hat{U}^\tw$ in Definition~\textnormal{\ref{def-moduli-utw},}
    provided that for each $\theta \in C^\sd (\calA),$
    there are only finitely many non-zero terms on each side
    with $\alpha_1 + \cdots + \alpha_n + \rho = \theta$.
\end{theorem}

\begin{proof}
    Consider the wall-crossing formulae~%
    \eqref{eq-wcf-star}--\eqref{eq-wcf-diamond}.
    By a combinatorial argument carried out
    in~\cite[Theorem~5.2 and~\S5.4]{Joyce2008IV}
    and in~\cite[Theorem~8.11]{Bu2023},
    respectively for the ordinary and self-dual cases,
    we obtain wall-crossing formulae
    \begin{align}
        \label{eq-wcf-delta}
        \delta_\alpha (\tilde{\tau}) & =
        \sum_{ \leftsubstack{
            \\[-2.5ex]
            & n > 0; \, 
            \alpha_1, \dotsc, \alpha_n \in C (\calA) \colon \\[-.5ex]
            & \alpha = \alpha_1 + \cdots + \alpha_n
        } } {} 
        S (\alpha_1, \dotsc, \alpha_n; \tau, \tilde{\tau}) \cdot
        \delta_{\alpha_1} (\tau) * \cdots *
        \delta_{\alpha_n} (\tau) ,
        \\[1ex]
        \label{eq-wcf-delta-sd}
        \delta^\sd_\theta (\tilde{\tau}) & =
        \sum_{ \leftsubstack{
            \\[-2.5ex]
            & n \geq 0; \, 
            \alpha_1, \dotsc, \alpha_n \in C (\calA); \,
            \rho \in \smash{C^\sd (\calA)} \colon \\[-.5ex]
            & \theta = \bar{\alpha}_1 + \cdots + \bar{\alpha}_n + \rho
        } } {} 
        S^\sd (\alpha_1, \dotsc, \alpha_n; \tau, \tilde{\tau}) \cdot
        \delta_{\alpha_1} (\tau) \diamond \cdots \diamond
        \delta_{\alpha_n} (\tau) \diamond 
        \delta^\sd_{\rho} (\tau) ,
    \end{align}
    where $S({\cdots})$ and $S^\sd ({\cdots})$ are coefficients
    defined in~\cite[Definition~8.8]{Bu2023}.

    Now, we expand the right-hand side of~\eqref{eq-ks-wcf-sd}
    using~\eqref{eq-wcf-delta}--\eqref{eq-wcf-delta-sd}.
    Let~$0$ denote the trivial stability condition,
    and apply~\cite[Lemma~11.8]{Bu2023} with $\tau, \tilde{\tau}, 0$
    in place of $\tau, \hat{\tau}, \tilde{\tau}$ there.
    By~\cite[Definition~8.8]{Bu2023}, we have
    \begin{align}
        S (\alpha_1, \dotsc, \alpha_n; \tau, 0)
        & = \begin{cases}
            1, & \tau (\alpha_1) > \cdots > \tau (\alpha_n), \\ 
            0, & \text{otherwise},
        \end{cases} \\
        S^\sd (\alpha_1, \dotsc, \alpha_n; \tau, 0)
        & = \begin{cases}
            1, & \tau (\alpha_1) > \cdots > \tau (\alpha_n) > 0, \\ 
            0, & \text{otherwise},
        \end{cases}
    \end{align}
    for any $\tau$, so that \eqref{eq-ks-wcf-sd}
    follows from~\cite[Lemma~11.8]{Bu2023}.
\end{proof}

\section{Self-dual quivers}

\label{sect-quiver}

In this section, we study a basic example
of the theory developed so far,
which is the self-dual category of representations of a self-dual quiver.
We verify that constructions of algebraic structures
in \S\ref{sect-vertex-moduli} apply to this case,
giving involutive graded vertex algebras and
involutive graded twisted modules,
and we state Theorem~\ref{thm-quiver-main},
which partially verifies our main conjecture,
Conjecture~\ref{conj-main}.

\subsection{Self-dual quivers}

We recall the notion of \emph{self-dual quivers}
from \cite[\S9]{Bu2023},
defined following Derksen--Weyman~\cite{DerksenWeyman2002}
and Young~\cite{Young2015,Young2020},
where the notion was called \emph{symmetric quivers}.
Representations of self-dual quivers
will serve as a basic example of our theory of
enumerative invariants in self-dual categories.

\begin{definition}
    \label{def-sd-quiver}
    A \emph{quiver} is a quadruple $Q = (Q_0, Q_1, s, t)$,
    where $Q_0, Q_1$ are finite sets,
    and $s, t \colon Q_1 \to Q_0$ are maps.
    We think of this data as a finite directed graph.

    A \emph{representation} of $Q$ over $\bbC$ is a set of data
    $E = ((E_i)_{i \in Q_0} \, , (e_a)_{a \in Q_1})$,
    where each $E_i$ is a finite-dimensional $\bbC$-vector space,
    and each $e_a \colon E_{s (a)} \to E_{t (a)}$ is a linear map.
    We will often use the lowercase of
    the same letter as the name of the representation
    to denote the structure maps.
    A \emph{morphism} of representations $h \colon E \to F$
    is a collection of maps $(h_i \colon E_i \to F_i)_{i \in Q_0}$\,,
    such that $f_a \circ h_{s (a)} = h_{t (a)} \circ e_a$
    for all $a \in Q_1$\,.
    We have the abelian category $\cat{Mod} (\bbC Q)$
    of representations of $Q$ over $\bbC$.

    The \emph{opposite quiver} of $Q$ is the quiver
    $Q^\op = (Q_0, Q_1, t, s)$,
    which is the quiver obtained from $Q$
    by inverting the directions of all arrows.

    A \emph{self-dual quiver}
    is a tuple $(Q, \sigma, u, v)$, where
    \begin{itemize}
        \item 
            $Q = (Q_0, Q_1, s, t)$ is a quiver.
        \item 
            $\sigma = (\sigma_i \colon Q_i \to Q_i)_{i=0,1}$
            consists of two maps,
            such that $\sigma_i^2 = \id_{\smash{Q_i}}$ for $i = 0, 1$,
            \,$\sigma_0 \circ s = t \circ \sigma_1$\,,
            and $\sigma_0 \circ t = s \circ \sigma_1$\,.
            This is called a \emph{contravariant involution} of~$Q$,
            and can be seen as an isomorphism $\sigma \colon Q \simto Q^\op$.
            We also write these maps as $(-)^\vee \colon Q_i \to Q_i$ for $i = 0, 1$.
        \item
            $u \colon Q_0 \to \{ \pm1 \}$ and $v \colon Q_1 \to \{ \pm1 \}$ are functions,
            such that $u_i = u_{\smash{i^\vee}}$ for all $i \in Q_0$\,,
            and $v_a \, v_{\smash{a^\vee}} = u_{s (a)} \, u_{t (a)}$
            for all $a \in Q_1$\,.
    \end{itemize}
    Given this data, one can define a self-dual structure on
    the abelian category $\cat{Mod} (\bbC Q)$,
    sending an object $E = ((E_i)_{i \in Q_0} \, , (e_a)_{a \in Q_1})$
    to the object $E^\vee = ((E \sss{i^\vee}{\vee})_{i \in Q_0} \, ,
    (v_a \cdot e \sss{a^\vee}{\vee})_{a \in Q_1})$.
    The identification $E^{\vee \vee} \simeq E$
    is given by
    $u_i \cdot \mathrm{ev} \sss{i}{-1} \colon E_i^{\vee \vee} \simeq E_i$\,,
    where $\mathrm{ev}_i \colon E_i \simeq E_i^{\vee \vee}$
    is the evaluation isomorphism.
    The self-dual objects in $\cat{Mod} (\bbC Q)$
    are called \emph{self-dual representations} of the self-dual quiver~$Q$.
\end{definition}

\begin{definition}
    Let $Q$ be a quiver. Define
    \[
        K (Q) = \bbZ^{Q_0}, \qquad
        C^\circ (Q) = \bbN^{Q_0} \subset K (Q), \qquad
        C (Q) = C^\circ (Q) \setminus \{ 0 \},
    \]
    whose elements are written as
    $\alpha = (\alpha_i)_{i \in Q_0}$\,,
    and called \emph{dimension vectors} of $Q$.
    The \emph{moduli stack of representations} of $Q$ is
    \begin{equation}
        \calM = \coprod_{\alpha \in C^\circ (Q)} \calM_\alpha \, ,
    \end{equation}
    where $\calM_\alpha = [V_\alpha / G_\alpha]$, with
    \begin{equation}
        \label{eq-m-quiv-expl}
        V_\alpha = \prod_{a \in Q_1} \Hom (\bbC^{\alpha_{s (a)}}, \bbC^{\alpha_{t (a)}}),
        \qquad
        G_\alpha = \prod_{i \in Q_0} \GL (\bbC^{\alpha_i}),
    \end{equation}
    and this can be extended to a $\bbC$-linear stack,
    as in~\cite[\S9.2]{Bu2023}.

    For a \emph{weak stability condition} $\tau$ on $\cat{Mod} (\bbC Q)$,
    as in~\cite[Definition~2.4]{Bu2023},
    which includes all stability conditions,
    and for any $\alpha \in C^\circ (Q)$, there is an open substack
    $\calM^\ss_\alpha (\tau) \subset \calM_\alpha$
    of $\tau$-semistable objects of class~$\alpha$,
    by~\cite[Proposition~9.9]{Bu2023}.

    We also have a derived moduli stack $\frM$ of objects
    in the $\bbC$-linear dg-category $\Db \cat{Mod} (\bbC Q)$,
    as in Toën--Vaquié \cite[Definition~3.33]{ToenVaquie2007},
    with a decomposition
    \begin{equation}
        \frM = \coprod_{\alpha \in K (Q)} \frM_\alpha \, .
    \end{equation}
    We have $\calM_\alpha \subset \frM_\alpha$
    as an open substack.

    Now, suppose that $(Q, \sigma, u, v)$ is a self-dual quiver.
    Define
    \begin{equation}
        C^\sd (Q) = \biggl\{ \,
            \theta \in C^\circ (Q) \biggm|
            \begin{aligned}[c]
                & \theta_i = \theta_{\smash{i^\vee}} \text{ for all } i, \\[-.5ex]
                & 2 \mid \theta_i \text{ if } i = i^\vee \text{ and } u_i = -1
            \end{aligned}
        \, \biggr\},
    \end{equation}
    and similarly $K^\sd (Q) \subset K (Q)$.
    The \emph{moduli stack of self-dual representations} of $Q$ is
    \begin{equation}
        \calM^\sd = \calM^{\bbZ_2}
        = \coprod_{\theta \in C^\sd (Q) \vphantom{^0}} \calM^\sd_\theta ,
    \end{equation}
    as in~\cite[\S9.2]{Bu2023},
    where $\calM^{\bbZ_2}$ denotes the homotopy fixed points
    of the $\bbZ_2$-action on $\calM$ coming from the self-dual structure of $Q$.
    The moduli stack $\calM^\sd_\theta$
    can be written more explicitly in~\eqref{eq-quiv-vsd}--\eqref{eq-quiv-gsd} below.

    For a \emph{self-dual weak stability condition} $\tau$ on $\cat{Mod} (\bbC Q)$,
    as in~\cite[Definition~3.12]{Bu2023},
    we have an open substack $\calM^\sdss_\theta (\tau) \subset \calM^\sd_\theta$
    of $\tau$-semistable self-dual objects of class~$\theta$.
    By Theorem~\ref{thm-ss-sd},
    this can be defined as $\calM^\ss (\tau)^{\bbZ_2} \cap \calM^\sd_\theta$,
    where $\calM^\ss (\tau)$ is the union of all
    $\calM^\ss_\alpha (\tau)$ in $\calM$,
    and we take the homotopy fixed points
    of the $\bbZ_2$-action coming from the self-dual structure.

    There is also a derived moduli stack
    \begin{equation}
        \frM^\sd = \frM^{\bbZ_2} 
        = \coprod_{\theta \in K^\sd (Q) \vphantom{^0}} \frM^\sd_\theta ,
    \end{equation}
    parametrizing self-dual objects in $\Db \cat{Mod} (\bbC Q)$.

    For $i \in Q_0$\,,
    there are universal vector bundles and perfect complexes
    \begin{equation}
        \label{eq-quiv-univ-vb}
        \begin{alignedat}{2}
            \calU_i & \in \cat{Vect} ( \calM ) , & \qquad
            \calV_i & \in \cat{Vect} ( \calM^\sd ) , \\
            \calU_i & \in \cat{Perf} ( \frM ) , & \qquad
            \calV_i & \in \cat{Perf} ( \frM^\sd ) .
        \end{alignedat}
    \end{equation}
    We have canonical isomorphisms $\calV_i^\vee \simeq \calV_{i^\vee}$
    for both $\calM^\sd$ and $\frM^\sd$.
\end{definition}

\begin{definition}
    Let $Q$ be a quiver.
    A \emph{stability function} on $Q$ is a function
    \[
        \tau \colon Q_0 \longrightarrow \bbQ .
    \]
    Given such a stability function, we can define
    a stability condition on $\cat{Mod} (\bbC Q)$,
    denoted by $\tau \colon C (Q) \longrightarrow \bbQ$,
    where $\bbQ$ is equipped with the usual order,
    by setting
    \begin{equation}
        \tau (\alpha) =
        \frac{\sum_{i \in Q_0} \tau (i) \, \alpha_i}{\sum_{i \in Q_0} \alpha_i}
    \end{equation}
    for all $\alpha \in C (Q)$.

    Now, suppose that $(Q, \sigma, u, v)$ is a self-dual quiver.
    We say that a stability function $\tau$ on $Q$ is \emph{self-dual},
    if $\tau \circ \sigma_0 = \tau \colon Q_0 \to \bbQ$.
    In this case, the corresponding stability condition on $\cat{Mod} (\bbC Q)$
    is a self-dual stability condition.
\end{definition}

\begin{notation}
    \label{ntn-sd-quiv}
    We fix a self-dual quiver $Q$.
    For convenience in notation,
    we choose total orders $\leq$ on $Q_0$ and $Q_1$\,,
    so that we can write
    \begin{align}
        Q_0 & = Q_0^+ \sqcup Q_0^- \sqcup Q_0^\tria \sqcup Q_0^{\tria\vee}, \\
        Q_1 & = Q_1^+ \sqcup Q_1^- \sqcup Q_1^\tria \sqcup Q_1^{\tria\vee}, 
    \end{align}
    where
    \begin{align}
        Q_0^\pm & = \{ i \in Q_0 \mid i = i^\vee, \, u_i = \pm 1 \}, \\
        Q_0^\tria & = \{ i \in Q_0 \mid i < i^\vee \}, \\
        Q_1^\pm & = \{ a \in Q_1 \mid a = a^\vee, \, u_{s (a)} \, v_a = \pm 1 \}, \\
        Q_1^\tria & = \{ a \in Q_1 \mid a < a^\vee \}.
    \end{align}
    This is purely for notational reasons,
    and the underlying theory will not depend on this choice.
\end{notation}

Now, using Notation~\ref{ntn-sd-quiv},
we may write $\calM^\sd_\theta = [V^\sd_\theta / G^\sd_\theta]$, where
\begin{align}
    \label{eq-quiv-vsd}
    V^\sd_\theta & =
    \prod_{a \in Q_1^\tria} {}
    \Hom (\bbC^{\theta_{s (a)}}, \bbC^{\theta_{t (a)}}) \times
    \prod_{a \in Q_1^+} {}
    \Sym^2 (\bbC^{\theta_{t (a)}}) \times
    \prod_{a \in Q_1^-} {}
    {\wedge}^2 (\bbC^{\theta_{t (a)}}), 
    \\
    \label{eq-quiv-gsd}
    G^\sd_\theta & =
    \prod_{i \in Q_0^\tria} {}
    \GL (\bbC^{\theta_i}) \times
    \prod_{i \in Q_0^+} {}
    \upO (\bbC^{\theta_i}) \times
    \prod_{i \in Q_0^-} {}
    \Sp (\bbC^{\theta_i}).
\end{align}

\begin{construction}
    \label{cons-quiv-va}
    \allowdisplaybreaks
    Let $Q$ be a self-dual quiver.
    We construct extra data on $\calM$ and $\calM^\sd$,
    satisfying conditions in
    Examples~\ref{eg-lin-cat} and \ref{eg-sd-lin-cat}.

    Consider the Ext-complex
    $\calExt \in \cat{Perf} (\calM \times \calM)$,
    whose fibre at $(E, F)$ gives the complex $\Ext^\bullet (E, F)$,
    as in~\cite[Construction~9.8]{Bu2023}.
    Define $\Thetaul \in \cat{Perf} (\calM \times \calM)$ by
    \begin{equation}
        \Thetaul = \calExt^\vee \oplus \sigma^* (\calExt),
    \end{equation}
    where $\sigma \colon \calM \times \calM \to \calM \times \calM$
    is the map swapping the two factors,
    and define complexes
    $\Thetaul[\dot] \in \cat{Perf} (\calM \times \calM^\sd)$ and
    $\Thetaul[\ddot] \in \cat{Perf} (\calM)$
    by~\eqref{eq-choice-theta-dot}--\eqref{eq-choice-theta-ddot}.
    Explicitly, we have
    \begin{align}
        \Thetaul & =
        \left( \ 
            \bigoplus_{a \in Q_1} \calU_{s (a)} \boxtimes
            \calU \sss{t (a)}{\vee}
            \overset{d}{\longrightarrow} 
            \biggl(
                \bigoplus_{i \in Q_0} \calU_i \boxtimes \calU_i^\vee
            \biggr)^{\oplus 2}
            \overset{d}{\longrightarrow} 
            \bigoplus_{a \in Q_1} \calU_{t (a)} \boxtimes
            \calU \sss{s (a)}{\vee}
        \ \right),
        \\[1ex]
        \Thetaul[\dot] & =
        \left( \ 
            \bigoplus_{a \in Q_1} \calU_{s (a)} \boxtimes 
            \calV \sss{t (a)}{\vee}
            \overset{d}{\longrightarrow}
            \biggl(
                \bigoplus_{i \in Q_0} \calU_i \boxtimes \calV_i^\vee
            \biggr)^{\oplus 2}
            \overset{d}{\longrightarrow} 
            \bigoplus_{a \in Q_1} \calU_{t (a)} \boxtimes
            \calV \sss{s (a)}{\vee}
        \ \right),
        \\[1ex]
        \Thetaul[\ddot] & =
        \left( \ \mathscalebox{.9}{
            \begin{aligned}[c]
                & \bigoplus_{a \in Q_1^\tria} {}
                \calU_{s (a)} \otimes \calU_{\smash{t (a)^\vee}} \\[-.5ex]
                & {} \oplus
                \bigoplus_{a \in Q_1^+} \Sym^2 (\calU_{s (a)}) \\[-.5ex]
                & {} \oplus
                \bigoplus_{a \in Q_1^-} {\wedge}^2 (\calU_{s (a)})
            \end{aligned}
            \ \overset{d}{\longrightarrow} \ 
            \left( \begin{aligned}[c]
                & \bigoplus_{i \in Q_0^\tria}
                \calU_i \otimes \calU_{\smash{i^\vee}} \\[-.5ex]
                & {} \oplus
                \bigoplus_{i \in Q_0^+} {\wedge}^2 (\calU_i) \\[-.5ex]
                & {} \oplus
                \bigoplus_{i \in Q_0^-} \Sym^2 (\calU_i)
            \end{aligned} \right)^{\oplus 2}
            \overset{d}{\longrightarrow} \ 
            \begin{aligned}[c]
                & \bigoplus_{a \in Q_1^\tria} {}
                \calU_{t (a)} \otimes \calU_{\smash{s (a)^\vee}} \\[-.5ex]
                & {} \oplus
                \bigoplus_{a \in Q_1^+} \Sym^2 (\calU_{t (a)}) \\[-.5ex]
                & {} \oplus
                \bigoplus_{a \in Q_1^-} {\wedge}^2 (\calU_{t (a)})
            \end{aligned}
        } \ \right),
    \end{align}
    where the amplitude shown here is $[-1, 1]$,
    and it is possible to write down the maps~$d$ explicitly,
    following~\cite[Construction~9.8]{Bu2023}.

    Next, we define functions $\varepsilon, \varepsilon^\sd$
    in Assumption~\ref{asn-h-sp}~\ref{itm-epsilon}
    and Assumption~\ref{asn-h-sp-mod}~\ref{itm-epsilon-sd}.
    For $\alpha, \beta \in C^\circ (Q)$ and $\theta \in C^\sd (Q)$, define
    \begin{align}
        \chi_Q (\alpha, \beta) & =
        \rank \calExt_{\alpha, \> \beta} \, , \\
        \dot{\chi}_Q (\alpha, \theta) & =
        \rank \calExt_{\alpha, \> j (\theta)} \, , \\
        \ddot{\chi}_Q (\alpha) & =
        \rank {} (\calExt_{\smash{\alpha, \alpha^\vee}})^{\bbZ_2},
    \end{align}
    where $j \colon C^\sd (Q) \hookrightarrow C^\circ (Q)$ is the inclusion,
    and $(-)^{\bbZ_2}$ means taking the fixed points of
    the $\bbZ_2$-action given by the self-dual structure.
    We have the relations
    \begin{align}
        \chi (\alpha, \beta) & =
        \chi_Q (\alpha, \beta) + \chi_Q (\beta, \alpha), \\
        \dot{\chi} (\alpha, \theta) & =
        \dot{\chi}_Q (\alpha, \theta) + \dot{\chi}_Q (\alpha^\vee, \theta), \\
        \ddot{\chi} (\alpha) & =
        \ddot{\chi}_Q (\alpha) + \ddot{\chi}_Q (\alpha^\vee),
    \end{align}
    where $\chi, \dot{\chi}, \ddot{\chi}$ are defined as the ranks of
    $\Thetaul, \Thetaul[\dot], \Thetaul[\ddot]$, respectively,
    as in Assumption~\ref{asn-h-sp}~\ref{itm-theta}
    and Assumption~\ref{asn-h-sp}~\ref{itm-theta-dots}.
    Then, as in Examples~\ref{eg-lin-cat} and~\ref{eg-sd-lin-cat}, we may take
    \begin{equation}
        \varepsilon (\alpha, \beta) =
        (-1)^{\chi_Q (\alpha, \> \beta)}, \qquad
        \varepsilon^\sd (\alpha, \beta) =
        (-1)^{\dot{\chi}_Q (\alpha, \> \theta) + \ddot{\chi}_Q (\alpha)}
    \end{equation}
    for all $\alpha, \beta \in C^\circ (Q)$ and $\theta \in C^\sd (Q)$.

    Finally, for Assumption~\ref{asn-h-sp}~\ref{itm-chi-check}
    and Assumption~\ref{asn-h-sp-mod}~\ref{itm-chi-ring}, we take
    \begin{equation}
        \hat{\chi} (\alpha) =
        \chi (\alpha, \alpha), \qquad
        \ring{\chi} (\theta) =
        \ddot{\chi} ( j (\theta) )
    \end{equation}
    for all $\alpha \in C^\circ (Q)$ and $\theta \in C^\sd (Q)$,
    as in Examples~\ref{eg-lin-cat} and~\ref{eg-sd-lin-cat}.
    One can verify using~\eqref{eq-m-quiv-expl}
    and~\eqref{eq-quiv-vsd}--\eqref{eq-quiv-gsd} that
    \begin{equation}
        \label{eq-dim-m-quiv}
        \dim_{\bbC} \calM_\alpha = -\frac{1}{2} \hat{\chi} (\alpha), \qquad
        \dim_{\bbC} \calM^\sd_\theta = -\frac{1}{2} \ring{\chi} (\theta).
    \end{equation}

    Now, we may apply Theorems~\ref{thm-va} and~\ref{thm-tw-mod}
    to obtain an involutive graded vertex algebra structure on $\hat{H}_* (\calM; \bbQ)$,
    and an involutive graded twisted module structure on $\ring{H}_* (\calM^\sd; \bbQ)$.

    Finally, Assumption~\ref{asn-rat-triv}
    follows from Joyce~\cite[\S6.4.1]{Joyce2021}
    verifying~\cite[Assumption~4.4~(g)]{Joyce2021}.
    Therefore, Theorem~\ref{thm-pl-lie-alg} gives
    an involutive graded Lie algebra structure on
    $\check{H}_* (\calM^\pl; \bbQ)$,
    and $\ring{H}_* (\calM^\sd; \bbQ)$
    is a graded twisted module for this graded Lie algebra.
\end{construction}

\begin{construction}
    \label{cons-quiv-va-der}
    In Construction~\ref{cons-quiv-va},
    we may replace $\calM$ and $\calM^\sd$
    by the derived moduli stacks $\frM$ and $\frM^\sd$,
    and replace all occurrences of $C^\circ (Q)$ and $C^\sd (Q)$
    by $K (Q)$ and $K^\sd (Q)$, respectively.

    We may construct the data
    $\Thetaul, \Thetaul[\dot], \Thetaul[\ddot]$, etc.,
    in essentially the same way,
    with the three-term complexes in Construction~\ref{cons-quiv-va}
    replaced by total complexes of double complexes supported on three columns.
    This gives rise to a graded vertex algebra structure on $\hat{H}_* (\frM; \bbQ)$,
    and a graded twisted module structure on $\ring{H}_* (\frM^\sd; \bbQ)$.

    However, as mentioned in Remark~\ref{rmk-rat-triv},
    Assumption~\ref{asn-rat-triv}~\ref{itm-rat-triv}
    may not hold for the class $0 \in K (Q)$.
    One solution to this problem is to choose
    $\frM_{\smash{>0}} = \coprod_{\smash{\alpha \in C (Q)}} \frM_\alpha$
    as in Remark~\ref{rmk-rat-triv},
    and to work with the graded Lie algebra
    $\check{H}_* (\frM^\pl_{\smash{>0}}; \bbQ)$ instead.
\end{construction}

\subsection{Enumerative invariants}
\label{sect-quiv-inv}

For self-dual quivers,
we can give a partial proof of the main conjecture,
Conjecture~\ref{conj-main}.
We state the main result as Theorem~\ref{thm-quiver-main} below.

\begin{definition}
    Let $Q$ be a quiver.
    We say that $Q$ has no \emph{oriented cycles},
    if we cannot find $n > 0$ and $a_1, \dotsc, a_n \in Q_1$\,,
    such that $t (a_i) = s (a_{i+1})$ for all $i = 1, \dotsc, n-1$,
    and $t (a_n) = s (a_1)$.
\end{definition}

\begin{proposition}
    \label{prop-quiver-properness}
    Let $Q$ be a self-dual quiver with no oriented cycles,
    and let $\tau$ be a self-dual stability function on $Q$.
    Let $\theta \in C^\sd (Q)$ be a class such that
    all semistable $\tau$-self-dual objects of class $\theta$ are $\tau$-stable.

    Then, $\smash{\calM^\sdss_\theta (\tau)}$
    is a smooth, proper Deligne--Mumford stack.
    In particular, it has a fundamental class
    \[
        \fund{\calM^\sdss_\theta (\tau)} \in
        H_{-\ring{\chi} (\theta)} (\calM^\sdss_\theta (\tau); \bbQ).
    \]
\end{proposition}

\begin{proof}
    The stack $\calM^\sdss_\theta$
    is smooth by its global quotient presentation,
    and is Deligne--Mumford by Theorem~\ref{thm-st-sd-decomp}.
    Therefore, it remains to show its properness.

    As shown by Young~\cite[Theorem~3.7]{Young2015},
    the $\tau$-semistable and $\tau$-stable loci
    in the space $V^\sd_\theta$ defined in~\eqref{eq-quiv-vsd}
    coincides with the semistable and stable loci
    for a certain geometric invariant theory~(GIT) stability,
    given by a linearization $\chi$ of the trivial line bundle on $V^\sd_\theta$.
    Let $M^\sdss_\theta (\tau) = V^\sd_\theta \dslashsub{\chi} G^\sd_\theta$
    denote this GIT quotient, which is a quasi-projective $\bbC$-variety.

    Since $Q$ has no oriented cycles, 
    if $\tau_0$ denotes the trivial stability condition on $Q$,
    then there is a unique \emph{$\tau_0$-polystable} self-dual object of class $\theta$.
    Here, a $\tau_0$-semistable self-dual representation
    being $\tau_0$-polystable means that
    it is isomorphic to the associated graded object
    of its $\tau_0$-Jordan--Hölder filtration.
    It follows that there is a unique $S$-equivalence class,
    in the sense of Young~\cite[\S3.2, p.~452]{Young2015},
    so that the affine GIT quotient $V^\sd_\theta \dslash G^\sd_\theta$
    is a point.

    As in King~\cite[Definition~2.1]{King1994},
    this implies that $M^\sdss_\theta (\tau)$ is a projective variety.
    By the Keel--Mori theorem~\cite{KeelMori1997},
    also stated in~\cite{Conrad2005},
    the projection $\calM^\sdss_\theta (\tau) \to M^\sdss_\theta (\tau)$
    is proper, which implies that $\calM^\sdss_\theta (\tau)$ is proper.
\end{proof}

\begin{theorem}
    \label{thm-quiver-main}
    Let $Q$ be a self-dual quiver with no oriented cycles.
    Then, there is a unique way to define homology classes
    \begin{align*}
        \inv{ \calM^\sdss_{\theta} (\tau) }
        & \in H_{-\ring{\chi} (\theta)} (\calM^\sd_{\theta}; \bbQ)
    \end{align*}
    for all self-dual stability conditions $\tau$ on $\cat{Mod} (\bbC Q),$
    and all classes $\theta \in C^\sd (Q),$
    satisfying the following conditions:

    \begin{enumerate}
        \item 
            \label{itm-quiver-fund}
            For any tame class $\theta \in C^\sd (Q)$
            in the sense of Definition~\textnormal{\ref{def-tame-class}}
            below, we have
            \begin{equation}
                \inv{ \calM^\sdss_{\theta} (\tau) } =
                \iota_* \fund{ \calM^\sdss_{\theta} (\tau) } \, ,
            \end{equation}
            where $\iota \colon \calM^\sdss_{\theta} (\tau)
            \hookrightarrow \calM^\sd_\theta$
            denotes the inclusion.
            
        \item 
            \label{itm-quiver-wcf}
            For any two self-dual stability conditions $\tau, \tilde{\tau}$ on $Q$,
            and any $\theta \in C^\sd (Q),$
            we have the wall-crossing formula
            \begin{align*}
                \hspace{2em} & \hspace{-2em}
                \inv{ \calM^\sdss_\theta (\tilde{\tau}) } = {}
                \numberthis
                \label{eq-wcf-quiver-sd}
                \\[.5ex] &
                \sum_{ \leftsubstack[8em]{
                    \\[-3ex]
                    & n \geq 0; \, m_1, \dotsc, m_n > 0; \\[-.5ex]
                    & \alpha_{1,1}, \dotsc, \alpha_{1,m_1}; \dotsc;
                    \alpha_{n,1}, \dotsc, \alpha_{n,m_n} \in C (Q); \,
                    \rho \in \smash{C^\sd (Q)} \colon \\[-.5ex]
                    & \theta = (\bar{\alpha}_{1,1} + \cdots + \bar{\alpha}_{1,m_1})
                    + \cdots + (\bar{\alpha}_{n,1} + \cdots + \bar{\alpha}_{n,m_n}) + \rho
                } } {} 
                \tilde{U}^\sd (\alpha_{1,1}, \dotsc, \alpha_{1,m_1}; \dotsc;
                \alpha_{n,1}, \dotsc, \alpha_{n,m_n}; \tau, \tilde{\tau}) \cdot {} 
                \\[.5ex] &
                \bigl[ \bigl[ \bigl[ 
                \inv{ \calM^\plss_{\alpha_{1,1}} (\tau) } \, , 
                \inv{ \calM^\plss_{\alpha_{1,2}} (\tau) } \bigr] , \dotsc \bigr] ,
                \inv{ \calM^\plss_{\alpha_{1,m_1}} (\tau) } \bigr] \heart \cdots \heart {}
                \\ 
                &
                \bigl[ \bigl[ \bigl[ 
                \inv{ \calM^\plss_{\alpha_{n,1}} (\tau) } \, , 
                \inv{ \calM^\plss_{\alpha_{n,2}} (\tau) } \bigr] , \dotsc \bigr] ,
                \inv{ \calM^\plss_{\alpha_{n,m_n}} (\tau) } \bigr] \heart 
                \inv{ \calM^\sdss_{\rho} (\tau) } \, .
            \end{align*}
            Here, the $\inv{\calM^\plss_{\alpha_{i, \> j}} (\tau)}$
            are Joyce's invariants defined
            in~\textnormal{\cite{GrossJoyceTanaka, Joyce2021}}.
    \end{enumerate}
\end{theorem}

The proof will be given in~\S\ref{sect-proof-quiver}.

\begin{remark}
    In Theorem~\ref{thm-quiver-main}~\ref{itm-quiver-fund},
    we hope to weaken the assumption that $\theta$ is a tame class
    to `all $\tau$-semistable self-dual objects of class $\theta$ are $\tau$-stable'.
    This still guarantees the existence of a fundamental class,
    by Proposition~\ref{prop-quiver-properness}.
\end{remark}

In the remainder of this subsection,
we formulate the definition of \emph{tame classes},
which are classes that we can prove
Theorem~\ref{thm-quiver-main}~\ref{itm-quiver-fund} for.

\begin{definition}
    \label{def-q-tilde}
    Let $Q$ be a self-dual quiver, and let $\theta \in C^\sd (Q)$.
    We define a new self-dual quiver $\tilde{Q}_\theta$ as follows.
    
    For convenience in notation, for each integer $n \geq 0$,
    write $\langle n \rangle$ for the set $\{ 1, \dotsc, n \}$,
    equipped with an involution
    $(-)^\vee \colon \langle n \rangle \to \langle n \rangle$,
    $j \mapsto n + 1 - j$.

    \begin{itemize}
        \item 
            Define the set of vertices to be
            $(\tilde{Q}_\theta)_0 = \{ (i, j) \mid i \in Q_0 \, , \ 
            j \in \langle \theta_i \rangle \}$.
            The involution is given by
            $(i, j)^\vee = (i^\vee, j^\vee)$.
            The vertex signs are $u_{(i, \> j)} = u_i$\,.
        \item
            Define the set of edges to be
            $(\tilde{Q}_\theta)_1 = \{ (a, j, k) \mid a \in Q_1 \, , \ 
            j \in \langle \theta_{s (a)} \rangle, \ 
            k \in \langle \theta_{t (a)} \rangle \}
            \setminus \{ (a, j, j^\vee) \mid a \in Q_1^-, \ 
            j \in \langle \theta_{s (a)} \rangle \}$.
            The source and target maps are
            $s (a, j, k) = (s (a), j)$ and
            $t (a, j, k) = (t (a), k)$.
            The involution is given by
            $(a, j, k)^\vee = (a^\vee, k^\vee, j^\vee)$.
            The edge signs are $v_{(a, \> j, \> k)} = v_a$\,.
    \end{itemize}

    In plain words, $\tilde{Q}_\theta$ is obtained from $Q$
    by splitting each vertex $i \in Q_0$ into $\theta_i$ new vertices,
    and splitting each arrow $a \in Q_1$ into new arrows connecting the new vertices,
    but we omit a new edge when $a \in Q_1^-$ and the new edge is self-dual.

    For a self-dual stability condition $\tau$ for $Q$,
    we define a self-dual stability condition $\tilde{\tau}$ for $\tilde{Q}_\theta$
    by setting
    \[
        \tilde{\tau} (\alpha) = \tau (\alpha_Q)
    \]
    for all $\alpha \in C (\tilde{Q}_\theta)$,
    where $\alpha_Q \in C (Q)$ is the class with
    $(\alpha_Q)_i = \sum_j \alpha_{(i, \> j)}$ for all $i \in Q_0$\,.
\end{definition}

\begin{definition}
    \label{def-tame-class}
    Let $Q$ be a self-dual quiver,
    and let $\tau$ be a self-dual stability condition for $Q$.
    We say that a class $\theta \in C^\sd (Q)$ is
    \emph{tame} with respect to $\tau$, if the following conditions are satisfied.
    \begin{enumerate}
        \item 
            \label{itm-tame-base}
            All $\tau$-semistable self-dual objects of class~$\theta$
            are $\tau$-stable.
        \item
            \label{itm-tame-extra}
            All $\tilde{\tau}$-semistable self-dual representations
            of~$\tilde{Q}_\theta$
            of class~$1 \in C^\sd (\tilde{Q}_\theta)$
            are $\tilde{\tau}$-stable.
            Here, $1$~is the class with
            $1_{(i, \> j)} = 1$ for all
            $(i, j) \in (\tilde{Q}_\theta)_0$\,.
    \end{enumerate}
    In particular, if $\theta_i \leq 1$ for all~$i$,
    then~\ref{itm-tame-extra} will not be needed.
\end{definition}

It is precisely the classes satisfying
Definition~\ref{def-tame-class}~\ref{itm-tame-base},
but not~\ref{itm-tame-extra}, that are causing problems for us,
in the sense that if we could also prove
Theorem~\ref{thm-quiver-main}~\ref{itm-quiver-fund} for these classes,
we would have a complete proof of Conjecture~\ref{conj-main} for self-dual quivers.

\begin{example}
    We discuss two examples demonstrating the notion of tame classes.

    \begin{enumerate}
        \item 
            Consider the self-dual quiver
            \[
                Q = \Biggl(
                    \begin{tikzcd}[column sep={1.2em,between origins}, row sep={1.2em,between origins}, nodes={inner sep=.2em}, baseline=-.5ex]
                        & \scriptstyle \bullet \ar[dr] & {} \\
                        \scriptstyle \bullet \ar[ur] \ar[dr] && \scriptstyle \bullet \\
                        & \scriptstyle \bullet \ar[ur]
                    \end{tikzcd}
                \Biggr),
            \]
            with the involution given by horizontal flipping,
            and use the self-dual stability function
            $\smash{\bigl( \begin{smallmatrix}
                & 0 \\[-.5ex] +1 && -1 \\[-.5ex] & 0
            \end{smallmatrix} \bigr)}$.
            Then the class
            $\theta = \smash{\bigl( \begin{smallmatrix}
                & 1 \\[-.5ex] 2 && 2 \\[-.5ex] & 1
            \end{smallmatrix} \bigr)} \in C^\sd (Q)$ is tame.

            Note that there are stable self-dual objects $E_1, E_2$ of classes
            $\smash{\bigl( \begin{smallmatrix}
                & 1 \\[-.5ex] 1 && 1 \\[-.5ex] & 0
            \end{smallmatrix} \bigr)}$
            and $\smash{\bigl( \begin{smallmatrix}
                & 0 \\[-.5ex] 1 && 1 \\[-.5ex] & 1
            \end{smallmatrix} \bigr)}$,
            and $E = E_1 \oplus E_2$ is a self-dual object of class $\theta$,
            but this does not contradict with
            Definition~\ref{def-tame-class}~\ref{itm-tame-base}
            by Theorem~\ref{thm-st-sd-decomp}.

        \item 
            Consider the self-dual quiver
            \[
                Q = \Bigl(
                    \begin{tikzcd}[column sep={2em,between origins}, nodes={inner sep=.3em}, baseline=-.5ex]
                        \scriptstyle \bullet 
                        & \scriptstyle \bullet \ar[l] \ar[r]
                        & \scriptstyle \bullet
                        & \scriptstyle \bullet \ar[l]
                    \end{tikzcd}
                \Bigr),
            \]
            with stability function $(-3, +2, -2, +3)$.
            Then the class $\theta = (1, 3, 3, 1) \in C^\sd (Q)$ is not tame.
            It satisfies Definition~\ref{def-tame-class}~\ref{itm-tame-base},
            as there are no semistable representations of class $\theta$.
            However, one can find
            semistable self-dual representations $\tilde{Q}_\theta$
            of class $1$ that are not stable.
    \end{enumerate}
\end{example}

\subsection{Morphisms of quivers}

We discuss a slightly unusual notion of morphisms of quivers,
introduced by Gross--Joyce--Tanaka~\cite[\S5.4]{GrossJoyceTanaka},
and we extend their definition to self-dual quivers.
Such morphisms interact well with enumerative invariants,
and will be useful in the proof of Theorem~\ref{thm-quiver-main}.

\begin{definition}
    \label{def-mor-quiv}
    Let $Q = (Q_0, Q_1, s, t)$
    and $Q' = (Q'_0, Q'_1, s', t')$
    be quivers.
    A \emph{morphism of quivers} $\lambda \colon Q \to Q'$ consists of the following data:

    \begin{itemize}
        \item 
            A map $\lambda \colon Q_0 \to Q'_0$\,.
        \item
            A subset $Q_1^\circ \subset Q_1$\,, and a map $\lambda \colon Q_1^\circ \to Q'_1$\,.
    \end{itemize}
    They should satisfy the following conditions:

    \begin{itemize}
        \item 
            For any $a \in Q_1^\circ$\,, we have $\lambda (s (a)) = s' (\lambda (a))$ and
            $\lambda (t (a)) = t' (\lambda (a))$.

        \item
            For any $i, j \in Q_0$ and $a' \in Q'_1$
            with $s' (a') = \lambda (i)$ and $t' (a') = \lambda (j)$,
            there exists a unique $a \in Q_1^\circ$ with $\lambda (a) = a'$.
    \end{itemize}
    If $\lambda \colon Q \to Q'$ and $\lambda' \colon Q' \to Q \second$
    are morphisms of quivers, then we have a \emph{composition}
    $\lambda' \circ \lambda \colon Q \to Q \second$.
\end{definition}

\begin{definition}
    \label{def-mor-quiv-fun}
    Let $\lambda \colon Q \to Q'$ be a morphism of quivers.
    Define a map
    \[
        \lambda_* \colon C^\circ (Q)
        \longrightarrow C^\circ (Q')
    \]
    by setting $\lambda_* (\alpha)_{i'} =
    \sum_{i \in \lambda^{-1} (i')} \alpha_i$
    for all $i' \in Q'_0$\,.
    Define a functor
    \[
        \lambda_* \colon
        \cat{Mod} (\bbC Q) \longrightarrow \cat{Mod} (\bbC Q')
    \]
    as follows.
    For each object $E = ((E_i)_{i \in Q_0} \, , (e_a)_{a \in Q_1}) \in \cat{Mod} (\bbC Q)$,
    define an object
    $\lambda_* (E) = ((E'_{\smash{i'}})_{i' \in Q'_0} \, ,
    (e'_{\smash{a'}})_{a' \in Q'_1}) \in \cat{Mod} (\bbC Q')$ by
    \[
        E'_{\smash{i'}} = \bigoplus_{i \in \lambda^{-1} (i')} E_i \, , \qquad
        e'_{\smash{a'}} = \sum_{a \in \lambda^{-1} (a')} e_a 
    \]
    for all $i' \in Q'_0$ and $a' \in Q'_1$\,.
    One can verify that this defines a morphism of moduli stacks
    \[
        F_\lambda \colon \calM \longrightarrow \calM' ,
    \]
    where $\calM, \calM'$ are moduli stacks of representations of $Q, Q'$,
    respectively.
    This also extends to a morphism of $\bbC$-linear stacks
    as in~\cite[\S9.2]{Bu2023}.
\end{definition}

\begin{definition}
    Let $\lambda \colon Q \to Q'$ be a morphism of quivers,
    and let $\tau$ be a stability condition on $\cat{Mod} (\bbC Q')$.
    Define a stability condition $\lambda^* (\tau)$ on $\cat{Mod} (\bbC Q)$,
    called the \emph{pullback} of $\tau$,
    by setting $\lambda^* (\tau) (\alpha) = \tau (\lambda_* (\alpha))$
    for all $\alpha \in C (Q)$.
\end{definition}

Next, we discuss \emph{morphisms of self-dual quivers}.

\begin{definition}
    \label{def-mor-quiv-sd}
    Let $Q = (Q_0, Q_1, s, t, \sigma, u, v)$
    and $Q' = (Q'_0, Q'_1, s', t', \sigma', u', v')$
    be self-dual quivers.
    A \emph{morphism of self-dual quivers} $\lambda \colon Q \to Q'$
    is a morphism of quivers,
    such that 
    
    \begin{itemize}
        \item 
            For all $i \in Q_0$\,, we have
            $\lambda (i^\vee) = \lambda (i)^\vee$
            and $u \sss{\lambda (i)}{\prime} = u_i$\,.
        \item
            For all $a \in Q_1^\circ$, we have $a^\vee \in Q_1^\circ$,
            $\lambda (a^\vee) = \lambda (a)^\vee$,
            and $v \sss{\lambda (a)}{\prime} = v_a$\,.
    \end{itemize}
    One can verify that the induced morphism of moduli stacks
    $\lambda_* \colon \calM^+ \to \calM'^+$
    preserves the self-dual structures,
    i.e.~induces self-dual functors after taking values on affine schemes.
\end{definition}

\begin{construction}
    \label{cons-quiv-omega}
    \allowdisplaybreaks
    Let $\lambda \colon Q \to Q'$
    be a morphism of self-dual quivers.
    Let $\calM, \calM'$ be the moduli stacks of representations of $Q, Q'$, and let
    \begin{align*}
        F_\lambda \colon \calM & \longrightarrow \calM', \\
        F^\sd_\lambda \colon \calM^\sd & \longrightarrow \calM'^\sd
    \end{align*}
    be the induced morphisms of moduli stacks.
    We define extra data for these morphisms,
    satisfying Assumptions~\ref{asn-h-sp-mor}
    and~\ref{asn-h-sp-mod-mor}, as follows.

    For $\alpha, \beta \in C (Q)$ and $\theta \in C^\sd (Q)$, define
    \begin{align}
        \Xi_{\alpha, \> \beta} & =
        \sum_{ \substack{ i, \> j \in Q_0 : \\ i \neq j, \, \lambda (i) = \lambda (j) } } {}
        \calU_i \boxtimes \calU_j^\vee 
        \ + \ 
        \sum_{ a \in Q_1 \setminus Q_1^\circ } {}
        \calU_{s (a)} \boxtimes \calU \sss{t (a)}{\vee} \, ,
        \\[.5ex]
        \hat{\Xi}_{\alpha} & =
        \sum_{ \substack{ i, \> j \in Q_0 : \\ i \neq j, \, \lambda (i) = \lambda (j) } } {}
        \calU_i \otimes \calU_j^\vee 
        \ + \ 
        \sum_{ a \in Q_1 \setminus Q_1^\circ } {}
        \calU_{s (a)} \otimes \calU \sss{t (a)}{\vee} \, ,
        \\[.5ex]
        \dot{\Xi}_{\alpha, \theta} & =
        \sum_{ \substack{ i, \> j \in Q_0 : \\ i \neq j, \, \lambda (i) = \lambda (j) } } {}
        \calU_i \boxtimes \calV_j^\vee 
        \ + \ 
        \sum_{ a \in Q_1 \setminus Q_1^\circ } {}
        \calU_{s (a)} \boxtimes \calV \sss{t (a)}{\vee} \, ,
        \\[.5ex]
        \ddot{\Xi}_{\alpha} & =
        \sum_{ \substack{
            i, \> j \in Q_0 : \\
            i < j, \, i \neq j^\vee, \\
            \lambda (i) = \lambda (j)^\vee
        } } {}
        \calU_i \otimes \calU_j
        \ + \ 
        \sum_{ \substack{
            i \in Q_0 : \,
            i \neq i^\vee, \\
            \lambda (i) \in Q'^+_0
        } } {}
        {\wedge}^2 (\calU_i)
        \ + \ 
        \sum_{ \substack{
            i \in Q_0 : \,
            i \neq i^\vee, \\
            \lambda (i) \in Q'^-_0
        } } {}
        \Sym^2 (\calU_i)
        \ + {}
        \notag \\*
        & \hspace{2em} 
        \sum_{ a \in Q_1^\tria \setminus Q_1^\circ } {} \hspace{-.5em}
        \calU_{s (a)} \otimes \calU_{t (a)^\vee}
        \ + 
        \sum_{ a \in Q_1^+ \setminus Q_1^\circ } {} \hspace{-.5em}
        \Sym^2 (\calU_{s (a)})
        \ + 
        \sum_{ a \in Q_1^- \setminus Q_1^\circ } {} \hspace{-.5em}
        {\wedge}^2 (\calU_{s (a)}) ,
        \\[.5ex]
        \ring{\Xi}_{\theta} & =
        \sum_{ \substack{
            i, \> j \in Q_0 : \\
            i < j, \, i \neq j^\vee, \\
            \lambda (i) = \lambda (j)^\vee
        } } {}
        \calV_i \otimes \calV_j
        \ + \ 
        \sum_{ \substack{
            i \in Q_0 : \,
            i \neq i^\vee, \\
            \lambda (i) \in Q'^+_0
        } } {}
        {\wedge}^2 (\calV_i)
        \ + \ 
        \sum_{ \substack{
            i \in Q_0 : \,
            i \neq i^\vee, \\
            \lambda (i) \in Q'^-_0
        } } {}
        \Sym^2 (\calV_i)
        \ + {}
        \notag \\*
        \label{eq-def-xi-ring}
        & \hspace{2em} 
        \sum_{ a \in Q_1^\tria \setminus Q_1^\circ } {} \hspace{-.5em}
        \calV_{s (a)} \otimes \calV_{t (a)^\vee}
        \ + 
        \sum_{ a \in Q_1^+ \setminus Q_1^\circ } {} \hspace{-.5em}
        \Sym^2 (\calV_{s (a)})
        \ + 
        \sum_{ a \in Q_1^- \setminus Q_1^\circ } {} \hspace{-.5em}
        {\wedge}^2 (\calV_{s (a)}) .
    \end{align}
    One can verify that this data satisfies
    Assumptions~\ref{asn-h-sp-mor} and~\ref{asn-h-sp-mod-mor}.

    By Theorems~\ref{thm-mor-va} and~\ref{thm-mor-tw-mod},
    we obtain maps
    \begin{alignat}{2}
        \Omega_\lambda \colon &&
        \hat{H}_* (\calM; \bbQ) & \longrightarrow
        \hat{H}_* (\calM'; \bbQ), \\
        \Omega^\pl_\lambda \colon &&
        \check{H}_* (\calM^\pl; \bbQ) & \longrightarrow
        \check{H}_* (\calM'^\pl; \bbQ), \\
        \Omega^\sd_\lambda \colon &&
        \ring{H}_* (\calM^\sd; \bbQ) & \longrightarrow
        \ring{H}_* (\calM'^\sd; \bbQ),
    \end{alignat}
    where $\Omega$ is a morphism of graded vertex algebras,
    $\Omega^\pl$ a morphism of graded Lie algebras,
    and $\Omega^\sd$ is compatible with the twisted module structures.
\end{construction}

\begin{theorem}
    \label{thm-omega-comp}
    Let $Q, Q', Q \second$ be quivers,
    and let $\lambda \colon Q \to Q',$
    $\lambda' \colon Q' \to Q \second$ be morphisms.
    Denote $\lambda \second = \lambda' \circ \lambda$.
    
    \begin{enumerate}
        \item
            The maps $\Omega_\lambda, \Omega_{\lambda'}, \Omega_{\lambda \second}$
            induced by $\lambda, \lambda', \lambda \second,$
            as in Theorem~\textnormal{\ref{thm-mor-va},} satisfy
            \begin{equation}
                \Omega_{\lambda \second} = 
                \Omega_{\lambda'} \circ \Omega_\lambda \, .
            \end{equation}
            In particular, we also have
            \begin{equation}
                \Omega^\pl_{\lambda \second} = 
                \Omega^\pl_{\lambda'} \circ \Omega^\pl_\lambda \, .
            \end{equation}
        \item
            Suppose, moreover, that $Q, Q', Q \second$ are self-dual quivers,
            and $\lambda, \lambda'$ are morphisms of self-dual quivers.
            Then the maps $\Omega_\lambda^\sd,
            \Omega^\sd_{\lambda \smash{'}}, \Omega^\sd_{\lambda \smash{\second}}$
            induced by $\lambda, \lambda', \lambda \second,$
            as in Theorem~\textnormal{\ref{thm-mor-tw-mod},} satisfy
            \begin{equation}
                \Omega^\sd_{\lambda \smash{\second}} = 
                \Omega^\sd_{\lambda \smash{'}} \circ \Omega^\sd_\lambda.
            \end{equation}
    \end{enumerate}
\end{theorem}

\begin{proof}
    One can verify that the data $\Xi, \hat{\Xi}, \dot{\Xi}, \ddot{\Xi}, \ring{\Xi}$
    chosen in Construction~\ref{cons-quiv-omega}
    satisfies the relations specified in
    Theorems~\ref{thm-comp-mor-va} and~\ref{thm-comp-mor-tw-mod},
    and the result follows from these theorems.
\end{proof}

\subsection{Computing the vertex algebra and twisted module}

Let $Q$ be a self-dual quiver.
We explicitly write down the graded vertex algebra and graded twisted module
defined in Construction~\ref{cons-quiv-va-der}.

\begin{theorem}
    \label{thm-quiv-ho-eq}
    We have homotopy equivalences
    \begin{align}
        \label{eq-ho-type-m}
        |\calM_\alpha| & \simeq
        \prod_{i \in Q_0} {} \BU (\alpha_i), \\
        \label{eq-ho-type-msd}
        |\calM^\sd_\theta| & \simeq
        \prod_{i \in Q_0^\tria} {} \BU (\theta_i) \times
        \prod_{i \in Q_0^+} {} \BO (\theta_i) \times
        \prod_{i \in Q_0^-} {} \BSp (\theta_i), \\
        \label{eq-ho-type-m-der}
        |\frM_\alpha| & \simeq
        \prod_{i \in Q_0} {} \BU, \\
        \label{eq-ho-type-msd-der}
        |\frM^\sd_\theta| & \simeq
        \prod_{i \in Q_0^\tria} {} \BU \times
        \prod_{i \in Q_0^+} {} |\OPerf_{\theta_i}| \times
        \prod_{i \in Q_0^-} {} |\SpPerf_{\theta_i}|.
    \end{align}
    Here, $\OPerf_{\smash{\theta_i}}$ and~$\SpPerf_{\smash{\theta_i}}$
    are the rank~$\theta_i$ components in
    the derived stacks $\OPerf$ and~$\SpPerf$
    in Definition~\textnormal{\ref{def-operf}}.
\end{theorem}

\begin{proof}
    Equations~\eqref{eq-ho-type-m}--\eqref{eq-ho-type-msd}
    follow from Example~\ref{eg-quot-st-top}
    and the explicit descriptions of $\calM_\alpha$
    and~$\calM^\sd_\theta$ as quotient stacks,
    as in \eqref{eq-m-quiv-expl} and~\eqref{eq-quiv-vsd}--\eqref{eq-quiv-gsd}.

    For~\eqref{eq-ho-type-m-der}--\eqref{eq-ho-type-msd-der},
    note that the homotopy type of $|\frM_\alpha|$
    does not depend on the edges of $Q$,
    since the closed substack consisting of objects
    $\smash{((E_i)_{i \in Q_0} \, , (e_a)_{a \in Q_1})}$
    with $e_a = 0$ for all $a$ is \Aone-homotopy equivalent to $\frM_\alpha$\,,
    since $\frM_\alpha$ can be retracted into this substack
    by scaling the $e_a$\,.
    We thus have an \Aone-homotopy equivalence
    $\frM_\alpha \simeq \smash{\prod_{i \in Q_0} \cat{Perf}_{\alpha_i}}$\,,
    where $\cat{Perf}_{\smash{\alpha_i}} \subset \cat{Perf}$
    is the rank~$\alpha_i$ component.
    Now, by~\eqref{eq-perf-bu} and the \Aone-homotopy invariance
    of the topological realization,
    we obtain~\eqref{eq-ho-type-m-der}.
    A similar argument establishes~\eqref{eq-ho-type-msd-der}.
\end{proof}

We would like to compute the homology of the spaces in
\eqref{eq-ho-type-m}--\eqref{eq-ho-type-msd-der}.
However, it is not straightforward to describe the spaces
$|\OPerf|$ and~$|\SpPerf|$, unless we know Conjecture~\ref{conj-homology}.
For this reason, we define an alternative version of $|\frM^\sd_\theta|$
which is easier to compute.

\begin{construction}
    \label{cons-msd-der-alt}
    Define a space $|\frM|^\sd$,
    as opposed to $|\frM^\sd|$, by
    \begin{equation}
        |\frM|^\sd = \coprod_{\theta \in K^\sd (Q) \vphantom{^0}} |\frM|^{\sd}_\theta,
    \end{equation}
    where
    \begin{equation}
        |\frM|^{\sd}_\theta =
        \prod_{i \in Q_0^\tria} {} \BU \times
        \prod_{i \in Q_0^+} {} \BO \times
        \prod_{i \in Q_0^-} {} \BSp.
    \end{equation}
    Note that the space $|\frM|^{\sd}_\theta$ does not depend on $\theta$.
    This space will be homotopy equivalent to $|\frM^\sd_\theta|$
    if Conjecture~\ref{conj-homology} is true.

    There is a natural map ${\oplus}^\sd \colon
    |\frM_\alpha| \times |\frM|^{\sd}_\theta
    \to |\frM|^{\sd}_{\bar{\alpha} + \theta}$\,,
    defined based on the maps $\oplus \colon \BU \times \BU \to \BU$,
    $\oplus^\sd \colon \BU \times \BO \to \BO$, and
    $\oplus^\sd \colon \BU \times \BSp \to \BSp$.
    This makes $|\frM|^\sd$ into an involutive module
    for the involutive commutative $H$-space $|\frM|$,
    in the sense of Definition~\ref{def-h-sp-invol}.

    Also, there are universal $K$-theory classes
    \begin{equation}
        \label{eq-quiv-univ-vb-alt}
        \calV_i \in K (|\frM|^\sd)
    \end{equation}
    for $i \in Q_0$\,, of rank $\theta_i$ on the component $|\frM|^\sd_\theta$\,.
    This is analogous to the universal perfect complexes $\calV_i$
    defined in~\eqref{eq-quiv-univ-vb},
    and can be constructed from universal $K$-theory classes
    on $\BU$, $\BO$, and $\BSp$,
    given by their natural maps to $\BU$.

    We can also define $\dot{\Theta} \in K (|\frM| \times |\frM|^\sd)$,
    similarly as in Constructions~\ref{cons-quiv-va}
    and~\ref{cons-quiv-va-der}, by
    \begin{equation}
        \dot{\Theta} =
        2 \sum_{i \in Q_0} \calU_i \boxtimes \calV \sss{i}{\vee}
        - \sum_{a \in Q_1} \calU_{s (a)} \boxtimes 
        \calV \sss{t (a)}{\vee}
        - \sum_{a \in Q_1} \calU_{t (a)} \boxtimes 
        \calV \sss{s (a)}{\vee} \, .
    \end{equation}
    We then choose other data including
    $\Theta, \ddot{\Theta}, \hat{\chi}, \ring{\chi}, \varepsilon, \varepsilon^\sd$
    in the same way as in Construction~\ref{cons-quiv-va-der}.

    Applying Theorem~\ref{thm-tw-mod},
    we obtain an involutive graded twisted module structure on
    the graded vector space~$\ring{H}_* (|\frM|^\sd; \bbQ)$,
    for the involutive graded vertex algebra $\hat{H}_* (\frM; \bbQ)$.
\end{construction}

\begin{remark}
    It is known that we have
    \begin{equation}
        \label{eq-bo-as-fixed-points}
        \BO \times \bbZ \simeq (\BU \times \bbZ)^{\bbZ_2},
    \end{equation} 
    where the right-hand side is the homotopy fixed points
    of the $\bbZ_2$-action given by taking the complex conjugate,
    which is equivalent to taking the dual vector bundle.
    A proof of the equivalence~\eqref{eq-bo-as-fixed-points} can be found in
    Dugger~\cite[Corollary~7.6]{Dugger2005}.

    By Bott periodicity,
    we have $\Omega^4 (\BO \times \bbZ) \simeq \BSp \times \bbZ$
    and $\Omega^4 (\BU \times \bbZ) \simeq \BU \times \bbZ$.
    Therefore, we may write
    \begin{equation}
        \BSp \times \bbZ \simeq (\BU \times \bbZ)^{\bbZ_2}
    \end{equation}
    for some $\bbZ_2$-action on $\BU \times \bbZ$.
    We expect that this coincides with the $\bbZ_2$-action
    given by taking the dual vector bundle,
    but introducing a sign $-1$ when identifying
    the double dual of a vector bundle with itself.

    We thus have $|\frM|^\sd \simeq |\frM|^{\bbZ_2}$,
    where the right-hand side denotes the homotopy $\bbZ_2$-fixed points.
\end{remark}

\begin{definition}
    For each $i \in Q_0$ and $k \geq 0$, define cohomology classes
    \begin{alignat*}{2}
        S \sss{i, \> k}{\alpha} & \in H^{2k} (\calM_\alpha; \bbQ), & \qquad
        S \sss{i, \> k}{\theta} & \in H^{2k} (\calM^\sd_\theta; \bbQ), \\
        S \sss{i, \> k}{\alpha} & \in H^{2k} (\frM_\alpha; \bbQ), & \qquad
        S \sss{i, \> k}{\theta} & \in H^{2k} (|\frM|^\sd_\theta; \bbQ)
    \end{alignat*}
    by
    \begin{equation}
        S \sss{i, \> k}{\alpha} = \ch_k (\calU_i), \qquad
        S \sss{i, \> k}{\theta} = \ch_k (\calV_i),
    \end{equation}
    where $\calU_i$ and $\calV_i$
    are as in~\eqref{eq-quiv-univ-vb} or~\eqref{eq-quiv-univ-vb-alt},
    depending on the base space.
    In particular, we have
    \begin{equation}
        S \sss{i, \> 0}{\alpha} = \alpha_i \cdot 1^\alpha , \qquad
        S \sss{i, \> 0}{\theta} = \theta_i \cdot 1^\theta ,
    \end{equation}
    where $1^\alpha$ and $1^\theta$
    are the unit elements in
    $H^0 (\calM_\alpha; \bbQ)$, $H^0 (\frM_\alpha; \bbQ)$,
    $H^0 (\calM^\sd_\theta; \bbQ)$, or $H^0 (|\frM|^\sd_\theta; \bbQ)$.
\end{definition}

\begin{theorem}
    \label{thm-quiv-cohom}
    We have isomorphisms of graded $\bbQ$-algebras
    \begin{align*}
        H^* (\calM_\alpha; \bbQ) \simeq
        \bbQ [ \, S \sss{i, \> k}{\alpha}
            & : i \in Q_0 \, , \, 0 < k \leq \alpha_i \, ], 
        \numberthis \\
        H^* (\calM^\sd_\theta; \bbQ) \simeq
        \bbQ [ \, S \sss{i, \> k}{\theta}
            & : i \in Q_0^{\mathrlap{\tria}} \, , \, 0 < k \leq \theta_i \, ; \\
            S \sss{i, \> 2k}{\theta}
            & : i \in Q_0^{\mathrlap{\pm}} \, , \, 0 < k \leq \theta_i / 2 \, ],
        \numberthis \\
        H^* (\frM_\alpha; \bbQ) \simeq
        \bbQ [ \, S \sss{i, \> k}{\alpha}
            & : i \in Q_0 \, , \, k > 0 \, ],
        \numberthis \\
        H^* (|\frM|^{\sd}_\theta; \bbQ) \simeq
        \bbQ [ \, S \sss{i, \> k}{\theta}
            & : i \in Q_0^{\mathrlap{\tria}} \, , \, k > 0 \, ; \\
            S \sss{i, \> 2k}{\theta}
            & : i \in Q_0^{\mathrlap{\pm}} \, , \, k > 0 \, ],
        \numberthis
    \end{align*}
    where $\deg S \sss{i, \> k}{\alpha} = \deg S \sss{i, \> k}{\theta} = 2k$.
\end{theorem}

\begin{proof}
    This follows from Theorem~\ref{thm-quiv-ho-eq},
    together with computations of
    $H^* (\BU (n); \bbQ)$,
    $H^* (\BO (n); \bbQ)$, and
    $H^* (\BSp (2n); \bbQ)$,
    as done in Borel~\cite[Examples~(1)--(4), p.\,200]{Borel1953}.
    The latter two isomorphisms follow from
    taking the colimits when $n \to \infty$,
    which commutes with taking homology, as in May~\cite[\S14.6]{May1999}.
\end{proof}

\begin{definition}
    \label{def-quiv-homol}
    By Theorem~\ref{thm-quiv-cohom}
    and by duality, we define isomorphisms of graded $\bbQ$-vector spaces
    \begin{align*}
        H_* (\calM_\alpha; \bbQ) \simeq
        \bbQ [ \, s \sss{i, \> k}{\alpha}
            & : i \in Q_0 \, , \, 0 < k \leq \alpha_i \, ],
        \numberthis \\
        H_* (\calM^\sd_\theta; \bbQ) \simeq
        \bbQ [ \, s \sss{i, \> k}{\theta}
            & : i \in Q_0^{\mathrlap{\tria}} \, , \, 0 < k \leq \theta_i \, ; \\
            s \sss{i, \> 2k}{\theta}
            & : i \in Q_0^{\mathrlap{\pm}} \, , \, 0 < k \leq \theta_i / 2 \, ],
        \numberthis \\
        H_* (\frM_\alpha; \bbQ) \simeq
        \bbQ [ \, s \sss{i, \> k}{\alpha}
            & : i \in Q_0 \, , \, k > 0 \, ],
        \numberthis \\
        H_* (\frM^\sd_\theta; \bbQ) \simeq
        \bbQ [ \, s \sss{i, \> k}{\theta}
            & : i \in Q_0^{\mathrlap{\tria}} \, , \, k > 0 \, ; \\
            s \sss{i, \> 2k}{\theta}
            & : i \in Q_0^{\mathrlap{\pm}} \, , \, k > 0 \, ],
        \numberthis
    \end{align*}
    by requiring that
    \begin{equation}
        \biggl( \prod_{i, \> k} {} (S \sss{i, \> k}{\alpha})^{m_{i, \> k}} \biggr) \cdot
        \biggl( \prod_{i, \> k} {} (s \sss{i, \> k}{\alpha})^{n_{i, \> k}} \biggr) 
        = \begin{cases}
            \prod_{i, \> k} m_{i, \> k} !,
            & m_{i, \> k} = n_{i, \> k} \text{ for all } i, k,
            \\ 0,
            & \text{otherwise,}
        \end{cases}
    \end{equation}
    where there are only finitely many $(i, k)$ with
    $m_{i, \> k} \neq 0$ or $n_{i, \> k} \neq 0$,
    and doing the same with $\theta$ in place of $\alpha$.

    This convention has the property that
    for any polynomial~$P$ in the variables~$S \sss{i, \> k}{\alpha}$
    and any polynomial~$q$ in the variables~$s \sss{i, \> k}{\alpha}$\,, we have
    \begin{equation}
        \label{eq-cap-equals-dds}
        q \cap P =
        P \Bigl( \frac{\partial}{\partial s \sss{i, \> k}{\alpha}} \Bigr)
        \bigl( q (s \sss{i, \> k}{\alpha}) \bigr),
    \end{equation}
    where the left-hand side is the cap product,
    and the right-hand side is a differential operator acting on a polynomial.
\end{definition}

When the classes $\alpha, \theta$ are understood,
we may sometimes omit the superscripts, and write
\begin{alignat*}{2}
    S_{i, \> k} & = S \sss{i, \> k}{\alpha} \, , & \qquad
    S_{i, \> k} & = S \sss{i, \> k}{\theta} \, , \\
    s_{i, \> k} & = s \sss{i, \> k}{\alpha} \, , & \qquad
    s_{i, \> k} & = s \sss{i, \> k}{\theta} \, .
\end{alignat*}

From now on, we only focus on the graded vertex algebra
$\hat{H}_* (\frM; \bbQ)$ and its graded twisted module
$\ring{H}_* (|\frM|^\sd; \bbQ)$,
defined in Constructions~\ref{cons-quiv-va-der}
and~\ref{cons-msd-der-alt},
as these have simpler explicit descriptions.

\begin{theorem}
    \label{thm-quiv-transl}
    Consider the graded vertex algebra $\hat{H}_* (\frM; \bbQ)$
    in Construction~\textnormal{\ref{cons-quiv-va-der}.}
    Then for $A \in H_* (\frM_\alpha; \bbQ),$ we have
    \begin{equation}
        \label{eq-quiv-transl}
        D (A) = \biggl(
            \sum_{i \in Q_0} \alpha_i \, s \sss{i, \> 1}{\alpha} +
            \sum_{ \substack{ i \in Q_0 \\ k > 0 } }
            s \sss{i, \> k + 1}{\alpha} \,
            \frac{\partial}{\partial s \sss{i, \> k}{\alpha}}
        \biggr) (A),
    \end{equation}
    where $D$ is the translation operator.
\end{theorem}

\begin{proof}
    \allowdisplaybreaks
    The universal complexes $\calU_i$ on $\calM_\alpha$
    have $\Gm$-weight $1$, so that
    \begin{align*}
        \odot^* (S \sss{i, \> k}{\alpha})
        & = \ch_k (L \boxtimes \calU_i)
        \\* & =
        \sum_{\ell = 0}^{k} \frac{1}{\ell!} \cdot
        T^{\ell} \boxtimes S \sss{i, \> k - \ell}{\alpha} \, ,
        \numberthis
    \end{align*}
    where $T$ is as in Definition~\ref{def-h-bu-1}.
    For a monomial
    $P = \prod_{i, \> k} {} (S \sss{i, \> k}{\alpha})^{m_{i, \> k}}
    \in H^* (\calM_\alpha; \bbQ)$,
    with $k > 0$ in each factor,
    we have
    \begin{align*}
        P \cdot D (A)
        & =
        \odot^* (P) \cdot (t \boxtimes A)
        \\* & =
        \prod_{i, \> k} {} 
        \biggl(
            \sum_{\ell = 0}^{k} \frac{1}{\ell!} \cdot
            T^{\ell} \boxtimes S \sss{i, \> k - \ell}{\alpha}
        \biggr)^{m_{i, \> k}}
        \cdot (t \boxtimes A)
        \\ & =
        \biggl[
            \sum_{i_0, \> k_0}
            m_{i_0, \> k_0} \cdot T \boxtimes 
            \biggl(
                (S \sss{i_0, \> k_0 - 1}{\alpha})^{m_{i_0, \> k_0 - 1}} \cdot
                \hspace{-1em} \prod_{(i, \> k) \neq (i_0, \> k_0)} \hspace{-1em} {}
                (S \sss{i, \> k}{\alpha})^{m_{i, \> k}}
            \biggr)
        \biggr]
        \cdot (t \boxtimes A)
        \\ & =
        \biggl(
            \sum_{ \substack{ i \in Q_0 \\ k > 0 } }
            S \sss{i, \> k-1}{\alpha} \,
            \frac{\partial}{\partial S \sss{i, \> k}{\alpha}}
        \biggr) (P) \cdot A
        \\ & =
        P \cdot \biggl(
            \sum_{i \in Q_0} \alpha_i \, s \sss{i, \> 1}{\alpha} +
            \sum_{ \substack{ i \in Q_0 \\ k > 0 } }
            s \sss{i, \> k + 1}{\alpha} \,
            \frac{\partial}{\partial s \sss{i, \> k}{\alpha}}
        \biggr) (A),
        \numberthis
    \end{align*}
    which proves~\eqref{eq-quiv-transl}.
\end{proof}

\begin{lemma}
    \label{lem-oplus-push}
    Let $A \in H_* (\frM_\alpha; \bbQ),$
    $B \in H_* (\frM_\beta; \bbQ),$ and
    $M \in H_* (|\frM|^\sd_\theta; \bbQ),$
    identified with polynomials in the variables
    $s \sss{i, \> k}{\alpha}\,,$
    $s \sss{i, \> k}{\beta}\,,$
    and $s \sss{i, \> k}{\theta}\,,$ respectively. Then
    \begin{align}
        \label{eq-oplus-push}
        ({\oplus}_{\alpha, \> \beta})_* (A \boxtimes B)
        & =
        \bigl(
            A (s \sss{i, \> k}{\alpha}) \cdot
            B (s \sss{i, \> k}{\prime \beta})
        \bigr) \ 
        \Big|_{ \, \leftsubstack[4em]{
            & s \sss{i, \> k}{\alpha} \mapsto s \sss{i, \> k}{\alpha + \beta} \\[.5ex]
            & s \sss{i, \> k}{\prime \beta} \mapsto s \sss{i, \> k}{\alpha + \beta}
        } } \ ,
        \\[1ex]
        \label{eq-oplus-push-sd}
        {\oplus}^\sd_* (A \boxtimes M)
        & =
        \bigl(
            A (s \sss{i, \> k}{\alpha}) \cdot
            M (s \sss{i, \> k}{\prime \theta})
        \bigr) \ 
        \Big|_{ \, \leftsubstack[4em]{
            & s \sss{i, \> k}{\alpha} \mapsto 2 s \sss{i, \> k}{\bar{\alpha} + \theta} \\[.5ex]
            & s \sss{i, \> k}{\prime \theta} \mapsto s \sss{i, \> k}{\bar{\alpha} + \theta}
        } } \ .
    \end{align}
    Here, it is understood that $s \sss{i, \> k}{\bar{\alpha} + \theta} = 0$
    when $i \in Q_0^\pm$ and $k$ is odd.
\end{lemma}

\begin{proof}
    We have
    \begin{align*}
        ({\oplus}_{\alpha, \> \beta})^* (S \sss{i, \> k}{\alpha + \beta})
        & =
        \ch_k (\calU_i \boxplus \calU_i)
        \\ & =
        1^{\alpha} \boxtimes
        S \sss{i, \> k}{\beta} +
        S \sss{i, \> k}{\alpha} \boxtimes
        1^{\beta} \, .
        \numberthis
    \end{align*}
    Therefore, for a monomial
    $P = \prod_{i, \> k} {} (S \sss{i, \> k}{\alpha + \beta})^{m_{i, \> k}}
    \in H^* (\calM_{\smash{\alpha + \beta}}; \bbQ)$,
    with $k > 0$ in each factor, we have
    \begin{align*}
        P \cdot ({\oplus}_{\alpha, \> \beta})_* (A \boxtimes B)
        & =
        \prod_{i, \> k} {}
        \bigl(
            1^{\alpha} \boxtimes
            S \sss{i, \> k}{\beta} +
            S \sss{i, \> k}{\alpha} \boxtimes
            1^{\beta}
        \bigr)^{m_{i, \> k}} \cdot
        (A \boxtimes B)
        \\ & =
        [1] \,
        \prod_{i, \> k} {} \Bigl(
            \frac{\partial}{\partial s_{i, \> k}}
        \Bigr)^{m_{i, \> k}} 
        \bigl(
            A (s_{i, \> k}) \,
            B (s_{i, \> k})
        \bigr) 
        \\ & =
        P \cdot \bigl(
            A (s \sss{i, \> k}{\alpha + \beta}) \,
            B (s \sss{i, \> k}{\alpha + \beta})
        \bigr) ,
        \numberthis
    \end{align*}
    where $[1]$ denotes taking the constant term,
    and the second step follows from~\eqref{eq-cap-equals-dds}
    and the rule of differentiating a product.
    This proves~\eqref{eq-oplus-push}.

    For~\eqref{eq-oplus-push-sd}, we have
    \begin{align*}
        ({\oplus} \sss{\alpha, \> \theta}{\sd})^*
        (S \sss{i, \> k}{\bar{\alpha} + \theta})
        & =
        \ch_k ((\calU_i \oplus \calU_i^\vee) \boxplus \calV_i)
        \\ & =
        1^{\alpha} \boxtimes
        S \sss{i, \> k}{\theta} +
        2 S \sss{i, \> k}{\alpha} \boxtimes
        1^{\theta} \, ,
        \numberthis
    \end{align*}
    with the assumption that $k$ is even.
    An analogous argument proves~\eqref{eq-oplus-push-sd}.
\end{proof}

\begin{lemma}
    \label{lem-quiv-ch-theta}
    For $i, i' \in Q_0$\,, define coefficients
    \begin{align}
        C_{i, \> i'} & =
        2 \delta_{i, \> i'}
        - \sum_{ \substack{ a \in Q_1 \\ a : \> i \to i' } } 1
        - \sum_{ \substack{ a \in Q_1 \\ a : \> i' \to i } } 1, \\
        \ddot{C}_{i} & =
        -u_i \cdot \biggl(
            2 \delta_{i, \> i^\vee}
            + \sum_{ \substack{ a \in Q_1^\pm \\ a : \> i \to i^\vee } } v_a
            + \sum_{ \substack{ a \in Q_1^\pm \\ a : \> i^\vee \to i } } v_a 
        \biggr) ,
    \end{align}
    where $\delta_{i, \> j}$ is the Kronecker delta.
    Then we have
    \begin{align}
        \label{eq-quiv-ch-theta}
        \ch (\Theta_{\alpha, \> \beta}) & =
        \sum_{ \substack{ i, \> i' \in Q_0 \\ k, \> k' \geq 0 } } {}
        (-1)^{k'} C_{i, \> i'} \,
        S \sss{i, \> k}{\alpha} \boxtimes S \sss{i', \> k'}{\beta} \, , \\
        \label{eq-quiv-ch-theta-dot}
        \ch (\dot{\Theta}_{\alpha, \> \theta}) & =
        \sum_{ \substack{ i, \> i' \in Q_0 \\ k, \> k' \geq 0 } } {}
        (-1)^{k'} C_{i, \> i'} \,
        S \sss{i, \> k}{\alpha} \boxtimes S \sss{i', \> k'}{\theta} \, , \\
        \label{eq-quiv-ch-theta-ddot}
        \ch (\ddot{\Theta}_{\alpha}) & =
        \frac{1}{2} \sum_{ \substack{ i, \> i' \in Q_0 \\ k, \> k' \geq 0 } } {}
        C_{i, \> i'^\vee} \,
        S \sss{i, \> k}{\alpha} \, S \sss{i', \> k'}{\alpha} +
        \frac{1}{2} \sum_{ \substack{ i \in Q_0 \\ k \geq 0 } } {}
        2^k \ddot{C}_i \, S \sss{i, \> k}{\alpha} \, ,
    \end{align}
    where $\Theta, \dot{\Theta}, \ddot{\Theta}$
    are $K$-theory classes defined in
    Construction~\textnormal{\ref{cons-quiv-va-der}}
    or Construction~\textnormal{\ref{cons-msd-der-alt}}.
\end{lemma}

\begin{proof}
    \allowdisplaybreaks
    We may use the explicit descriptions of
    $\Thetaul, \Thetaul[\dot], \Thetaul[\ddot]$
    in Construction~\ref{cons-quiv-va},
    since after translating those expressions into $K$-theory,
    the relations between $\Theta, \dot{\Theta}, \ddot{\Theta}$
    and $\calU_i\,, \calV_i$ are also valid for
    Constructions~\ref{cons-quiv-va-der} and~\ref{cons-msd-der-alt}.

    We have
    \begin{align*}
        \numberthis
        \ch (\Theta_{\alpha, \> \beta}) & =
        2 \sum_{ \substack{ i \in Q_0 \\ k, \> k' \geq 0 } } {}
        (-1)^{k'} S \sss{i, \> k}{\alpha} \boxtimes S \sss{i, \> k'}{\beta}
        - \sum_{ \substack{ a \in Q_1 \\ k, \> k' \geq 0 } } {}
        (-1)^{k'} S \sss{s (a), \> k}{\alpha} \boxtimes S \sss{t (a), \> k'}{\beta}
        \\*[-3ex] & \hspace{16em}
        - \sum_{ \substack{ a \in Q_1 \\ k, \> k' \geq 0 } } {}
        (-1)^{k'} S \sss{t (a), \> k}{\alpha} \boxtimes S \sss{s (a), \> k'}{\beta} \ ,
        \\[1ex]
        \numberthis
        \ch (\dot{\Theta}_{\alpha, \> \theta}) & =
        2 \sum_{ \substack{ i \in Q_0 \\ k, \> k' \geq 0 } } {}
        (-1)^{k'} S \sss{i, \> k}{\alpha} \boxtimes S \sss{i, \> k'}{\theta}
        - \sum_{ \substack{ a \in Q_1 \\ k, \> k' \geq 0 } } {}
        (-1)^{k'} S \sss{s (a), \> k}{\alpha} \boxtimes S \sss{t (a), \> k'}{\theta}
        \\*[-3ex] & \hspace{16em}
        - \sum_{ \substack{ a \in Q_1 \\ k, \> k' \geq 0 } } {}
        (-1)^{k'} S \sss{t (a), \> k}{\alpha} \boxtimes S \sss{s (a), \> k'}{\theta} \ ,
        \\[1ex]
        \numberthis
        \ch (\ddot{\Theta}_{\alpha}) & =
        2 \sum_{ \substack{ i \in Q_0^\tria \\ k, \> k' \geq 0 } } {}
        S \sss{i, \> k}{\alpha} \, S \sss{i^\vee \! , \> k'}{\alpha} +
        \sum_{ i \in Q_0^\pm } {} \biggl(
            \sum_{k, \> k' \geq 0} S \sss{i, \> k}{\alpha} \, S \sss{i, \> k'}{\alpha}
            - u_i \cdot \sum_{ k \geq 0 } 2^k S \sss{i, \> k}{\alpha}
        \biggr)
        \\*[-1ex] & \hspace{-2em}
        - \sum_{ \substack{ a \in Q_1^\tria \\ k, \> k' \geq 0 } } {}
        S \sss{s (a), \> k}{\alpha} \, S \sss{t (a)^\vee \! , \> k'}{\alpha}
        - \frac{1}{2} \sum_{ a \in Q_1^\pm } {} \biggl(
            \sum_{k, \> k' \geq 0} S \sss{s (a), \> k}{\alpha} \, S \sss{s (a), \> k'}{\alpha} +
            u_{s (a)} \, v_a \cdot \sum_{ k \geq 0 } 2^k S \sss{s (a), \> k}{\alpha}
        \biggr)
        \\*[-1ex] & \hspace{-2em}
        - \sum_{ \substack{ a \in Q_1^\tria \\ k, \> k' \geq 0 } } {}
        S \sss{t (a), \> k}{\alpha} \, S \sss{s (a)^\vee \! , \> k'}{\alpha}
        - \frac{1}{2} \sum_{ a \in Q_1^\pm } {} \biggl(
            \sum_{k, \> k' \geq 0} S \sss{t (a), \> k}{\alpha} \, S \sss{t (a), \> k'}{\alpha} +
            u_{t (a)} \, v_a \cdot \sum_{ k \geq 0 } 2^k S \sss{t (a), \> k}{\alpha}
        \biggr) ,
    \end{align*}
    where we used Lemma~\ref{lem-ch-sym}
    to compute the Chern characters of $\Sym^2$ and $\wedge^2$.
    These are equivalent to the expressions
    \eqref{eq-quiv-ch-theta}--\eqref{eq-quiv-ch-theta-ddot}.
\end{proof}

\begin{remark}
    The matrix $(C_{i, \> i'})_{i, \> i' \in Q_0}$
    can be seen as a generalized Cartan matrix.
    In particular, if $Q$ is a simply laced Dynkin quiver,
    that is, the underlying undirected graph of $Q$
    is a Dynkin diagram of type A/D/E,
    then $(C_{i, \> i'})$ coincides with the Cartan matrix of the Dynkin diagram.
\end{remark}

\begin{theorem}
    Consider the graded vertex algebra $\hat{H}_* (\frM; \bbQ)$
    in Construction~\textnormal{\ref{cons-quiv-va}}.

    Let $A \in H_* (\frM_\alpha; \bbQ)$ and
    $B \in H_* (\frM_\beta; \bbQ),$
    identified with polynomials in the variables
    $s \sss{i, \> k}{\alpha}$ and $s \sss{i, \> k}{\beta}\,,$ respectively. 
    Then we have
    \begin{align*}
        \numberthis
        & Y (A, z) \, B =
        (-1)^{\chi_Q (\alpha, \beta)} \cdot
        z^{\chi (\alpha, \beta)} \cdot 
        \exp \biggl[ z
            \biggl( \,
                \sum_{ \nicesubstack{ i \in Q_0 } } {}
                \alpha_i \, s \sss{i, \> 1}{\alpha} +
                \sum_{ \leftsubstack[2em]{
                    & i \in Q_0 \> , \\[-.5ex]
                    & k > 0
                } } {}
                s \sss{i, \> k + 1}{\alpha} \,
                \frac{\partial}{\partial s \sss{i, \> k}{\alpha}} 
            \biggr)
        \biggr] \circ {}
        \\*[-1ex]
        & \hspace{2em}
        \exp \Biggl[ \ 
            \sum_{ \leftsubstack[3em]{
                & i, i' \in Q_0 \> , \\[-.5ex]
                & k, k' \geq 0 \colon \\[-.5ex]
                & k + k' > 0
            } } {}
            (-1)^{k-1} (k + k' - 1)! \, z^{-(k+k')} \cdot C_{i, \> i'} \,
            \partial \sss{i, \> k}{\> \alpha} \,
            \partial \sss{i', \> k'}{\> \prime \> \beta}
        \Biggr] 
        \\*[-6ex]
        & \hspace{21em}
        \bigl(
            A (s \sss{i, \> k}{\alpha}) \cdot
            B (s \sss{i, \> k}{\prime \beta})
        \bigr) \ 
        \Bigg|_{ \, \leftsubstack[4em]{
            & s \sss{i, \> k}{\alpha} \mapsto s \sss{i, \> k}{\alpha + \beta} \\[.5ex]
            & s \sss{i, \> k}{\prime \beta} \mapsto s \sss{i, \> k}{\alpha + \beta}
        } } \ ,
    \end{align*}
    where
    \[
        \partial \sss{i, \> k}{\> \alpha} = \begin{cases}
            \alpha_i \, , & k = 0, \\
            \partial / \partial s \sss{i, \> k}{\alpha} \, , & k > 0,
        \end{cases} \qquad
        \partial \sss{i, \> k}{\> \prime \beta} = \begin{cases}
            \beta_i \, , & k = 0, \\
            \partial / \partial s \sss{i, \> k}{\prime \beta} \, , & k > 0.
        \end{cases}
    \]
\end{theorem}

\begin{proof}
    This follows from the construction of the vertex algebra
    in Theorem~\ref{thm-va},
    together with Theorem~\ref{thm-quiv-transl},
    Lemma~\ref{lem-oplus-push},
    and~\eqref{eq-quiv-ch-theta} in Lemma~\ref{lem-quiv-ch-theta},
    providing explicit expressions for
    the translation operator~$D$,
    the pushforward along the direct sum map, 
    and the Chern character of $\Theta$,
    as well as Lemma~\ref{lem-ch-to-c} translating
    from the Chern character to the Chern class.
\end{proof}

\begin{theorem}
    Consider the graded vertex algebra $\hat{H}_* (\frM; \bbQ)$
    in Construction~\textnormal{\ref{cons-quiv-va},}
    and its graded twisted module $\ring{H}_* (|\frM|^\sd; \bbQ)$
    in Construction~\textnormal{\ref{cons-msd-der-alt}}.

    Let $A \in H_* (\frM_\alpha; \bbQ)$ and
    $M \in H_* (\frM^\sd_\theta; \bbQ)$, 
    identified with polynomials in the variables
    $s \sss{i, \> k}{\alpha}$ and $s \sss{i, \> k}{\theta}\,,$ respectively. 
    Then we have
    \begin{align*}
        \numberthis
        & Y^\sd (A, z) \, M =
        (-1)^{\chi^\sd_Q (\alpha, \theta)} \cdot
        z^{\dot{\chi} (\alpha, \theta)} \cdot
        (2z)^{\ddot{\chi} (\alpha)} \cdot {}
        \\*
        & \hspace{2em}
        \exp \biggl[ z
            \biggl( \,
                \sum_{ \nicesubstack{ i \in Q_0 } } {}
                \alpha_i \, s \sss{i, \> 1}{\alpha} +
                \sum_{ \leftsubstack[2em]{
                    & i \in Q_0 \> , \\[-.5ex]
                    & k > 0
                } } {}
                s \sss{i, \> k + 1}{\alpha} \,
                \frac{\partial}{\partial s \sss{i, \> k}{\alpha}}
            \biggr)
        \biggr] \circ {}
        \\*[-1ex]
        & \hspace{2em}
        \exp \Biggl[ \ 
            \sum_{ \leftsubstack[3em]{
                & i, i' \in Q_0 \> , \\[-.5ex]
                & k, k' \geq 0 \colon \\[-.5ex]
                & k + k' > 0
            } } {}
            (-1)^{k-1} (k + k' - 1)! \, z^{-(k+k')} \cdot C_{i, \> i'} \,
            \partial \sss{i, \> k}{\> \alpha} \,
            \partial \sss{i', \> k'}{\> \prime \theta}
        \\*
        & \hspace{4em}
            {} + \frac{1}{2}
            \sum_{ \leftsubstack[3em]{
                & i, i' \in Q_0 \> , \\[-.5ex]
                & k, k' \geq 0 \colon \\[-.5ex]
                & k + k' > 0
            } } {}
            (-1)^{k+k'-1} \, 2^{-(k+k')} (k + k' - 1)! \, z^{-(k+k')} \cdot C_{i, \> i'^\vee} \,
            \partial \sss{i, \> k}{\> \alpha} \,
            \partial \sss{i', \> k'}{\> \alpha}
        \\*[-1ex]
        & \hspace{4em}
            {} + \frac{1}{2}
            \sum_{ \leftsubstack[2em]{
                & i \in Q_0 \> , \\[-.5ex]
                & k > 0
            } } {}
            (-1)^{k-1} (k - 1)! \, z^{-k} \cdot \ddot{C}_i \,
            \partial \sss{i, \> k}{\> \alpha}
        \Biggr]
        \bigl(
            A (s \sss{i, \> k}{\alpha}) \cdot
            M (s \sss{i, \> k}{\prime \theta})
        \bigr) \ 
        \Bigg|_{ \, \leftsubstack[4em]{
            & s \sss{i, \> k}{\alpha} \mapsto 2 s \sss{i, \> k}{\bar{\alpha} + \theta} \\[.5ex]
            & s \sss{i, \> k}{\prime \theta} \mapsto s \sss{i, \> k}{\bar{\alpha} + \theta}
        } } \ ,
    \end{align*}
    where
    \[
        \partial \sss{i, \> k}{\> \alpha} = \begin{cases}
            \alpha_i \, , & k = 0, \\
            \partial / \partial s \sss{i, \> k}{\alpha} \, , & k > 0,
        \end{cases}
        \qquad
        \partial \sss{i, \> k}{\> \prime \theta} = \begin{cases}
            \theta_i \, , & k = 0, \\
            \partial / \partial s \sss{i, \> k}{\prime \theta} \, ,
            & k > 0, \, i \in Q_0^\tria \cup Q_0^\pm \, , \\
            (-1)^{k} \, \partial / \partial s \sss{i^\vee \! , \> k}{\prime \theta} \, ,
            & k > 0, \, i \in Q_0^{\tria\vee} \, .
        \end{cases}
    \]
\end{theorem}

\begin{proof}
    This follows from the construction of the twisted module
    in Theorem~\ref{thm-tw-mod},
    together with Theorem~\ref{thm-quiv-transl},
    Lemma~\ref{lem-oplus-push},
    and~\eqref{eq-quiv-ch-theta-dot}--\eqref{eq-quiv-ch-theta-ddot},
    providing explicit expressions for
    the translation operator~$D$,
    the pushforward along $\oplus^\sd$, 
    and the Chern character of $\dot{\Theta}$ and $\ddot{\Theta}$,
    as well as Lemma~\ref{lem-ch-to-c} translating
    from the Chern character to the Chern class.
\end{proof}

\section{Orthogonal and symplectic bundles}

\label{sect-principal}

In this section, we apply our theory to the self-dual category of
vector bundles on a smooth, projective variety,
whose self-dual objects are orthogonal or symplectic principal bundles.
This is also discussed in~\cite*[\S10]{Bu2023}.

We have not yet been able to obtain any results on constructing
the enumerative invariants in Conjecture~\ref{conj-main}.
However, the results in \S\ref{sect-vertex-moduli}
do apply in this case,
and we obtain vertex algebra and twisted module structures
on the homology of the relevant moduli stacks.

We explain two different settings where such algebraic structures can be obtained.
Namely, we have the \emph{algebraic} version,
coming from moduli stacks of algebraic vector bundles and principal bundles,
and the \emph{topological} version,
coming from moduli spaces of topological vector bundles and principal bundles.

In this section, for a vector bundle or a perfect complex $E \to X$,
we will always denote by $E^\vee$ the \emph{usual} dual of $E$,
given by $\calHom (E, \calO_X)$ or $\RcalHom (E, \calO_X)$,
rather than its dual under a possibly different self-dual structure,
which will occur below.

\subsection{The algebraic version}

For a $\bbC$-variety $X$ and a reductive $\bbC$-group $G$,
we have the notion of a \emph{principal $G$-bundle} on $X$,
defined in Serre~\cite[\S2.2]{Serre1958},
and discussed in \cite[Definition~10.1]{Bu2023}.

When $G$ is the orthogonal group $\upO (n)$
and the symplectic group $\Sp (2n)$,
principal $G$-bundles can be identified with
vector bundles equipped with symmetric or antisymmetric
non-degenerate bilinear forms,
as in~\cite[\S5.6]{Serre1958} or~\cite[Definition~10.2]{Bu2023}.
Therefore, we may identify orthogonal or symplectic bundles
with self-dual objects in the category of vector bundles on $X$.

\begin{construction}
    \label{cons-princ}
    Let $X$ be a connected, smooth, projective $\bbC$-variety,
    and let $\epsilon = \pm 1$ be a sign.
    Choose a line bundle $L \to X$,
    which we often take to be the trivial line bundle $\calO_X$\,.

    Consider the $\bbC$-linear category
    $\cat{Vect} (X)$ of vector bundles on $X$,
    equipped with a self-dual structure given by
    \begin{equation}
        \calHom (-, L) \colon \cat{Vect} (X) \longsimto \cat{Vect} (X),
    \end{equation}
    with the natural isomorphism $\eta \colon D^\op \circ D \simTo \id$
    given by
    \begin{equation}
        \eta_E = \epsilon \cdot \mathrm{ev}_E^{-1} \colon
        E^{\vee\vee} \longsimto E,
    \end{equation}
    where $\mathrm{ev}_E \colon E \simto E^{\vee\vee}$ is the evaluation isomorphism.
    When $L \simeq \calO_X$\,,
    the self-dual objects of $\cat{Vect} (X)$
    can be identified with orthogonal or symplectic bundles on $X$,
    depending on whether $\epsilon$ is equal to $+1$ or $-1$.

    We have a self-dual $\bbC$-linear stack $\calM^+$
    as the moduli stack of $\cat{Vect} (X)$,
    as in~\cite[Construction~10.5]{Bu2023}.
    Let $\calM$ be its underlying stack, which is an algebraic stack,
    and it is equipped with an involution given by the self-dual structure of $\calM^+$.
    Let $\calM^\sd = \calM^{\bbZ_2}$ be the homotopy fixed points
    of the $\bbZ_2$-action given by this involution.

    In the case when $L \simeq \calO_X$,
    we write $\calM^{\upO}$ and $\calM^{\Sp}$
    for the stack $\calM^\sd$ when $\epsilon = +1$ and $-1$, respectively.
    These can be interpreted as moduli stacks of
    orthogonal and symplectic principal bundles on $X$, respectively.

    We may choose a derived enhancement of $\calM$ to be
    \begin{equation}
        \calMul = \coprod_{n \geq 0} \cat{Map} (X, [*/\GL (n)]),
    \end{equation}
    where $\cat{Map} (-, -)$ denotes the derived mapping stack,
    as in~\cite[\S2.2.6.3]{ToenVezzosi2008}.
    
    As in Toën--Vaquié~\cite[Definition~3.28]{ToenVaquie2007},
    we have a derived stack
    \begin{equation}
        \frM = \cat{Map} (X, \cat{Perf})
    \end{equation}
    of perfect complexes on $\frM$.
    Then $\calMul \subset \frM$
    is an open substack.
    By~\cite[Corollary~3.29]{ToenVaquie2007},
    we see that $\calMul$ is smooth when $X$ is a smooth, projective curve,
    in which case $\calM = \calMul$,
    and that $\calMul$ is quasi-smooth when $X$ is a smooth, projective surface.

    Since taking the tangent complex
    commutes with taking homotopy limits,
    we see that $\calMul^\sd$ is smooth when $X$ is a smooth, projective curve,
    in which case $\calM^\sd = \calMul^\sd$,
    and that $\calMul$ is quasi-smooth when $X$ is a smooth, projective surface.

    Next, we define extra data on the moduli stacks $\calM, \calM^\sd$,
    as in Examples~\ref{eg-lin-cat} and~\ref{eg-sd-lin-cat},
    so that Theorems~\ref{thm-va} and~\ref{thm-tw-mod}
    define vertex algebra and twisted module structures
    on the homology of $\calM$ and $\calM^\sd$.

    As in Remarks~\ref{rmk-choice-ca} and~\ref{rmk-choice-csda},
    we choose an involutive abelian monoid $C^\circ (X)$
    and an involutive module $C^\sd (X)$,
    as alternatives to $\pi_0 (\calM)$ and $\pi_0 (\calM^\sd)$.
    Define $C^\circ (X) \subset H^* (X; \bbQ)$
    to be the set of classes that can be written as the
    Chern character of a vector bundle.
    Define $C^\sd (X) \subset C^\circ (X)$
    to be the set of classes $\alpha$ with $\alpha^\vee = \alpha$.
    Here, $\alpha^\vee$ is the class given by
    $\smash{(\alpha^\vee)_j} = \smash{(-1)^j \alpha_j}$ for all $j$,
    where $\smash{\alpha_j}$ is the part of $\alpha$ living in $H^{2j} (X; \bbQ)$.
    We have natural quotient maps
    $\pi_0 (\calM) \to C^\circ (X)$ and
    $\pi_0 (\calM^\sd) \to C^\sd (X)$.

    Let $\calU \to X \times \calM$ be the universal vector bundle.
    Define the Ext-complex
    $\calExt \in \cat{Perf} (\calM \times \calM)$ by
    \begin{equation}
        \label{eq-def-calext}
        \calExt = (\pr_{2, 3})_* \bigl(
            \pr_{1, 2}^* (\calU^\vee) \otimes \pr_{1, 3}^* (\calU)
        \bigr),
    \end{equation}
    where the $\pr_{i, \> j}$ are projections from
    $X \times \calM \times \calM$ to its factors,
    and we use the derived pushforward functor.
    Define $\calExt[\dot] \in \cat{Perf} (\calM \times \calM^\sd)$ by
    \begin{equation}
        \label{eq-def-calext-dot}
        \calExt[\dot] = (\id \times j)^* (\calExt),
    \end{equation}
    where $j \colon \calM^\sd \to \calM$ is the projection,
    and define $\calExt[\ddot] \in \cat{Perf} (\calM)$ by
    \begin{equation}
        \label{eq-def-calext-ddot}
        \calExt[\ddot] = \begin{cases}
            (\pr_2)_* \bigl( \Sym^2 (\calU^\vee) \bigr), & \epsilon = +1, \\
            (\pr_2)_* \bigl( {\wedge}^2 (\calU^\vee) \bigr), & \epsilon = -1,
        \end{cases}
    \end{equation}
    where $\pr_2 \colon X \times \calM \to \calM$ is the projection.
    
    As in Examples~\ref{eg-lin-cat} and~\ref{eg-sd-lin-cat},
    define perfect complexes
    \begin{alignat}{2}
        \label{eq-def-thetaul-princ-first}
        \Thetaul & = \calExt^\vee \oplus \sigma^* (\calExt)
        && \in \cat{Perf} (\calM \times \calM), \\
        \Thetaul[\dot] & = \calExt[\dot]^\vee + (I \times \id_{\calM^\sd})^* (\calExt[\dot])
        && \in \cat{Perf} (\calM \times \calM^\sd), \\
        \label{eq-def-thetaul-princ-last}
        \Thetaul[\ddot] & = \calExt[\ddot]^\vee + I^* (\calExt[\ddot])
        && \in \cat{Perf} (\calM),
    \end{alignat}
    where $\sigma \colon \calM \times \calM \to \calM \times \calM$
    swaps the two factors, and $I \colon \calM \simto \calM$ is the involution
    given by the self-dual structure.

    Next, we define functions $\varepsilon, \varepsilon^\sd$
    in Assumption~\ref{asn-h-sp}~\ref{itm-epsilon}
    and Assumption~\ref{asn-h-sp-mod}~\ref{itm-epsilon-sd}.
    For $\alpha, \beta \in C^\circ (X)$ and $\theta \in C^\sd (X)$, define
    \begin{align}
        \label{eq-def-chi-x-first}
        \chi_X (\alpha, \beta) & =
        \rank \calExt_{\alpha, \> \beta} \, , \\
        \dot{\chi}_X (\alpha, \theta) & =
        \rank \calExt[\dot]_{\alpha, \> \theta} \, , \\
        \label{eq-def-chi-x-last}
        \ddot{\chi}_X (\alpha) & =
        \rank \calExt[\ddot]_{\alpha} \, .
    \end{align}
    Explicitly, by the Riemann--Roch formula, we have
    \begin{align}
        \label{eq-expl-chi-x-first}
        \chi_X (\alpha, \beta) & \textstyle =
        \int_X \alpha^\vee \cdot \beta \cdot \td (X), \\
        \dot{\chi}_X (\alpha, \theta) & \textstyle =
        \int_X \alpha^\vee \cdot \theta \cdot \td (X), \\
        \label{eq-expl-chi-x-last}
        \ddot{\chi}_X (\alpha) & =
        \begin{cases}
            \int_X \Sym^2 (\alpha^\vee) \cdot \td (X), & \epsilon = +1, \\
            \int_X {\wedge}^2 (\alpha^\vee) \cdot \td (X), & \epsilon = -1.
        \end{cases}
    \end{align}
    for $\alpha, \beta \in C^\circ (X)$ and $\theta \in C^\sd (X)$,
    also regarded as elements of $H^* (X; \bbQ)$,
    where $\td (X)$ is the Todd class of~$X$.
    We have the relations
    \begin{align}
        \label{eq-chi-chi-x-first}
        \chi (\alpha, \beta) & =
        \chi_X (\alpha, \beta) + \chi_X (\beta, \alpha), \\
        \dot{\chi} (\alpha, \theta) & =
        \dot{\chi}_X (\alpha, \theta) + \dot{\chi}_X (\alpha^\vee, \theta), \\
        \label{eq-chi-chi-x-last}
        \ddot{\chi} (\alpha) & =
        \ddot{\chi}_X (\alpha) + \ddot{\chi}_X (\alpha^\vee),
    \end{align}
    where $\chi, \dot{\chi}, \ddot{\chi}$ are defined as the ranks of
    $\Thetaul, \Thetaul[\dot], \Thetaul[\ddot]$, respectively,
    as in Assumption~\ref{asn-h-sp}~\ref{itm-theta}
    and Assumption~\ref{asn-h-sp}~\ref{itm-theta-dots}.
    Then, as in Examples~\ref{eg-lin-cat} and~\ref{eg-sd-lin-cat}, we may take
    \begin{equation}
        \label{eq-epsilon-princ}
        \varepsilon (\alpha, \beta) =
        (-1)^{\chi_X (\alpha, \> \beta)}, \qquad
        \varepsilon^\sd (\alpha, \beta) =
        (-1)^{\dot{\chi}_X (\alpha, \> \theta) + \ddot{\chi}_X (\alpha)}
    \end{equation}
    for all $\alpha, \beta \in C^\circ (X)$ and $\theta \in C^\sd (X)$.
    
    Finally, for Assumption~\ref{asn-h-sp}~\ref{itm-chi-check}
    and Assumption~\ref{asn-h-sp-mod}~\ref{itm-chi-ring}, we take
    \begin{equation}
        \label{eq-chi-hat-ring-princ}
        \hat{\chi} (\alpha) =
        \chi (\alpha, \alpha), \qquad
        \ring{\chi} (\theta) =
        \ddot{\chi} ( j (\theta) )
    \end{equation}
    for all $\alpha \in C^\circ (X)$ and $\theta \in C^\sd (X)$,
    as in Examples~\ref{eg-lin-cat} and~\ref{eg-sd-lin-cat}.

    Now, we may apply Theorems~\ref{thm-va} and~\ref{thm-tw-mod}
    to obtain an involutive graded vertex algebra structure on $\hat{H}_* (\calM; \bbQ)$,
    and an involutive graded twisted module structure on $\ring{H}_* (\calM^\sd; \bbQ)$.
\end{construction}

\begin{remark}
    \label{rmk-princ-choice-lb}
    In Construction~\ref{cons-princ},
    regarding the choice of the line bundle $L$,
    if two line bundles $L_1, L_2$ differ by a square of a line bundle,
    i.e.~we have $L_2 \simeq L_1 \otimes (L')^{\otimes 2}$
    for some line bundle $L'$,
    then choosing $L = L_1$ or $L = L_2$
    actually produces the same result.
    This is because we have an equivalence of self-dual categories
    \begin{equation}
        {-} \otimes L' \colon
        \bigl( \cat{Vect} (X), \calHom (-, L_1), \epsilon \bigr) \longsimto
        \bigl( \cat{Vect} (X), \calHom (-, L_2), \epsilon \bigr).
    \end{equation}
    For this reason, instead of choosing $L$,
    we only need to choose an element of the group
    \[
        \mathrm{Pic} (X) / 2 \mathrm{Pic} (X).
    \]
\end{remark}

We summarize Construction~\ref{cons-princ}
in the special case when $L \simeq \calO_X$
as follows.

\begin{theorem}
    Let $X$ be a connected, smooth, projective $\bbC$-variety.
    Let $\calM$ be the moduli stack of vector bundles on $X,$
    and let $\calM^{\upO}, \calM^{\Sp}$
    be the moduli stacks of orthogonal and symplectic bundles on $X,$
    respectively.

    \begin{enumerate}
        \item 
            There is an involutive graded vertex algebra structure on the graded vector space
            \begin{equation}
                \hat{H}_* (\calM; \bbQ) =
                \bigoplus_{\alpha \in C^\circ (X)}
                H_{* - \hat{\chi} (\alpha)} (\calM_\alpha; \bbQ),
            \end{equation}
            with $\hat{\chi} (\alpha) = 2 \int_X \alpha^\vee \cdot \alpha \cdot \td (X)$
            for all $\alpha \in C^\circ (X)$.

        \item
            There is an involutive graded twisted module structure
            for the involutive graded vertex algebra $\hat{H}_* (\calM; \bbQ)$
            on the graded vector space 
            \begin{equation}
                \ring{H}_* (\calM^{\upO}; \bbQ) =
                \bigoplus_{\theta \in C^\sd (X) \vphantom{^0}}
                H_{* - \ring{\chi} (\theta)} (\calM^\upO_\theta; \bbQ),
            \end{equation}
            with $\ring{\chi} (\theta) = 2 \int_X \Sym^2 (\theta) \cdot \td (X)$
            for all $\theta \in C^\sd (X)$.

        \item
            There is an involutive graded twisted module structure
            for the involutive graded vertex algebra $\hat{H}_* (\calM; \bbQ)$
            on the graded vector space 
            \begin{equation}
                \ring{H}_* (\calM^{\Sp}; \bbQ) =
                \bigoplus_{\theta \in C^\sd (X) \vphantom{^0}}
                H_{* - \ring{\chi} (\theta)} (\calM^{\smash{\Sp}}_\theta; \bbQ),
            \end{equation}
            with $\ring{\chi} (\theta) = 2 \int_X {\wedge}^2 (\theta) \cdot \td (X)$
            for all $\theta \in C^\sd (X)$.
    \end{enumerate}
\end{theorem}

Instead of the moduli stack $\calM$ of vector bundles on $X$ considered above,
we could also consider the full moduli stack $\frM$
of perfect complexes on $X$.
This will give us a different vertex algebra
and a different twisted module.
We sketch the construction below.

\begin{construction}
    \label{cons-princ-perf}
    Let $X$ be a connected, smooth, projective $\bbC$-variety,
    and let $\epsilon = \pm 1$ be a sign.
    Let $\cat{Perf} (X) \simeq \Db \cat{Coh} (X)$
    be the $\bbC$-linear stable $\infty$-category
    of perfect complexes on $X$.
    Choose an object $L \in \cat{Perf} (X)$
    that is \emph{invertible} with respect to the derived tensor product.
    Equivalently~\citestacks{0FPG},
    $L[n]$ is a line bundle for some $n \in \bbZ$.
    We typically choose $L = \calO_X$\,, the trivial line bundle.

    Consider the self-dual structure on $\cat{Perf} (X)$ given by
    \begin{equation}
        \RcalHom (-, L) \colon \cat{Perf} (X)
        \longsimto \cat{Perf} (X),
    \end{equation}
    with a sign $\epsilon$ introduced when identifying
    the double dual of a complex with itself.

    As in Toën--Vaquié~\cite[Definition~3.28]{ToenVaquie2007},
    we have a derived stack
    \begin{equation}
        \frM = \cat{Map} (X, \cat{Perf})
    \end{equation}
    of perfect complexes on $\frM$.
    This can be seen as a moduli stack of objects in $\cat{Perf} (X)$.
    The self-dual structure on $\cat{Perf} (X)$
    can be extended to an involution $I \colon \frM \simto \frM$.
    Let $\frM^\sd = \frM^{\bbZ_2}$ be the homotopy fixed points
    of this $\bbZ_2$-action, which can be seen as
    a moduli stack of self-dual objects in $\cat{Perf} (X)$.
    As in Construction~\ref{cons-princ} above,
    $\frM$ and $\frM^\sd$ are smooth when $X$ is a curve,
    and quasi-smooth when $X$ is a surface.

    Let $K (X), K^\sd (X) \subset H^* (X; \bbQ)$
    be subgroups spanned by the subsets
    $C^\circ (X), C^\sd (X)$ defined in Construction~\ref{cons-princ}, respectively.
    These are temporary notation
    and are only effective in this construction.

    Define perfect complexes
    \[
        \calExt \in \cat{Perf} (\frM \times \frM), \qquad
        \calExt[\dot] \in \cat{Perf} (\frM \times \frM^\sd), \qquad
        \calExt[\ddot] \in \cat{Perf} (\frM)
    \]
    by the formulae~\eqref{eq-def-calext}--\eqref{eq-def-calext-ddot},
    with $\frM, \frM^\sd$ in place of $\calM, \calM^\sd$ there,
    and with the universal perfect complex $\calU \in \cat{\Perf} (X \times \frM)$
    in place of the universal vector bundle there.
    Then, define perfect complexes
    \[
        \Thetaul \in \cat{Perf} (\frM \times \frM), \qquad
        \Thetaul[\dot] \in \cat{Perf} (\frM \times \frM^\sd), \qquad
        \Thetaul[\ddot] \in \cat{Perf} (\frM)
    \]
    by \eqref{eq-def-thetaul-princ-first}--\eqref{eq-def-thetaul-princ-last}.

    Proceeding as in Construction~\ref{cons-princ},
    choosing extra data~$\varepsilon, \varepsilon^\sd, \hat{\chi}, \ring{\chi}$
    as in~\eqref{eq-def-chi-x-first}--\eqref{eq-chi-hat-ring-princ},
    with the classes $\alpha, \beta, \theta$, etc.\ 
    now lying in $K (X), K^\sd (X)$ instead of $C^\circ (X), C^\sd (X)$,
    all requirements of Assumptions~\ref{asn-h-sp} and~\ref{asn-h-sp-mod}
    are satisfied.

    Now, applying Theorems~\ref{thm-va} and~\ref{thm-tw-mod},
    we obtain an involutive graded vertex algebra structure
    on $\hat{H}_* (\frM; \bbQ)$,
    and an involutive graded twisted module structure
    on $\ring{H}_* (\frM^\sd; \bbQ)$.

    Moreover, applying Theorems~\ref{thm-mor-va} and~\ref{thm-mor-tw-mod},
    with the extra data $\Xi, \hat{\Xi}, \dot{\Xi}, \ddot{\Xi}, \ring{\Xi}$
    taken to be zero, we obtain morphisms
    \begin{align}
        i_* \colon \hat{H}_* (\calM; \bbQ)
        & \longrightarrow \hat{H}_* (\frM; \bbQ), \\
        i^\sd_* \colon \ring{H}_* (\calM^\sd; \bbQ)
        & \longrightarrow \ring{H}_* (\frM^\sd; \bbQ)
    \end{align}
    that respect the vertex algebra and twisted module structures,
    where $i \colon \calM \hookrightarrow \frM$
    and $i^\sd \colon \calM^\sd \hookrightarrow \frM^\sd$
    are the inclusions.
\end{construction}

\subsection{The topological version}

Instead of the moduli stacks $\calM, \calM^\sd, \frM, \frM^\sd$
considered in the previous subsection,
we can consider topological versions of these moduli spaces,
and obtain vertex algebras and twisted modules from them.
The topological versions have the advantage that
it is easier to write down their homology, and hence,
the algebraic structures that we obtain.

\begin{construction}
    \label{cons-princ-top}
    Let $X$ be a connected compact Kähler manifold,
    and choose a sign $\epsilon = \pm 1$.
    Define topological spaces
    \begin{equation}
        M = [X, \BU \times \bbZ], \qquad
        M^{\upO} = [X, \BO \times \bbZ], \qquad
        M^{\Sp} = [X, \BSp \times \bbZ], 
    \end{equation}
    where $[-, -]$ denotes the mapping space of topological spaces.
    In fact, these are spaces computing the topological
    $K$-theory, $\mathit{KO}$-theory, and $\mathit{KSp}$-theory of $X$,
    respectively.
    We may regard them as completions of moduli spaces of
    complex vector bundles, complex orthogonal vector bundles,
    and complex symplectic vector bundles on $X$,
    respectively.

    By Bott periodicity, we have
    \begin{equation}
        \Omega^2 M \simeq M, \qquad
        \Omega^4 M^{\upO} \simeq M^{\Sp}, \qquad
        \Omega^4 M^{\Sp} \simeq M^{\upO}.
    \end{equation}

    For ease in notation, write
    \[
        M^\sd = \begin{cases}
            M^{\upO}, & \epsilon = +1, \\
            M^{\Sp}, & \epsilon = -1,
        \end{cases}
        \qquad
        K^\sd (X) = \begin{cases}
            \mathit{KO} (X), & \epsilon = +1, \\
            \mathit{KSp} (X), & \epsilon = -1.
        \end{cases}
    \]
    We thus have $\pi_0 (M) \simeq K (X)$
    and $\pi_0 (M^\sd) \simeq K^\sd (X)$.

    Let $j \colon M^\sd \to M$ be the forgetful map,
    induced by morphisms of topological groups $\upO \to \upU$ and $\Sp \to \upU$
    sending a complex orthogonal or symplectic matrix to its underlying matrix,
    where we used the equivalence $\GL \simeq \upU$.

    We have a universal $K$-theory class $U \in K (X \times M)$.

    Define a class $\itExt \in K (M \times M)$ by
    \begin{equation}
        \itExt = (\pr_{2, 3})_* \bigl(
            \pr_{1, 2}^* (U)^\vee \otimes \pr_{1, 3}^* (U)
        \bigr),
    \end{equation}
    where the $\pr_{i, \> j}$ are projections from $X \times M \times M$
    to its factors, and
    $(\pr_{2, 3})_*$ is the Gysin map, as in Karoubi~\cite[\S5.23]{Karoubi1978}.
    Define $\itExt[\dot] \in K (M \times M^\sd)$ by
    \begin{equation}
        \itExt[\dot] = (\id \times j)^* (\itExt),
    \end{equation}
    and define $\itExt[\ddot] \in K (M)$ by
    \begin{equation}
        \itExt[\ddot] = \begin{cases}
            (\pr_2)_* \bigl( \Sym^2 (U^\vee) \bigr), & \epsilon = +1, \\
            (\pr_2)_* \bigl( {\wedge}^2 (U^\vee) \bigr), & \epsilon = -1.
        \end{cases}
    \end{equation}
    Define the data $\Theta, \dot{\Theta}, \ddot{\Theta}$
    as in Assumption~\ref{asn-h-sp}~\ref{itm-theta}
    and Assumption~\ref{asn-h-sp-mod}~\ref{itm-theta-dots} by
    \begin{alignat}{2}
        \Theta & = \itExt^\vee + \sigma^* (\itExt)
        && \in K (M \times M), \\
        \dot{\Theta} & = \itExt[\dot]^\vee + (I \times \id_{M^\sd})^* (\itExt[\dot])
        && \in K (M \times M^\sd), \\
        \ddot{\Theta} & = \itExt[\ddot]^\vee + I^* (\itExt[\ddot])
        && \in K (M),
    \end{alignat}
    where $\sigma \colon M \times M \to M \times M$ exchanges the two factors,
    and $I \colon M \to M$ is the involution given by taking the dual.

    By a version of the Riemann--Roch formula for topological $K$-theory,
    as in Baum--Fulton--MacPherson~\cite[Theorem~4.1]{Baum1979},
    the formulae for $\chi, \dot{\chi}, \ddot{\chi}$
    in Construction~\ref{cons-princ} also apply to the topological version.
    Namely, for $\alpha, \beta \in K (X)$
    and $\theta \in K^\sd (X)$, if we define
    \begin{align}
        \chi_X (\alpha, \beta) & =
        \rank \itExt_{\alpha, \> \beta} \, , \\
        \dot{\chi}_X (\alpha, \theta) & =
        \rank \itExt[\dot]_{\alpha, \> \theta} \, , \\
        \ddot{\chi}_X (\alpha) & =
        \rank \itExt[\ddot]_{\alpha} \, ,
    \end{align}
    then $\chi_X, \dot{\chi}_X, \ddot{\chi}_X$ are given precisely
    by~\eqref{eq-expl-chi-x-first}--\eqref{eq-expl-chi-x-last},
    by the topological Riemann--Roch formula cited above.
    The relations~\eqref{eq-chi-chi-x-first}--\eqref{eq-chi-chi-x-last}
    also hold in this case.

    Choose the data $\varepsilon, \varepsilon^\sd, \hat{\chi}, \ring{\chi}$
    as in~\eqref{eq-epsilon-princ}--\eqref{eq-chi-hat-ring-princ},
    but with the classes $\alpha, \beta, \theta$
    now living in $K (X), K^\sd (X)$.

    Now, we may apply Theorems~\ref{thm-va} and~\ref{thm-tw-mod}
    to obtain an involutive graded vertex algebra structure on $\hat{H}_* (M; \bbQ)$,
    and an involutive graded twisted module structure on $\ring{H}_* (M^\sd; \bbQ)$.
\end{construction}

\begin{remark}
    In Construction~\ref{cons-princ-top},
    it is also possible to twist the self-dual structure
    by a line bundle $L \to X$,
    or more generally, a class $L \in K (X)$ with $L \cdot L^\vee = 1$,
    as in Constructions~\ref{cons-princ} and~\ref{cons-princ-perf}.
    The space $M^\sd$ can be defined as the homotopy fixed points
    of the $\bbZ_2$-action on $M$, twisted by $L$.
    This will produce possibly different twisted modules
    for the vertex algebra $\hat{H}_* (M; \bbQ)$.
\end{remark}

\begin{remark}
    In Construction~\ref{cons-princ-top},
    we could also consider the non-completed moduli spaces
    \[
        N = \biggl[ X, \coprod_{n \geq 0} \BU (n) \biggr], \qquad
        N^{\upO} = \biggl[ X, \coprod_{n \geq 0} \BO (n) \biggr], \qquad
        N^{\Sp} = \biggl[ X, \coprod_{n \geq 0} \BSp (2n) \biggr].
    \]
    However, we choose to formulate the completed version,
    as we have an explicit description in Theorem~\ref{thm-coh-top-expl} below.
\end{remark}

\begin{remark}
    By a theorem of Gross~\cite[Theorem~4.15]{Gross2019},
    if $X$ is a connected, smooth, projective $\bbC$-variety
    satisfying a condition called being of \emph{class~D},
    which includes curves, surfaces, as well as
    many higher dimensional varieties,
    then we have a homotopy equivalence
    \[
        |\frM| \simeq M,
    \]
    with $\frM$ as in Construction~\ref{cons-princ-perf},
    and $M$ as in Construction~\ref{cons-princ-top}.
    Consequently, for varieties of class~D,
    Constructions~\ref{cons-princ-perf} and~\ref{cons-princ-top}
    produce the same vertex algebras.
    However, it is unclear whether the equivalence
    $|\frM^\sd| \simeq M^\sd$ holds for these varieties, or whether 
    Constructions~\ref{cons-princ-perf} and~\ref{cons-princ-top}
    produce the same twisted module.
\end{remark}

Next, we give an explicit description of the
vertex algebras and twisted modules
obtained in Construction~\ref{cons-princ-top}.

\begin{definition}
    Let $X$ be a connected compact Kähler manifold of complex dimension~$n$,
    and use the notations in Construction~\ref{cons-princ-top}.

    For each $0 \leq k \leq n$,
    choose a basis $(e_{j, \> k})_{j = 1, \dotsc, b_k}$ for $H_k (X; \bbQ)$,
    where $b_k$ is the $k$-th Betti number of $X$.
    Assume that $e_{1, \> 0}$ is the class of a point,
    and that $e_{1, \> n}$ is the fundamental class of $X$.
    Write $J = \{ (j, k) \mid 0 \leq k \leq n, \, 1 \leq j \leq b_k \}$.

    For $\alpha \in K (X)$, $(j, k) \in J$, 
    and an integer $\ell \geq k/2$, define a cohomology class
    $S \sss{j, \> k, \> \ell}{\alpha} \in H^{2 \ell - k} (M_\alpha; \bbQ)$ by
    \begin{equation}
        S \sss{j, \> k, \> \ell}{\alpha} =
        \ch_\ell (U_\alpha) \setminus e_{j, \> k} \, ,
    \end{equation}
    where $U_\alpha \in K (X \times M_\alpha)$ is the universal class,
    classified by the restriction of the evaluation map
    $\ev \colon X \times M \to \BU \times \bbZ$ to $X \times M_\alpha$\,,
    and ${\setminus} \colon H^* (X \times M; \bbQ) \otimes
    H_* (X; \bbQ) \to H^* (M; \bbQ)$
    is the \emph{slant product}, as in, for example,
    Hatcher~\cite[\S3.B]{Hatcher2002}.

    Similarly, for $\theta \in K^\sd (X)$, $(j, k) \in j$, 
    and an integer $\ell \geq k/4$, define a cohomology class
    $S \sss{j, \> k, \> 2 \ell}{\theta} \in H^{4 \ell - k} (M^\sd_\theta; \bbQ)$ by
    \begin{equation}
        S \sss{j, \> k, \> 2 \ell}{\theta} =
        \ch_{2 \ell} (V_\theta) \setminus e_{j, \> k} \, ,
    \end{equation}
    where $V_\theta \in K^\sd (X \times M^\sd_\theta)$ is the universal class,
    classified by the restriction of the evaluation map
    $\ev \colon X \times M^\sd \to (\BO \times \bbZ$ or $\BSp \times \bbZ)$
    to $X \times M^\sd_\theta$\,.
    Note that in this case, the odd Chern characters are always zero,
    since the class $V_\theta$ is self-dual.
\end{definition}

\begin{theorem}
    \label{thm-coh-top-expl}
    In Construction~\textnormal{\ref{cons-princ-top},}
    we have isomorphisms of graded $\bbQ$-algebras
    \begin{align}
        \label{eq-cohom-gl-expl}
        H^* (M_\alpha; \bbQ) & \simeq
        \bbQ [S \sss{j, \> k, \> \ell}{\alpha} : (j, k) \in J, \, 2\ell > k], \\
        \label{eq-cohom-princ-expl}
        H^* (M^\sd_\theta; \bbQ) & \simeq
        \bbQ [S \sss{j, \> k, \> 2\ell}{\theta} : (j, k) \in J, \, 4\ell > k].
    \end{align}
\end{theorem}

\begin{proof}
    The isomorphism~\eqref{eq-cohom-gl-expl} was proved in
    Gross~\cite[Theorem~4.15]{Gross2019},
    using tools from rational homotopy theory.
    We sketch the main ingredients of his proof,
    and modify it to prove~\eqref{eq-cohom-princ-expl}.

    The proof uses a theorem of Milnor--Moore~%
    \cites[Appendix]{MilnorMoore1965}[Theorem~9.2.5]{MayPonto2012},
    which states that if $G$ is a connected commutative $H$-group,
    then there is an isomorphism of graded Hopf algebras
    \begin{equation}
        H_* (G; \bbQ) \simeq \Sym (\pi_* (G) \otimes \bbQ),
    \end{equation}
    where $\pi_* (G) \otimes \bbQ$ is the graded $\bbQ$-vector space
    consisting of the homotopy groups of $G$,
    and $\Sym$ denotes the graded symmetric algebra,
    where odd variables anti-commute.

    We also use a description of the rational homotopy type of the mapping space by
    Brown--Szczarba~\cite[Corollary~1.6]{BrownSzczarba1997},
    which implies that
    \begin{equation}
        \pi_n (M_0) \otimes \bbQ \simeq
        \bigoplus_{k \geq 0} {}
        (\pi_{n+k} (\BU) \otimes \bbQ)
        \otimes H^{k} (X; \bbQ)
    \end{equation}
    for all $n > 0$,
    where $M_0$ is the identity component of $M$.
    We abbreviate this as
    \begin{equation}
        \pi_* (M_0) \otimes \bbQ \simeq
        (\pi_* (\BU) \otimes \bbQ)
        \otimes H^{-*} (X; \bbQ),
    \end{equation}
    but the reader should keep in mind that
    the non-positively graded pieces on the right-hand side are discarded.
    This, together with the Milnor--Moore theorem, gives
    \begin{equation}
        H_* (M_0; \bbQ) \simeq
        \Sym \Bigl( (\pi_* (\BU) \otimes \bbQ)
        \otimes H^{-*} (X; \bbQ) \Bigr).
    \end{equation}

    Now, using a description of the evaluation map
    $\ev \colon M_0 \times X \to \BU$
    in \cites[Theorem~1.1]{BuijsMurillo2006}[Proposition~2.5]{Gross2019},
    one can conclude that the induced map on cohomology
    \begin{equation}
        \ev^* \colon 
        \Sym \Bigl( (\pi_* (\BU) \otimes \bbQ)^\vee \Bigr)
        \longrightarrow
        H^* (X; \bbQ) \otimes
        \Sym \Bigl( (\pi_* (\BU) \otimes \bbQ)^\vee
        \otimes H_{-*} (X; \bbQ) \Bigr)
    \end{equation}
    is given by
    \begin{equation}
        \gamma \longmapsto \sum_{(j, \> k) \in J}
        e \sss{j, \> k}{\vee} \otimes (\gamma \otimes e_{j, \> k})
    \end{equation}
    for $\gamma \in (\pi_* (\BU) \otimes \bbQ)^\vee$.
    This shows that if we define classes
    $\tilde{S} {}\sss{j, \> k, \> \ell}{\alpha} =
    c_\ell (U_\alpha) \setminus \smash{e_{j, \> k}}$\,,
    using the Chern class~$c_\ell$ instead of the Chern character~$\ch_\ell$\,,
    then $\tilde{S} {}\sss{j, \> k, \> \ell}{0}$
    corresponds to the element $c_\ell \otimes \smash{e_{j, \> k}}$
    in the above presentation.
    Since the Chern classes~$c_\ell$ span $(\pi_* (\BU) \otimes \bbQ)^\vee$,
    we see that $H^* (M_0; \bbQ)$ is
    freely generated by the variables $\tilde{S} {}\sss{j, \> k, \> \ell}{0}$\,.
    Then, by observing that
    $S \sss{j, \> k, \> \ell}{0} -
    ((-1)^{\ell-1}/\ell!) \, \tilde{S} {}\sss{j, \> k, \> \ell}{0}$
    is a polynomial in the variables $\tilde{S} {}\sss{j', \> k', \> \ell'}{0}$
    for $\ell' < \ell$,
    one concludes that $H^* (M_0; \bbQ)$ is
    freely generated by the variables $S \sss{j, \> k, \> \ell}{0}$\,.

    Finally, for general $\alpha \neq 0$, there is a canonical isomorphism
    $M_0 \simto M_\alpha$\,, given by the $H$-group structure of $M$,
    which identifies the class $S \sss{j, \> k, \> \ell}{0}$ with
    $S \sss{j, \> k, \> \ell}{\alpha}$\,.
    This proves~\eqref{eq-cohom-gl-expl}.

    The above argument can be adapted to the case of $M^\sd$
    with little modification.
    For convenience, we choose $\epsilon = +1$,
    but the case $\epsilon = -1$ is analogous.
    First, by the result of Brown--Szczarba~\cite{BrownSzczarba1997}
    mentioned above, we have
    \begin{equation}
        \pi_* (M^\sd_0) \otimes \bbQ \simeq
        (\pi_* (\BO) \otimes \bbQ)
        \otimes H^{-*} (X; \bbQ).
    \end{equation}
    Then, the evaluation map
    $\ev \colon M^\sd_0 \times X \to \BO$
    induces a map on cohomology
    \begin{equation}
        \ev^* \colon 
        \Sym \Bigl( (\pi_* (\BO) \otimes \bbQ)^\vee \Bigr)
        \longrightarrow
        H^* (X; \bbQ) \otimes
        \Sym \Bigl( (\pi_* (\BO) \otimes \bbQ)^\vee
        \otimes H_{-*} (X; \bbQ) \Bigr)
    \end{equation}
    given by
    \begin{equation}
        \gamma \longmapsto \sum_{(j, \> k) \in J}
        e \sss{j, \> k}{\vee} \otimes (\gamma \otimes e_{j, \> k})
    \end{equation}
    for $\gamma \in (\pi_* (\BO) \otimes \bbQ)^\vee$.
    This means that if we define classes
    $\tilde{S} {}\sss{j, \> k, \> 2 \ell}{\theta} =
    p_\ell (V_\theta) \setminus \smash{e_{j, \> k}}$\,,
    using the Pontryagin class~$p_\ell$
    instead of the Chern character~$\ch_{2 \ell}$\,,
    then $\tilde{S} {}\sss{j, \> k, \> 2 \ell}{0}$
    corresponds to the element $p_\ell \otimes \smash{e_{j, \> k}}$
    in the above presentation.
    The Pontryagin classes $p_\ell$ span $(\pi_* (\BO) \otimes \bbQ)^\vee$,
    so that $H^* (M^\sd_0; \bbQ)$ is
    freely generated by the variables $\tilde{S} {}\sss{j, \> k, \> 2 \ell}{0}$\,.
    But $S \sss{j, \> k, \> 2 \ell}{0} -
    ((-1)^{\ell-1} / (2 \ell)!) \, \tilde{S} {}\sss{j, \> k, \> 2 \ell}{0}$
    is a polynomial in the variables~$\tilde{S} {}\sss{j', \> k', \> 2 \ell'}{0}$
    for $\ell' < \ell$,
    so that $H^* (M^\sd_0; \bbQ)$ is
    freely generated by the variables $S \sss{j, \> k, \> 2 \ell}{0}$\,.

    Finally, for general $\theta \neq 0$, there is a canonical isomorphism
    $M^\sd_0 \simto M^\sd_\theta$\,, given by the $H$-group structure of $M^\sd$,
    which identifies the class $S \sss{j, \> k, \> \ell}{0}$ with
    $S \sss{j, \> k, \> \ell}{\theta}$\,.
    This proves~\eqref{eq-cohom-princ-expl}.
\end{proof}

The following is a comparison between the algebraic and topological versions
of the vertex algebras and twisted modules.

\begin{proposition}
    Let $X$ be a connected, smooth, projective $\bbC$-variety.
    Then there are morphisms
    \begin{align}
        \hat{H}_* (\frM; \bbQ)
        & \longrightarrow \hat{H}_* (M; \bbQ), \\
        \ring{H}_* (\frM^\sd; \bbQ)
        & \longrightarrow \ring{H}_* (M^\sd; \bbQ)
    \end{align}
    that respect the vertex algebra and twisted module structures.
\end{proposition}

\begin{proof}
    We have a natural map $|\frM| \to M$ given by~\eqref{eq-map-kalg-ktop}.
    We also have the map $|\frM^\sd| \to M^\sd$,
    induced by the $\bbZ_2$-equivariant map $|\frM^\sd| \to |\frM| \to M$.
    Using Baum--Fulton--MacPherson~\cite[Theorem~4.1]{Baum1979}
    comparing pushforwards in topological and algebraic $K$-theory,
    one can deduce that the maps $|\frM| \to M$ and $|\frM^\sd| \to M^\sd$
    given above are compatible with the data
    $\Theta, \dot{\Theta}, \ddot{\Theta}$, etc.
    Therefore, we can apply Theorems~\ref{thm-mor-va} and~\ref{thm-mor-tw-mod},
    with the extra data $\Xi, \hat{\Xi}, \dot{\Xi}, \ddot{\Xi}, \ring{\Xi}$
    taken to be zero, obtaining the desired morphisms.
\end{proof}

\clearpage
\appendix

\section[Proof of results in \S\ref*{sect-vertex-moduli}]{Proof of results in \S\ref{sect-vertex-moduli}}

\label{sect-proofs-va}

\subsection[Proof of Theorem~\ref*{thm-va}]{Proof of Theorem~\ref{thm-va}}
\label{sect-proof-va}

We present a proof of Theorem~\ref{thm-va},
which states Joyce's vertex algebra construction.
We follow ideas from~\cite{JoyceHall},
but we use an alternative argument,
so that it generalizes to prove
other results in~\S\ref{sect-vertex-moduli}.

Suppose that we are given data satisfying Assumption~\ref{asn-h-sp}.
For convenience, write
\[
    V = \hat{H}_* (X; \bbQ), \qquad
    V_\alpha = \hat{H}_* (X_\alpha; \bbQ)
\]
for $\alpha \in \pi_0 (X)$. Write
\[
    D_\alpha \colon V_\alpha \longrightarrow V_\alpha
\]
for the restriction of the translation operator $D$ in~\eqref{eq-constr-transl}
to the subspace $V_\alpha$\,.

\begin{definition}
    Let $n \geq 0$ be an integer. Define a map
    \[
        X_n (-, \dotsc, -; z_1, \dotsc, z_n) \colon
        V^{\otimes n} \longrightarrow
        V \llbr z_1, \dotsc, z_n \rrbr [(z_i - z_j)^{-1} : 1 \leq i < j \leq n]
    \]
    as follows.
    For elements $A_i \in H_{a_i} (X_{\alpha_i}; \bbQ)$, for $i = 1, \dotsc, n$, define
    \begin{align*}
        \numberthis
        \hspace{2em} & \hspace{-2em} \smash[b]{
        X_n (A_1, \dotsc, A_n; z_1, \dotsc, z_n)= 
        } \\[\LSs]
        & \smash{
        (-1)^{ \sum_{1 \leq i < j \leq n} a_i \hat{\chi} (\alpha_j) } 
        \varepsilon (\alpha_1, \dotsc, \alpha_n) \cdot 
        (\oplus_{\alpha_1, \> \dotsc \> , \> \alpha_n})_* \circ
        (\upe^{z_1 D_{\alpha_1}} \otimes \cdots \otimes \upe^{z_n D_{\alpha_n}})
        } \\[\LS]
        & \hspace{2em} \smash[t]{
        \biggl[
            (A_1 \boxtimes \cdots \boxtimes A_n) \cap
            \prod_{1 \leq i < j \leq n} {} \bigl(
                (z_i - z_j)^{\chi (\alpha_i, \> \alpha_j)} \,
                c_{1/(z_i - z_j)} (\Theta_{\alpha_i, \> \alpha_j})
            \bigr)
        \biggr]. }
    \end{align*}
\end{definition}

\begin{lemma}
    \label{lem-cd-ezd-chern}
    Let $Y$ be a topological space,
    and let $E \in K (X_{\alpha_1} \times \cdots \times X_{\alpha_n} \times Y)$
    be of rank $r$ and $\upU (1)$-weight $(k_1, \dotsc, k_n)$. Write
    \begin{align*}
        W
        & = H_* (Y; \bbQ),
        \\
        V_{\alpha_1, \> \dotsc \> , \> \alpha_n}
        & =
        H_* (X_{\alpha_1}; \bbQ) \otimes
        \cdots \otimes
        H_* (X_{\alpha_n}; \bbQ) \otimes W.
    \end{align*}
    Then we have a commutative diagram
    \begin{equation}
    \label{eq-cd-ezd-chern}
    \begin{tikzcd}
        V_{\alpha_1, \> \dotsc \> , \> \alpha_n}
        \ar[r, equals]
        \ar[d, "{\cap \, \tilde{y}^r \, c_{1/\tilde{y}} (E)}"']
        & V_{\alpha_1, \> \dotsc \> , \> \alpha_n}
        \ar[d, "{ \upe^{z_1 D_{\alpha_1}} \, \otimes \, \cdots \, \otimes \, \upe^{z_n D_{\alpha_n}} \, \otimes \, \id_W }"]
        \\ \phantom{ [\tilde{y}^{\pm 1}] }
        V_{\alpha_1, \> \dotsc \> , \> \alpha_n} [\tilde{y}^{\pm 1}]
        \ar[d, "{ \upe^{z_1 D_{\alpha_1}} \, \otimes \, \cdots \, \otimes \, \upe^{z_n D_{\alpha_n}} \, \otimes \, \id_W }"']
        & V_{\alpha_1, \> \dotsc \> , \> \alpha_n} 
        \mathrlap{ \llbr z_i \rrbr }
        \ar[d, "{\cap \, y^r \, c_{1/y} (E)}"]
        \\ \phantom{ \llbr z_i \rrbr [\tilde{y}^{\pm 1}] }
        V_{\alpha_1, \> \dotsc \> , \> \alpha_n} 
        \llbr z_i \rrbr [\tilde{y}^{\pm 1}]
        \ar[d, hook]
        & V_{\alpha_1, \> \dotsc \> , \> \alpha_n} 
        \mathrlap{ [y^{\pm 1}] \llbr z_i \rrbr } \ar[d, hook]
        \\ \phantom{ \llparen y, z_i \rrparen }
        V_{\alpha_1, \> \dotsc \> , \> \alpha_n} 
        \llparen y, z_i \rrparen
        \ar[r, hook, "\iota_{y, \> \{ z_i \}}"]
        & V_{\alpha_1, \> \dotsc \> , \> \alpha_n} 
        \mathrlap{ \llparen y \rrparen \llparen z_i \rrparen \ , }
    \end{tikzcd}
    \end{equation} 
    where $\tilde{y} = y + \sum_i k_i z_i \, ,$
    and $\iota_{y, \{ z_i \}}$ is the map in Definition~\textnormal{\ref{def-iota}.}
\end{lemma}

\begin{proof}
    \allowdisplaybreaks
    The upper-left vertical map in~\eqref{eq-cd-ezd-chern} is well-defined,
    since for any element $A \in V_{\alpha_1, \> \dotsc \> , \> \alpha_n}$\,,
    we have $A \cap c_i (E) = 0$ for $i \gg 0$, for degree reasons.
    Similarly, the middle right map is also well-defined.
    
    Let
    \[
        \odot \colon
        \BU (1)^n \times X_{\alpha_1} \times \cdots \times X_{\alpha_n} \times Y
        \longrightarrow X_{\alpha_1} \times \cdots \times X_{\alpha_n} \times Y
    \]
    be the map given by the natural $\BU (1)$-actions on each $X_{\alpha_i}$\,.
    Then the composition of the left vertical arrows in~\eqref{eq-cd-ezd-chern}
    is equal to
    \begin{equation}
        \sum_{a_1, \> \dotsc \> , \> a_n \geq 0}
        \frac{1}{a_1 ! \cdots a_n !} \cdot
        z\sss{1}{a_1} \cdots z\sss{n}{a_n} \cdot \odot_* \bigl(
            t\sss{1}{a_1} \otimes \cdots \otimes t\sss{n}{a_n} \otimes
            ({-} \cap \tilde{y}^r \, c_{1/\tilde{y}} (E))
        \bigr),
    \end{equation}
    and the composition of the right vertical arrows in~\eqref{eq-cd-ezd-chern}
    is equal to
    \begin{equation}
        \sum_{a_1, \> \dotsc \> , \> a_n \geq 0}
        \frac{1}{a_1 ! \cdots a_n !} \cdot
        z\sss{1}{a_1} \cdots z\sss{n}{a_n} \cdot
        \odot_* (t\sss{1}{a_1} \otimes \cdots \otimes t\sss{n}{a_n} \otimes -)
        \cap y^r \, c_{1/y} (E).
    \end{equation}
    Therefore, it is enough to prove that
    \begin{multline}
        \sum_{a_1, \> \dotsc \> , \> a_n \geq 0}
            \frac{1}{a_1 ! \cdots a_n !} \cdot
            z\sss{1}{a_1} t\sss{1}{a_1} \otimes \cdots \otimes
            z\sss{n}{a_n} t\sss{n}{a_n} \otimes
            ({-} \cap \tilde{y}^r \, c_{1/\tilde{y}} (E)) = \\
        \sum_{a_1, \> \dotsc \> , \> a_n \geq 0} {}
            \frac{1}{a_1 ! \cdots a_n !} \cdot
            ( z\sss{1}{a_1} t\sss{1}{a_1} \otimes \cdots \otimes
            \cdots z\sss{n}{a_n} t\sss{n}{a_n} \otimes -)
            \cap y^r \, c_{1/y} (\odot^* (E))
    \end{multline}
    up to expanding $\tilde{y}^{-1}$ in terms of $y$,
    since applying $\odot_*$ gives the desired equality.
    Looking at the coefficient of a term
    $t\sss{1}{a_1} \otimes \cdots \otimes t\sss{n}{a_n}$,
    we see that it suffices to prove that
    \begin{equation}
        \label{eq-pf-ezd-chern-1}
        z\sss{1}{a_1} \cdots z\sss{n}{a_n} \cdot
        \tilde{y}^r \, c_{1/\tilde{y}} (E) =
        \exp \biggl( \sum_i z_i \, \frac{\partial}{\partial T_i} \biggr) \,
        \Bigl(
            T\sss{1}{a_1} \cdots T\sss{n}{a_n} \cdot
            y^r \, c_{1/y} (\odot^* (E))
        \Bigr) \,
        \bigg|_{T_i = 0}
    \end{equation}
    in cohomology.

    Using
    $\odot^* (E) = L\sss{1}{\otimes k_1} \boxtimes \cdots
    \boxtimes L\sss{n}{\otimes k_n} \boxtimes E$,
    and using Lemma~\ref{lem-ch-to-c},
    the right-hand side of~\eqref{eq-pf-ezd-chern-1} equals
    \begin{align*}
        \hspace{1em} & \hspace{-1em}
        T\sss{1}{a_1} \cdots T\sss{n}{a_n} \, y^r \cdot \exp \Biggl[ \,
            \sum_{j > 0} {}
            (-1)^{j-1} (j - 1)! \, y^{-j}
            \ch_j \bigl( L\sss{1}{\otimes k_1} \boxtimes \cdots
            \boxtimes L\sss{n}{\otimes k_n} \boxtimes E \bigr)
        \Biggr] \, \Bigg|_{T_i = z_i} \\
        & = 
        z\sss{1}{a_1} \cdots z\sss{n}{a_n} \, y^r \cdot \exp \Biggl[ \,
            \sum_{ \leftsubstack[4em]{
                & j_0, \> \dotsc \> , \> j_n \geq 0, \ j > 0 : \\[-.5ex]
                & j = j_0 + \cdots + j_n
            } } {}
            (-1)^{j-1} (j - 1)! \, y^{-j} \cdot
            \biggl( \prod_{i = 1}^{n} \frac{(k_i z_i)^{j_i}}{j_i !} \biggr)
            \cdot \ch_{j_0} (E)
        \Biggr] \\
        & = 
        z\sss{1}{a_1} \cdots z\sss{n}{a_n} \, y^r \cdot \exp \Biggl[ \,
            \sum_{j > 0} {}
            \frac{(-1)^{j-1}}{j} \cdot y^{-j} \cdot
            \biggl( \sum_{i=1}^n k_i z_i \biggr)^j \cdot r \\
        && \mathllap{
            {} + \sum_{j_0 > 0, \ j' \geq 0} {}
            (-1)^{j_0 + j' - 1} \cdot \frac{(j_0 + j' - 1)!}{j'!} \cdot 
            y^{-j_0} \cdot
            \biggl( y^{-1} \sum_{i=1}^n k_i z_i \biggr)^{j'} \cdot
            \ch_{j_0} (E)
        \Biggr] } \\
        & = 
        z\sss{1}{a_1} \cdots z\sss{n}{a_n} \, y^r \cdot \exp \Biggl[
            r \cdot \log \biggl( 1 + y^{-1} \sum_{i=1}^n k_i z_i \biggr) \\*[-1ex]
        && \mathllap{
            \sum_{j_0 > 0} {} (-1)^{j_0 - 1} (j_0 - 1)! \, y^{-j_0} \cdot
            \biggl( 1 + y^{-1} \sum_{i=1}^n k_i z_i \biggr)^{-j_0} \cdot
            \ch_{j_0} (E)
        \Biggr] } \\
        & = 
        z\sss{1}{a_1} \cdots z\sss{n}{a_n} \, \tilde{y}^r \cdot \exp \Biggl[ \,
            \sum_{j_0 > 0} {} (-1)^{j_0 - 1} (j_0 - 1)! \, \tilde{y}^{-j_0} \cdot
            \ch_{j_0} (E)
        \Biggr] \ ,
        \numberthis
    \end{align*}
    which agrees with the left-hand side of~\eqref{eq-pf-ezd-chern-1}.
\end{proof}

\begin{lemma}
    \label{lem-pf-va-main}
    The operations $X_n$ satisfy the conditions in
    Theorem~\textnormal{\ref{thm-va-alt-def}.}
\end{lemma}

\begin{proof}
    For Theorem~\ref{thm-va-alt-def}~\ref{itm-va-multi-unit},
    we have by definition
    \begin{equation}
        X_1 (A; z) = \upe^{z D} (A).
    \end{equation}
    Theorem~\ref{thm-va-alt-def}~\ref{itm-va-multi-equivar}
    also follows from the definitions.

    For Theorem~\ref{thm-va-alt-def}~\ref{itm-va-multi-assoc},
    suppose that $A_i \in V_{\alpha_i}$ and
    $B_i \in V_{\beta_i}$ for classes
    $\alpha_i, \beta_i \in \pi_0 (X)$.
    Suppose that $\deg A_i = a_i$ and $\deg B_i = b_i$ for all $i$.
    Write $\alpha_\ell = \beta_1 + \cdots + \beta_n$\,, and write
    $a_\ell = b_1 + \cdots + b_n$\,. Set
    \begin{align*}
        N & =
        A_1 \otimes \cdots \otimes A_{\ell-1} \otimes
        B_1 \otimes \cdots \otimes B_{n} \otimes
        A_{\ell+1} \otimes \cdots \otimes A_m
        \\ & \in
        V_{\alpha_1} \otimes \cdots \otimes V_{\alpha_{\ell-1}} \otimes
        V_{\beta_1} \otimes \cdots \otimes V_{\beta_n} \otimes
        V_{\alpha_{\ell+1}} \otimes \cdots \otimes V_{\alpha_m} \, .
    \end{align*}
    Expanding the left-hand side of~\eqref{eq-va-multi-assoc}, we obtain
    \begin{align*}
        \hspace{2em} & \hspace{-2em} \smash[b]{
        (-1)^{
            \sum_{1 \leq i < j \leq m} a_i \hat{\chi} (\alpha_j) \, + \,
            \sum_{1 \leq i < j \leq n} b_i \hat{\chi} (\beta_j)
        } \cdot
        \varepsilon (\alpha_1, \dotsc, \alpha_m) \cdot
        \varepsilon (\beta_1, \dotsc, \beta_n) \cdot {}
        } \\[\LSs]
        & \smash{
        \bigl( \oplus_{\alpha_1, \> \dotsc \> , \> \alpha_m} \bigr)_* \circ
        \bigl( \upe^{y_1 D_{\alpha_1}} \otimes \cdots \otimes
        \upe^{y_m D_{\alpha_m}} \bigr)
        } \\[\LS]
        & \smash{
        \Biggl\{
            \Biggl[
                \biggl(
                    \id_{V_{\alpha_1} \otimes \cdots \otimes V_{\alpha_{\ell-1}}} \otimes
                    \Bigl(
                        \bigl( \oplus_{\beta_1, \> \dotsc \> , \> \beta_n} \bigr)_* \circ
                        \bigl(
                            \upe^{z_1 D_{\beta_1}} \otimes \cdots
                            \otimes \upe^{z_n D_{\beta_n}}
                        \bigr)
                    \Bigr) \otimes
                    \id_{V_{\alpha_{\ell+1}} \otimes \cdots \otimes V_{\alpha_m}}
                \biggr)
        } \\[\LSl]
        & \hspace{3em} \smash{
                \biggl(
                    N \cap
                    \prod_{1 \leq i < j \leq n} {} \bigl(
                        (z_i - z_j)^{\chi (\beta_i, \> \beta_j)} \,
                        c_{1/(z_i - z_j)} (\Theta_{\beta_i, \> \beta_j})
                    \bigr)
                \biggr)
            \Biggr]
        } \\[\LSl]
        & \hspace{1.5em} \smash[t]{
            \mathllap{ {} \cap \biggl( \hspace{-.3em} }
                \LAP{2em}{ \prod_{1 \leq i < j \leq m} }
                \bigl(
                    (y_i - y_j)^{\chi (\alpha_i, \> \alpha_j)} \,
                    c_{1/(y_i - y_j)} (\Theta_{\alpha_i, \> \alpha_j})
                \bigr)
            \biggr) 
        \Biggr\} . }
        \numberthis \label{eq-pf-va-1}
    \end{align*}
    Here, it is understood that the classes
    $\smash{\Theta_{\alpha_i, \> \alpha_j}}$\,,
    $\smash{\Theta_{\beta_i, \> \beta_j}}$
    are pulled back to the space
    $\smash{X_{\alpha_1}} \times \cdots \times \smash{X_{\alpha_{\ell-1}}}
    \times \smash{X_{\beta_1}} \times \cdots \times \smash{X_{\beta_n}}
    \times \smash{X_{\alpha_{\ell+1}}} \times \cdots \times \smash{X_{\alpha_m}}$
    using appropriate projections.
    In particular, for classes involving $\alpha_\ell$\,,
    this uses the map $\oplus_{\smash{\beta_1, \> \dotsc, \> \beta_n}}$ as well.

    Using the relations~\eqref{eq-theta-linear-first}--\eqref{eq-theta-linear-second}
    and~\eqref{eq-theta-sym}
    to expand the terms involving $\alpha_\ell$\,,
    the last line of~\eqref{eq-pf-va-1} becomes
    \begin{align*}
        & \hspace{1.5em} \smash[b]{
            \mathllap{ {\cdots} \cap \biggl( \hspace{-.3em} }
                \LAP{2em}{ \prod_{ \substack{ 1 \leq i < j \leq m \\ i, j \neq \ell} } }
                \bigl(
                    (y_i - y_j)^{\chi (\alpha_i, \> \alpha_j)} \,
                    c_{1/(y_i - y_j)} (\Theta_{\alpha_i, \> \alpha_j})
                \bigr) \cdot {}
        } \\[\LSl]
        & \hspace{4.5em} \smash[t]{
                \LAP{2em}{ \prod_{ \substack{ 1 \leq i \leq m, \ i \neq \ell \\ 1 \leq j \leq n } } }
                \bigl(
                    (y_i - y_\ell)^{\chi (\alpha_i, \> \beta_j)} \,
                    c_{1/(y_i - y_\ell)} (\Theta_{\alpha_i, \> \beta_j})
                \bigr)
            \biggr) . }
        \numberthis
    \end{align*}
    Next, we apply Lemma~\ref{lem-cd-ezd-chern} to move
    all the Chern classes into the inner layer $(N \cap {\cdots})$.
    We may also merge the two direct sum operations.
    This gives
    \begin{align*}
        \hspace{2em} & \hspace{-2em} \smash[b]{
        (-1)^{
            \sum_{1 \leq i < j \leq m} a_i \hat{\chi} (\alpha_j) \, + \,
            \sum_{1 \leq i < j \leq n} b_i \hat{\chi} (\beta_j)
        } \cdot
        \varepsilon (\alpha_1, \dotsc, \alpha_m) \cdot
        \varepsilon (\beta_1, \dotsc, \beta_n) \cdot {}
        } \\*[\LSs]
        & \smash{
        \bigl(
            \oplus_{\alpha_1, \> \dotsc \> , \> \alpha_{\ell-1}, \> 
            \beta_1, \> \dotsc \> , \> \beta_n, \,
            \alpha_{\ell+1}, \> \dotsc \> , \> \alpha_m}
        \bigr)_* \circ \Bigl(
            \upe^{y_1 D_{\alpha_1}} \otimes \cdots \otimes
            \upe^{y_{\ell-1} D_{\alpha_{\ell-1}}} \otimes {}
        } \\[\LSs]
        & \hspace{2.5em} \smash{
            \upe^{(y_\ell \, + \, z_1) D_{\beta_1}} \otimes \cdots \otimes
            \upe^{(y_\ell \, + \, z_n) D_{\beta_n}} \otimes
            \upe^{y_{\ell+1} D_{\alpha_{\ell+1}}} \otimes \cdots \otimes
            \upe^{y_m D_{\alpha_m}}
        \Bigr)
        } \\[\LS]
        & \hspace{2.5em} \smash{
        \mathllap{ \Biggl\{
            N \cap \Biggl( \hspace{-.3em} }
                \LAP{2em}{ \prod_{1 \leq i < j \leq n} }
                \bigl(
                    (z_i - z_j)^{\chi (\beta_i, \> \beta_j)} \,
                    c_{1/(z_i - z_j)} (\Theta_{\beta_i, \> \beta_j})
                \bigr) \cdot {}
        } \\[\LSl]
        & \hspace{2.5em} \smash{
                \LAP{2em}{ \prod_{ \substack{ 1 \leq i < j \leq m \\ i, j \neq \ell} } }
                \bigl(
                    (y_i - y_j)^{\chi (\alpha_i, \> \alpha_j)} \,
                    c_{1/(y_i - y_j)} (\Theta_{\alpha_i, \> \alpha_j})
                \bigr) \cdot {}
        } \\[\LSl]
        & \hspace{2.5em} \smash[t]{
                \LAP{2em}{ \prod_{ \substack{ 1 \leq i \leq m, \ i \neq \ell \\ 1 \leq j \leq n } } }
                \bigl(
                    (y_i - y_\ell - z_j)^{\chi (\alpha_i, \> \beta_j)} \,
                    c_{1/(y_i - y_\ell - z_j)} (\Theta_{\alpha_i, \> \beta_j})
                \bigr)
            \Biggr) 
        \Biggr\} , }
        \numberthis
    \end{align*}
    which matches the right-hand side of~\eqref{eq-va-multi-assoc}.
    Note that we have
    \begin{equation}
        \varepsilon (\alpha_1, \dotsc, \alpha_m) \cdot
        \varepsilon (\beta_1, \dotsc, \beta_n) =
        \varepsilon (\alpha_1, \dotsc, \alpha_{\ell-1},
            \beta_1, \dotsc, \beta_n,
            \alpha_{\ell+1}, \dotsc, \alpha_m) ,
    \end{equation}
    which can be deduced from~\eqref{eq-epsilon-assoc}
    via a combinatorial argument.
\end{proof}

\begin{bproof}[Proof of Theorem~\ref{thm-va}]
    By Lemma~\ref{lem-pf-va-main} and Theorem~\ref{thm-va-alt-def},
    we have a graded vertex algebra structure on $V$.
    It remains to verify~\eqref{eq-constr-unit}--\eqref{eq-constr-va}.

    For~\eqref{eq-constr-unit},
    we have by definition $X_0 (\, ; \,) = [0]$,
    so~\eqref{eq-constr-unit} follows from~\eqref{eq-va-alt-unit}.
    For~\eqref{eq-constr-transl}, we have by definition
    \begin{equation}
        X_1 (A; z) = \upe^{z D} (A)
    \end{equation}
    for all $A \in V$, so~\eqref{eq-constr-transl}
    follows from~\eqref{eq-va-alt-transl}.
    Similarly, \eqref{eq-constr-va} follows from~\eqref{eq-va-alt-mult}
    and the definition of $X_2$\,.
\end{bproof}

\subsection[Proof of Theorem~\ref*{thm-tw-mod}]{Proof of Theorem~\ref{thm-tw-mod}}
\label{sect-proof-tw-mod}

Suppose that we are given data satisfying Assumption~\ref{asn-h-sp-mod}.
As in \S\ref{sect-proof-va}, write
\[
    V = \hat{H}_* (X; \bbQ), \qquad
    V_\alpha = \hat{H}_* (X_\alpha; \bbQ)
\]
for $\alpha \in \pi_0 (X)$. In addition, we write
\[
    W = \ring{H}_* (X^\sd; \bbQ), \qquad
    W_\theta = \ring{H}_* (X^\sd_\theta; \bbQ)
\]
for $\theta \in \pi_0 (X^\sd)$.

\begin{definition}
    Let $n \geq 0$ be an integer. Define a map
    \[
        X^\sd_n (-, \dotsc, -; z_1, \dotsc, z_n) \colon
        V^{\otimes n} \otimes W \longrightarrow
        W \llbr z_1, \dotsc, z_n \rrbr [z_i^{-1}, (z_i \pm z_j)^{-1}]
    \]
    as follows,
    where on the right-hand side, we invert all $z_i$ for $1 \leq i \leq n$,
    and $(z_i \pm z_j)$ for $1 \leq i < j \leq n$.
    For $A_i \in H_{a_i} (X_{\alpha_i}; \bbQ)$, where $i = 1, \dotsc, n$,
    and for $M \in H_m (X^\sd_\theta; \bbQ)$, define
    \begin{align*}
        \numberthis
        \hspace{2em} & \hspace{-2em} \smash[b]{
        X_n^\sd (A_1, \dotsc, A_n, M; z_1, \dotsc, z_n) = 
        } \\*[\LSs]
        & \smash{
        (-1)^{
            \sum_{1 \leq i < j \leq n} a_i \hat{\chi} (\alpha_j) \, + \,
            \sum_{1 \leq i \leq n} a_i \ring{\chi} (\theta)
        } \cdot
        \varepsilon^\sd (\alpha_1, \dotsc, \alpha_n, \theta) \cdot {}
        } \\*[\LSs]
        & \smash{
        (\oplus_{\alpha_1, \> \dotsc \> , \> \alpha_n, \> \theta}^\sd)_* \circ
        (\upe^{z_1 D_{\alpha_1}} \otimes \cdots \otimes \upe^{z_n D_{\alpha_n}} \otimes \id_{W_\theta}) \, \biggl[
        (A_1 \boxtimes \cdots \boxtimes A_n \boxtimes M) \cap {}
        } \\*[\LS]
        & \hspace{2em} \smash{
        \mathrlap{ \biggl( }
        \LAP{2em}{ \prod_{1 \leq i < j \leq n} }
            \bigl(
                (z_i - z_j)^{\chi (\alpha_i, \> \alpha_j)} \,
                c_{1/(z_i - z_j)} (\Theta_{\alpha_i, \> \alpha_j})
            \bigr) \cdot {}
        } \\*[\LS]
        & \hspace{2em} \smash{
            \LAP{2em}{ \prod_{1 \leq i < j \leq n} }
            \bigl(
                (z_i + z_j)^{\chi (\alpha_i, \> \alpha_j^\vee)} \,
                c_{1/(z_i + z_j)} (\Theta_{\alpha_i, \> \alpha_j^\vee}) 
            \bigr) \cdot {}
        } \\*[\LS]
        & \hspace{2em} \smash[t]{
            \LAP{2em}{ \prod_{1 \leq i \leq n} }
            \bigl(
                z\sss{i}{\dot{\chi} (\alpha_i, \> \theta)} \,
                c_{1/z_i} (\dot{\Theta}_{\alpha_i, \> \theta}) 
            \bigr) \cdot
            \LAP{1.2em}{ \prod_{1 \leq i \leq n} }
            \bigl(
                (2 z_i)^{\ddot{\chi} (\alpha_i)} \,
                c_{1/(2z_i)} (\ddot{\Theta}_{\alpha_i})
            \bigr)
        \biggr) \biggr]. }
    \end{align*}
\end{definition}

\begin{lemma}
    The operations $X_n^\sd$ satisfy the relations
    \textnormal{\eqref{eq-tw-mod-assoc-x1}--\eqref{eq-tw-mod-assoc-x2}.}
\end{lemma}

\begin{proof}
    \allowdisplaybreaks
    For~\eqref{eq-tw-mod-assoc-x1},
    suppose that $A_i \in V_{\alpha_i}$\,, $B_i \in V_{\beta_i}$\,,
    and $M \in W_\theta$\,, for classes
    $\alpha_i, \beta_i \in \pi_0 (X)$ and $\theta \in \pi_0 (X^\sd)$.
    Suppose that $\deg A_i = a_i$ and $\deg B_i = b_i$ for all $i$.
    Write $\alpha_\ell = \beta_1 + \cdots + \beta_n$\,, and write
    $a_\ell = b_1 + \cdots + b_n$\,. Set
    \begin{align*}
        N & =
        A_1 \otimes \cdots \otimes A_{\ell-1} \otimes
        B_1 \otimes \cdots \otimes B_{n} \otimes
        A_{\ell+1} \otimes \cdots \otimes A_m \otimes M
        \\ & \in
        V_{\alpha_1} \otimes \cdots \otimes V_{\alpha_{\ell-1}} \otimes
        V_{\beta_1} \otimes \cdots \otimes V_{\beta_n} \otimes
        V_{\alpha_{\ell+1}} \otimes \cdots \otimes V_{\alpha_m} \otimes
        W_{\theta}.
    \end{align*}
    Expanding the left-hand side of~\eqref{eq-tw-mod-assoc-x1}, we obtain
    \begin{align*}
        \hspace{1em} & \hspace{-1em} \smash[b]{
        (-1)^{
            \sum_{1 \leq i < j \leq m} a_i \hat{\chi} (\alpha_j) \, + \,
            \sum_{1 \leq i < j \leq n} b_i \hat{\chi} (\beta_j) \, + \,
            \sum_{1 \leq i \leq m} a_i \ring{\chi} (\theta)
        } \cdot
        \varepsilon^\sd (\alpha_1, \dotsc, \alpha_m, \theta) \cdot
        \varepsilon (\beta_1, \dotsc, \beta_n) \cdot {}
        } \\[\LSs]
        & \smash{
        \bigl( \oplus \sss{\alpha_1, \> \dotsc \> , \> \alpha_m, \> \theta}{\sd} \bigr)_* \circ
        \bigl( \upe^{y_1 D_{\alpha_1}} \otimes \cdots \otimes
        \upe^{y_m D_{\alpha_m}} \otimes \id_{W_\theta} \bigr)
        } \\[\LS]
        & \smash{
        \Biggl\{
            \Biggl[
                \biggl(
                    \id_{V_{\alpha_1} \otimes \cdots \otimes V_{\alpha_{\ell-1}}} \otimes
                    \Bigl(
                        \bigl( \oplus_{\beta_1, \> \dotsc \> , \> \beta_n} \bigr)_* \circ
                        \bigl(
                            \upe^{z_1 D_{\beta_1}} \otimes \cdots
                            \otimes \upe^{z_n D_{\beta_n}}
                        \bigr)
                    \Bigr) \otimes
                    \id_{V_{\alpha_{\ell+1}} \otimes \cdots \otimes V_{\alpha_m} \otimes W_\theta}
                \biggr)
        } \\[\LSl]
        & \hspace{3em} \smash{
                \biggl(
                    N \cap
                    \prod_{1 \leq i < j \leq n} {} \bigl(
                        (z_i - z_j)^{\chi (\beta_i, \> \beta_j)} \,
                        c_{1/(z_i - z_j)} (\Theta_{\beta_i, \> \beta_j})
                    \bigr)
                \biggr)
            \Biggr]
        } \\[\LSl]
        & \hspace{1.5em} \smash{
            \mathllap{ {} \cap \biggl( \hspace{-.3em} }
                \LAP{2em}{ \prod_{1 \leq i < j \leq m} }
                \bigl(
                    (y_i - y_j)^{\chi (\alpha_i, \> \alpha_j)} \,
                    c_{1/(y_i - y_j)} (\Theta_{\alpha_i, \> \alpha_j})
                \bigr) \cdot {}
        } \\[\LSl]
        & \hspace{1.5em} \smash{
                \LAP{2em}{ \prod_{1 \leq i < j \leq m} }
                \bigl(
                    (y_i + y_j)^{\chi (\alpha_i, \> \alpha_j^\vee)} \,
                    c_{1/(y_i + y_j)} (\Theta_{\alpha_i, \> \alpha_j^\vee}) 
                \bigr) \cdot {}
        } \\[\LSl]
        & \hspace{1.5em} \smash[t]{
                \LAP{2em}{ \prod_{1 \leq i \leq m} }
                \bigl(
                    y\sss{i}{\dot{\chi} (\alpha_i, \> \theta)} \,
                    c_{1/y_i} (\dot{\Theta}_{\alpha_i, \> \theta}) 
                \bigr) \cdot
                \prod_{1 \leq i \leq m} {} \bigl(
                    (2 y_i)^{\ddot{\chi} (\alpha_i)} \,
                    c_{1/(2y_i)} (\ddot{\Theta}_{\alpha_i})
                \bigr)
            \biggr) 
        \Biggr\} . }
        \numberthis
    \end{align*}
    Using the relations
    \eqref{eq-theta-linear-first}--\eqref{eq-theta-linear-second},
    \eqref{eq-theta-sym},
    \eqref{eq-theta-dot-linear-first},
    \eqref{eq-theta-ddot-linear}
    to expand the terms involving~$\alpha_\ell$\,, this becomes
    \begin{align*}
        \hspace{1em} & \hspace{-1em} \smash[b]{
        (-1)^{
            \sum_{1 \leq i < j \leq m} a_i \hat{\chi} (\alpha_j) \, + \,
            \sum_{1 \leq i < j \leq n} b_i \hat{\chi} (\beta_j) \, + \,
            \sum_{1 \leq i \leq m} a_i \ring{\chi} (\theta)
        } \cdot
        \varepsilon^\sd (\alpha_1, \dotsc, \alpha_m, \theta) \cdot
        \varepsilon (\beta_1, \dotsc, \beta_n) \cdot {}
        } \\*[\LSs]
        & \smash{
        \bigl( \oplus \sss{\alpha_1, \> \dotsc \> , \> \alpha_m, \> \theta}{\sd} \bigr)_* \circ
        \bigl( \upe^{y_1 D_{\alpha_1}} \otimes \cdots \otimes
        \upe^{y_m D_{\alpha_m}} \otimes \id_{W_\theta} \bigr)
        } \\[\LS]
        & \smash{
        \Biggl\{
            \Biggl[
                \biggl(
                    \id_{V_{\alpha_1} \otimes \cdots \otimes V_{\alpha_{\ell-1}}} \otimes
                    \Bigl(
                        \bigl( \oplus_{\beta_1, \> \dotsc \> , \> \beta_n} \bigr)_* \circ
                        \bigl(
                            \upe^{z_1 D_{\beta_1}} \otimes \cdots
                            \otimes \upe^{z_n D_{\beta_n}}
                        \bigr)
                    \Bigr) \otimes
                    \id_{V_{\alpha_{\ell+1}} \otimes \cdots \otimes V_{\alpha_m} \otimes W_\theta}
                \biggr)
        } \\[\LSl]
        & \hspace{3em} \smash{
                \biggl(
                    N \cap
                    \prod_{1 \leq i < j \leq n} {} \bigl(
                        (z_i - z_j)^{\chi (\beta_i, \> \beta_j)} \,
                        c_{1/(z_i - z_j)} (\Theta_{\beta_i, \> \beta_j})
                    \bigr)
                \biggr)
            \Biggr]
        } \\[\LSl]
        & \hspace{1.5em} \smash{
            \mathllap{ {} \cap \Biggl( \hspace{-.3em} }
                \LAP{2em}{ \prod_{ \substack{ 1 \leq i < j \leq m \\ i, j \neq \ell} } }
                \bigl(
                    (y_i - y_j)^{\chi (\alpha_i, \> \alpha_j)} \,
                    c_{1/(y_i - y_j)} (\Theta_{\alpha_i, \> \alpha_j})
                \bigr) \cdot {}
        } \\[\LSl]
        & \hspace{1.5em} \smash{
                \LAP{2em}{ \prod_{ \substack{ 1 \leq i \leq m, \ i \neq \ell \\ 1 \leq j \leq n } } }
                \bigl(
                    (y_i - y_\ell)^{\chi (\alpha_i, \> \beta_j)} \,
                    c_{1/(y_i - y_\ell)} (\Theta_{\alpha_i, \> \beta_j})
                \bigr) \cdot {}
        } \\[\LSl]
        & \hspace{1.5em} \smash{
                \LAP{2em}{ \prod_{ \substack{ 1 \leq i < j \leq m \\ i, j \neq \ell} } }
                \bigl(
                    (y_i + y_j)^{\chi (\alpha_i, \> \alpha_j^\vee)} \,
                    c_{1/(y_i + y_j)} (\Theta_{\alpha_i, \> \alpha_j^\vee}) 
                \bigr) \cdot {}
        } \\[\LSl]
        & \hspace{1.5em} \smash{
                \LAP{2em}{ \prod_{ \substack{ 1 \leq i \leq m, \ i \neq \ell \\ 1 \leq j \leq n } } }
                \bigl(
                    (y_i + y_\ell)^{\chi (\alpha_i, \> \beta_j^\vee)} \,
                    c_{1/(y_i + y_\ell)} (\Theta_{\alpha_i, \> \beta_j^\vee})
                \bigr) \cdot {}
        } \\[\LSl]
        & \hspace{1.5em} \smash{
                \LAP{2em}{ \prod_{ 1 \leq i \leq m, \ i \neq \ell } }
                \bigl(
                    y\sss{i}{\dot{\chi} (\alpha_i, \> \theta)} \,
                    c_{1/y_i} (\dot{\Theta}_{\alpha_i, \> \theta}) 
                \bigr) \cdot 
                \prod_{ 1 \leq j \leq n } {}
                \bigl(
                    y\sss{\ell}{\dot{\chi} (\beta_j, \> \theta)} \,
                    c_{1/y_\ell} (\dot{\Theta}_{\beta_j, \> \theta}) 
                \bigr) \cdot {}
        } \\[\LSl]
        & \hspace{1.5em} \smash{
                \LAP{2em}{ \prod_{ 1 \leq i \leq m, \ i \neq \ell } }
                \bigl(
                    (2 y_i)^{\ddot{\chi} (\alpha_i)} \,
                    c_{1/(2y_i)} (\ddot{\Theta}_{\alpha_i})
                \bigr) \cdot {}
        } \\[\LSl]
        & \hspace{1.5em} \smash{
                \LAP{2em}{ \prod_{ 1 \leq i < j \leq n } }
                \bigl(
                    (2 y_\ell)^{\chi (\beta_i, \> \beta_j^\vee)} \,
                    c_{1/(2y_\ell)} (\Theta_{\beta_i, \> \beta_j^\vee})
                \bigr) \cdot {}
        } \\*[\LSl]
        & \hspace{1.5em} \smash[t]{
                \LAP{2em}{ \prod_{ 1 \leq j \leq n } }
                \bigl(
                    (2 y_\ell)^{\ddot{\chi} (\beta_j)} \,
                    c_{1/(2y_\ell)} (\ddot{\Theta}_{\beta_j})
                \bigr)
            \Biggr) 
        \Biggr\} . }
        \numberthis
    \end{align*}
    Next, we apply Lemma~\ref{lem-cd-ezd-chern} to move
    all the Chern classes into the inner layer $(N \cap {\cdots})$.
    We may also merge the two direct sum operations.
    This gives
    \begin{align*}
        \hspace{2em} & \hspace{-2em} \smash[b]{
        (-1)^{
            \sum_{1 \leq i < j \leq m} a_i \hat{\chi} (\alpha_j) \, + \,
            \sum_{1 \leq i < j \leq n} b_i \hat{\chi} (\beta_j) \, + \,
            \sum_{1 \leq i \leq m} a_i \ring{\chi} (\theta)
        } \cdot
        \varepsilon^\sd (\alpha_1, \dotsc, \alpha_m, \theta) \cdot
        \varepsilon (\beta_1, \dotsc, \beta_n) \cdot {}
        } \\*[\LSs]
        & \smash{
        \bigl(
            \oplus \sss{\alpha_1, \> \dotsc \> , \> \alpha_{\ell-1}, \> 
            \beta_1, \> \dotsc \> , \> \beta_n, \,
            \alpha_{\ell+1}, \> \dotsc \> , \> \alpha_m, \> \theta}{\sd}
        \bigr)_* \circ \Bigl(
            \upe^{y_1 D_{\alpha_1}} \otimes \cdots \otimes
            \upe^{y_{\ell-1} D_{\alpha_{\ell-1}}} \otimes {}
        } \\[\LSs]
        & \hspace{2.5em} \smash{
            \upe^{(y_\ell \, + \, z_1) D_{\beta_1}} \otimes \cdots \otimes
            \upe^{(y_\ell \, + \, z_n) D_{\beta_n}} \otimes
            \upe^{y_{\ell+1} D_{\alpha_{\ell+1}}} \otimes \cdots \otimes
            \upe^{y_m D_{\alpha_m}} \otimes \id_{W_\theta}
        \Bigr)
        } \\[\LS]
        & \hspace{2.5em} \smash{
        \mathllap{ \Biggl\{
            N \cap \Biggl( \hspace{-.3em} }
                \LAP{2em}{ \prod_{1 \leq i < j \leq n} }
                \bigl(
                    (z_i - z_j)^{\chi (\beta_i, \> \beta_j)} \,
                    c_{1/(z_i - z_j)} (\Theta_{\beta_i, \> \beta_j})
                \bigr) \cdot {}
        } \\[\LSl]
        & \hspace{2.5em} \smash{
                \LAP{2em}{ \prod_{ \substack{ 1 \leq i < j \leq m \\ i, j \neq \ell} } }
                \bigl(
                    (y_i - y_j)^{\chi (\alpha_i, \> \alpha_j)} \,
                    c_{1/(y_i - y_j)} (\Theta_{\alpha_i, \> \alpha_j})
                \bigr) \cdot {}
        } \\[\LSl]
        & \hspace{2.5em} \smash{
                \LAP{2em}{ \prod_{ \substack{ 1 \leq i \leq m, \ i \neq \ell \\ 1 \leq j \leq n } } }
                \bigl(
                    (y_i - y_\ell - z_j)^{\chi (\alpha_i, \> \beta_j)} \,
                    c_{1/(y_i - y_\ell - z_j)} (\Theta_{\alpha_i, \> \beta_j})
                \bigr) \cdot {}
        } \\[\LSl]
        & \hspace{2.5em} \smash{
                \LAP{2em}{ \prod_{ \substack{ 1 \leq i < j \leq m \\ i, j \neq \ell} } }
                \bigl(
                    (y_i + y_j)^{\chi (\alpha_i, \> \alpha_j^\vee)} \,
                    c_{1/(y_i + y_j)} (\Theta_{\alpha_i, \> \alpha_j^\vee}) 
                \bigr) \cdot {}
        } \\[\LSl]
        & \hspace{2.5em} \smash{
                \LAP{2em}{ \prod_{ \substack{ 1 \leq i \leq m, \ i \neq \ell \\ 1 \leq j \leq n } } }
                \bigl(
                    (y_i + y_\ell + z_j)^{\chi (\alpha_i, \> \beta_j^\vee)} \,
                    c_{1/(y_i + y_\ell + z_j)} (\Theta_{\alpha_i, \> \beta_j^\vee})
                \bigr) \cdot {}
        } \\[\LSl]
        & \hspace{2.5em} \smash{
                \LAP{2em}{ \prod_{ 1 \leq i \leq m, \ i \neq \ell } }
                \bigl(
                    y\sss{i}{\dot{\chi} (\alpha_i, \> \theta)} \,
                    c_{1/y_i} (\dot{\Theta}_{\alpha_i, \> \theta}) 
                \bigr) \cdot {}
        } \\[\LSl]
        & \hspace{2.5em} \smash{
                \LAP{2em}{ \prod_{ 1 \leq j \leq n } }
                \bigl(
                    (y_\ell + z_j)^{\dot{\chi} (\beta_j, \> \theta)} \,
                    c_{1/(y_\ell + z_j)} (\dot{\Theta}_{\beta_j, \> \theta}) 
                \bigr) \cdot {}
        } \\[\LSl]
        & \hspace{2.5em} \smash{
                \LAP{2em}{ \prod_{ 1 \leq i \leq m, \ i \neq \ell } }
                \bigl(
                    (2 y_i)^{\ddot{\chi} (\alpha_i)} \,
                    c_{1/(2 y_i)} (\ddot{\Theta}_{\alpha_i})
                \bigr) \cdot {}
        } \\[\LSl]
        & \hspace{2.5em} \smash{
                \LAP{2em}{ \prod_{ 1 \leq i < j \leq n } }
                \bigl(
                    (2 y_\ell + z_i + z_j)^{\chi (\beta_i, \> \beta_j^\vee)} \,
                    c_{1/(2 y_\ell + z_i + z_j)} (\Theta_{\beta_i, \> \beta_j^\vee})
                \bigr) \cdot {}
        } \\*[\LSl]
        & \hspace{2.5em} \smash[t]{
                \LAP{2em}{ \prod_{ 1 \leq j \leq n } }
                \bigl(
                    (2 y_\ell + 2 z_j)^{\ddot{\chi} (\beta_j)} \,
                    c_{1/(2 y_\ell + 2 z_j)} (\ddot{\Theta}_{\beta_j})
                \bigr)
            \Biggr) 
        \Biggr\} , }
        \numberthis
    \end{align*}
    which matches the right-hand side of~\eqref{eq-tw-mod-assoc-x1}.
    Note that we have
    \begin{equation}
        \varepsilon^\sd (\alpha_1, \dotsc, \alpha_m, \theta) \cdot
        \varepsilon (\beta_1, \dotsc, \beta_n) =
        \varepsilon^\sd (\alpha_1, \dotsc, \alpha_{\ell-1},
            \beta_1, \dotsc, \beta_n,
            \alpha_{\ell+1}, \dotsc, \alpha_m, \theta) ,
    \end{equation}
    which can be deduced from
    \eqref{eq-epsilon-assoc} and~\eqref{eq-epsilon-sd-assoc}
    via a combinatorial argument.

    To verify \eqref{eq-tw-mod-assoc-x2},
    suppose that $A_i \in V_{\alpha_i}$\,, $B_i \in V_{\beta_i}$\,,
    and $M \in W_\theta$\,, for classes
    $\alpha_i, \beta_i \in \pi_0 (X)$ and $\rho \in \pi_0 (X^\sd)$.
    Suppose that $\deg A_i = a_i$ and $\deg B_i = b_i$ for all $i$.
    Write $\theta = \bar{\beta}_1 + \cdots + \bar{\beta}_n + \rho$\,. Set
    \begin{align*}
        N & =
        A_1 \otimes \cdots \otimes A_m \otimes
        B_1 \otimes \cdots \otimes B_n \otimes M
        \\ & \in
        V_{\alpha_1} \otimes \cdots \otimes V_{\alpha_m} \otimes
        V_{\beta_1} \otimes \cdots \otimes V_{\beta_n} \otimes
        W_{\rho} \, .
    \end{align*}
    Expanding the left-hand side of~\eqref{eq-tw-mod-assoc-x2}, we obtain
    \begin{align*}
        \hspace{2em} & \hspace{-2em} \smash[b]{
        (-1)^{
            \sum_{1 \leq i < j \leq m} a_i \hat{\chi} (\alpha_j) \, + \,
            \sum_{1 \leq i < j \leq n} b_i \hat{\chi} (\beta_j) \, + \,
            \sum_{1 \leq i \leq m} a_i \ring{\chi} (\theta) \, + \,
            \sum_{1 \leq i \leq n} b_i \ring{\chi} (\rho)
        } \cdot {}
        } \\*[\LSs] 
        & \smash{
        \varepsilon^\sd (\alpha_1, \dotsc, \alpha_m, \theta) \cdot
        \varepsilon^\sd (\beta_1, \dotsc, \beta_n, \rho) \cdot {}
        } \\*[\LSs]
        & \smash{
        \bigl( \oplus \sss{\alpha_1, \> \dotsc \> , \> \alpha_m, \> \theta}{\sd} \bigr)_* \circ
        \bigl( \upe^{y_1 D_{\alpha_1}} \otimes \cdots \otimes
        \upe^{y_m D_{\alpha_m}} \otimes
        \id_{V_{\beta_1} \otimes \cdots \otimes V_{\beta_n} \otimes W_\theta} \bigr)
        } \\[\LS]
        & \smash{
        \Biggl\{
            \Biggl[
                \biggl(
                    \id_{V_{\alpha_1} \otimes \cdots \otimes V_{\alpha_m}} \otimes
                    \Bigl(
                        \bigl( \oplus \sss{\beta_1, \> \dotsc \> , \> \beta_n, \> \rho}{\sd} \bigr)_* \circ
                        \bigl(
                            \upe^{z_1 D_{\beta_1}} \otimes \cdots
                            \otimes \upe^{z_n D_{\beta_n}} \otimes \id_{W_\rho}
                        \bigr)
                    \Bigr)
                \biggr)
        } \\[\LSl]
        & \hspace{4.5em} \smash{
                \mathllap{ \Biggl(
                    N \cap \biggl( \hspace{-.3em} }
                        \LAP{2em}{ \prod_{1 \leq i < j \leq n} }
                        \bigl(
                            (z_i - z_j)^{\chi (\beta_i, \> \beta_j)} \,
                            c_{1/(z_i - z_j)} (\Theta_{\beta_i, \> \beta_j})
                        \bigr) \cdot {}
        } \\[\LSl]
        & \hspace{4.5em} \smash{
                        \LAP{2em}{ \prod_{1 \leq i < j \leq n} }
                        \bigl(
                            (z_i + z_j)^{\chi (\beta_i, \> \beta_j^\vee)} \,
                            c_{1/(z_i + z_j)} (\Theta_{\beta_i, \> \beta_j^\vee}) 
                        \bigr) \cdot {}
        } \\[\LSl]
        & \hspace{4.5em} \smash{
                        \LAP{2em}{ \prod_{1 \leq i \leq n} } \bigl(
                            z\sss{i}{\dot{\chi} (\beta_i, \> \rho)} \,
                            c_{1/z_i} (\dot{\Theta}_{\beta_i, \> \rho}) 
                        \bigr) \cdot
                        \prod_{1 \leq i \leq n} {} \bigl(
                            (2 z_i)^{\ddot{\chi} (\beta_i)} \,
                            c_{1/(2z_i)} (\ddot{\Theta}_{\beta_i})
                        \bigr)
                    \biggr) 
                \Biggr)
            \Biggr]
        } \\[\LSl]
        & \hspace{1.5em} \smash{
            \mathllap{ {} \cap \biggl( \hspace{-.3em} }
                \LAP{2em}{ \prod_{1 \leq i < j \leq m} }
                \bigl(
                    (y_i - y_j)^{\chi (\alpha_i, \> \alpha_j)} \,
                    c_{1/(y_i - y_j)} (\Theta_{\alpha_i, \> \alpha_j})
                \bigr) \cdot {}
        } \\*[\LSl]
        & \hspace{1.5em} \smash{
                \LAP{2em}{ \prod_{1 \leq i < j \leq m} }
                \bigl(
                    (y_i + y_j)^{\chi (\alpha_i, \> \alpha_j^\vee)} \,
                    c_{1/(y_i + y_j)} (\Theta_{\alpha_i, \> \alpha_j^\vee}) 
                \bigr) \cdot {}
        } \\*[\LSl]
        & \hspace{1.5em} \smash[t]{
                \LAP{2em}{ \prod_{1 \leq i \leq m} }
                \bigl(
                    y\sss{i}{\dot{\chi} (\alpha_i, \> \theta)} \,
                    c_{1/y_i} (\dot{\Theta}_{\alpha_i, \> \theta}) 
                \bigr) \cdot
                \prod_{1 \leq i \leq m} {} 
                \bigl(
                    (2 y_i)^{\ddot{\chi} (\alpha_i)} \,
                    c_{1/(2y_i)} (\ddot{\Theta}_{\alpha_i})
                \bigr)
            \biggr) 
        \Biggr\} . }
        \numberthis
    \end{align*}
    Using the relation~\eqref{eq-theta-dot-linear-second}, we may rewrite
    \begin{align*}
        &  \smash[b]{
        \prod_{1 \leq i \leq m} {} \bigl(
            y\sss{i}{\dot{\chi} (\alpha_i, \> \theta)} \,
            c_{1/y_i} (\dot{\Theta}_{\alpha_i, \> \theta}) 
        \bigr) = {}
        } \\*[\LSl]
        & \hspace{2em} \smash{
        \LAP{2em}{ \prod_{ \substack{ 1 \leq i \leq m \\ 1 \leq j \leq n } } }
        \bigl(
            y\sss{i}{\chi (\alpha_i, \> \beta_j)}
            c_{1/y_i} (\Theta_{\alpha_i, \> \beta_j}) 
        \bigr) \cdot
        \prod_{ \substack{ 1 \leq i \leq m \\ 1 \leq j \leq n } } {}
        \bigl(
            y\sss{i}{\chi (\alpha_i, \> \beta_j^\vee)}
            c_{1/y_i} (\Theta_{\alpha_i, \> \beta_j^\vee}) 
        \bigr) \cdot {}
        } \\*[\LSl]
        & \hspace{2em} \smash[t]{
        \LAP{2em}{ \prod_{ 1 \leq i \leq m } }
        \bigl(
            y\sss{i}{\dot{\chi} (\alpha_i, \> \rho)}
            c_{1/y_i} (\dot{\Theta}_{\alpha_i, \> \rho}) 
        \bigr) . }
    \end{align*}
    Using this, and applying Lemma~\ref{lem-cd-ezd-chern} to move
    all the Chern classes into the inner layer $(N \cap {\cdots})$,
    we obtain
    \begin{align*}
        \hspace{2em} & \hspace{-2em} \smash[b]{
        (-1)^{
            \sum_{1 \leq i < j \leq m} a_i \hat{\chi} (\alpha_j) \, + \,
            \sum_{1 \leq i < j \leq n} b_i \hat{\chi} (\beta_j) \, + \,
            \sum_{ \substack{ 
                \scriptscriptstyle 1 \leq i \leq m \\
                \scriptscriptstyle 1 \leq j \leq n
            } } a_i \hat{\chi} (\beta_j) \, + \,
            \sum_{1 \leq i \leq m} a_i \ring{\chi} (\rho) \, + \,
            \sum_{1 \leq i \leq n} b_i \ring{\chi} (\rho)
        } \cdot {}
        } \\*[\LSs]
        & \smash{
        \varepsilon^\sd (\alpha_1, \dotsc, \alpha_m, \theta) \cdot
        \varepsilon^\sd (\beta_1, \dotsc, \beta_n, \rho) \cdot {}
        } \\*[\LSs]
        & \smash{ \bigl( \oplus \sss{\alpha_1, \> \dotsc \> , \> \alpha_m, \,
            \beta_1, \> \dotsc \> , \> \beta_n, \> \rho}{\sd} \bigr)_* \circ
        \Bigl( \upe^{y_1 D_{\alpha_1}} \otimes \cdots \otimes
            \upe^{y_m D_{\alpha_m}} \otimes
            \upe^{z_1 D_{\beta_1}} \otimes \cdots
            \otimes \upe^{z_n D_{\beta_n}} \otimes \id_{W_\rho} \Bigr)
        } \\[\LS]
        & \hspace{2.5em} \smash{
        \mathllap{ \Biggl\{ N \cap 
            \biggl( \hspace{-.3em} }
                \LAP{2em}{ \prod_{1 \leq i < j \leq n} }
                \bigl(
                    (z_i - z_j)^{\chi (\beta_i, \> \beta_j)} \,
                    c_{1/(z_i - z_j)} (\Theta_{\beta_i, \> \beta_j})
                \bigr) \cdot {}
        } \\[\LSl]
        & \hspace{2.5em} \smash{
                \LAP{2em}{ \prod_{1 \leq i < j \leq n} }
                \bigl(
                    (z_i + z_j)^{\chi (\beta_i, \> \beta_j^\vee)} \,
                    c_{1/(z_i + z_j)} (\Theta_{\beta_i, \> \beta_j^\vee}) 
                \bigr) \cdot {}
        } \\[\LSl]
        & \hspace{2.5em} \smash{
                \LAP{2em}{ \prod_{1 \leq i \leq n} }
                \bigl(
                    z\sss{i}{\dot{\chi} (\beta_i, \> \rho)} \,
                    c_{1/z_i} (\dot{\Theta}_{\beta_i, \> \rho}) 
                \bigr) \cdot
                \prod_{1 \leq i \leq n} {} \bigl(
                    (2 z_i)^{\ddot{\chi} (\beta_i)} \,
                    c_{1/(2z_i)} (\ddot{\Theta}_{\beta_i})
                \bigr) \cdot {}
        } \\[\LSl]
        & \hspace{2.5em} \smash{
                \LAP{2em}{ \prod_{1 \leq i < j \leq m} }
                \bigl(
                    (y_i - y_j)^{\chi (\alpha_i, \> \alpha_j)} \,
                    c_{1/(y_i - y_j)} (\Theta_{\alpha_i, \> \alpha_j})
                \bigr) \cdot {}
        } \\[\LSl]
        & \hspace{2.5em} \smash{
                \LAP{2em}{ \prod_{1 \leq i < j \leq m} }
                \bigl(
                    (y_i + y_j)^{\chi (\alpha_i, \> \alpha_j^\vee)} \,
                    c_{1/(y_i + y_j)} (\Theta_{\alpha_i, \> \alpha_j^\vee}) 
                \bigr) \cdot {}
        } \\[\LSl]
        & \hspace{2.5em} \smash{
                \LAP{2em}{ \prod_{ \substack{ 1 \leq i \leq m \\ 1 \leq j \leq n } } }
                \bigl(
                    (y_i - z_j)^{\chi (\alpha_i, \> \beta_j)}
                    c_{1/(y_i - z_j)} (\Theta_{\alpha_i, \> \beta_j}) 
                \bigr) \cdot {}
        } \\[\LSl]
        & \hspace{2.5em} \smash{
                \LAP{2em}{ \prod_{ \substack{ 1 \leq i \leq m \\ 1 \leq j \leq n } } }
                \bigl(
                    (y_i + z_j)^{\chi (\alpha_i, \> \beta_j^\vee)}
                    c_{1/(y_i + z_j)} (\Theta_{\alpha_i, \> \beta_j^\vee}) 
                \bigr) \cdot {}
        } \\[\LSl]
        & \hspace{2.5em} \smash[t]{
                \LAP{2em}{ \prod_{ 1 \leq i \leq m } }
                \bigl(
                    y\sss{i}{\dot{\chi} (\alpha_i, \> \rho)}
                    c_{1/y_i} (\dot{\Theta}_{\alpha_i, \> \rho}) 
                \bigr) \cdot
                \prod_{1 \leq i \leq m} {} \bigl(
                    (2 y_i)^{\ddot{\chi} (\alpha_i)} \,
                    c_{1/(2y_i)} (\ddot{\Theta}_{\alpha_i})
                \bigr)
            \biggr) 
        \Biggr\} , }
        \numberthis
    \end{align*}
    which matches the right-hand side of~\eqref{eq-tw-mod-assoc-x2}.
    Note that we have
    \begin{equation}
        \varepsilon^\sd (\alpha_1, \dotsc, \alpha_m, \theta) \cdot
        \varepsilon^\sd (\beta_1, \dotsc, \beta_n, \rho) =
        \varepsilon^\sd (\alpha_1, \dotsc, \alpha_m, 
            \beta_1, \dotsc, \beta_n, \rho),
    \end{equation}
    which can be deduced from
    \eqref{eq-epsilon-assoc} and~\eqref{eq-epsilon-sd-assoc}
    via a combinatorial argument.
\end{proof}

\subsection[Proof of Theorem~\ref*{thm-mor-va}]{Proof of Theorem~\ref{thm-mor-va}}
\label{sect-proof-mor-va}

Write
\begin{alignat*}{2}
    V & = \hat{H}_* (X; \bbQ), & \qquad
    V_\alpha & = \hat{H}_* (X_\alpha; \bbQ), \\
    V' & = \hat{H}_* (X'; \bbQ), & \qquad
    V'_{\smash{\alpha'}} & = \hat{H}_* (X'_{\smash{\alpha'}}; \bbQ)
\end{alignat*}
for $\alpha \in \pi_0 (X)$ and $\alpha' \in \pi_0 (X')$.

\begin{lemma}
    \label{lem-pf-mor-va}
    For any $\alpha, \beta \in \pi_0 (X),$
    and any $A \in V_\alpha$ and $B \in V_\beta$\,, we have
    \begin{align}
        & \smash[b]{
        (\oplus_{\alpha, \> \beta})_* \circ (\upe^{z D_\alpha} \otimes \id_{V_\beta})
        \Bigl[ 
            (A \boxtimes B) \cap 
            \bigl( z^{\chi (\alpha, \> \beta)} \, c_{1/z} (\Theta_{\alpha, \> \beta}) \bigr)
        \Bigr]
        \cap c_{\hat{\xi} (\alpha + \beta)} (\hat{\Xi} \sss{\alpha + \beta}{\vee})
        } \notag \\[\LS]
        = {}
        & \smash{
        (-1)^{\xi (\alpha, \> \beta)} \cdot
        (\oplus_{\alpha, \> \beta})_* \circ (\upe^{z D_\alpha} \otimes \id_{V_\beta}) \biggl\{
            \Bigl[
                \bigl( A \cap c_{\hat{\xi} (\alpha)} (\hat{\Xi} \sss{\alpha}{\vee}) \bigr) \boxtimes
                \bigl( B \cap c_{\hat{\xi} (\beta)} (\hat{\Xi} \sss{\beta}{\vee}) \bigr)
            \Bigr]
        } \notag \\*[\LS]
        \label{eq-lem-pf-mor-va}
        & \hspace{2em} \smash[t]{
            {} \cap \biggl[ 
                z^{\chi (\alpha, \> \beta)} \,
                c_{1/z} (\Theta_{\alpha, \> \beta}) \cdot
                z^{\xi (\alpha, \> \beta)} \,
                c_{1/z} (\Xi_{\alpha, \> \beta}) \cdot
                z^{\xi (\beta, \> \alpha)} \,
                c_{1/z} (\sigma \sss{\alpha, \> \beta}{*} (\Xi \sss{\beta, \> \alpha}{\vee})) 
            \biggr]
        \biggr\} , }
    \end{align}
    where $\sigma_{\alpha, \> \beta} \colon X_\alpha \times X_\beta
    \to X_\beta \times X_\alpha$ exchanges the two factors.
\end{lemma}

\begin{proof}
    \allowdisplaybreaks
    For convenience in notation, we write
    $\Xi_{\beta, \> \alpha}$\,, etc.,
    for their pullbacks to $X_\alpha \times X_\beta$\,.
    We have
    \begin{align*}
        & \smash[b]{
        (\oplus_{\alpha, \> \beta})_* \circ (\upe^{z D_\alpha} \otimes \id_{V_\beta})
        \Bigl[ 
            (A \boxtimes B) \cap 
            \bigl( z^{\chi (\alpha, \> \beta)} \, c_{1/z} (\Theta_{\alpha, \> \beta}) \bigr)
        \Bigr]
        \cap (-y)^{\hat{\xi} (\alpha + \beta)} \, c_{1/y} (\hat{\Xi}_{\alpha + \beta})
        } \\[\LS]
        = {}
        & \smash{
        (-1)^{\hat{\xi} (\alpha + \beta)} \cdot (\oplus_{\alpha, \> \beta})_* \biggl\{
            (\upe^{z D_\alpha} \otimes \id_{V_\beta})
            \Bigl[ 
                (A \boxtimes B) \cap 
                \bigl( z^{\chi (\alpha, \> \beta)} \, c_{1/z} (\Theta_{\alpha, \> \beta}) \bigr)
            \Bigr]
        } \\*[\LS]
        & \hspace{1em} \smash{
            {} \cap \Bigl[
                y^{\hat{\xi} (\alpha)} \, c_{1/y} (\hat{\Xi}_\alpha) \cdot
                y^{\hat{\xi} (\beta)} \, c_{1/y} (\hat{\Xi}_\beta) \cdot
                y^{\xi (\alpha, \> \beta)} \, c_{1/y} (\Xi_{\alpha, \> \beta}) \cdot
                y^{\xi (\beta, \> \alpha)} \, c_{1/y} (\Xi_{\beta, \> \alpha})
            \Bigr]
        \biggr\}
        } \\[\LS]
        = {}
        & \smash{
        (-1)^{\hat{\xi} (\alpha + \beta)} \cdot 
        (\oplus_{\alpha, \> \beta})_* \circ (\upe^{z D_\alpha} \otimes \id_{V_\beta}) 
        } \\*[\LS]
        & \hspace{1em} \smash{
        \biggl\{
            (A \boxtimes B)
            {} \cap \Bigl[
                z^{\chi (\alpha, \> \beta)} \, c_{1/z} (\Theta_{\alpha, \> \beta}) \cdot
                y^{\hat{\xi} (\alpha)} \, c_{1/y} (\hat{\Xi}_\alpha) \cdot
                y^{\hat{\xi} (\beta)} \, c_{1/y} (\hat{\Xi}_\beta) \cdot {}
        } \\*[\LS]
        & \hspace{3em} \smash{
                (y + z)^{\xi (\alpha, \> \beta)} \, c_{1/(y + z)} (\Xi_{\alpha, \> \beta}) \cdot
                (y - z)^{\xi (\beta, \> \alpha)} \, c_{1/(y - z)} (\Xi_{\beta, \> \alpha})
            \Bigr]
        \biggr\}
        } \\[\LS]
        = {} & \smash{
        (-1)^{\xi (\alpha, \> \beta)} \cdot
        (\oplus_{\alpha, \> \beta})_* \circ (\upe^{z D_\alpha} \otimes \id_{V_\beta}) 
        } \\*[\LS]
        & \hspace{1em} \smash{
        \biggl\{
            \Bigl[
                \bigl( A \cap (-y)^{\hat{\xi} (\alpha)} \, c_{1/y} (\hat{\Xi}_\alpha) \bigr) \boxtimes
                \bigl( B \cap (-y)^{\hat{\xi} (\beta)} \, c_{1/y} (\hat{\Xi}_\beta) \bigr)
            \Bigr]
        } \\*[\LS]
        & \hspace{3em} \smash{
            {} \cap \biggl[ 
                z^{\chi (\alpha, \> \beta)} \, c_{1/z} (\Theta_{\alpha, \> \beta}) \cdot {}
        } \\*[\LS]
        & \hspace{6em} \smash[t]{
                (y + z)^{\xi (\alpha, \> \beta)} \, c_{1/(y + z)} (\Xi_{\alpha, \> \beta}) \cdot
                (z - y)^{\xi (\beta, \> \alpha)} \, c_{1/(y - z)} (\Xi_{\beta, \> \alpha})
            \biggr]
        \biggr\} , }
        \numberthis
    \end{align*}
    where the first step uses Assumption~\ref{asn-h-sp-mor}~\ref{itm-xi-check},
    and the second step uses Lemma~\ref{lem-cd-ezd-chern}.
    The third step uses~\eqref{eq-xi-check-rk} to adjust the signs.

    The above equation is, strictly speaking, an equality in 
    $V_{\smash{\alpha + \beta}} \llparen y \rrparen \llparen z \rrparen$,
    where we apply the inclusion
    $\iota_{y, z} \colon V_{\smash{\alpha + \beta}} \llparen y, z \rrparen
    \hookrightarrow V_{\smash{\alpha + \beta}} \llparen y \rrparen \llparen z \rrparen$
    to the right-hand side.
    However, note that the left-hand side lies in the subspace
    $V_{\smash{\alpha + \beta}} \llparen z \rrparen [y]$,
    since by Assumption~\ref{asn-h-sp-mor}~\ref{itm-xi},
    the class $\hat{\Xi}_{\smash{\alpha + \beta}}$ comes from a vector bundle,
    so that its Chern classes above $\hat{\xi} (\alpha + \beta)$ are zero.
    Therefore, the right-hand side must also lie in $V_{\smash{\alpha + \beta}} \llparen z \rrparen [y]$.
    This allows us to take $y = 0$ on both sides,
    giving the desired equality~\eqref{eq-lem-pf-mor-va}.
\end{proof}

\begin{lemma}
    \label{lem-omega-compat-y}
    For any $\alpha, \beta \in \pi_0 (X),$
    and any $A \in V_\alpha$ and $B \in V_\beta$\,, we have
    \begin{equation}
        \Omega (Y (A, z) \, B) = Y' (\Omega (A), z) \, \Omega (B),
    \end{equation}
    where $Y, Y'$ are structure maps for the vertex algebra structures on $V, V',$
    and $\Omega$ is the map given in Theorem~\textnormal{\ref{thm-mor-va}}.
\end{lemma}

\begin{proof}
    \allowdisplaybreaks
    For $A \in H_a (X_\alpha; \bbQ)$ and $B \in H_b (X_\beta; \bbQ)$,
    we have
    \begin{align*}
        & \smash[b]{
        \Omega (Y (A, z) \, B)
        } \\[\LS]
        = {}
        & \smash{
        (F_{\alpha + \beta})_* \biggl\{
            (-1)^{a \hat{\chi} (\beta)} \varepsilon (\alpha, \beta) \cdot
            (\oplus_{\alpha, \> \beta})_* \circ (\upe^{z D_\alpha} \otimes \id_{V_\beta})
        } \\*[\LS]
        & \hspace{2em} \smash{
            \Bigl[ 
                (A \boxtimes B) \cap 
                \bigl( z^{\chi (\alpha, \> \beta)} \, c_{1/z} (\Theta_{\alpha, \> \beta}) \bigr)
            \Bigr]
            \cap c_{\hat{\xi} (\alpha + \beta)} (\hat{\Xi} \sss{\alpha + \beta}{\vee})
        \biggr\}
        } \\[\LS]
        = {}
        & \smash{
        (-1)^{\xi (\alpha, \> \beta) + a \hat{\chi} (\beta)} \varepsilon (\alpha, \beta) \cdot
        (F_{\alpha + \beta})_* \Biggl\{
            (\oplus_{\alpha, \> \beta})_* \circ (\upe^{z D_\alpha} \otimes \id_{V_\beta})
        } \\*[\LS]
        & \hspace{2em} \smash{
            \biggl\{
                \Bigl[
                    \bigl( A \cap c_{\hat{\xi} (\alpha)} (\hat{\Xi} \sss{\alpha}{\vee}) \bigr) \boxtimes
                    \bigl( B \cap c_{\hat{\xi} (\beta)} (\hat{\Xi} \sss{\beta}{\vee}) \bigr)
                \Bigr]
        } \\*[\LS]
        & \hspace{4em} \smash{
                {} \cap \biggl[ 
                    z^{\chi (\alpha, \> \beta)} \,
                    c_{1/z} (\Theta_{\alpha, \> \beta}) \cdot
                    z^{\xi (\alpha, \> \beta)} \,
                    c_{1/z} (\Xi_{\alpha, \> \beta}) \cdot
                    z^{\xi (\beta, \> \alpha)} \,
                    c_{1/z} (\Xi \sss{\beta, \> \alpha}{\vee}) 
                \biggr]
            \biggr\}
        \Biggr\}
        } \\[\LS]
        = {}
        & \smash{
        (-1)^{a \hat{\chi}' (F (\beta))} \varepsilon' (F (\alpha), F (\beta)) \cdot
        (\oplus \sss{\alpha, \> \beta}{\prime})_* \circ
        (\upe^{z D \sss{\alpha}{\prime}} \otimes \id_{V \sss{\beta}{\prime}})
        } \\*[\LS]
        & \hspace{2em} \smash{
            \biggl\{
            \Bigl[
                (F_\alpha)_* \bigl( A \cap c_{\hat{\xi} (\alpha)}
                (\hat{\Xi} \sss{\alpha}{\vee}) \bigr) \boxtimes
                (F_\beta)_* \bigl( B \cap c_{\hat{\xi} (\beta)}
                (\hat{\Xi} \sss{\beta}{\vee}) \bigr)
            \Bigr] \cap \bigl(
                z^{\chi' (\alpha, \> \beta)} \,
                c_{1/z} (\Theta \sss{\alpha, \> \beta}{\prime})
            \bigr)
        \biggr\}
        } \\[\LS]
        = {}
        & \smash[t]{
        Y' (\Omega (A), z) \, \Omega (B), }
        \numberthis
    \end{align*}
    where the first step uses the definitions of~$\Omega$ and~$Y$,
    the second uses Lemma~\ref{lem-pf-mor-va},
    the third uses
    Assumption~\ref{asn-h-sp-mor}~\ref{itm-xi-compat-theta}--\ref{itm-xi-compat-epsilon},
    and the last step uses the definitions of~$\Omega$ and~$Y'$.
\end{proof}

\begin{bproof}[Proof of Theorem~\ref{thm-mor-va}]
    We have to show that the map
    \[
        \Omega \colon V \longrightarrow V'
    \]
    is a morphism of graded vertex algebras.
    By~\eqref{eq-chi-check-diff},
    the map $\Omega$ is compatible with the gradings of $V$ and $V'$.
    By Lemma~\ref{lem-omega-compat-y}, $\Omega$ is compatible with the vertex operations $Y$ and $Y'$.
    
    To see that $\Omega (1) = 1$, note that~\eqref{eq-chi-check-diff} forces
    $\hat{\xi} (0) = 0$, so that $\hat{\Xi}_0 = 0$.

    Finally, using Lemma~\ref{lem-va-facts}~\ref{itm-va-az1},
    it follows from Lemma~\ref{lem-omega-compat-y}
    that $\Omega \circ D = D' \circ \Omega$.
\end{bproof}

\subsection[Proof of Theorem~\ref*{thm-comp-mor-va}]{Proof of Theorem~\ref{thm-comp-mor-va}}
\label{sect-proof-comp-mor-va}

First, let us verify Assumption~\ref{asn-h-sp-mor}
for $F \second$, equipped with the data $\Xi \second$ and $\hat{\Xi} \second$
given in the statement of the theorem.
Verifying Assumption~\ref{asn-h-sp-mor}~\ref{itm-xi}--\ref{itm-xi-check}
is straightforward.
For Assumption~\ref{asn-h-sp-mor}~\ref{itm-xi-compat-theta},
we have
\begin{align*}
    & \smash[b]{
    (F \sss{\alpha}{\dprime} \times F \sss{\beta}{\dprime})^* 
    (\Theta \sss{F \second (\alpha), \> F \second (\beta)}{\dprime}) 
    } \\*[\LS]
    = {} & \smash{
    (F_\alpha \times F_\beta)^* \Bigl(
        \Theta \sss{F (\alpha), \> F (\beta)}{\prime} +
        \Xi \sss{F (\alpha), \> F (\beta)}{\prime} +
        \sigma \sss{F (\alpha), \> F (\beta)}{*}
        (\Xi \sss{F (\beta), \> F (\alpha)}{\prime \vee})
    \Bigr)
    } \\*[\LS]
    = {} & \smash{
    \Bigl(
        \Theta_{\alpha, \> \beta} +
        \Xi_{\alpha, \> \beta} +
        \sigma \sss{\alpha, \> \beta}{*} (\Xi \sss{\beta, \> \alpha}{\vee})    
    \Bigr) +
    (F_\alpha \times F_\beta)^* \Bigl(
        \Xi \sss{F (\alpha), \> F (\beta)}{\prime} +
        \sigma \sss{F (\alpha), \> F (\beta)}{*}
        (\Xi \sss{F (\beta), \> F (\alpha)}{\prime \vee})
    \Bigr)
    } \\*[\LS]
    = {} & \smash[t]{
    \Theta_{\alpha, \> \beta} +
    \Xi \sss{\alpha, \> \beta}{\dprime} +
    \sigma \sss{\alpha, \> \beta}{*}
    (\Xi \sss{\beta, \> \alpha}{\dprime \vee}) .
    }
    \numberthis
\end{align*}
Assumption~\ref{asn-h-sp-mor}~\ref{itm-xi-compat-epsilon}
follows from the relation
\begin{equation}
    \xi \second (\alpha, \beta) =
    \xi (\alpha, \beta) + \xi' (F (\alpha), F (\beta)).
\end{equation}

Finally, to see that
\[
    \Omega \second = \Omega' \circ \Omega,
\] 
we observe that for any $A \in \hat{H}_* (X_\alpha; \bbQ)$, we have
\begin{align*}
    \Omega' \circ \Omega (A) = {} & \smash[b]{
    F'_* \Bigl(
        F_* (A \cap c_{\hat{\xi} (\alpha)}
        (\hat{\Xi} \sss{\alpha}{\vee}))
        \cap c_{\hat{\xi}' (F (\alpha))}
        (\hat{\Xi} \sss{F (\alpha)}{\prime \vee})
    \Bigr)
    } \\*[\LS]
    = {} & \smash{
    F'_* \circ F_* \Bigl[
        A \cap \Bigl( c_{\hat{\xi} (\alpha)}
        (\hat{\Xi} \sss{\alpha}{\vee})
        \cdot c_{\hat{\xi}' (F (\alpha))}
        (F^* (\hat{\Xi} \sss{F (\alpha)}{\prime \vee}) ) \Bigr)
    \Bigr]
    } \\*[\LS]
    = {} & \smash{
    F \second _* \Bigl(
        A \cap c_{\hat{\xi} (\alpha) + \hat{\xi}' (F (\alpha))}
        (\hat{\Xi} \sss{\alpha}{\vee} +
        F^* (\hat{\Xi} \sss{F (\alpha)}{\prime \vee}))  
    \Bigr)
    } \\*[\LS]
    = {} & \smash{
    \Omega \second (A),
    }
    \numberthis
\end{align*}
where in the third step,
we used the condition in Assumption~\ref{asn-h-sp-mor}~\ref{itm-xi}
that the classes $\hat{\Xi}_\alpha$ and $\hat{\Xi} \sss{F (\alpha)}{\prime}$
come from vector bundles,
so their Chern classes above their ranks are zero.
\qed

\subsection[Proof of Theorem~\ref*{thm-mor-tw-mod}]{Proof of Theorem~\ref{thm-mor-tw-mod}}
\label{sect-proof-mor-tw-mod}

\begin{lemma}
    \label{lem-pf-mor-tw-mod}
    For any $\alpha \in \pi_0 (X),$ $\theta \in \pi_0 (X^\sd),$
    and any $A \in V_\alpha$ and $M \in W_\theta$\,, we have
    \begin{align}
        & \smash[b]{
        (\oplus \sss{\alpha, \> \theta}{\sd})_* \circ (\upe^{z D_\alpha} \otimes \id_{W_\theta})
        \Bigl[ 
            (A \boxtimes M) 
        } \notag \\*[\LS]
        & \hspace{3em} \smash{
            {} \cap 
            \Bigl(
                z^{\dot{\chi} (\alpha, \> \theta)} \,
                c_{1/z} (\dot{\Theta}_{\alpha, \> \theta}) \cdot
                (2z)^{\ddot{\chi} (\alpha)} \,
                c_{1/(2z)} (\pr_1^* (\ddot{\Theta}_{\alpha}))
            \Bigr)
        \Bigr]
        \cap c_{\ring{\xi} (\bar{\alpha} + \theta)} (\ring{\Xi} \sss{\bar{\alpha} + \theta}{\vee})
        } \notag \\[\LS]
        = {}
        & \smash{
        (-1)^{\dot{\xi} (\alpha, \> \theta) + \ddot{\xi} (\alpha)} \cdot
        (\oplus \sss{\alpha, \> \theta}{\sd})_* \circ (\upe^{z D_\alpha} \otimes \id_{W_\theta})
        \biggl\{
            \Bigl[
                \bigl( A \cap c_{\hat{\xi} (\alpha)} (\hat{\Xi} \sss{\alpha}{\vee}) \bigr) \boxtimes
                \bigl( M \cap c_{\ring{\xi} (\theta)} (\ring{\Xi} \sss{\theta}{\vee}) \bigr)
            \Bigr]
        } \notag \\*[\LS]
        & \hspace{1em} \smash{
            {} \cap \biggl[ 
                z^{\dot{\chi} (\alpha, \> \theta)} \,
                c_{1/z} (\dot{\Theta}_{\alpha, \> \theta}) \cdot
                (2z)^{\ddot{\chi} (\alpha)} \,
                c_{1/(2z)} (\pr_1^* (\ddot{\Theta}_{\alpha})) \cdot {}
        } \notag \\*[\LS]
        & \hspace{3em} \smash{
                z^{\dot{\xi} (\alpha, \> \theta)} \,
                c_{1/z} (\dot{\Xi}_{\alpha, \> \theta}) \cdot
                (2z)^{\ddot{\xi} (\alpha)} \,
                c_{1/(2z)} (\pr_1^* (\ddot{\Xi}_{\alpha})) \cdot {}
        } \notag \\*[\LS]
        & \hspace{3em} \smash[t]{
                z^{\dot{\xi} (\alpha^\vee, \> \theta)} \,
                c_{1/z} \bigl( (I_\alpha \times \id_{X^\sd_\theta})^*
                (\dot{\Xi} \sss{\alpha^\vee, \> \theta}{\vee}) \bigr) \cdot
                (2z)^{\ddot{\xi} (\alpha^\vee)} \,
                c_{1/(2z)} ( \pr_1^* I \sss{\alpha}{*} (\ddot{\Xi} \sss{\alpha^\vee}{\vee}) )
            \biggr]
        \biggr\} , }
        \label{eq-lem-pf-mor-tw-mod}
    \end{align}
    where $\pr_1 \colon X_\alpha \times X^\sd_\theta \to X_\alpha$
    is the projection, and
    $I_\alpha \colon X_\alpha \to X_{\smash{\alpha^\vee}}$ 
    is the involution of $X$.
\end{lemma}

\begin{proof}
    \allowdisplaybreaks
    For convenience in notation, we write
    $\ddot{\Theta}_{\alpha}$\,, etc.,
    for their pullbacks to $X_\alpha \times X^\sd_\theta$.
    We have
    \begin{align*}
        & \smash[b]{
        (\oplus \sss{\alpha, \> \theta}{\sd})_* \circ (\upe^{z D_\alpha} \otimes \id_{W_\theta})
        \Bigl[ 
            (A \boxtimes M) \cap 
            \bigl( z^{\dot{\chi} (\alpha, \> \theta)} \, c_{1/z} (\dot{\Theta}_{\alpha, \> \theta}) \cdot
            (2z)^{\ddot{\chi} (\alpha)} \, c_{1/(2z)} (\ddot{\Theta}_{\alpha}) \bigr)
        \Bigr]
        } \\*[\LS]
        & \hspace{3em} \smash{
        {} \cap (-y)^{\ring{\xi} (\bar{\alpha} + \theta)} \,
        c_{1/y} (\ring{\Xi}_{\bar{\alpha} + \theta})
        } \\[\LS]
        = {}
        & \smash{
        (-1)^{\ring{\xi} (\bar{\alpha} + \theta)} \cdot
        (\oplus \sss{\alpha, \> \theta}{\sd})_* \biggl\{
            (\upe^{z D_\alpha} \otimes \id_{W_\theta})
            \Bigl[ 
                (A \boxtimes M)
        } \\*[\LS]
        & \hspace{3em} \smash{
            {} \cap \Bigl(
                    z^{\dot{\chi} (\alpha, \> \theta)} \,
                    c_{1/z} (\dot{\Theta}_{\alpha, \> \theta}) \cdot
                    (2z)^{\ddot{\chi} (\alpha)} \,
                    c_{1/(2z)} (\ddot{\Theta}_{\alpha})
                \Bigr)
            \Bigr]
        } \\*[\LS]
        & \hspace{1em} \smash{
            {} \cap \Bigl[
                y^{\hat{\xi} (\alpha)} \, c_{1/y} (\hat{\Xi}_\alpha) \cdot
                y^{\ring{\xi} (\theta)} \, c_{1/y} (\ring{\Xi}_\theta) \cdot
                y^{\dot{\xi} (\alpha, \> \theta)} \, c_{1/y} (\dot{\Xi}_{\alpha, \> \theta}) \cdot
                y^{\ddot{\xi} (\alpha)} \, c_{1/y} (\ddot{\Xi}_{\alpha}) \cdot {}
        } \\*[\LS]
        & \hspace{3em} \smash{
                y^{\dot{\xi} (\alpha^\vee, \> \theta)} \,
                c_{1/y} (\dot{\Xi}_{\smash{\alpha^\vee, \> \theta}}) \cdot
                y^{\ddot{\xi} (\alpha^\vee)} \,
                c_{1/y} (\ddot{\Xi}_{\smash{\alpha^\vee}})
            \Bigr]
        \biggr\}
        } \\[\LS]
        = {}
        & \smash{
        (-1)^{\ring{\xi} (\bar{\alpha} + \theta)} \cdot
        (\oplus \sss{\alpha, \> \theta}{\sd})_* \circ
        (\upe^{z D_\alpha} \otimes \id_{W_\theta})
        \biggl\{
            (A \boxtimes M)
        } \\*[\LS]
        & \hspace{1em} \smash{
            {} \cap \Bigl[
                z^{\dot{\chi} (\alpha, \> \theta)} \,
                c_{1/z} (\dot{\Theta}_{\alpha, \> \theta}) \cdot
                (2z)^{\ddot{\chi} (\alpha)} \,
                c_{1/(2z)} (\ddot{\Theta}_{\alpha}) \cdot
                y^{\hat{\xi} (\alpha)} \, c_{1/y} (\hat{\Xi}_\alpha) \cdot
                y^{\ring{\xi} (\theta)} \, c_{1/y} (\ring{\Xi}_\theta) \cdot {}
        } \\*[\LS]
        & \hspace{3em} \smash{
                (y + z)^{\dot{\xi} (\alpha, \> \theta)} \,
                c_{1/(y + z)} (\dot{\Xi}_{\alpha, \> \theta}) \cdot
                (y + 2z)^{\ddot{\xi} (\alpha)} \,
                c_{1/(y + 2z)} (\ddot{\Xi}_{\alpha}) \cdot {}
        } \\*[\LS]
        & \hspace{3em} \smash{
                (y - z)^{\dot{\xi} (\alpha^\vee, \> \theta)} \,
                c_{1/(y - z)} (\dot{\Xi}_{\smash{\alpha^\vee, \> \theta}}) \cdot
                (y - 2z)^{\ddot{\xi} (\alpha^\vee)} \,
                c_{1/(y - 2z)} (\ddot{\Xi}_{\smash{\alpha^\vee}})
            \Bigr]
        \biggr\}
        } \\[\LS]
        = {}
        & \smash{
        (-1)^{\dot{\xi} (\alpha, \> \theta) + \ddot{\xi} (\alpha)} \cdot
        (\oplus \sss{\alpha, \> \theta}{\sd})_* \circ
        (\upe^{z D_\alpha} \otimes \id_{W_\theta})
        } \\*[\LS]
        & \hspace{1em} \smash{
        \biggl\{
            \Bigl[
                \bigl( A \cap (-y)^{\hat{\xi} (\alpha)} \, c_{1/y} (\hat{\Xi}_\alpha) \bigr) \boxtimes
                \bigl( M \cap (-y)^{\ring{\xi} (\theta)} \, c_{1/y} (\ring{\Xi}_\theta) \bigr)
            \Bigr]
        } \\*[\LS]
        & \hspace{3em} \smash{
            {} \cap \biggl[ 
                z^{\dot{\chi} (\alpha, \> \theta)} \,
                c_{1/z} (\dot{\Theta}_{\alpha, \> \theta}) \cdot
                (2z)^{\ddot{\chi} (\alpha)} \,
                c_{1/(2z)} (\ddot{\Theta}_{\alpha}) \cdot {}
        } \\*[\LS]
        & \hspace{6em} \smash{
                (y + z)^{\dot{\xi} (\alpha, \> \theta)} \,
                c_{1/(y + z)} (\dot{\Xi}_{\alpha, \> \theta}) \cdot
                (y + 2z)^{\ddot{\xi} (\alpha)} \,
                c_{1/(y + 2z)} (\ddot{\Xi}_{\alpha}) \cdot {}
        } \\*[\LS]
        & \hspace{6em} \smash[t]{
                (z - y)^{\dot{\xi} (\alpha^\vee, \> \theta)} \,
                c_{1/(y - z)} (\dot{\Xi}_{\smash{\alpha^\vee, \> \theta}}) \cdot
                (2z - y)^{\ddot{\xi} (\alpha^\vee)} \,
                c_{1/(y - 2z)} (\ddot{\Xi}_{\smash{\alpha^\vee}})
            \Bigr]
        \biggr\} , }
        \numberthis
    \end{align*}
    where the first step uses Assumption~\ref{asn-h-sp-mod-mor}~\ref{itm-xi-ring},
    and the second step uses Lemma~\ref{lem-cd-ezd-chern}.
    The third step adjusts the signs using~\eqref{eq-xi-ring-rk}.

    An analogous argument as in the proof of Lemma~\ref{lem-pf-mor-va}
    shows that one can take $y = 0$ on both sides,
    using the condition in Assumption~\ref{asn-h-sp-mod-mor}~\ref{itm-xi-dots}
    that $\ring{\Xi}$ is the class of a vector bundle.
    This gives the desired equality~\eqref{eq-lem-pf-mor-tw-mod}.
\end{proof}

\begin{lemma}
    \label{lem-omega-sd-compat-y}
    For any $\alpha \in \pi_0 (X),$ $\theta \in \pi_0 (X^\sd),$
    and any $A \in V_\alpha$ and $M \in W_\theta$\,, we have
    \begin{equation}
        \Omega^\sd (Y^\sd (A, z) \, M) = Y'^\sd (\Omega (A), z) \, \Omega^\sd (M),
    \end{equation}
    where $Y^\sd, Y'^\sd$ are structure maps for the twisted module structures on $W, W',$
    and $\Omega, \Omega^\sd$ are the maps given in
    Theorems~\textnormal{\ref{thm-mor-va}} and~\textnormal{\ref{thm-mor-tw-mod}}.
\end{lemma}

\begin{proof}
    \allowdisplaybreaks
    For $A \in H_a (X_\alpha; \bbQ)$ and $M \in H_m (X^\sd_\theta; \bbQ)$,
    we have
    \begin{align*}
        & \smash[b]{
        \Omega^\sd (Y^\sd (A, z) \, M)
        } \\[\LS]
        = {}
        & \smash{
        (F \sss{\bar{\alpha} + \theta}{\sd})_* \biggl\{
            (-1)^{a \ring{\chi} (\theta)} \varepsilon^\sd (\alpha, \theta) \cdot
            (\oplus \sss{\alpha, \> \theta}{\sd})_* \circ (\upe^{z D_\alpha} \otimes \id_{W_\theta})
        } \\*[\LS]
        & \hspace{2em} \smash{
            \Bigl[ 
                (A \boxtimes M) \cap 
                \Bigl(
                    z^{\dot{\chi} (\alpha, \> \theta)} \,
                    c_{1/z} (\dot{\Theta}_{\alpha, \> \theta}) \cdot
                    (2z)^{\ddot{\chi} (\alpha)} \,
                    c_{1/(2z)} (\ddot{\Theta}_{\alpha})
                \Bigr)
            \Bigr]
            \cap c_{\ring{\xi} (\bar{\alpha} + \theta)} (\ring{\Xi} \sss{\bar{\alpha} + \theta}{\vee})
        \biggr\}
        } \\[\LS]
        = {}
        & \smash{
        (-1)^{\dot{\xi} (\alpha, \> \theta) + \ddot{\xi} (\alpha) + a \ring{\chi} (\theta)}
        \varepsilon^\sd (\alpha, \theta) \cdot
        (F \sss{\bar{\alpha} + \theta}{\sd})_* \Biggl\{
            (\oplus \sss{\alpha, \> \theta}{\sd})_* \circ (\upe^{z D_\alpha} \otimes \id_{W_\theta})
        } \\*[\LS]
        & \hspace{2em} \smash{
            \biggl\{
                \Bigl[
                    \bigl( A \cap c_{\hat{\xi} (\alpha)} (\hat{\Xi} \sss{\alpha}{\vee}) \bigr) \boxtimes
                    \bigl( M \cap c_{\ring{\xi} (\theta)} (\ring{\Xi} \sss{\theta}{\vee}) \bigr)
                \Bigr]
        } \\*[\LS]
        & \hspace{4em} \smash{
                {} \cap \biggl[ 
                    z^{\dot{\chi} (\alpha, \> \theta)} \,
                    c_{1/z} (\dot{\Theta}_{\alpha, \> \theta}) \cdot
                    (2z)^{\ddot{\chi} (\alpha)} \,
                    c_{1/(2z)} (\ddot{\Theta}_{\alpha}) \cdot {}
        } \\*[\LS]
        & \hspace{7em} \smash{
                    z^{\dot{\xi} (\alpha, \> \theta)} \,
                    c_{1/z} (\dot{\Xi}_{\alpha, \> \theta}) \cdot
                    (2z)^{\ddot{\xi} (\alpha)} \,
                    c_{1/(2z)} (\ddot{\Xi}_{\alpha}) \cdot {}
        } \\*[\LS]
        & \hspace{7em} \smash{
                    z^{\dot{\xi} (\alpha^\vee, \> \theta)} \,
                    c_{1/z} ( \dot{\Xi} \sss{\alpha^\vee, \> \theta}{\vee}) \cdot
                    (2z)^{\ddot{\xi} (\alpha^\vee)} \,
                    c_{1/(2z)} (\ddot{\Xi} \sss{\alpha^\vee}{\vee}) 
                \biggr]
            \biggr\}
        \Biggr\}
        } \\[\LS]
        = {}
        & \smash{
        (-1)^{a \ring{\chi}' (F^\sd (\theta))}
        \varepsilon'^\sd (F (\alpha), F^\sd (\theta)) \cdot
        (\oplus \sss{\alpha, \> \theta}{\prime \sd})_* \circ
        (\upe^{z D \sss{\alpha}{\prime}} \otimes \id_{W \sss{\theta}{\prime}})
        } \\*[\LS]
        & \hspace{2em} \smash{
            \biggl\{
            \Bigl[
                (F_\alpha)_* \bigl( A \cap c_{\hat{\xi} (\alpha)}
                (\hat{\Xi} \sss{\alpha}{\vee}) \bigr) \boxtimes
                (F^\sd_\theta)_* \bigl( M \cap c_{\ring{\xi} (\theta)}
                (\ring{\Xi} \sss{\theta}{\vee}) \bigr)
            \Bigr]
        } \\*[\LS]
        & \hspace{4em} \smash{
            {} \cap \Bigl(
                z^{\dot{\chi}' (F (\alpha), \> F^\sd (\theta))} \,
                c_{1/z} (\dot{\Theta} \sss{F (\alpha), \> F^\sd (\theta)}{\prime}) \cdot
                (2z)^{\ddot{\chi}' (F (\alpha))} \,
                c_{1/(2z)} (\ddot{\Theta} \sss{F (\alpha)}{\prime})
            \Bigr)
        \biggr\}
        } \\[\LS]
        = {}
        & \smash[t]{
        Y'^\sd (\Omega (A), z) \, \Omega^\sd (M), }
        \numberthis
    \end{align*}
    where the first step uses the definitions of~$\Omega^\sd$ and~$Y^\sd$,
    the second uses Lemma~\ref{lem-pf-mor-tw-mod},
    the third uses
    Assumption~\ref{asn-h-sp-mod-mor}~\ref{itm-xi-dots-compat-theta}--\ref{itm-xi-dots-compat-epsilon},
    and the last step uses the definitions of~$\Omega$, $\Omega^\sd$, and~$Y'^\sd$.
\end{proof}

\begin{bproof}[Proof of Theorem~\ref{thm-mor-tw-mod}]
    By~\eqref{eq-chi-ring-diff}, the map
    \[
        \Omega^\sd \colon W \longrightarrow W'
    \]
    is compatible with the gradings of $W$ and $W'$.
    Lemma~\ref{lem-omega-sd-compat-y} shows 
    the commutativity of the diagram~\eqref{eq-cd-omega-sd}.
\end{bproof}

\subsection[Proof of Theorem~\ref*{thm-comp-mor-tw-mod}]{Proof of Theorem~\ref{thm-comp-mor-tw-mod}}
\label{sect-proof-comp-mor-tw-mod}

First, let us verify Assumption~\ref{asn-h-sp-mod-mor}
for $F^{\dprime \sd}$, equipped with the data
$\dot{\Xi} \second, \ddot{\Xi} \second, \ring{\Xi} \second$
given in the statement of the theorem.
Assumption~\ref{asn-h-sp-mod-mor}~\ref{itm-f-invol}
follows from Theorem~\ref{thm-comp-mor-va}.
Verifying Assumption~\ref{asn-h-sp-mod-mor}~\ref{itm-fsd}--\ref{itm-xi-ring}
is straightforward.
For Assumption~\ref{asn-h-sp-mod-mor}~\ref{itm-xi-dots-compat-theta},
we have
\begin{align*}
    & \smash[b]{
    (F \sss{\alpha}{\dprime} \times F \sss{\theta}{\dprime \sd})^* 
    (\dot{\Theta} \sss{F \second (\alpha), \> F^{\dprime \sd} (\theta)}{\dprime}) 
    } \\*[\LS]
    = {} & \smash{
    (F_\alpha \times F \sss{\theta}{\sd})^* \Bigl(
        \dot{\Theta} \sss{F (\alpha), \> F^\sd (\theta)}{\prime} +
        \dot{\Xi} \sss{F (\alpha), \> F^\sd (\theta)}{\prime} +
        (I \sss{F (\alpha)}{\prime} \times \id_{X _{F^\sd (\theta)}^{\prime \sd}})^*
        (\dot{\Xi} \sss{F (\alpha)^\vee, \> F^\sd (\theta)}{\prime \vee})
    \Bigr)
    } \\*[\LS]
    = {} & \smash{
    \dot{\Theta}_{\alpha, \> \theta} +
    \dot{\Xi}_{\alpha, \> \theta} +
    (I_\alpha \times \id_{X^\sd_\theta})^*
    (\dot{\Xi} \sss{\alpha^\vee, \> \theta}{\vee}) + {}
    } \\*[\LS]
    & \smash{
    \hspace{6em}
    (F_\alpha \times F^\sd_\theta)^* \Bigl(
        \dot{\Xi} \sss{F (\alpha), \> F^\sd (\theta)}{\prime} +
        (I \sss{F (\alpha)}{\prime} \times \id_{X _{F^\sd (\theta)}^{\prime \sd}})^*
        (\dot{\Xi} \sss{F (\alpha)^\vee, \> F^\sd (\theta)}{\prime \vee})
    \Bigr)
    } \\*[\LS]
    = {} & \smash[t]{
    \dot{\Theta}_{\alpha, \> \theta} +
    \dot{\Xi} \sss{\alpha, \> \theta}{\dprime} +
    (I_\alpha \times \id_{X^\sd_\theta})^*
    (\dot{\Xi} \sss{\alpha^\vee, \> \theta}{\dprime \vee}) ,
    }
    \numberthis
\end{align*}
where $I \sss{F (\alpha)}{\prime} \colon
X \sss{F (\alpha)}{\prime} \simto X \sss{F (\alpha)^\vee}{\prime}$
denotes the involution, and
\begin{align*}
    & \smash[b]{
    (F \sss{\alpha}{\dprime})^* 
    (\ddot{\Theta} \sss{F \second (\alpha)}{\dprime}) 
    } \\*[\LS]
    = {} & \smash{
    F_\alpha^* \Bigl(
        \ddot{\Theta} \sss{F (\alpha)}{\prime} +
        \ddot{\Xi} \sss{F (\alpha)}{\prime} +
        (I \sss{F (\alpha)}{\prime})^* (\ddot{\Xi} \sss{F (\alpha)^\vee}{\prime \vee})
    \Bigr)
    } \\*[\LS]
    = {} & \smash{
    \ddot{\Theta}_{\alpha} +
    \ddot{\Xi}_{\alpha} +
    I_\alpha^* (\ddot{\Xi} \sss{\alpha^\vee}{\vee}) +
    F_\alpha^* \Bigl(
        \ddot{\Xi} \sss{F (\alpha)}{\prime} +
        (I \sss{F (\alpha)}{\prime})^*
        (\ddot{\Xi} \sss{F (\alpha)^\vee}{\prime \vee})
    \Bigr)
    } \\*[\LS]
    = {} & \smash[t]{
    \ddot{\Theta}_{\alpha} +
    \ddot{\Xi} \sss{\alpha}{\dprime} +
    I_\alpha^* (\ddot{\Xi} \sss{\alpha^\vee}{\dprime \vee}) ,
    }
    \numberthis
\end{align*}
and also
\begin{align*}
    \ring{\chi} \second (F^{\dprime \sd} (\theta))
    & = \ring{\chi}' (F'^\sd (\theta)) + 2 \ring{\xi}' (F^\sd (\theta)) \\
    & = \ring{\chi} (\theta) + 2 \ring{\xi} (\theta) + 2 \ring{\xi}' (F^\sd (\theta)) \\
    & = \ring{\chi} (\theta) + 2 \ring{\xi} \second (\theta).
    \numberthis
\end{align*}
Assumption~\ref{asn-h-sp-mod-mor}~\ref{itm-xi-dots-compat-epsilon}
follows from the relations
\begin{align}
    \dot{\xi} \second (\alpha, \theta) & =
    \dot{\xi} (\alpha, \theta) + \dot{\xi}' (F (\alpha), F^\sd (\theta)), \\
    \ddot{\xi} \second (\alpha) & =
    \ddot{\xi} (\alpha) + \ddot{\xi}' (F (\alpha)).
\end{align}

Finally, to see that
\[
    \Omega^{\dprime \sd} = \Omega'^\sd \circ \Omega^\sd,
\] 
we observe that for any $M \in \ring{H}_* (X^\sd_\theta; \bbQ)$, we have
\begin{align*}
    \Omega'^\sd \circ \Omega^\sd (M) = {} & \smash[b]{
    (F'^\sd)_* \Bigl(
        F^\sd_* (M \cap c_{\ring{\xi} (\theta)}
        (\ring{\Xi} \sss{\theta}{\vee}))
        \cap c_{\ring{\xi}' (F^\sd (\theta))}
        (\ring{\Xi} \sss{F^\sd (\theta)}{\prime \vee})
    \Bigr)
    } \\*[\LS]
    = {} & \smash{
    (F'^\sd)_* \circ F^\sd_* \Bigl[
        M \cap \Bigl( c_{\ring{\xi} (\theta)}
        (\ring{\Xi} \sss{\theta}{\vee})
        \cdot c_{\ring{\xi}' (F^\sd (\theta))}
        ((F^\sd)^* (\ring{\Xi} \sss{F^\sd (\theta)}{\prime \vee}) ) \Bigr)
    \Bigr]
    } \\*[\LS]
    = {} & \smash{
    (F^{\dprime \sd}) _* \Bigl(
        M \cap c_{\ring{\xi} (\theta) \> + \> \ring{\xi}' (F^\sd (\theta))}
        (\ring{\Xi} \sss{\theta}{\vee} +
        F^* (\ring{\Xi} \sss{F^\sd (\theta)}{\prime \vee}))  
    \Bigr)
    } \\*[\LS]
    = {} & \smash{
    \Omega^{\dprime \sd} (M),
    }
    \numberthis
\end{align*}
where in the third step,
we used the condition in Assumption~\ref{asn-h-sp-mod-mor}~\ref{itm-xi-dots}
that the classes $\ring{\Xi}_\theta$ and $\ring{\Xi} \sss{F^\sd (\theta)}{\prime}$
come from vector bundles,
so their Chern classes above their ranks are zero.
\qed

\section[Proof of Theorem~\ref*{thm-quiver-main}]{Proof of Theorem~\ref{thm-quiver-main}}

\label{sect-proof-quiver}

We use an argument similar to that of
Gross--Joyce--Tanaka~\cite[\S6]{GrossJoyceTanaka}.

Throughout this section, we fix a self-dual quiver $Q$
with no oriented cycles.
In addition to Notation~\ref{ntn-sd-quiv},
we also introduce the following notions.

\begin{definition}
    \label{def-quiv-basic}
    Let $\alpha \in C (Q)$.
    The \emph{support} of $\alpha$, denoted by $\supp \alpha$,
    is the maximal sub-quiver of $Q$
    whose vertices are those $i \in Q_0$ with $\alpha_i \neq 0$.

    We say that $\alpha$ is \emph{binary},
    if $\alpha_i \in \{ 0, 1 \}$ for all $i \in Q_0$\,.

    We say that $\alpha$ is \emph{primitive},
    if $\alpha = \delta_i$ for some $i \in Q_0$\,.
    Here, $\delta_i \in C (Q)$ denotes the class with
    $(\delta_i)_i = 1$ and $(\delta_i)_j = 0$ for all $j \neq i$.

    Likewise, for $\theta \in C^\sd (Q)$,
    its \emph{support}, denoted by $\supp \theta$,
    is the maximal sub-quiver of $Q$
    whose vertices are those $i \in Q_0$ with $\theta_i \neq 0$.
    This sub-quiver inherits a self-dual structure from that of $Q$.

    We say that $\theta$ is \emph{binary},
    if $\theta_i \in \{ 0, 1 \}$ for all $i \in Q_0$\,.

    We say that $\theta$ is \emph{primitive},
    if $\theta = \delta_{i_1} + \cdots + \delta_{i_m}$
    for some $m \geq 0$ and distinct elements $i_1, \dotsc, i_m \in Q_0^+$.
    In this case, we write $|\theta| = m$.

    Note that the class $0 \in C^\circ (Q)$ is not primitive,
    but the class $0 \in C^\sd (Q)$ is primitive.
\end{definition}

\subsection{Defining invariants}

We define the invariants $\smash{\inv{\calM^\sdss_{\theta} (\tau)}}$
by first considering \emph{increasing} stability functions,
whose semistable moduli spaces are easy to describe,
and then using the wall-crossing formula~\eqref{eq-wcf-quiver-sd}
to define invariants for all stability conditions.

\begin{definition}
    A stability function $\mu$ for $Q$ is \emph{increasing},
    if for any edge $i \to j$ in $Q$,
    we have $\mu (i) < \mu (j)$.
    
    Note that since it is assumed that $Q$ has no oriented loops,
    there always exists an increasing self-dual stability function for $Q$.
\end{definition}

\begin{lemma}
    \label{lem-increasing-msdss}
    Let $\mu$ be an increasing self-dual stability function for $Q$.
    Let $\alpha \in C (Q)$ and $\theta \in C^\sd (Q)$.
    \begin{enumerate}
        \item
            \label{itm-increasing-mplss}
            If $\alpha$ is primitive, then
            \[
                \calM^\plss_\alpha (\mu) \simeq *
            \]
            is a single point.
            Otherwise, we have either
            $\calM^\plss_\alpha (\mu) \simeq \varnothing,$
            or $\dim \calM^\plss_\alpha (\mu) < 0$.

        \item 
            \label{itm-increasing-msdss}
            If $\theta$ is primitive, then
            \[
                \calM^\sdss_{\theta} (\mu) \simeq
                [*/\bbZ \sss{2}{|\theta|}].
            \]
            Otherwise, we have either
            $\calM^\sdss_{\theta} (\mu) \simeq \varnothing,$
            or
            $\dim \calM^\sdss_{\theta} (\mu) < 0$.
    \end{enumerate}
\end{lemma}

\begin{proof}
    For~\ref{itm-increasing-mplss},
    see Gross--Joyce--Tanaka~\cite[Proposition~5.6]{GrossJoyceTanaka}.

    For~\ref{itm-increasing-msdss}, 
    let us first suppose that $\theta$ is primitive.
    Then $\supp \theta$ has no edges,
    since if there exist two vertices $i, j \in Q_0^+$ and an edge $a \in Q_1$ from $i$ to $j$,
    then the edge $a^\vee$ would go from $j$ to $i$,
    forming an oriented cycle with $a$, which is impossible.
    Therefore, there is a unique self-dual object of class $\theta$,
    with automorphism group $\upO (1)^{|\theta|} \simeq \bbZ \sss{2}{|\theta|}$.

    If $\theta$ is not primitive,
    let $i \in (\supp \theta)_0$ be a vertex with $\mu (i)$ maximal.
    In particular, as $\mu$ is increasing,
    this implies that there are no outgoing arrows from $i$ in $\supp \theta$.
    If $\mu (i) > 0$,
    then every self-dual object of class $\theta$
    will be destabilized by a non-zero subobject supported on the vertex $i$,
    and hence $\calM^\sdss_\theta (\mu) \simeq \varnothing$.
    Otherwise, we have $\mu (i) = 0$
    for all $i \in (\supp \theta)_0$\,.
    As in the previous case,
    we see that $\supp \theta$ has no edges,
    and there is at most one self-dual object of class $\theta$, up to an isomorphism.
    Therefore, $\calM^\sdss_\theta (\mu)$ is either empty, or a product of stacks of the form
    $[* / \GL (n)]$,
    $[* / \upO (n)]$, or $[* / \Sp (2n)]$,
    for $n > 0$.
    We see that $\dim \calM^\sdss_\theta (\mu) \geq 0$ precisely when all these factors
    are $[* / \upO (1)]$,
    i.e.\ when $\theta_i = 1$ and $i \in Q_0^+$ for all $i \in (\supp \theta)_0$\,.
\end{proof}

In particular, for an increasing stability function~$\mu$, we have
\begin{equation}
    \label{eq-quiv-incr-inv-pl}
    \inv{ \calM^\plss_{\alpha} (\mu) } =
    \begin{cases}
        1^\pl_\alpha,
        & \alpha \text{ primitive,} \\
        0,
        & \text{otherwise},
    \end{cases}
\end{equation}
as in Gross--Joyce--Tanaka~\cite[Equation~(5.8)]{GrossJoyceTanaka},
where $1^\pl_\alpha \in H_0 (\calM^\pl_\alpha; \bbQ)$
is the class of a point.

\begin{theorem}
    \label{thm-quiver-uniqueness}
    There is a unique way to define classes
    $\inv{ \calM^\sdss_{\theta} (\tau) }$
    for all self-dual stability conditions $\tau$ for $Q$
    and all classes $\theta \in C^\sd (Q),$
    satisfying Theorem~\textnormal{\ref{thm-quiver-main}~\ref{itm-quiver-wcf},}
    together with the following condition:

    \begin{enumerate}
        \item 
            For an increasing stability function $\mu$
            and a class $\theta \in C^\sd (Q)$,
            \begin{equation}
                \label{eq-quiv-incr-inv-sd}
                \inv{ \calM^\sdss_{\theta} (\mu) } =
                \begin{cases}
                    2^{-|\theta|} \cdot 1^\sd_\theta,
                    & \theta \text{ primitive,} \\
                    0,
                    & \text{otherwise,}
                \end{cases}
            \end{equation}
            where $1^\sd_\theta \in H_0 (\calM^\sd_\theta; \bbQ)$
            is the class of a point.
    \end{enumerate}
    In particular, the uniqueness part of
    Theorem~\textnormal{\ref{thm-quiver-main}} holds.
\end{theorem}

\begin{proof}
    Let $\mu$ be an increasing stability function.
    Then all the classes $\inv{ \calM^\ss_{\alpha} (\mu) }$
    and $\inv{ \calM^\sdss_{\theta} (\mu) }$
    are already defined.
    For any other stability condition $\tau$,
    the condition in Theorem~\ref{thm-quiver-main}~\ref{itm-quiver-wcf}
    forces $\inv{ \calM^\sdss_{\theta} (\tau) }$
    to be given by~\eqref{eq-wcf-quiver-sd},
    with $\mu, \tau$ in place of $\tau, \tilde{\tau}$ there.
    This proves the uniqueness part.

    Explicitly, substituting~%
    \eqref{eq-quiv-incr-inv-pl}--\eqref{eq-quiv-incr-inv-sd}
    into~\eqref{eq-wcf-quiver-sd}, we obtain
    \begin{align*}
        \hspace{2em} & \hspace{-2em}
        \inv{ \calM^\sdss_\theta (\tau) } = {}
        \\[.5ex] &
        \sum_{ \leftsubstack[8em]{
            \\[-3ex]
            & n \geq 0; \, m_1, \dotsc, m_n > 0; \\[-.5ex]
            & \alpha_{1,1}, \dotsc, \alpha_{1,m_1}; \dotsc;
            \alpha_{n,1}, \dotsc, \alpha_{n,m_n} \in C (Q)\text{ primitive}; \\[-.5ex]
            & \rho \in \smash{C^\sd (Q)} \text{ primitive} \colon \\[-.5ex]
            & \theta = (\bar{\alpha}_{1,1} + \cdots + \bar{\alpha}_{1,m_1})
            + \cdots + (\bar{\alpha}_{n,1} + \cdots + \bar{\alpha}_{n,m_n}) + \rho
        } } {} 
        2^{-|\rho|} \cdot
        \tilde{U}^\sd (\alpha_{1,1}, \dotsc, \alpha_{1,m_1}; \dotsc;
        \alpha_{n,1}, \dotsc, \alpha_{n,m_n}; \mu, \tau) \cdot {} 
        \\[.5ex]
        & \hspace{2em}
        [ [ [ 
        1^\pl_{\alpha_{1,1}} \, , 
        1^\pl_{\alpha_{1,2}} ] , \dotsc ] ,
        1^\pl_{\alpha_{1,m_1}} ] \heart \cdots \heart
        [ [ [ 
        1^\pl_{\alpha_{n,1}} \, , 
        1^\pl_{\alpha_{n,2}} ] , \dotsc ] ,
        1^\pl_{\alpha_{n,m_n}} ] \heart 
        1^\sd_{\rho} \, .
        \numberthis
        \label{eq-def-quiver-inv}
    \end{align*}

    Next, we verify that~\eqref{eq-wcf-quiver-sd}
    remains true for any self-dual stability conditions $\tau, \tilde{\tau}$,
    if we define the invariants by~\eqref{eq-def-quiver-inv}.
    For this purpose, we rewrite~\eqref{eq-def-quiver-inv}
    using~\eqref{eq-wcf-diamond}, to obtain
    \begin{align}
        \label{eq-def-quiver-inv-diamond}
        & \inv{ \calM^\sdss_\theta (\tau) } =
        \sum_{ \leftsubstack[6em]{
            \\[-3ex]
            & n \geq 0; \, m_1, \dotsc, m_n > 0; \\[-.5ex]
            & \alpha_1, \dotsc, \alpha_n \in C (Q)\text{ primitive}; \\[-.5ex]
            & \rho \in \smash{C^\sd (Q)} \text{ primitive} \colon \\[-.5ex]
            & \theta = \bar{\alpha}_1 + \cdots + \bar{\alpha}_n + \rho
        } } {} 
        2^{-|\rho|} \cdot
        U^\sd (\alpha_1, \dotsc, \alpha_n; \mu, \tau) \cdot
        1^\pl_{\alpha_1} \diamond \cdots \diamond 1^\pl_{\alpha_n}
        \diamond 1^\sd_{\rho}
        \raisetag{2ex}
    \end{align}
    in the space $U^\tw$ in Definition~\ref{def-moduli-utw}.
    Now, we may apply~\cite[Lemma~11.8]{Bu2023}
    with $\mu, \tau, \tilde{\tau}$ in place of
    $\tau, \hat{\tau}, \tilde{\tau}$ there.
    This shows that if we expand the
    right-hand side of~\eqref{eq-wcf-diamond}
    using~\eqref{eq-def-quiver-inv-diamond},
    the result is the same as~\eqref{eq-def-quiver-inv-diamond}
    with $\tilde{\tau}$ in place of $\tau$.
    In other words, invariants defined by~\eqref{eq-def-quiver-inv},
    then the wall-crossing formula~\eqref{eq-wcf-quiver-sd}
    is always satisfied.

    Finally, we need to show that~\eqref{eq-quiv-incr-inv-sd}
    remains true for any other increasing
    self-dual stability function~$\tilde{\mu}$.
    We follow an analogous argument as in
    Gross--Joyce--Tanaka~\cite[Proposition~6.1]{GrossJoyceTanaka},
    and fix a class~$\theta$.
    Suppose that \eqref{eq-quiv-incr-inv-sd} fails for $\tilde{\mu}$.
    For $t \in [0, 1]$, define an $\bbR$-valued
    increasing self-dual stability function
    $\mu^t = t \, \tilde{\mu} + (1 - t) \, \mu$.
    Since every stability function \emph{dominates}
    any other stability function in a small neighbourhood,
    in the sense of~\cite[Definition~3.5]{GrossJoyceTanaka},
    we can find $t', t \second \in [0, 1]$ with
    $\mu' = \mu^{\smash{t'}}$ and $\mu \second = \mu^{\smash{t \second}}$,
    such that $\mu'$ either dominates or is dominated by $\mu \second$,
    and \eqref{eq-quiv-incr-inv-sd} holds for $\mu'$ but fails for $\mu \second$.
    Now, consider \eqref{eq-def-quiver-inv-diamond}
    with $\mu', \mu \second$ in place of $\mu, \tau$.
    If $\mu'$ dominates $\mu \second$, then any non-zero term
    on the right-hand side of \eqref{eq-def-quiver-inv-diamond}
    would have $\mu' (\alpha_1) = \cdots = \mu' (\alpha_n) = 0$,
    so there are no edges joining the supports of $\alpha_1, \dotsc, \alpha_n, \rho$.
    If $\mu \second$ dominates $\mu'$, a similar statement is true.
    But in this case,
    any term on the right-hand side of~\eqref{eq-def-quiver-inv}
    involving the classes $\alpha_1, \dotsc, \alpha_n, \rho$ is zero
    unless $n = 0$,
    by the definition of the operation~$\heart$.
    In other words, if we consider
    \eqref{eq-def-quiver-inv} with $\mu', \mu \second$ in place of $\mu, \tau$,
    then the only possible non-zero term on the right-hand side
    is the term with $n = 0$ and $\rho = \theta$,
    which is only present when $\theta$ is primitive.
    This establishes~\eqref{eq-quiv-incr-inv-sd} for $\mu \second$,
    which contradicts with the assumption that it fails.
    This completes the proof.
\end{proof}

\subsection{Invariants and morphisms of quivers}

Recall the slightly unusual notion of a \emph{morphism of quivers}
from Definition~\ref{def-mor-quiv}.
Gross--Joyce--Tanaka~\cite[Theorem~5.11]{GrossJoyceTanaka}
showed that such morphisms of quivers
interact well with the invariants,
and we state their result as follows.

\begin{theorem}
    \label{thm-omega-inv}
    Let $\lambda \colon Q \to Q'$
    be a morphism of quivers,
    where $Q, Q'$ have no oriented cycles.
    Let $\tau'$ be a weak stability condition for $Q',$
    such that $\tau = \lambda^* (\tau')$
    is a weak stability condition for $Q$.

    Then for all $\alpha \in C (Q),$ we have
    \begin{equation}
        \label{eq-omega-inv}
        \prod_{i \in Q_0} \alpha_i ! \cdot
        \Omega^\pl_\lambda \bigl( \inv{ \calM^\ss_{\alpha} (\tau) } \bigr) 
        =
        \prod_{i' \in Q'_0} \alpha \sss{i'}{\prime} ! \cdot
        \inv{ \calM'^\ss_{\smash{\alpha'}} (\tau') } \, ,
    \end{equation}
    where $\alpha' = \lambda_* (\alpha),$
    and $\Omega^\pl_\lambda$ is the map defined in
    Construction~\textnormal{\ref{cons-quiv-omega}}.
\end{theorem}

We prove the following analogous result,
which relates the self-dual invariants
when we have a morphism of self-dual quivers,
as in Definition~\ref{def-mor-quiv-sd}.

\begin{theorem}
    \label{thm-omega-sd-inv}
    Let $\lambda \colon Q \to Q'$
    be a morphism of self-dual quivers,
    where $Q, Q'$ have no oriented cycles.
    Let $\tau'$ be a self-dual weak stability condition for $Q',$
    such that $\tau = \lambda^* (\tau')$
    is a self-dual weak stability condition for $Q$.

    Then for all $\theta \in C^\sd (Q),$ we have
    \begin{multline}
        \label{eq-omega-sd-inv}
        \prod_{i \in Q_0^\tria} \theta_i ! \cdot
        \prod_{i \in Q_0^\pm} 
        2^{\tilde{\theta}_i} \, \tilde{\theta}_i ! \cdot
        \Omega^\sd_\lambda \bigl( \inv{ \calM^\sdss_{\theta} (\tau) } \bigr)
        = \\[-1ex]
        \prod_{i' \in Q'^\tria_0} \theta'_{\smash{i'}} ! \cdot
        \prod_{i' \in Q'^\pm_0}
        2^{ \tilde{\theta}'_{\smash{i'}} } \, \tilde{\theta}'_{\smash{i'}} ! \cdot
        \inv{ \calM'^\sdss_{\theta \smash{'}} (\tau') } \, ,
    \end{multline}
    where $\theta' = \lambda_* (\theta),$ and
    \begin{align*}
        \tilde{\theta}_i & =
        \Bigl\lfloor \frac{\theta_i}{2} \Bigr\rfloor ,
        \\
        \tilde{\theta}'_{\smash{i'}} & = \frac{1}{2} \Bigl(
            \theta \sss{i'}{\prime}
            - \# \bigl\{ i \in Q_0^\pm \bigm|
            \lambda (i) = i', \ \theta_i \textnormal{ odd} \bigr\}
        \Bigr)
    \end{align*}
    for $i \in Q^\pm_0$ and $i' \in Q'^\pm_0$.
    Note that $\tilde{\theta}_i$ and $\tilde{\theta}'_{\smash{i'}}$ 
    are always non-negative integers.

    In particular, if $\lambda$ induces a bijection on vertices, then
    \begin{equation}
        \Omega^\sd_\lambda \bigl( \inv{ \calM^\sdss_{\theta} (\tau) } \bigr)
        = \inv{ \calM'^\sdss_{\theta \smash{'}} (\tau') } \, .
    \end{equation}
\end{theorem}

\begin{proof}
    \allowdisplaybreaks
    We proceed in three steps, where we prove the result in special cases.

    \begin{enumerate}
        \item
            \label{itm-omega-sd-case-i}
            Both $\tau$ and $\tau'$ are given by increasing stability functions.
    \end{enumerate}
    
    In this case, by Lemma~\ref{lem-increasing-msdss},
    if $\theta$ is not primitive,
    then both sides of~\eqref{eq-omega-sd-inv} are zero.
    Therefore, we can assume that $\theta$ is primitive.
    We split into two cases.
    
    If $\theta' = \lambda_* (\theta)$ is binary,
    then it is also primitive,
    and $\ring{\xi} (\theta) = 0$,
    so that $\ring{\Xi}_\theta = 0$,
    and $\Omega^\sd_\lambda$ is just a pushforward map on $H_* (\calM^\sd_{\theta}; \bbQ)$.
    We thus have
    \begin{equation}
        \Omega^\sd_\lambda \bigl( \inv{ \calM^\sdss_{\theta} (\tau) } \bigr) =
        \Omega^\sd_\lambda \bigl( 2^{-|\theta|} \cdot 1^\sd_\theta \bigr) =
        2^{-|\theta'|} \cdot 1^\sd_{\theta'} =
        \inv{ \calM'^\sdss_{\theta \smash{'}} (\tau') } \, ,
    \end{equation}
    where $1^\sd_\theta \in H_0 (\calM^\sd_\theta; \bbQ)$
    and $1^\sd_{\theta'} \in H_0 (\calM'^\sd_{\theta'}; \bbQ)$
    are classes of a point, 
    and we used that $|\theta| = |\theta'|$.
    This proves~\eqref{eq-omega-sd-inv}.
    
    If $\theta'$ is not binary, then one necessarily has
    $\ring{\xi} (\theta) > 0$ and $\dim \calM'^\sdss_{\theta \smash{'}} (\tau') < 0$,
    as $Q'$ cannot have any edge between two vertices in $Q'^+_0$\,.
    This shows that both sides of~\eqref{eq-omega-sd-inv} are zero.
    
    \begin{enumerate}[resume]
        \item
            \label{itm-omega-sd-case-ii}
            There exists an increasing self-dual stability function
            $\mu'$ for~$Q'$,
            such that $\mu = \lambda^* (\mu)$ is an increasing stability function
            for~$Q$.
    \end{enumerate}
    
    In this case, in the space~$U^\tw \bigl( \check{H}_* (\calM'^\pl; \bbQ); \,
    \ring{H}_* (\calM'^\sd; \bbQ) \bigr)$,
    as in Definition~\ref{def-moduli-utw}, we have
    \begin{align*}
        & \Omega^\sd_\lambda \bigl( \inv{ \calM^\sdss_{\theta} (\tau) } \bigr)
        \\*[2ex]
        = {} & 
        \sum_{ \leftsubstack[6em]{
            & \alpha_1, \dotsc, \alpha_n \in C (Q), \,
            \rho \in C^\sd (Q) \colon \\[-.5ex]
            & \theta = \bar{\alpha}_1 + \cdots + \bar{\alpha}_n + \rho 
        } } {}
        U^\sd (\alpha_1, \dotsc, \alpha_n; \mu, \tau) \cdot {}
        \\*[1ex]
        & \hspace{3em}
        \Omega^\pl_\lambda \bigl( \inv{ \calM^\ss_{\alpha_1} (\mu) } \bigr) \diamond
        \cdots \diamond
        \Omega^\pl_\lambda \bigl( \inv{ \calM^\ss_{\alpha_n} (\mu) } \bigr) \diamond
        \Omega^\sd_\lambda \bigl( \inv{ \calM^\sdss_{\rho} (\mu) } \bigr)
        \\[2ex]
        = {}
        &
        \sum_{ \leftsubstack[6em]{
            & \alpha_1, \dotsc, \alpha_n \in C (Q), \,
            \rho \in C^\sd (Q) \colon \\[-.5ex]
            & \theta = \bar{\alpha}_1 + \cdots + \bar{\alpha}_n + \rho . \\[-.5ex]
            & \text{Write } \alpha \sss{j}{\prime} = \lambda_* (\alpha_j), \,
            \rho' = \lambda_* (\rho)
        } } {}
        U^\sd (\alpha_1, \dotsc, \alpha_n; \mu, \tau) \cdot {}
        \\*
        & \hspace{3em}
        \prod_{j=1}^{n} {} \left(
            \frac{ 
                \prod \limits _{i' \in Q'_0} (\alpha'_j)_{i'} !
            }{ 
                \prod \limits _{i \in Q_0} (\alpha_j)_i !
                \vphantom{2^{\tilde{\rho}_i}}
            } 
        \right) \cdot
        \frac{ 
            \prod \limits _{i' \in Q'^\tria_0} \rho \sss{i'}{\prime} ! \cdot
            \prod \limits _{i' \in Q'^\pm_0} 2^{\tilde{\rho} \sss{i'}{\prime}}
            \tilde{\rho} \sss{i'}{\prime} !
        }{ 
            \prod \limits _{i \in Q^\tria_0} \rho_i ! \cdot
            \prod \limits _{i \in Q^\pm_0} 2^{\tilde{\rho}_i} \tilde{\rho}_i !
        } \cdot {}
        \\*[1ex]
        & \hspace{3em}
        \inv{ \calM'^\ss_{\alpha \sss{1}{\prime}} (\mu') } \diamond \cdots \diamond
        \inv{ \calM'^\ss_{\alpha \sss{n}{\prime}} (\mu') } \diamond
        \inv{ \calM'^\sdss_{\rho'} (\mu') }
        \\[2ex]
        = {} &
        \sum_{ \leftsubstack[6em]{
            & \alpha \sss{1}{\prime}, \dotsc, \alpha \sss{n}{\prime} \in C (Q'), \,
            \rho' \in C^\sd (Q') \colon \\[-.5ex]
            & \theta' = \bar{\alpha} \sss{1}{\prime} + \cdots +
            \bar{\alpha} \sss{n}{\prime} + \rho'
        } } {}
        U^\sd (\alpha'_1, \dotsc, \alpha'_n; \mu', \tau') \cdot {}
        \\*
        & \hspace{1em}
        \sum_{ \leftsubstack[6em]{
            & \alpha_1, \dotsc, \alpha_n \in C (Q), \,
            \rho \in C^\sd (Q) \colon \\[-.5ex]
            & \theta = \bar{\alpha}_1 + \cdots + \bar{\alpha}_n + \rho, \\[-.5ex]
            & \alpha \sss{j}{\prime} = \lambda_* (\alpha_j), \,
            \rho' = \lambda_* (\rho)
        } } \hspace{4em} \left[
        \prod_{j=1}^{n} {} \left(
            \frac{ 
                \prod \limits _{i' \in Q'_0} (\alpha'_j)_{i'} !
            }{ 
                \prod \limits _{i \in Q_0} (\alpha_j)_i !
                \vphantom{2^{\tilde{\rho}_i}}
            } 
        \right) \cdot
        \frac{ 
            \prod \limits _{i' \in Q'^\tria_0} \rho \sss{i'}{\prime} ! \cdot
            \prod \limits _{i' \in Q'^\pm_0} 2^{\tilde{\rho} \sss{i'}{\prime}}
            \tilde{\rho} \sss{i'}{\prime} !
        }{ 
            \prod \limits _{i \in Q^\tria_0} \rho_i ! \cdot
            \prod \limits _{i \in Q^\pm_0} 2^{\tilde{\rho}_i} \tilde{\rho}_i !
        } \right] \cdot {}
        \\*[1ex]
        & \hspace{3em}
        \inv{ \calM'^\ss_{\alpha \sss{1}{\prime}} (\mu') } \diamond \cdots \diamond
        \inv{ \calM'^\ss_{\alpha \sss{n}{\prime}} (\mu') } \diamond
        \inv{ \calM'^\sdss_{\rho'} (\mu') }
        \\[2ex]
        = {} &
        \sum_{ \leftsubstack[6em]{
            & \alpha \sss{1}{\prime}, \dotsc, \alpha \sss{n}{\prime} \in C (Q'), \,
            \rho' \in C^\sd (Q') \colon \\[-.5ex]
            & \theta' = \bar{\alpha} \sss{1}{\prime} + \cdots +
            \bar{\alpha} \sss{n}{\prime} + \rho'
        } } {}
        U^\sd (\alpha'_1, \dotsc, \alpha'_n; \mu', \tau') \cdot
        \frac{ 
            \prod \limits _{i' \in Q'^\tria_0} \theta \sss{i'}{\prime} ! \cdot
            \prod \limits _{i' \in Q'^\pm_0} 2^{\tilde{\theta} \sss{i'}{\prime}}
            \tilde{\theta} \sss{i'}{\prime} !
        }{ 
            \prod \limits _{i \in Q^\tria_0} \theta_i ! \cdot
            \prod \limits _{i \in Q^\pm_0} 2^{\tilde{\theta}_i} \tilde{\theta}_i !
        } \cdot {}
        \\*[1ex]
        & \hspace{3em}
        \inv{ \calM'^\ss_{\alpha \sss{1}{\prime}} (\mu') } \diamond \cdots \diamond
        \inv{ \calM'^\ss_{\alpha \sss{n}{\prime}} (\mu') } \diamond
        \inv{ \calM'^\sdss_{\rho'} (\mu') } \ ,
        \numberthis
        \label{eq-pf-omega-sd-ii}
    \end{align*}
    where $\tilde{\rho}_i$ and $\tilde{\rho} \sss{i'}{\prime}$
    are defined analogously to
    $\tilde{\theta}_i$ and $\tilde{\theta} \sss{i'}{\prime}$\,,
    but with $\rho$ in place of $\theta$.
    The first step uses~\eqref{eq-wcf-diamond},
    and the diagram~\eqref{eq-cd-omega-pl-sd}
    in Theorem~\ref{thm-mor-tw-mod}.
    The second step uses step~\ref{itm-omega-sd-case-i}
    and Theorem~\ref{thm-omega-inv}.
    In the third step, we rearrange the sums, and we use that
    $U^\sd (\alpha_1, \dotsc, \alpha_n; \mu, \tau) =
    U^\sd (\alpha'_1, \dotsc, \alpha'_n; \mu', \tau')$.
    The final step involves a combinatorial identity,
    which we explain as follows.

    Fix $\theta' = \bar{\alpha}'_1 + \cdots + \bar{\alpha}'_n + \rho$
    and $\theta$ with $\lambda_* (\theta) = \theta'$,
    as in the equation.
    For each $i \in Q_0$\,, let $T_i$ be a set with $|T_i| = \theta_i$\,,
    and let $T = \coprod_i T_i$\,.
    For $i' \in Q'_0$\,, define
    $T \sss{i'}{\prime} = \coprod_{\smash{i \in \lambda^{-1} (i')}} T_i$\,,
    so that $|T \sss{i'}{\prime}| = \theta \sss{i'}{\prime}$\,.
    Choose a $\bbZ_2$-action on $T$, given by a map $(-)^\vee \colon T \to T$, 
    such that $T_i^\vee = T_{i^\vee}$ for all $i \in Q_0$\,,
    and for any $i \in Q^\pm_0$,
    the map $(-)^\vee$ fixes at most one element of $T_i$\,.
    Let $T_i^\circ \subset T_i$ be the $\bbZ_2$-fixed subset.
    The fraction in the last step of~\eqref{eq-pf-omega-sd-ii}
    is equal to the number of ways to choose a decomposition
    \begin{equation}
        \label{eq-pf-omega-sd-decomp-t}
        T = \coprod_{i \in Q_0} S_i
    \end{equation}
    with $S_i \subset T \sss{\lambda (i)}{\prime}$\,,
    $S_i^\vee = S_{i^\vee}$\,,
    and $|S_i| = \theta_i$ for all $i \in Q_0$\,,
    and with $T_i^\circ \subset S_i$ for all $i \in Q^\pm_0$.

    On the other hand, in the second last step of~\eqref{eq-pf-omega-sd-ii},
    the coefficient counts the same number of decompositions by 
    fixing a different decomposition
    \[
        T = A_1 \sqcup A_1^\vee \sqcup \cdots \sqcup A_n \sqcup A_n^\vee \sqcup B,
    \]
    with $|A_j \cap T \sss{i'}{\prime}| = (\alpha'_j)_{i'}$ and
    $|B \cap T \sss{i'}{\prime}| = \rho \sss{i'}{\prime}$ for all $i' \in Q'_0$\,,
    and going through all possibilities of
    $(\alpha_j)_i = |A_j \cap S_i|$
    and $\rho_i = |B \cap S_i|$,
    where $i \in Q_0$\,, and the $S_i$ are subsets
    in the decomposition~\eqref{eq-pf-omega-sd-decomp-t}.
    We show that the number of decompositions~\eqref{eq-pf-omega-sd-decomp-t}
    corresponding to some fixed $\alpha_j$ and $\rho$
    is equal to the inner summand
    in the second last step of~\eqref{eq-pf-omega-sd-ii}.
    Indeed, each decomposition~\eqref{eq-pf-omega-sd-decomp-t}
    corresponds to a way to choose decompositions
    \[
        A_j = \coprod_{i \in Q_0} A_{j, \> i} \, , \qquad
        B = \coprod_{i \in Q_0} B_i \, ,
    \]
    with $A_{j, \> i} \subset T \sss{\lambda (i)}{\prime}$\,,
    $|A_{j, \> i}| = (\alpha_j)_i$\,,
    $B_i \subset T \sss{\lambda (i)}{\prime}$\,,
    $B_i^\vee = B_{i^\vee}$\,,
    and $|B_i| = \rho_i$ for all $i \in Q_0$\,,
    and with $T_i^\circ \subset B_i$ for all $i \in Q^\pm_0$.
    The number of such choices can be seen to be
    equal to the inner summand
    in the second last step of~\eqref{eq-pf-omega-sd-ii}.

    Finally, applying~\eqref{eq-wcf-diamond} for $Q'$,
    with $\mu', \tau'$ in place of $\tau, \tilde{\tau}$ there,
    we deduce~\eqref{eq-omega-sd-inv} from~\eqref{eq-pf-omega-sd-ii}.
    
    \begin{enumerate}[resume]
        \item
            \label{itm-omega-sd-case-iii}
            The general case.
    \end{enumerate}

    Define an auxiliary self-dual quiver $Q \second$ from $Q$
    by removing all edges that are not in $Q_1^\circ$,
    where $Q_1^\circ \subset Q_1$ is the subset used to define $\lambda$,
    as in Definition~\ref{def-mor-quiv}.
    Let $\lambda' \colon Q \to Q \second$ be the morphism
    removing the edges,
    and let $\lambda \second \colon Q \second \to Q'$
    be the morphism induced from $\lambda$, so that
    $\lambda = \lambda \second \circ \lambda'$.
    We claim that both $\lambda'$ and $\lambda \second$
    satisfy the condition in step~\ref{itm-omega-sd-case-ii}.
    Indeed, every increasing self-dual stability function on $Q$
    comes from an increasing self-dual stability function on $Q \second$,
    and every increasing self-dual stability function on $Q'$
    induces one on $Q \second$.

    Now, by step~\ref{itm-omega-sd-case-ii}, we obtain
    \begin{align*}
        &
        \prod_{i \in Q_0^\tria} \theta_i ! \cdot
        \prod_{i \in Q_0^\pm} 
        2^{\tilde{\theta}_i} \, \tilde{\theta}_i ! \cdot
        \Omega^\sd_\lambda \bigl( \inv{ \calM^\sdss_{\theta} (\tau) } \bigr)
        \\*
        = &
        \prod_{i \in Q_0^\tria} \theta_i ! \cdot
        \prod_{i \in Q_0^\pm} 
        2^{\tilde{\theta}_i} \, \tilde{\theta}_i ! \cdot
        \Omega^\sd_{\lambda \smash{\second}}
        \bigl( \inv{ \calM^{\dprime \sd}_{\theta} (\tau \second) } \bigr)
        \\*
        = &
        \prod_{i' \in Q'^\tria_0} \theta'_{\smash{i'}} ! \cdot
        \prod_{i' \in Q'^\pm_0}
        2^{ \tilde{\theta}'_{\smash{i'}} } \, \tilde{\theta}'_{\smash{i'}} ! \cdot
        \inv{ \calM'^\sdss_{\theta \smash{'}} (\tau') } \, ,
        \numberthis
    \end{align*}
    where $\tau \second = (\lambda \second)^* (\tau)$,
    and we used Theorem~\ref{thm-omega-comp}.
\end{proof}

\subsection{Fundamental classes for binary classes}

We prove Theorem~\ref{thm-quiver-main}~\ref{itm-quiver-fund},
the compatibility between the invariants
and the fundamental classes of moduli stacks when the latter exist,
in the special case when $\theta$ is a \emph{binary} class,
in the sense of Definition~\ref{def-quiv-basic}.

\begin{lemma}
    \label{lem-tree-neg-dim}
    Let $\theta \in C^\sd (Q)$ be a binary class,
    and assume that all components of $\supp \theta$
    are self-dual and are trees. Then,
    \[
        \dim \calM^\sd_\theta \leq 0.
    \]
\end{lemma}

\begin{proof}
    We may assume that $\supp \theta$ is connected,
    since in any case, $\calM^\sd_\theta$ is always a product
    of such moduli stacks for classes~$\theta$ with connected support.

    Let $Q' = \supp \theta$.
    By~\eqref{eq-dim-m-quiv}, we have
    \begin{align*}
        \dim_{\bbC} \calM^\sd_\theta
        & = - \frac{1}{2} \ring{\chi}_Q (\theta)
        = |Q'^+_1| + |Q'^\tria_1| - |Q'^\tria_0|
        \\ & =
        \frac{1}{2} \bigl( |Q'^+_0| + |Q'^+_1| - |Q'^-_1| - 1 \bigr),
        \numberthis
    \end{align*}
    where we used that $|Q'_0| = |Q'_1| + 1$ since $Q'$ is a tree,
    and that $Q'^-_0 = \varnothing$ for dimension reasons.

    However, we have $|Q'^+_0| \leq 1$, since if it had two distinct elements,
    then any path connecting them would form a cycle with its dual path,
    contradicting the fact that $Q'$ is a tree.
    A similar argument shows that $|Q'^\pm_1| \leq 1$.
    This means that $\dim \calM^\sd_\theta \leq 0$.
\end{proof}

\begin{lemma}
    \label{lem-quiv-comp-sd}
    Let $\theta \in C^\sd (Q)$.
    If $\supp \theta$ has a connected component that is not self-dual,
    then there are no $\tau$-stable self-dual representations of class~$\theta$
    for any self-dual stability condition~$\tau$.
\end{lemma}

\begin{proof}
    Let $\alpha \in C (Q)$ be the restriction of~$\theta$
    to one of the components of $\supp \theta$,
    such that $\alpha^\vee \neq \alpha$.
    Given a self-dual stability condition~$\tau$,
    we may assume that $\tau (\alpha) \geq 0$,
    since otherwise, we may replace $\alpha$ by $\alpha^\vee$.
    Then, for any self-dual representation~$E$ of class~$\theta$,
    we have the restriction $E_\alpha$ of $E$ to~$\supp \alpha$,
    which is a representation of class~$\alpha$.
    Then $E_\alpha \subset E$ is an isotropic subobject
    with $\tau (E_\alpha) = \tau (\alpha) \geq 0$.
\end{proof}

\begin{lemma}
    \label{lem-quiv-va-on-units}
    \allowdisplaybreaks
    Let $\alpha, \beta \in C (Q)$ be binary classes,
    such that $\supp \alpha$ and $\supp \beta$ are trees,
    $\alpha + \beta$ is binary,
    and all connected components of $\supp (\alpha + \beta)$ are trees.
    
    Let $1^\pl_\alpha \in H_0 (\calM^\pl_\alpha; \bbQ)$
    and $1^\pl_\beta \in H_0 (\calM^\pl_\beta; \bbQ)$
    be classes of a point. Then,
    \begin{equation}
        \label{eq-quiv-binary-lie}
        [1^\pl_\alpha \, , 1^\pl_\beta]
        = \left\{ \, \begin{alignedat}[c]{3}
            & {-1^\pl_{\alpha+\beta}} \, , \quad
            && \chi_Q (\alpha, \beta) = -1, \ 
            && \chi_Q (\beta, \alpha) = 0, \\
            & 1^\pl_{\alpha+\beta} \, , \quad
            && \chi_Q (\alpha, \beta) = 0, \ 
            && \chi_Q (\beta, \alpha) = -1, \\
            & 0, \quad
            && \chi_Q (\alpha, \beta) = 0, \ 
            && \chi_Q (\beta, \alpha) = 0,
        \end{alignedat} \right.
    \end{equation}
    and these cover all the possibilities.
    
    Similarly, let $\alpha \in C (Q)$ and $\theta \in C^\sd (Q)$
    be binary classes,
    such that $\supp \alpha$ is a tree, $\bar{\alpha} + \theta$ is binary,
    and all connected components of $\supp (\bar{\alpha} + \theta)$ are trees.
    
    Let $1^\pl_\alpha \in H_0 (\calM^\pl_\alpha; \bbQ)$
    and $1^\sd_\theta \in H_0 (\calM^\sd_\theta; \bbQ)$
    be classes of a point. Then,
    \begin{equation}
        \label{eq-quiv-binary-heart}
        1^\pl_\alpha \heart 1^\sd_\theta 
        = \left\{ \, \begin{alignedat}[c]{5}
            & {-1^\sd_{\bar{\alpha} + \theta}} \, , \quad
            && \dot{\chi}_Q (\alpha, \theta) = -1, \ 
            && \ddot{\chi}_Q (\alpha) = 0, \ 
            && \dot{\chi}_Q (\alpha^\vee, \theta) = 0, \ 
            && \ddot{\chi}_Q (\alpha^\vee) = 0, \\
            & {-\tfrac{1}{2}} \cdot 1^\sd_{\bar{\alpha} + \theta} \, , \quad
            && \dot{\chi}_Q (\alpha, \theta) = 0, \ 
            && \ddot{\chi}_Q (\alpha) = -1, \ 
            && \dot{\chi}_Q (\alpha^\vee, \theta) = 0, \ 
            && \ddot{\chi}_Q (\alpha^\vee) = 0, \\
            & 1^\sd_{\bar{\alpha} + \theta} \, , \quad
            && \dot{\chi}_Q (\alpha, \theta) = 0, \ 
            && \ddot{\chi}_Q (\alpha) = 0, \ 
            && \dot{\chi}_Q (\alpha^\vee, \theta) = -1, \ 
            && \ddot{\chi}_Q (\alpha^\vee) = 0, \\
            & \tfrac{1}{2} \cdot 1^\sd_{\bar{\alpha} + \theta} \, , \quad
            && \dot{\chi}_Q (\alpha, \theta) = 0, \ 
            && \ddot{\chi}_Q (\alpha) = 0, \ 
            && \dot{\chi}_Q (\alpha^\vee, \theta) = 0, \ 
            && \ddot{\chi}_Q (\alpha^\vee) = -1, \\
            & 0, \quad
            && \dot{\chi}_Q (\alpha, \theta) = 0, \ 
            && \ddot{\chi}_Q (\alpha) = 0, \ 
            && \dot{\chi}_Q (\alpha^\vee, \theta) = 0, \ 
            && \ddot{\chi}_Q (\alpha^\vee) = 0,
        \end{alignedat} \right.
    \end{equation}
    and these cover all the possibilities.
\end{lemma}

\begin{proof}
    \allowdisplaybreaks
    For~\eqref{eq-quiv-binary-lie},
    see Gross--Joyce--Tanaka~\cite[Equation~(6.9)]{GrossJoyceTanaka}.
    To see that these cover all the possibilities,
    note that under our assumptions,
    $-\chi_Q (\alpha, \beta)$
    is the number of edges from $\supp \alpha$ to $\supp \beta$,
    and similarly for $-\chi_Q (\beta, \alpha)$.
    Therefore, since we assumed that all components
    of $\supp (\alpha + \beta)$ are trees,
    at most one of these two numbers is non-zero,
    and it has to be~$1$.

    To prove~\eqref{eq-quiv-binary-heart}, 
    we have by definition
    \begin{align*}
        1^\pl_\alpha \heart 1^\sd_\theta 
        & = \sum_{ \substack{
            i, j \geq 0: \\
            i + j \geq \chi^\sd (\alpha, \> \theta) + 1
        } }
        \frac{(-1)^{\chi^\sd_Q (\alpha, \> \theta)} \cdot
        2^{\ddot{\chi} (\alpha)}}
        {2^j \cdot (i + j - \chi^\sd (\alpha, \theta) - 1)!} \cdot {}
        \\* & \hspace{3em}
        (\oplus^\sd_{\alpha, \> \theta})_* \circ
        (D^{i + j - \chi^\sd (\alpha, \> \theta) - 1} \otimes \id)
        \bigl(
            (1_\alpha \boxtimes 1^\sd_\theta) \cap
            (c_i (\dot{\Theta}_{\alpha, \> \theta}) \cdot c_j (\ddot{\Theta}_\alpha))
        \bigr)
        \\[2ex]
        & = \begin{cases}
            \displaystyle
            \frac{(-1)^{\chi^\sd_Q (\alpha, \> \theta)} \cdot
            2^{\ddot{\chi} (\alpha)}}
            {(-\chi^\sd (\alpha, \theta) - 1)! \vphantom{^0}} \cdot
            (\oplus^\sd_{\alpha, \> \theta})_*
            \bigl( D^{- \chi^\sd (\alpha, \> \theta) - 1} (1_\alpha) 
            \boxtimes 1^\sd_\theta \bigr),
            & \chi^\sd (\alpha, \theta) \leq -1,
            \\ 0,
            & \chi^\sd (\alpha, \theta) \geq 0,
        \end{cases}
        \numberthis
    \end{align*}
    where $1_\alpha \in H_0 (\calM_\alpha; \bbQ)$
    is the class of a point.
    The last step is because
    the cap product is zero unless $i = j = 0$.
    Note that this expression is valid for any~$\alpha$ and $\theta$.

    For $\alpha, \theta$ satisfying our assumptions,
    note that $-\dot{\chi}_Q (\alpha, \theta)$
    is the number of edges from $\supp \alpha$ to $\supp \theta$,
    which is equal to the number of edges
    from $\supp \theta$ to $\supp \alpha^\vee$.
    $-\ddot{\chi}_Q (\alpha)$ is the number of edges
    from $\supp \alpha$ to $\supp \alpha^\vee$.
    Similar descriptions apply to
    $-\dot{\chi}_Q (\alpha^\vee, \theta)$ and 
    $-\ddot{\chi}_Q (\alpha^\vee)$.
    Therefore, by our assumption that all components of
    $\supp (\bar{\alpha} + \theta)$ are trees,
    at most one of these four numbers is non-zero,
    and it has to be~$1$.
    
    Now, \eqref{eq-quiv-binary-heart} follows from the calculation above.
\end{proof}

\begin{lemma}
    Let $\theta \in C^\sd (Q)$ be a binary class,
    and assume that all components of $\supp \theta$
    are self-dual and are trees.

    Then, the classes $\inv{\calM^\sdss_{\theta} (\tau)}$
    defined in Theorem~\textnormal{\ref{thm-quiver-uniqueness}}
    satisfy Theorem~\textnormal{\ref{thm-quiver-main}~\ref{itm-quiver-fund}}.
\end{lemma}

\begin{proof}
    \allowdisplaybreaks
    If $\dim \calM^\sd_\theta < 0$,
    the lemma is trivially true,
    since in this case, $\inv{\calM^\sdss_{\theta} (\tau)}$
    lies in a negative degree homology group,
    and must be zero.
    Therefore, by Lemma~\ref{lem-tree-neg-dim},
    we may assume that $\dim \calM^\sd_\theta = 0$.

    Let $r_\theta$ be the number of connected components of $\supp \theta$.

    Consider the open substack $\calM^\sdss_\theta (\tau) \subset \calM^\sd_\theta$.
    It is a quotient stack $[V^\sdss_\theta (\tau) / G^\sd_\theta]$,
    with $G^\sd_\theta$ as in~\eqref{eq-quiv-gsd},
    and $V^\sdss_\theta (\tau) \subset V^\sd_\theta$
    is an open subset, so it is either empty or connected.
    Since we have `stable~$=$ semistable' for the class~$\theta$ as an assumption,
    by Theorem~\ref{thm-st-sd-decomp}, we have
    \[
        \calM^\sdss_\theta (\tau) \simeq
        [*/\bbZ_2^{\smash{r_\theta}}] \text{ or } \varnothing.
    \]
    Note that the stabilizer group has to be $\bbZ_2^{\smash{r_\theta}}$,
    since by Lemma~\ref{lem-quiv-comp-sd},
    the decomposition in Theorem~\ref{thm-st-sd-decomp}
    must have precisely $r_\theta$~terms.
    Consequently, we have
    \begin{equation}
        \label{eq-quiv-fund-singleton}
        (\iota^\sd_\theta)_* \fund{\calM^\sdss_\theta (\tau)} = \begin{cases}
            2^{-r_\theta} \cdot 1_\theta^\sd \, ,
            & \calM^\sdss_\theta (\tau) \neq \varnothing, \\
            0, & \calM^\sdss_\theta (\tau) = \varnothing,
        \end{cases}
    \end{equation}
    where
    $\iota^\sd_\theta \colon \calM^\sdss_\theta (\tau)
    \hookrightarrow \calM^\sd_\theta$
    is the inclusion map.
    
    On the other hand, consider a version of the
    \emph{naïve Donaldson--Thomas invariants}
    $\upJ^\sharps_\alpha (\tau) \in \bbQ$ and
    $\upJ^\sharpsd_\theta (\tau) \in \bbQ$,
    defined in~\cite[Definition~8.31]{Bu2023}
    for $\alpha \in C (Q)$ and $\theta \in C^\sd (Q)$.
    The invariants~$\upJ^\sharps_\alpha (\tau)$
    were considered in Joyce~\cite{Joyce2008IV},
    while the invariants~$\upJ^\sharpsd_\theta (\tau)$
    were constructed by the author in~\cite{Bu2023}.

    By the definition of these invariants,
    for the class~$\theta$ considered above, we have
    \begin{equation}
        \label{eq-quiv-dt-binary}
        \upJ^\sharpsd_\theta (\tau) = \begin{cases}
            1, & \calM^\sdss_\theta (\tau) \neq \varnothing, \\
            0, & \calM^\sdss_\theta (\tau) = \varnothing.
        \end{cases}
    \end{equation}
    By~\cite[Theorem~8.33]{Bu2023},
    these invariants satisfy the wall-crossing formula
    \begin{align*}
        \upJ^\sharpsd_\theta (\tau)
        & = 
        \sum_{ \leftsubstack[7.5em]{
            \\[-3ex]
            & n \geq 0; \, m_1, \dotsc, m_n > 0; \\[-.5ex]
            & \alpha_{1,1}, \dotsc, \alpha_{1,m_1}; \dotsc;
            \alpha_{n,1}, \dotsc, \alpha_{n,m_n} \in C (Q); \,
            \rho \in \smash{C^\sd (Q)} \colon \\[-.5ex]
            & \theta = (\bar{\alpha}_{1,1} + \cdots + \bar{\alpha}_{1,m_1})
            + \cdots + (\bar{\alpha}_{n,1} + \cdots + \bar{\alpha}_{n,m_n}) + \rho
        } } {}
        \tilde{U}^\sd (\alpha_{1,1}, \dotsc, \alpha_{1,m_1}; \dotsc;
        \alpha_{n,1}, \dotsc, \alpha_{n,m_n}; \mu, \tau) \cdot {}
        \\*[1ex]
        & \hspace{1em}
        \tilde{\chi}^\sd (\alpha_{1,1}, \dotsc, \alpha_{1,m_1}; \dotsc;
        \alpha_{n,1}, \dotsc, \alpha_{n,m_n}; \rho) \cdot {}
        \\*[1ex] 
        & \hspace{3em}
        \bigl( \upJ^\sharps_{\alpha_{1,1}} (\mu) \cdots
        \upJ^\sharps_{\alpha_{1,m_1}} (\mu) \bigr) \cdots
        \bigl( \upJ^\sharps_{\alpha_{n,1}} (\mu) \cdots
        \upJ^\sharps_{\alpha_{n,m_n}} (\mu) \bigr) \cdot
        \upJ^\sharpsd_\rho (\mu)
        \numberthis
        \label{eq-quiv-dt-wcf}
    \end{align*}
    for two self-dual stability conditions~$\mu, \tau$,
    where $\tilde{\chi}^\sd ({\cdots})$ are coefficients defined there.

    Now, consider binary classes $\alpha \in C (Q')$ and $\rho \in C^\sd (Q')$
    satisfying the conditions of Lemma~\ref{lem-quiv-va-on-units}.
    Let $\mu$ be an increasing stability function.
    Then $\calM^\plss_\alpha (\mu)$ is a point if $\alpha$ is primitive,
    or empty otherwise.
    Similarly, $\calM^\sdss_\rho (\mu)$ is $[*/\bbZ_2^{\smash{|\rho|}}]$
    if $\rho$ is primitive, or empty otherwise. We thus have
    \begin{align}
        \label{eq-dt-incr}
        \upJ^\sharps_\alpha (\mu) 
        & = \begin{cases}
            1, & \alpha \text{ primitive}, \\
            0, \hspace{1em} & \text{otherwise},
        \end{cases}
        \\
        \label{eq-dt-sd-incr}
        \upJ^\sharpsd_\rho (\mu) 
        & = \begin{cases}
            1, & \rho \text{ primitive}, \\
            0, \hspace{1em} & \text{otherwise}.
        \end{cases}
    \end{align}

    Consider the involutive Lie algebra $\tilde{A} (Q; \bbQ)$
    and its twisted module $\tilde{A}^\sd (Q; \bbQ)$,
    defined in~\cite[\S8.3]{Bu2023}
    with $\calA = \cat{Mod} (\bbC Q)$ and $\Omega = \bbQ$.
    They are spanned by elements $\omega_\alpha$ and $\omega^\sd_\rho$\,,
    respectively, with $\alpha \in C^\circ (Q)$ and $\rho \in C^\sd (Q)$.
    Note that our $\chi_Q (\alpha, \beta)$ and $\chi^\sd_Q (\alpha, \rho)$
    are denoted by $\chi (\alpha, \beta)$ and $\chi^\sd (\rho, \alpha^\vee)$ there.
    Now, for classes $\alpha, \beta, \rho$
    satisfying the conditions of Lemma~\ref{lem-quiv-va-on-units},
    we have by definition
    \begin{align}
        \label{eq-quiv-binary-lie-motivic}
        [\omega_\alpha \, , \omega_\beta] & =
        \begin{cases}
            -\omega_{\alpha + \beta} \, , &
            \chi_Q (\alpha, \beta) = -1, \ 
            \chi_Q (\beta, \alpha) = 0, \\
            \omega_{\alpha + \beta} \, , &
            \chi_Q (\alpha, \beta) = 0, \ \phantom{-}
            \chi_Q (\beta, \alpha) = -1, \\
            0, &
            \chi_Q (\alpha, \beta) = 0, \ \phantom{-}
            \chi_Q (\beta, \alpha) = 0,
        \end{cases}
        \\
        \label{eq-quiv-binary-heart-motivic}
        \omega_\alpha \heart \omega^\sd_\rho & =
        \begin{cases}
            -\omega^\sd_{\bar{\alpha} + \rho} \, , &
            \chi_Q^\sd (\alpha, \rho) = -1, \ 
            \chi_Q^\sd (\alpha^\vee, \rho) = 0, \\
            \omega^\sd_{\bar{\alpha} + \rho} \, , &
            \chi_Q^\sd (\alpha, \rho) = 0, \ \phantom{-}
            \chi_Q^\sd (\alpha^\vee, \rho) = -1, \\
            0, &
            \chi_Q^\sd (\alpha, \rho) = 0, \ \phantom{-}
            \chi_Q^\sd (\alpha^\vee, \rho) = 0.
        \end{cases}
    \end{align}
    The coefficients $\tilde{\chi}^\sd ({\cdots})$
    in~\eqref{eq-quiv-dt-wcf} are defined so that
    \begin{multline}
        {} [ [ \omega_{\alpha_{1,1}}, \dotsc ] ,
        \omega_{\alpha_{1,m_1}} ] \heart \cdots \heart
        [ [ \omega_{\alpha_{n,1}}, \dotsc ] ,
        \omega_{\alpha_{n,m_n}} ] \heart 
        \omega^\sd_{\rho} =
        \\
        \tilde{\chi}^\sd (\alpha_{1,1}, \dotsc, \alpha_{1,m_1}; \dotsc;
            \alpha_{n,1}, \dotsc, \alpha_{n,m_n}; \rho) \cdot
        \omega^\sd_{
            (\bar{\alpha}_{1,1} + \cdots + \bar{\alpha}_{1,m_1}) + \cdots +
            (\bar{\alpha}_{n,1} + \cdots + \bar{\alpha}_{n,m_n}) + \rho
        }
    \end{multline}
    for $\alpha_{i, \> j} \in C^\circ (Q)$ and $\rho \in C^\sd (Q)$,
    as in~\cite[Lemma~8.24]{Bu2023}.
    Therefore,
    using~\eqref{eq-quiv-binary-lie}--\eqref{eq-quiv-binary-heart}
    and~\eqref{eq-quiv-binary-lie-motivic}--\eqref{eq-quiv-binary-heart-motivic},
    we also have
    \begin{align*}
        &{} [ [ 1_{\alpha_{1,1}}, \dotsc ] ,
        1_{\alpha_{1,m_1}} ] \heart \cdots \heart
        [ [ 1_{\alpha_{n,1}}, \dotsc ] ,
        1_{\alpha_{n,m_n}} ] \heart 
        1^\sd_{\rho} =
        \\ & \hspace{2em}
        2^{r_\rho - r_\theta} \cdot
        \tilde{\chi}^\sd (\alpha_{1,1}, \dotsc, \alpha_{1,m_1}; \dotsc;
            \alpha_{n,1}, \dotsc, \alpha_{n,m_n}; \rho) \cdot {}
        \\ & \hspace{16em}
        1^\sd_{
            (\bar{\alpha}_{1,1} + \cdots + \bar{\alpha}_{1,m_1}) + \cdots +
            (\bar{\alpha}_{n,1} + \cdots + \bar{\alpha}_{n,m_n}) + \rho
        }
        \numberthis
        \label{eq-quiv-binary-heart-multi}
    \end{align*}
    for choices of $\alpha_{i, \> j} \in C (Q)$ and $\rho \in C^\sd (Q)$
    that appear in non-zero terms of~\eqref{eq-quiv-dt-wcf},
    with~$\mu$ increasing.
    Here, $r_\rho$ is the number of connected components of $\supp \rho$,
    and $r_\theta$ is that of $\supp \theta$.
    The factor $2^{r_\rho - r_\theta}$
    comes from the factors $1/2$ in~\eqref{eq-quiv-binary-heart},
    as each occurrence of such a factor results in
    the number of components increasing by one,
    which can be seen from the proof of Lemma~\ref{lem-quiv-va-on-units}.

    Now, consider the formula~\eqref{eq-def-quiver-inv}
    defining the invariants.
    Substituting in~\eqref{eq-quiv-binary-heart-multi}, we obtain
    \begin{align*}
        \inv{ \calM^\sdss_\theta (\tau) } = {} &
        \sum_{ \leftsubstack[8em]{
            \\[-3ex]
            & n \geq 0; \, m_1, \dotsc, m_n > 0; \\[-.5ex]
            & \alpha_{1,1}, \dotsc, \alpha_{1,m_1}; \dotsc;
            \alpha_{n,1}, \dotsc, \alpha_{n,m_n} \in C (Q)\text{ primitive}; \\[-.5ex]
            & \rho \in \smash{C^\sd (Q)} \text{ primitive} \colon \\[-.5ex]
            & \theta = (\bar{\alpha}_{1,1} + \cdots + \bar{\alpha}_{1,m_1})
            + \cdots + (\bar{\alpha}_{n,1} + \cdots + \bar{\alpha}_{n,m_n}) + \rho
        } } {} 
        \tilde{U}^\sd (\alpha_{1,1}, \dotsc, \alpha_{1,m_1}; \dotsc;
        \alpha_{n,1}, \dotsc, \alpha_{n,m_n}; \mu, \tau) \cdot {} 
        \\[.5ex] & \hspace{2em}
        2^{-r_\theta} \cdot 
        \tilde{\chi}^\sd (\alpha_{1,1}, \dotsc, \alpha_{1,m_1}; \dotsc;
            \alpha_{n,1}, \dotsc, \alpha_{n,m_n}; \rho) \cdot
        1^\sd_{\theta} \, .
        \numberthis
    \end{align*}
    Comparing this with~\eqref{eq-quiv-dt-wcf},
    with~\eqref{eq-dt-incr}--\eqref{eq-dt-sd-incr} substituted in,
    we see that
    \begin{equation}
        \inv{ \calM^\sdss_\theta (\tau) } =
        2^{-r_\theta} \cdot \upJ^\sharpsd_\theta (\tau) \cdot 1^\sd_\theta \, .
    \end{equation}
    This, together with~\eqref{eq-quiv-fund-singleton}
    and~\eqref{eq-quiv-dt-binary},
    completes the proof.
\end{proof}

\begin{lemma}
    \label{lem-omega-sd-fund}
    Let $\lambda \colon Q \to Q'$ be a morphism of self-dual quivers
    that induces a bijection on vertices.
    Fix a class $\theta \in C^\sd (Q)$.
    Let $\tau$ be a self-dual stability condition on $Q$, such that
    all $\tau$-semistable self-dual representations of $Q$ of class $\theta$ are $\tau$-stable.
    Then the same is true for $Q'$, and we have
    \begin{align}
        \label{eq-omega-sd-inv-sub}
        \Omega^\sd_\lambda \bigl( \inv{ \calM^\sdss_\theta (\tau) } \bigr) 
        & = \inv{ \calM'^\sdss_\theta (\tau') } \, ,
        \\ \label{eq-omega-sd-fund}
        \Omega^\sd_\lambda \bigl( (\iota^\sd_\theta)_* \fund{ \calM^\sdss_\theta (\tau) } \bigr) 
        & = (\iota'^\sd_\theta)_* \fund{ \calM'^\sdss_\theta (\tau') } \, ,
    \end{align}
    where $\Omega^\sd_\lambda$ is the map in
    Construction~\textnormal{\ref{cons-quiv-omega},}
    and the maps $\iota^\sd_\theta \colon \calM^\sdss_\theta (\tau)
    \hookrightarrow \calM^\sd_\theta$
    and $\iota'^\sd_\theta \colon \calM'^\sdss_\theta (\tau')
    \hookrightarrow \calM'^\sd_\theta$
    are the inclusions.
\end{lemma}

\begin{proof}
    Note that~\eqref{eq-omega-sd-inv-sub} is a special case of Theorem~\ref{thm-omega-sd-inv}.

    For~\eqref{eq-omega-sd-fund}, consider the morphisms
    \begin{alignat*}{2}
        j \colon &&
        \calM'^\sd_\theta & \longrightarrow
        \calM^\sd_\theta \, , \\
        p \colon &&
        \calM^\sd_\theta & \longrightarrow
        \calM'^\sd_\theta \, , 
    \end{alignat*}
    where $j$ is the inclusion given by assigning zero maps to edges in $Q_1 \setminus Q'_1$\,,
    and $p$ is the projection given by restriction from $Q$ to $Q'$.
    We have $p \circ j \simeq \id_{\smash{\calM'^\sd_\theta}}$\,.
    
    Consider the vector bundle
    $\ring{\Xi}_\theta \in \cat{Vect} (\calM^\sd_\theta)$
    defined as in~\eqref{eq-def-xi-ring},
    but with `$\oplus$' in place of `$+$'. In this case, we have
    \begin{equation}
        \ring{\Xi} \sss{\theta}{\vee} \simeq
        \bigoplus_{ a \in Q_1^\tria \setminus Q'_1 } {} \hspace{-.5em}
        \calV \sss{s (a)}{\vee} \otimes \calV_{t (a)}
        \ \oplus 
        \bigoplus_{ a \in Q_1^+ \setminus Q'_1 } {} \hspace{-.5em}
        \Sym^2 (\calV_{t (a)})
        \ \oplus 
        \bigoplus_{ a \in Q_1^- \setminus Q'_1 } {} \hspace{-.5em}
        {\wedge}^2 (\calV_{t (a)}) .
    \end{equation}
    There is a canonical section $s$ of $\ring{\Xi} \sss{\theta}{\vee}$\,, given by
    \begin{equation}
        s \colon 
        \bigl( (E_i, \phi_i)_{i \in Q_0} \, , \, (e_a)_{a \in Q_1} \bigr)
        \longmapsto
        (e_a)_{a \in Q_1 \setminus Q'_1} \, ,
    \end{equation}
    and $j (\calM'^\sd_\theta) = s^{-1} (0) \subset \calM^\sd_\theta$
    is the zero locus of $s$.
    Moreover, we have
    $j (\calM'^\sdss_\theta) = \calM^\sdss_\theta \cap s^{-1} (0)$, and
    $j (\calM'^\sdst_\theta) = \calM^\sdst_\theta \cap s^{-1} (0)$,
    where `st' means `stable'.
    This verifies that $\calM'^\sdss_\theta = \calM'^\sdst_\theta$,
    so it is a smooth, proper Deligne--Mumford stack
    which is a closed substack of $\calM^\sdss_\theta = \calM^\sdst_\theta$.
    
    One can verify that the section $s$
    is transverse to the zero section, so that
    \begin{equation}
        (j|_{\calM'^\sdss_\theta})_* \bigl( \fund{\calM'^\sdss_\theta} \bigr) =
        \fund{\calM^\sdss_\theta} \cap
        c_{\ring{\xi} (\theta)} (\ring{\Xi} \sss{\theta}{\vee} |_{\calM^\sdss_\theta}).
    \end{equation}
    Applying $p_* \circ (\iota'^\sd_\theta)_*$ to both sides, we obtain~\eqref{eq-omega-sd-fund}.
\end{proof}

\begin{lemma}
    \label{lem-quiv-main-binary}
    Let $\theta \in C^\sd (Q)$ be a binary class,
    and assume that all components of $\supp \theta$
    are self-dual.

    Then, the classes $\inv{\calM^\sdss_{\theta} (\tau)}$
    defined in Theorem~\textnormal{\ref{thm-quiver-uniqueness}}
    satisfy Theorem~\textnormal{\ref{thm-quiver-main}~\ref{itm-quiver-fund}}.
\end{lemma}

\begin{proof}
    If $\ring{\chi}_Q (\theta) > 0$, then the result is trivial,
    since $\inv{\calM^\sdss_{\theta} (\tau)} \in
    \smash{H_{-2 \ring{\chi}_Q (\theta)}} (\calM^\sd_\theta)$ and must be zero,
    and so is $\fund{\calM^\sdss_{\theta} (\tau)}$\,.

    We use induction on $-\ring{\chi}_Q (\theta)$.
    For convenience, we assume that $\supp \theta = Q$,
    since otherwise, we could replace $Q$ by $\supp \theta$.
    
    For $a \in Q_1^\tria$\,, let $Q^a$
    be the self-dual quiver obtained from $Q$
    by removing the edges $a$ and $a^\vee$.
    Let $\calM^a$, $\calM^{\smash{a, \> \sd}}$, etc.,
    denote the moduli stacks associated with $Q^a$.

    We use Definition~\ref{def-quiv-homol} to make the identifications
    \[
        H_* (\calM^\sd_\theta; \bbQ) \simeq
        H_* (\calM^{\smash{a, \> \sd}}_\theta; \bbQ) \simeq
        \bbQ [ s_{i, \> 1} : i \in Q_0^\tria ].
    \]
    The map
    \[
        \Omega^\sd_{a} \colon 
        H_* (\calM^\sd_\theta; \bbQ) \longrightarrow
        H_{*-2} (\calM^{\smash{a, \> \sd}}_\theta; \bbQ)
    \]
    is given by taking the cap product with the class
    \begin{alignat*}{2}
        c_{\ring{\xi} (\theta)} (\ring{\Xi} \sss{\theta}{\vee})
        & =
        c_1 ( \calV \sss{s (a)}{\vee} \otimes \calV_{\smash{t (a)}} )\\
        & =
        S \sss{t (a), \> 1}{\theta} -
        S \sss{s (a), \> 1}{\theta} 
        && \in
        H^2 (\calM^\sd_\theta; \bbQ),
        \numberthis
    \end{alignat*}
    and can be written as
    \begin{equation}
        \Omega^\sd_{a} =
        \frac{\partial}{\partial s \sss{t (a), \> 1}{\theta} \vphantom{^\theta}} -
        \frac{\partial}{\partial s \sss{s (a), \> 1}{\theta} \vphantom{^\theta}} \, .
    \end{equation}
    It follows that the map
    \begin{equation}
        \Omegaul^\sd = \bigoplus_{a \in Q_1^\tria} \Omega^\sd_{a} \colon \ 
        H_* (\calM^\sd_\theta; \bbQ) \longrightarrow
        \bigoplus_{a \in Q_1^\tria}
        H_{*-2} (\calM^{\smash{a, \> \sd}}_\theta; \bbQ)
    \end{equation}
    is injective.
    This is because the polynomial ring
    $\bbQ [ s_{i, \> 1} : i \in Q_0^\tria ]$
    is generated by the elements
    $s \sss{t (a), \> 1}{\theta} - s \sss{s (a), \> 1}{\theta}$
    for $a \in Q_1^\tria$\,.
    Note that this uses the fact that all components of $\supp \theta$
    are self-dual.

    However, by Lemma~\ref{lem-omega-sd-fund},
    and by the induction hypothesis applied to the self-dual quivers $Q^a$,
    under the map $\Omegaul^\sd$,
    the classes $\inv{ \calM^\sdss_\theta (\tau) }$ and
    $\iota_* \fund{ \calM^\sdss_\theta (\tau) }$
    are mapped to the same element,
    where $\iota$ is the inclusion into $\calM^\sd_\theta$.
    Since $\Omegaul^\sd$ is injective,
    it follows that $\inv{ \calM^\sdss_\theta (\tau) } =
    \iota_* \fund{ \calM^\sdss_\theta (\tau) }$.
\end{proof}

\begin{lemma}
    \label{lem-quiv-inv-disconn}
    Let $\theta \in C^\sd (Q)$ be a binary class.
    Suppose that $\supp \theta$ has a component that is not self-dual. Then
    \begin{align*}
        \inv{ \calM^\sdss_\theta (\tau) } = 0,
    \end{align*}
    and the classes $\inv{\calM^\sdss_{\theta} (\tau)}$
    satisfy Theorem~\textnormal{\ref{thm-quiver-main}~\ref{itm-quiver-fund}}.
\end{lemma}

\begin{proof}
    Consider the formula~\eqref{eq-def-quiver-inv},
    and consider a non-zero term in it.
    Let $\alpha_{i, \> j}$ and $\rho$ be the classes in that term.

    As in the proof of~\cite[Corollary~6.7]{GrossJoyceTanaka},
    since the element
    \[
        A_i = [ [ [ 
        1^\pl_{\alpha_{i,1}} \, , 
        1^\pl_{\alpha_{i,2}} ] , \dotsc ] ,
        1^\pl_{\alpha_{i,m_i}} ]
    \]
    is assumed to be non-zero,
    we must have that $\supp \alpha_i$ is connected,
    where $\alpha_i = \alpha_{i, \> 1} + \cdots + \alpha_{i, \> m_i}$\,.
    
    This means that there must exist $1 \leq i \leq n$,
    such that the three parts
    \[
        \supp \alpha_i \, , \qquad
        \supp \alpha_i^\vee, \qquad
        \supp (\bar{\alpha}_{i+1} + \cdots + \bar{\alpha}_n + \rho)
    \]
    are not joined by any edges.
    Indeed, a non-self-dual component of $\supp \theta$
    must contain $\supp \alpha_i$ or $\supp \alpha_i^\vee$ for some $i$,
    since $\rho$ does not have any such components as it is primitive.
    We can then just choose the largest $i$ with this property.
    But this means that
    \[
        A_i \heart \bigl( A_{i+1} \heart \cdots
        \heart A_n \heart 1^\sd_\rho \bigr) = 0,
    \]
    since if we write
    $\sigma = \bar{\alpha}_{i+1} + \cdots + \bar{\alpha}_n + \rho$,
    then $\chi^\sd (\alpha, \sigma) = 0$,
    and there are no negative powers of $z$
    in $Y^\sd (A, z) M$ if $A \in H_* (\calM_{\alpha_i}; \bbQ)$
    and $M \in H_* (\calM^\sd_\sigma; \bbQ)$.

    For the final statement,
    by Lemma~\ref{lem-quiv-comp-sd}
    and the assumption that all semistable self-dual objects of class~$\theta$
    are stable, we have $\calM^\ss_\theta (\tau) = \varnothing$.
\end{proof}

Combining Lemmas~\ref{lem-quiv-main-binary} and~\ref{lem-quiv-inv-disconn},
we obtain the following.

\begin{theoremq}
    Let $\theta \in C^\sd (Q)$ be any binary class.
    Then the classes $\inv{\calM^\sdss_{\theta} (\tau)}$
    defined in Theorem~\textnormal{\ref{thm-quiver-uniqueness}}
    satisfy Theorem~\textnormal{\ref{thm-quiver-main}~\ref{itm-quiver-fund}}.
\end{theoremq}

\subsection{Fundamental classes for tame classes}

We can now finish the proof of Theorem~\ref{thm-quiver-main}.
We follow the approach of Gross--Joyce--Tanaka~\cite[\S6.5]{GrossJoyceTanaka},
making use of a result of Martin~\cite[Theorem~B$'$, p.~13]{Martin2000},
which we state below.

\begin{situation}
    \label{sit-moment-map}
    Let $(X, \omega)$ be a symplectic manifold,
    with a Hamiltonian action by a compact, connected Lie group~$G$,
    with Lie algebra~$\frg$.
    Let $T \subset G$ be a maximal torus, with Lie algebra $\frt \subset \frg$.
    Let
    \[
        \mu_G \colon X \longrightarrow \frg^\vee, \qquad
        \mu_T \colon X \longrightarrow \frt^\vee
    \]
    be the moment maps.
    Suppose that $\mu_G$ and $\mu_T$ both have~$0$ as a regular value,
    and that the manifolds $\mu_G^{-1} (0)$ and $\mu_T^{-1} (0)$ are compact.
    Let
    \[
        X \dslash G = \mu_G^{-1} (0) / G, \qquad
        X \dslash T = \mu_T^{-1} (0) / T
    \]
    be the symplectic quotients, which are compact symplectic orbifolds.

    There is the notion of a \emph{lift} of a class
    $a \in H^* (X \dslash G; \bbQ)$ to $\tilde{a} \in H^* (X \dslash T; \bbQ)$.
    This is defined by considering the compact orbifold
    $Y = \mu_G^{-1} (0) / T$, and we say that $\tilde{a}$ is a \emph{lift} of $a$,
    if $\pi^* (a) = i^* (\tilde{a})$,
    where $\pi \colon Y \to X \dslash G$ is the projection,
    and $i \colon Y \to X \dslash T$ is the inclusion.

    Let $\frg = \frt \oplus \frm$ be the $T$-invariant splitting.
    Let $E \to X \dslash T$ be the complex vector bundle given by
    $E = (\mu_T^{-1} (0) \times \frm_{\bbC}) / T$,
    where $\frm_{\bbC} = \frm \otimes_{\bbR} \bbC$
    is seen as a complex representation of~$T$.
\end{situation}

\begin{theorem}[Martin]
    \label{thm-martin}
    In Situation~\textnormal{\ref{sit-moment-map},}
    for any class $a \in H^* (X \dslash G; \bbQ)$,
    with a lift $\tilde{a} \in H^* (X \dslash T; \bbQ)$, we have
    \begin{equation}
        \int_{X \dslash G} a =
        \frac{1}{|W|} \cdot
        \frac{o_T (\mu_T^{-1} (0))}{o_G (\mu_G^{-1} (0))} \cdot
        \int_{X \dslash T} \tilde{a} \cup c_{\mathrm{top}} (E),
    \end{equation}
    where $W$ is the Weyl group of $G$,
    and $o_G (-)$ denotes the number of elements of $G$
    that act trivially on $(-)$.
\end{theorem}

\begin{theorem}
    The classes $\inv{\calM^\sdss_{\theta} (\tau)}$
    defined in Theorem~\textnormal{\ref{thm-quiver-uniqueness}}
    satisfy Theorem~\textnormal{\ref{thm-quiver-main}~\ref{itm-quiver-fund}}
    when $\theta$ is tame.
\end{theorem}

\begin{proof}
    Define subsets $Q \sss{0}{\mathrm{B}}, Q \sss{0}{\mathrm{D}} \subset Q_0^+$
    as follows.
    Let $Q \sss{0}{\mathrm{B}} \subset Q_0^+$
    be the set of vertices~$i$ with $\theta_i$ odd,
    and let $Q \sss{0}{\mathrm{D}} \subset Q_0^+$
    be the set of vertices~$i$ with $\theta_i > 0$ even.

    Write $V = V^\sd_\theta$ and $\bar{G} = G^\sd_\theta$,
    as defined in~\eqref{eq-quiv-vsd}--\eqref{eq-quiv-gsd}, so that
    $\calM^\sd_\theta \simeq [V / \bar{G}]$.
    Let $G \subset \bar{G}$ be the identity component,
    and let $\frg$ be the Lie algebra of~$G$. We may identify
    \begin{equation}
        G \simeq
        \prod_{i \in Q_0^\tria} {}
        \GL (\theta_i; \bbC) \times
        \prod_{i \in Q_0^+} {}
        \SO (\theta_i; \bbC) \times
        \prod_{i \in Q_0^-} {}
        \Sp (\theta_i; \bbC).
    \end{equation}
    Let $\bar{G}_{\bbR} \subset \bar{G}$ and $G_{\bbR} \subset G$
    be the maximal compact subgroups,
    with Lie algebra $\frg_{\bbR}$\,.
    The quotient group $\pi_0 (\bar{G}) = \bar{G} / G$
    is isomorphic to~$\bbZ \sss{2}{|Q_0^{\mathrm{B}}| + |Q_0^{\mathrm{D}}|}$.

    Let $\bar{T} \subset \bar{G}$ be the subgroup
    consisting of the diagonal matrices.
    Here, we use the following convention.
    For the groups $\upO (n; \bbC)$ and $\Sp (n; \bbC)$,
    we use a basis $(e_i)_{1 \leq i \leq n}$ of the underlying vector space
    such that the orthogonal or symplectic form is given by
    $\langle e_i, e_j \rangle = \pm \delta_{i+j, \> n+1}$\,,
    with the sign being negative only for the symplectic group and when $i > j$.
    Explicitly, we have
    \begin{equation}
        \bar{T} \simeq
        \prod_{i \in Q_0^\tria} {}
        (\Gm)^{\theta_i} \times
        \prod_{\substack{i \in Q_0^{\pm} \\
            \mathclap{\theta_i = 2 \tilde{\theta}_i}}} {}
        (\Gm)^{\tilde{\theta}_i} \times
        \prod_{\substack{i \in Q_0^+ \\
            \mathclap{\theta_i = 2 \tilde{\theta}_i + 1}}} {}
        ((\Gm)^{\tilde{\theta}_i} \times \bbZ_2),
    \end{equation}
    where the second and third products involve the vertices~$i$
    with $\theta_i$ even or odd, respectively.
    Let $\bar{T}_{\bbR} = \bar{T} \cap G_{\bbR}$ be its compact real form.

    Let~$T \subset \bar{T}$ be the identity component of $\bar{T}$,
    which is a maximal torus of $G$,
    and let $\frt \subset \frg$ be its Lie algebra.
    Let $T_{\bbR} = T \cap G_{\bbR}$ be the corresponding
    maximal torus of~$G_{\bbR}$\,,
    with Lie algebra $\frt_{\bbR} \subset \frg_{\bbR}$\,.
    The quotient group $\pi_0 (\bar{T}) = \bar{T}/T$ is isomorphic to~%
    $\bbZ \sss{2}{|Q_0^{\mathrm{B}}|}$.

    Consider the self-dual quiver $\tilde{Q}_\theta$
    in Definition~\ref{def-q-tilde}, and let
    $\smash{\tilde{\calM}^\sd_1}$ be the moduli stack of
    self-dual representations of $\tilde{Q}_\theta$
    of class $1 \in C^\sd (\tilde{Q}_\theta)$,
    as in Definition~\ref{def-tame-class}.
    Then we can make the identification
    \[
        \tilde{\calM}^\sd_1 \simeq [V/\bar{T}].  
    \]

    Let $\lambda \colon \tilde{Q}_\theta \to Q$
    be the morphism of quivers sending
    $(i, j) \in (\tilde{Q}_\theta)_0$ to $i \in Q_0$\,,
    and $(a, i, j) \in (\tilde{Q}_\theta)_1$ to $a \in Q_1$\,.
    Then the morphism
    $F^\sd_\lambda \colon \tilde{\calM}^\sd_1 \to \calM^\sd_\theta$
    induced by $\lambda$ coincides with the projection
    $[V/\bar{T}] \to [V/\bar{G}]$.

    By the comparison of symplectic quotients
    and GIT quotients in Kirwan~\cite[\S\S8--9]{Kirwan1984},
    there is a moment map $\mu_{G_{\bbR}} \colon V \to \frg \sss{\bbR}{\vee}$\,,
    and the symplectic quotient $\mu \sss{G_{\bbR}}{-1} (0) / G_{\bbR}$
    can be identified with the GIT quotient $V^\ss (\tau) / G$,
    where $V^\ss (\tau) \subset V$ is the semistable locus
    for the linearization chosen as in the proof of
    Proposition~\ref{prop-quiver-properness}.
    (We omit the square bracket for a quotient stack
    when it can be identified with an orbifold.)

    By Proposition~\ref{prop-quiver-properness},
    the orbifold $V^\ss (\tau) / \bar{G}$ is compact,
    and since we have a principal $Q$-bundle
    $V^\ss (\tau) / G \to V^\ss (\tau) / \bar{G}$,
    we see that $V^\ss (\tau) / G$ is also compact,
    and so is $\mu \sss{G_{\bbR}}{-1} (0)$.
    To see that $0$ is a regular value of $\mu_{G_{\bbR}}$\,,
    for any $x \in \mu \sss{G_{\bbR}}{-1} (0)$,
    the tangent map $d (\mu_{G_{\bbR}})_x$ is non-degenerate,
    since under the symplectic form, it is dual to the infinitesimal action
    $\frg_{\bbR} \to T_x V$, which is injective by Theorem~\ref{thm-st-sd-decomp},
    as no continuous stabilizer groups can occur.

    Similarly, there is a moment map
    $\mu_{T_{\bbR}} \colon V \to \frt \sss{\bbR}{\vee}$\,,
    and the symplectic quotient $\mu \sss{T_{\bbR}}{-1} (0) / T_{\bbR}$
    can be identified with the quotient $V^\ss (\tilde{\tau}) / T$,
    where $\tilde{\tau}$ is the stability condition for $\tilde{Q}_\theta$
    induced from~$\tau$.
    Moreover, $\mu \sss{T_{\bbR}}{-1} (0)$ is compact,
    and $0$~is a regular value of~$\mu_{T_{\bbR}}$\,.

    Now, we may apply Theorem~\ref{thm-martin} with
    the group~$G_{\bbR}$ acting on $X = V$.
    Let $a \in H^* ([V/\bar{G}]; \bbQ)$ be a class.
    Let $b \in H^* (V^\ss (\tau) / G; \bbQ)$
    and $\tilde{b} \in H^* (V^\ss (\tilde{\tau}) / T; \bbQ)$
    be the pullback of $a$ along the natural morphisms.
    Then $\tilde{b}$ is a lift of $b$
    in the sense of Situation~\ref{sit-moment-map},
    and Theorem~\ref{thm-martin} implies that
    \begin{equation}
        \label{eq-quiv-martin}
        \int_{V \dslash G} b =
        \frac{1}{|W|} \cdot
        \int_{V \dslash T}
        \tilde{b} \cup c_{\mathrm{top}} (E),
    \end{equation}
    where $W$ is the Weyl group of~$G$.
    The factors $o_G ({\cdots})$ and $o_T ({\cdots})$
    in Theorem~\ref{thm-martin} are equal,
    since elements of $G$ that act trivially must be
    diagonal matrices with $\pm 1$ on the diagonal,
    which can be deduced from Theorem~\ref{thm-st-sd-decomp}.

    Now, let $c \in H^* (V^\ss (\tau) / \bar{G}; \bbQ)$
    and $\tilde{c} \in H^* (V^\ss (\tilde{\tau}) / \bar{T}; \bbQ)$
    be the pullback of $a$ along the natural morphisms.
    As we will see shortly, $E$ descends to a vector bundle 
    $\bar{E} \to [V / \bar{T}]$.
    Using the principal $\pi_0 (\bar{G})$-bundle
    $V^\ss (\tau) / G \to V^\ss (\tau) / \bar{G}$
    and the principal $\pi_0 (\bar{T})$-bundle
    $V^\ss (\tilde{\tau}) / T \to V^\ss (\tilde{\tau}) / \bar{T}$,
    from~\eqref{eq-quiv-martin}, we obtain
    \begin{equation}
        \label{eq-quiv-martin-non-conn}
        \int_{\calM^\sdss_\theta (\tau)} c =
        \frac{1}{\prod \limits_{i \in Q_0^\tria} \theta_i ! \cdot
        \prod \limits_{i \in Q_0^\pm}
        2^{\tilde{\theta}_i} \, \tilde{\theta}_i !} \cdot
        \int_{\tilde{\calM}^\sdss_1 (\tilde{\tau})}
        \tilde{c} \cup c_{\mathrm{top}} (\bar{E}),
    \end{equation}
    where $\tilde{\theta}_i = \lfloor \theta_i / 2 \rfloor$.
    Here, we used the fact that
    \begin{equation}
        \prod_{i \in Q_0^\tria} \theta_i ! \cdot
        \prod_{i \in Q_0^\pm}
        2^{\tilde{\theta}_i} \, \tilde{\theta}_i !
        =
        2^{|Q_0^{\mathrm{D}}|} \cdot |W|
        =
        \frac{|\pi_0 (\bar{G})|}{|\pi_0 (\bar{T})|} \cdot |W|,
    \end{equation}
    since the Weyl groups of
    $\GL (n; \bbC)$, $\SO (2n+1; \bbC)$, $\Sp (2n; \bbC)$, $\SO (2n; \bbC)$
    have sizes $n!$, $2^n \, n!$, $2^n \, n!$, and $2^{n-1} \, n!$, respectively,
    except for $\SO (0)$, whose Weyl group is trivial.

    We may rewrite~\eqref{eq-quiv-martin-non-conn} as
    \begin{multline*}
        \prod_{i \in Q_0^\tria} \theta_i ! \cdot
        \prod_{i \in Q_0^\pm} 2^{\tilde{\theta}_i} \, \tilde{\theta}_i ! \cdot \bigl(
            a \cdot \iota_* \fund{\calM^\sdss_\theta (\tau)}
        \bigr) = {}
        \\[-2ex]
        a \cdot \Bigl(
            (F^\sd_\lambda)_* \bigl(
                \tilde{\iota}_*
                \fund{\tilde{\calM}^\sdss_1 (\tilde{\tau})}
            \bigr)
            \cap c_{\mathrm{top}} (\bar{E})
        \Bigr),
        \numberthis
    \end{multline*}
    where $\iota \colon V^\ss (\tau) / \bar{G}
    \hookrightarrow [V / \bar{G}]$
    and $\tilde{\iota} \colon V^\ss (\tilde{\tau}) / \bar{T}
    \hookrightarrow [V / \bar{T}]$
    are the inclusions.

    By the definition of~$E$, one can show that
    $E$ descends to the vector bundle
    $\bar{E} = \tilde{\iota}^* (\ring{\Xi}_1)^\vee
    \to \tilde{\calM}^\sdss_1 (\tilde{\tau})$,
    where $\ring{\Xi}_1$
    was defined in Construction~\ref{cons-quiv-omega} as a $K$-theory class,
    but can be interpreted as a vector bundle.
    
    Since the choice of~$a$ was arbitrary, 
    by the definition of the map $\Omega^\sd_\lambda$\,,
    we have
    \begin{equation}
        \prod_{i \in Q_0^\tria} \theta_i ! \cdot
        \prod_{i \in Q_0^\pm} 2^{\tilde{\theta}_i} \, \tilde{\theta}_i ! \cdot
        \iota_* \fund{\calM^\sdss_\theta (\tau)}
        = \Omega^\sd_\lambda \bigl(
            \tilde{\iota}_*
            \fund{\tilde{\calM}^\sdss_1 (\tilde{\tau})}
        \bigr).
    \end{equation}
    Now, by Lemma~\ref{lem-quiv-main-binary},
    we have $\tilde{\iota}_* \fund{\tilde{\calM}^\sdss_1 (\tilde{\tau})}
    = \inv{\tilde{\calM}^\sdss_1 (\tilde{\tau})}$\,.
    By Theorem~\ref{thm-omega-sd-inv}, we see that
    \begin{equation}
        \iota_* \fund{\calM^\sdss_\theta (\tau)} =
        \inv{\calM^\sdss_\theta (\tau)} \, ,
    \end{equation}
    as desired.
\end{proof}

This concludes the proof of Theorem~\ref{thm-quiver-main}.

\newpage
\phantomsection
\addcontentsline{toc}{section}{References}
\sloppy
\renewcommand*{\bibfont}{\normalfont\small}
\printbibliography

\par\noindent\rule{0.38\textwidth}{0.4pt}
{\par\noindent\small
\hspace*{2em}Mathematical Institute, University of Oxford, Oxford OX2 6GG, United Kingdom.\\[-2pt]
\hspace*{2em}Email: \texttt{bu@maths.ox.ac.uk}
}

\end{document}